\documentclass[10pt,twoside,a4paper]{article}

\usepackage{authblk}

\usepackage{geometry}
\usepackage[utf8]{inputenc}
\usepackage{hyperref}

\usepackage{todonotes}

\usepackage[toc,title,page]{appendix}
\usepackage{hyperref}

\usepackage{bbm}
\usepackage{mathtools}
\usepackage{amsmath,amsfonts,amssymb,amsthm}
\usepackage{stackengine} 

\newcommand{\tildhatY}[2]{%
\widetilde{\widehat{Y}\mathrlap{^#1}}_{#2}
}

\newcommand{\tildhatX}[2]{%
	\widetilde{\widehat{X}\mathrlap{^#1}}_{#2}
}

\newcommand{\tildhatUT}[1]{%
	\widetilde{\widehat{U}\mathrlap{^\top}}_{#1}
}

\newcommand{\tildhatsigma}[3]{%
	\widetilde{\widehat{\sigma}\mathrlap{^#1}}_{#2, #3}
}

\usepackage{enumitem}

\usepackage{algorithm}
\usepackage{algpseudocode}

\newtheorem{Theorem}{Theorem}[section]

\newtheorem{Proposition}[Theorem]{Proposition}
\newtheorem{Lemma}[Theorem]{Lemma}
\newtheorem{Assumption}[]{Assumption}
\newtheorem{Corollary}[Theorem]{Corollary}
\newtheorem{Remark}[Theorem]{Remark}

\usepackage[a-1b]{pdfx}

\newcommand{\sbigotimes}{%
	\mathop{\mathchoice{\textstyle\bigotimes}{\bigotimes}{\bigotimes}{\bigotimes}}%
}

\usepackage{graphicx}
\usepackage{subcaption}

\usepackage[maxbibnames=50]{biblatex}
\definecolor{green(colorwheel)(x11green)}{rgb}{0.16, 0.5, 0.0}
\usepackage[normalem]{ulem}
\definecolor{electricultramarine}{rgb}{0.25, 0.0, 1.0}

\newcommand*{\fzcst}[1]{\relax\ifmmode\text{\textcolor{green(colorwheel)(x11green)}{\sout{\ensuremath{#1}}}}\else\textcolor{green(colorwheel)(x11green)}{\sout{#1}}\fi}

\newcommand*{\ykcst}[1]{\relax\ifmmode\text{\textcolor{blue}{\sout{\ensuremath{#1}}}}\else\textcolor{blue}{\sout{#1}}\fi}
\definecolor{ikb}{rgb}{0.0, 0.18, 0.65}

\title{Numerical Methods for Dynamical Low-Rank \\Approximations
of Stochastic Differential Equations \\
Part II: Stochastic discretization}
\author[1]{Yoshihito Kazashi}
\author[2]{Fabio Nobile}
\author[2]{Fabio Zoccolan}

\affil[1]{Department of Mathematics,
	University of Manchester, Oxford Road, Manchester, M13 9PL, UK. email: y.kazashi@manchester.ac.uk}
\affil[2]{Institut de Mathématiques, École Polytechnique Fédérale de Lausanne, 1015 Lausanne, Switzerland. email: fabio.nobile@epfl.ch, fabio.zoccolan@epfl.ch}
\date{}

\begin{document}

\maketitle


\begin{abstract}
	\noindent 
	This paper is the follow up (Part II) of \cite{kazashi2025dynamicalpartI}, which proposed and analyzed time discretization schemes for the Dynamical Low-Rank Approximation (DLRA) of high-dimensional systems of Stochastic Differential Equations (SDEs). In this work, we discuss and analyze the corresponding stochastic discretization. Specifically, we consider Monte Carlo type methods to approximate the stochastic space and all related quantities of interest, leading to produce noisy interacting particle systems whose error analysis is non standard.
	
	Assuming subgaussian tails of the initial condition, in the case of elliptic diffusion, we provide convergence results for the Monte Carlo approximation of the DLR-PS-SDE integrator presented in Part I. The latter requires inversion of the Gram matrix associated to the stochastic modes. When the corresponding empirical Gramian in the Monte Carlo approximation remains full rank over the evolution, we prove a convergence rate close to the usual Monte Carlo one, i.e.\ $O(\frac{1}{\sqrt{M}})$, where $M$ is the sample size. If a regularization of this matrix is employed at each time-step, the rate is close to $O(\frac{1}{\sqrt[3]{M}})$. We also investigate intermediate options where the regularization is activated only when the smallest eigenvalue of the empirical covariance falls below a prescribed tolerance, showing convergence rates in between the one mentioned above.

    Numerical simulations will support this theoretical analysis.
\end{abstract}	
\section{Introduction}
This work concerns the stochastic error analysis of the so-called Dynamical Low-Rank Approximation (DLRA) for Stochastic Differential Equations (SDEs). In detail, here we study the Monte Carlo discretization of the involved stochastic space and how this procedure behaves in terms of propagation of errors. Together with the first part of this project \cite{kazashi2025dynamicalpartI}, which focuses on time integration, the present work introduces the first fully discrete formulation of the DLRA for SDEs in the literature.

There are different ways to toggle the DLRA: here we consider the Dynamically Orthogonal (DO) framework. The DO approximation is a time-dependent reduced order model made by two components, a deterministic basis and a stochastic one, both allowed to evolve in time. Precisely, considering a generic stochastic space $\Omega$, then for any random realization $\omega \in \Omega$ the DO solution is built as the following linear combination of $k$ members:
\begin{equation}\label{eq: DO}
	X(t,\omega) = \sum_{i = 1}^{k} U^{i}(t)Y^{i}(t,\omega), \quad t \geq 0,
\end{equation}
where $\{U^{i}\}_{i=1,\dots,k}$ are the deterministic modes and $\{Y^{i}\}_{i=1,\dots,k}$ is the stochastic basis, respectively, and $k$ the rank of the approximation. Usually, the equations that $U$ and $Y$ satisfy are strongly coupled and their well-posedness is absolutely not trivial. The existence and uniqueness of solutions for the framework \eqref{eq: DO}  in the context of SDEs was first studied in \cite{kazashi2025dynamical}. 

In the latest years, DLRA started to be increasingly and extensively applied in the context of uncertainty quantification, above all in the setting of random partial differential equations and SDEs (see e.g., \cite{carrel2024randomised,sondergaard2013dataI,sondergaard2013dataII,kazashi2021stability,kazashi2021existence,lu2021bayesian,musharbash2015error,musharbash2018dual,musharbash2020symplectic,schmidt2023rank,Vidlickovathesis}). Concerning its simulation in these contexts, reduction of dimensionality can be exploited with respect to the spatial dimension, or to the stochastic one, or to both. Regarding the SDE setting, in many applications, such as filtering of noisy chaotic dynamics, one wants to generate a huge number of SDE paths in order to well estimate the distribution of the studied system. However, assessing this task via full-order model algorithms is often unpractical, as the computational time and storage can be easily out of machine availability. In this sight, DLRA is definitely cheaper from the computational perspective and can show great accuracy when dealing with low-dimensional problems.

A preliminary discretization study has been developed in \cite{kazashi2025dynamicalpartI}, where all the approximation analyses were made in time, leaving untouched the stochastic space. These results were derived under the assumptions of Lipschitzianity and linear-growth bound of the drift and the diffusion of the involved SDE. Three methods were presented: the \textit{DLR Euler-Maruyama} (DLR EM), the \textit{DLR Projector Splitting for SDEs} (DLR PS SDE), and the \textit{DLR Projector Splitting for EM} (DLR PS EM). Under a condition of well-approximability of the dynamics onto its tangent space, similar to the usual one present in the literature (see e.g., \cite{kieri2016discretized,kieri2019projection,Vidlickovathesis}), results of in-time convergence were also provided for all three algorithms. More in details, the DLR EM requires the time step to be smaller than the $k$-th singular value of the Gramian of the stochastic basis to assure boundedness of its second moment. On the other hand, this phenomenon is not observed for the other two algorithms, thanks to the their structure of projected methods. Moreover, DLR PS SDE guarantees error estimates independent of this singular value.

Due to the coupling nature of the DO equations, how to discretize the stochastic framework and to understand how this procedure affects the final error are absolutely not a trivial issue. Therefore, studying the Monte Carlo approximation is of primary interest in the DLRA community.

For standard Itô SDEs, the Monte Carlo procedure consists in generating $M$ independent realizations of the studied equations, making their trajectory evolve through a time integration scheme. For sequential methods, for instance like the Euler-Maruyama \cite{kloeden1992stochastic}, each path of the SDE usually evolves independently from the others. Therefore, to compute the final statistical error, one has just to average all the simulated paths out at final time, obtaining the convergence rate for the root mean square error of $O(\frac{1}{\sqrt{M}})$ thanks to the law of large numbers \cite{kloeden1992stochastic}.

Unlike standard Itô SDEs, the Monte Carlo discretization of the DLRA presents several challenges and, therefore, its analysis is not standard. Indeed, first the DO equations~\cite{kazashi2025dynamicalpartI} are characterized by the presence of the inverse of a Gram matrix and of an expectation in the equation of the deterministic modes. This former term is directly affected by the Monte Carlo procedure and it is not straightforward how the stochastic error would propagate in this case, knowning that a priori the smallest singular value of the Gramian could directly influence the error analysis, similarly to the treatment of Part I. Then, the presence of the expectation in the modes makes the stochastic discretization for DLRA resembles the one of McKean-Vlasov SDEs. Indeed, each of the aforementioned algorithms produces an interacting noisy particle systems as output. Another additional difficulty is erased by the discretization of the deterministic modes $U$: under the Monte Carlo method, the numerical solution of the modes is not longer deterministic with respect to the general probability space. It implies that also the Monte Carlo error analysis concerning the stochastic basis cannot be directly pursued via the standard one for SDEs. Furthermore, all realizations of the stochastic basis are computed using the same deterministic mode, implying that these realizations are not independent. As such, an additional difficulty is erased concerning bounds for the singular values of the Gramian, because the majority of results concerning concentration inequalities required the independence of the realizations (see e.g., \cite{tropp2012user,tropp2015introduction}).

In this paper, we present a self-contained analysis of the Monte Carlo error for the DLR PS SDE. We will denote by $U_n$, $Y_n$ the deterministic and stochastic bases, respectively, obtained by discretize-in-time-only algorithms at time $t_n$, whereas $\widehat{U}_n$, $\widehat{Y}_n$ will indicate the respective ones derived when also the stochastic discretization is employed. First, in Section \ref{sec: problem setting} we provide the general setting of this work. Section \ref{sec: time integration} briefly shows an overview of the computational procedures introduced in \cite{kazashi2025dynamicalpartI} with the related numerical results, whereas in Section \ref{sec: stochastic discretization} the stochastic discretization framework is illustrated. Therein we presented these algorithms under the Monte Carlo discretization and regularization of the Gramian. In Sections \ref{sec: notation} the main notation is explained, where we provide a preliminary discussion on the involved quantities concerning the DLRA.

Section \ref{sec: stoc proj} is devoted to the numerical study of the DLR Projector Splitting for SDEs. Therein we provide some lower bounds on the smallest singular values of the Gramian: these estimates will turn out to be essential to estimate the final error. Subsequently, we derive useful bounds on the moments of the numerical solutions and Lipschitz-type bound between the involved projectors. Following these results, we propose several versions of the DLR PS SDE, which differentiate according to which kind of regularization is employed. Different regularization strategies imply different Monte Carlo converge rates. Moreover, in our error estimate we will not consider a stochastic discretized version of {\cite[Relation (52)]{kazashi2025dynamicalpartI}, i.e.\ the well-approximability assumption of the coefficients of the SDE onto the tangent space of low-dimensional processes. Indeed, in this case, such a condition would employ the discretized expectation and, hence, it would involve a direct dependence on the parameters of the stochastic discretization. Therefore, error estimates might be affected by the smallest singular value of the Gramian of the fully discretized stochastic basis. In case of elliptic diffusion and subgaussian tails of the initial condition, DLR PS SDE is convergent with rate close to the standard Monte-Carlo-type error. For the case of degenerate diffusion, we propose a regularization of it which turns out to a convergent algorithm.

This kind of analysis cannot be converted straightforwardly into the case of the DLR EM and of the DLR PS EM (see Section \ref{sec: DLR EM & KNV}) and, hence, their stochastic convergence remains an open question (even though numerical experiments seems promising). Concerning the former, by being a fully explicit method, it is not immediate to derive bounds of the moments of its numerical solution. About the latter, the presence of the noise increment inside the expectation of the DO equations introduces an irksome additional error term with respect to the analysis of the DLR PS SDE. 

Finally, numerical simulations are shown in Section \ref{sec: numerical experiments}. For the sake of readability, we deferred the proofs of intermediary results presented throughout the paper to the Appendix at the end of this document.

\section{Problem setting and time integration}\label{sec: problem setting}
The SDE framework considered in this work is the same as the one described in our companion paper Part I \cite{kazashi2025dynamicalpartI}. For the sake of completeness, we briefly recall it here along with the three DLRA time-discretization algorithms proposed in \cite{kazashi2025dynamicalpartI}. We then introduce the stochastic discretization, which is the focus of this work. 

Let $\left( \Omega, \mathcal{F}, \mathbb{P}, (\mathcal{F}_t)_{t \geq 0} \right)$ be a filtered complete probability space with the usual conditions \cite[Remark 6.24]{schilling2021brownian}.
For $T>0$, we aim to approximate a $d$-dimensional SDE of the following type
\begin{equation}\label{eq:SDE-diff}
	\mathrm{d}X^{\mathrm{true}}(t) = a(t,X^{\mathrm{true}}_t)\mathrm{d}t+b(t,X^{\mathrm{true}}_t)\mathrm{d}W_t, \quad \text{ for all } t \in [0, T], \quad X^{\mathrm{true}}(0) = X^{\mathrm{true}}_0,
\end{equation}
whose solution at time $t$ is a random vector $X^{\mathrm{true}}(t)=\left(X^{\mathrm{true}}_1(t), \ldots, X^{\mathrm{true}}_d(t)\right)^{\top}$ and where we denote by $W$ a real $m$-dimensional $(\mathcal{F}_t)$-Brownian motion,
i.e., $W(t) = \left(W_1(t), \ldots , W_m(t) \right)^{\top}.$ We assume that $\{W(t)\}_{t\in [0,T]}$ and $X_0$ are independent.
Furthermore, in what follows, we will always work under the following Assumptions.
\begin{Assumption}[Lipschitz coefficients]\label{ass: Lip}
	The drift $a\colon[0,\infty)\times\mathbb{R}^{d}\to\mathbb{R}^{d}$ 
	and the diffusion $b\colon[0,\infty)\times\mathbb{R}^{d}\to\mathbb{R}^{d\times m}$
	are measurable between the Borel fields $\left([0,\infty)\times\mathbb{R}^{d};\mathcal{B}([0,\infty)\times\mathbb{R}^{d})\right)$ and $\left(\mathbb{R}^{d};\mathcal{B}(\mathbb{R}^{d})\right)$, $\left(\mathbb{R}^{d\times m};\mathcal{B}(\mathbb{R}^{d \times m})\right)$, respectively, and Lipschitz continuous with respect to the second variable, uniformly in the first one:
	\begin{equation}
		\begin{cases}
			\!\!\!\! & |a(s,x)-a(s,y)|\leq C_{\mathrm{Lip}}|x-y|, \quad \text{ for all } x,y \in \mathbb{R}^d, \ s \in [0,\infty),\\
			\!\!\!\! & \|b(s,x)-b(s,y)\|_{\mathrm{F}}\leq C_{\mathrm{Lip}}|x-y|, \quad \text{ for all } x,y \in \mathbb{R}^d, \ s \in [0,\infty),
		\end{cases}\label{eq:lip}
	\end{equation}
	for some constant $C_{\mathrm{Lip}}>0$, where $|\cdot|$ and $\|\cdot\|_{\mathrm{F}}$ are the usual Euclidean norm on $\mathbb{R}^{d}$ and the Frobenius one on $\mathbb{R}^{d \times m}$, respectively. 
\end{Assumption}
\begin{Assumption}[Linear growth bound]\label{ass: linear-growth bound}
	The drift $a$
	and the diffusion $b$ fulfil the following linear-growth condition:
	\begin{equation}
		|a(s,x)|^{2}+\|b(s,x)\|_{\mathrm{F}}^{2}\leq C_{\mathrm{lgb}}(1+|x|^{2}), \quad \text{ for all } x \in \mathbb{R}^d, \ s \in [0,\infty),\label{eq:lin-growth}
	\end{equation}
	for some constant $C_{\mathrm{lgb}}>0$.
\end{Assumption}
\begin{Assumption}[Square integrable initial condition]\label{ass: initial value}
	\begin{equation}
		X^{\mathrm{true}}_0 \mbox{ is } \mathcal{F}_0\mbox{-measurable and satisfies } \mathbb{E}[|X^{\mathrm{true}}_0|^2] < +\infty.
	\end{equation}
\end{Assumption}
These assumptions assure the well-posedness of \eqref{eq:SDE-diff} in the sense of strong solutions; see for example \cite[Theorem 21.13]{schilling2021brownian}. 

Our dynamical low-rank approximation is defined by the pair $\left(U,Y\right)$ of processes, with $U_t=(U_t^1,\dots, U_t^k)^{\top} \in \mathbb{R}^{k \times d}$ absolutely continuous and orthogonal for all $t\in [0,T]$, $Y_t$ $k$-dimensional and $\mathcal{F}_t$-adapted solving the so-called DO equations~\cite{kazashi2025dynamical,sapsis2009dynamically}:
\begin{equation}\label{DLR-conditions}
	\begin{aligned}
		C_{Y_t}\dot{U}_t & = \mathbb{E}[Y_t a(t,U_{t}^{\top}Y_t)^{\top}]\left(I_{d \times d} - P^{\text{row} }_{U_t} \right), \\
	\mathrm{d}Y_t  & = U_t a(t,U_{t}^{\top}Y_t) \mathrm{d}t + U_t b(t,U_{t}^{\top}Y_t)\mathrm{d}W_t,
	\end{aligned}
\end{equation}
where $C_{Y_t}:= \mathbb{E}[Y_tY_{t}^{\top}]$ is the Gramian, $P^{\text{row}}_{U_t}$ is the projector onto the vector space $\operatorname{span}\{U^{1}_t, \ldots, U^{k}_t\} \subset \mathbb{R}^{d}$, where $U^{i}_t$ is the $i$-th row of  $U_t$, and reduces to $P^{\text{row}}_{U_t}:= U_t^{\top}U_t$ when $U_t$ has orthonormal rows.  The Dynamical Low-Rank approximation $X$ is then defined by the product $X=U^{\top}Y$. 
Given a suitable rank-$k$ approximation $X_0 :=U_0^{\top}Y_0$, with $U_0$ orthogonal, of $X_0^{\text{true}}$, whose rank is assumed to be at least $k$, 
a unique strong DO solution $(U,Y)$ exists under Assumptions 1--3 \cite{kazashi2025dynamical}, at least for short times. For the sake of readability, we will drop the superscript $\mathrm{row}$ and, write henceforth $P^{\text{row}}_{U_t} = P_{U_t}$.

Lastly, we say that a process $Z \in L^2(\Omega,\mathbb{R}^d)$ is of rank $k$ if $\mathrm{rank}(\mathbb{E}[ZZ^{\top}])=k$ (for more details see e.g. \cite[Section 2.2]{kazashi2025dynamical}).

\subsection{Time integration schemes}\label{sec: time integration}
This section briefly reviews the time integration schemes that we introduced in Part I of this work \cite{kazashi2025dynamicalpartI}.
Let us consider a partition 
\[
\Delta := \{t_n : 0 = t_0 < t_1 < \ldots < t_{N-1} < t_{N} = T\}
\]
of \([0,T]\), with \(\Delta t_n = t_{n+1} - t_n\). We seek reasonable approximations \((U_n)_n\) and \((Y_n)_n\) of \((U(t_n))_n\) and \((Y(t_n))_n\), respectively,  from which we build then an approximate DLRA solution \(X_n\) at time \(t_n\) as the product of these two basis approximations, i.e.,
\[
X_n = U_n^{\top} Y_n \approx X(t_n).
\]
From now on, for the sake of notation, we will write ${a}_n := a(t_n, {X}_n)$, ${b}_n := b(t_n, {X}_n)$, the drift and diffusion computed with the time-discretized numerical solution at time $t_n$, respectively.

Let $\Delta W_n: = W(t_{n+1})-W(t_{n}) \sim \mathcal{N}(0,\Delta t_{n}I_{m \times m})$ be the Brownian increment at time $t_n$, which satisfies the basic properties
$\mathbb{E}[\Delta W_n] = 0,$ $\mathbb{E}[\Delta W_n \Delta W_n^{\top}] = \Delta t_nI_{m \times m},$ for all $n$ and $\mathbb{E}[\Delta W_n \Delta W_m^{\top}] = 0$, for all $m \neq n$.

All the time marching schemes proposed in \cite{kazashi2025dynamicalpartI} start at  $X_0=U_0^{\top}Y_0$.
\begin{itemize}
	\item The \textit{DLR Euler-Maruyama} approximation consists of the following update
	\begin{equation}\label{eq: DLR EM U & Y}
		\begin{aligned}
			C_{Y_n} \widetilde{U}_{n+1} & = C_{Y_n} U_n + \mathbb{E}\left[Y_n a_n^{\top}\right]\left(I_{d \times d} - P_{U_n} \right) \Delta t_n, \\
			\widetilde{Y}_{n+1} &= Y_n + U_n a_n \Delta t_n +  U_n b_n \Delta W_n,
		\end{aligned}
	\end{equation}
	where $C_{Y_n} := \mathbb{E}[Y_nY_n^{\top}]$, and $ P_{U_n}$ is the projector onto the rows of $U_n$, 
	which reduces to $P_{U_n} = U_n^{\top}U_n$ when $U_n$ has orthogonal rows.
	The final update $(U_{n+1}, Y_{n+1})$ with $U_{n+1}$ orthogonal is then constructed by QR decomposition, namely $(U_{n+1}^{\top}, R_{n+1}) = \texttt{QR}(\widetilde{U}_{n+1}^{\top})$ and $Y_{n+1} = R_{n+1} \widetilde{Y}_{n+1}$, where \texttt{QR} is the QR decomposition, that for a rectangular matrix $A \in \mathbb{R}^{n \times k}$ returns $(Q,R) =\texttt{QR}(A)$ with $Q \in \mathbb{R}^{n \times k}$ with orthogonal columns and $R \in \mathbb{R}^{k \times k}$. Clearly, one has $U_{n+1}^{\top}Y_{n+1}= \widetilde{U}_{n+1}^{\top}\widetilde{Y}_{n+1}$.
	We can also derive an update equation for $X_n$ assuming $C_{Y_n}$ invertible, which reads as (see \cite{kazashi2025dynamicalpartI} for more details)
		\begin{equation}\label{eq: DLR EM X}
		\begin{aligned}
			X_{n+1} := & U_{n+1}^{\top}Y_{n+1} =  \widetilde{U}_{n+1}^{\top} \widetilde{Y}_{n+1}\\
			=& X_n  + P_{U_n^{\top}Y_n}[a_n]\Delta t_n+ P_{U_n} b_n \Delta W_n \\
			&+\left(I_{d \times d} - P_{U_n} \right) \mathbb{E}\left[a_n Y_n ^{\top}\right]C^{-1}_{Y_n} \left[U_n a_n (\Delta t_n)^2 + U_n b_n (\Delta W_n) \Delta t_n \right], \\
		\end{aligned}
	    \end{equation}
	    where $P_{U_n^{\top}Y_n}$ denotes the projection onto the tangent space of rank-$k$ processes, defined as 
	    \begin{equation}\label{eq: disc proj}
	    	P_{U_n^{\top}Y_n}[v] := U^{\top}_n U_nv + (I_{d \times d}-U^{\top}_n U_n)\mathbb{E}[vY_n^{\top}]C^{-1}_{Y_n}Y_n.
	    \end{equation}
Assuming that the continuous Gramian 
$C_{Y_t} := \mathbb{E}[Y_tY_t^{\top}]$ and the discretized one $C_{Y_n} := \mathbb{E}[Y_nY_n^{\top}]$ are always of full-rank, convergence in time at the order of $O((\Delta t)^{1/2})$ is proved in \cite{kazashi2025dynamicalpartI}, but with hidden constant inversely proportional to the smallest singular values of $C_{Y_t}$ and $C_{Y_n}$. Therefore, if $C_{Y_t}$ or $C_{Y_n}$ become ``nearly singular'' during the time trajectory, then the error bound might blow up.

\item The \textit{DLR Projector Splitting for SDEs} algorithm approximates the stochastic and the deterministic bases in the following staggered way
\begin{equation}\label{eq: Stoc Proj U & Y}
	\begin{aligned}
		\widetilde{Y}_{n+1} &= Y_n + U_n a_n \Delta t_n +  U_n b_n  \Delta W_n\\
		C_{\widetilde{Y}_{n+1}} \widetilde{U}_{n+1} &=  C_{\widetilde{Y}_{n+1}} U_n + \mathbb{E}\left[\widetilde{Y}_{n+1} \left(a_n^{\top}\Delta t_n\right)\right]\left(I_{d \times d} - P_{U_n} \right),
	\end{aligned}
\end{equation}
followed by an orthonormalization step
$({U}_{n+1}^{\top}, R_{n+1}) = \texttt{QR}(\widetilde{U}_{n+1}^{\top})$ and $Y_{n+1} = R_{n+1} \widetilde{Y}_{n+1}$. 
The corresponding DLR approximation is then given by
\begin{equation} \label{eq: Stoc Proj X}
	\begin{aligned}
		X_{n+1} &= {U}_{n+1}^{\top}{Y}_{n+1} =\widetilde{U}_{n+1}^{\top}\widetilde{Y}_{n+1} \\
		& = X_n + P_{{U}^{\top}_n\widetilde{Y}_{n+1}}[ a_n \Delta t_n] + P_{{U}_n} [b_n] \Delta W_n,
	\end{aligned}
\end{equation} 
Results in Part I show that this method converges in time with rate
$O((\Delta t)^{1/2})$, up to a saturation level given by the projection error of the coefficients of the SDE onto the tangent space to the manifold of rank-$k$ processes (see \eqref{eps bound projection} below); in contrast to the DLR EM scheme, the error and stability estimates do not exhibit any explicit dependence on the smallest singular values of the Gramians appearing in the algorithm.

\item The \textit{DLR Projector Splitting for EM} algorithm is defined by the following step
\begin{equation}\label{eq: KNV Proj Split U & Y}
	\begin{aligned}
		\widetilde{Y}_{n+1} &= Y_n + U_n a_n \Delta t_n +  U_n b_n  \Delta W_n\\
		C_{\widetilde{Y}_{n+1}} \widetilde{U}_{n+1} &=  C_{\widetilde{Y}_{n+1}} U_n + \mathbb{E}\left[\widetilde{Y}_{n+1} \left(a_n^{\top}\Delta t_n+ (\Delta W_n)^{\top}b_n^{\top}\right) \right]\left(I_{d \times d} - P_{U_n} \right),
	\end{aligned}
\end{equation}
again followed by an orthonormalization step,
$({U}_{n+1}^{\top}, R_{n+1}) = \texttt{QR}(\widetilde{U}_{n+1}^{\top})$ and $Y_{n+1} = R_{n+1} \widetilde{Y}_{n+1}$. 
The corresponding update for the DLR approximation is then given by
\begin{equation} \label{eq: KNV Proj Split X}
	\begin{aligned}
		X_{n+1} &= {U}_{n+1}^{\top}{Y}_{n+1} =\widetilde{U}_{n+1}^{\top}\widetilde{Y}_{n+1} \\
		& = X_n + P_{{U}^{\top}_n\widetilde{Y}_{n+1}}[ a_n \Delta t_n+ b_n \Delta W_n].
	\end{aligned}
\end{equation} 
Results in Part I show that this method also converges in time with rate
$O((\Delta t)^{1/2})$ up to the same saturation level as the previous scheme. In contrast to the DLR EM scheme, the numerical stability does not exhibit any explicit dependence on the smallest singular values of the Gramians appearing in the time step. On the other hand, our error analysis involves constants that blow up when the smallest singular value goes to zero, although this degeneracy was not observed in our numerical results, indicating a possible lack of sharpness of our theoretical results.
\end{itemize}

In the convergence analysis of the time integration schemes proposed in~\cite{kazashi2025dynamicalpartI}, the final error is bounded by the sum of three contributions, namely, the error on the initial condition, the Euler-Maruyama error, and the error due to the low-rank approximation (hereafter called ``modelling error"), which is characterized as
\begin{equation}\label{eps bound projection}
		\varepsilon^2 := \sup_{\substack{t \in [0,T] \\ Z \in L^2(\Omega,\mathbb{R}^d) \text{ of rank } k \\
			\text{ with } ·\mathbb{E}[|Z|^2] \leq \tilde{K}}} \max\{\mathbb{E}[|a(t,Z) -  P_{Z}a(t, Z) |^2], \ \mathbb{E}[\| b(t, Z) - \mathcal{P}_{\mathcal{U}(Z)}b(t, Z)\|^2_{\mathrm{F}} ] \},
\end{equation}
where $\tilde{K}:=\min \{ K > 0 \text{ such that } \max\{\mathbb{E}[\sup\limits_{0 \leq t \leq T} |X_t|^2], \mathbb{E}[\max\limits_{0 \leq n \leq N} |X_n|^2] \} \leq K \}$, $P_{Z}[ \ \cdot\ ]:= (I_{d\times d}-\mathcal{P}_{\mathcal{U}(Z)})[\ \cdot \ ]\mathcal{P}_{\mathcal{Y}(Z)}+\mathcal{P}_{\mathcal{U}(Z)}[\ \cdot \ ]$ is the orthogonal projector onto the tangent space to the manifold of random vectors of rank $k$ in the point $Z$ (as in \eqref{eq: disc proj}), whereas $P_\mathcal{U(Z)}$ and $\mathcal{P}_{\mathcal{Y}(Z)}$ denote the projection onto range and corange of $Z$, respectively. We recall that error bounds involving similar definitions of the modelling error have been derived in other notable works concerning dynamical low-rank approximations for deterministic or random dynamical systems, see e.g., \cite{koch2007dynamical,kieri2016discretized,kieri2019projection,nobile2026high}.
\subsection{Stochastic discretization}\label{sec: stochastic discretization}
The algorithms recalled in Section \ref{sec: time integration} assume that the expectation is computable exactly, which is generally not the case in practice. In what follows, we will refer to the solution $X_n=U_n^{\top}Y_n$ obtained by any of the algorithms of the previous section as ``semi-discrete". In this paper, we discuss the discretization of all the stochastic quantities present in \eqref{eq: DLR EM U & Y}, \eqref{eq: Stoc Proj U & Y}, and \eqref{eq: KNV Proj Split U & Y} via the Monte Carlo method. The corresponding fully-discrete versions of the three algorithms are detailed in this section and given in the Algorithms \ref{alg: DLR EM SDE algorithm}, \ref{alg: Stoc Proj algorithm}, and \ref{alg: Eva Proj Splitt SDE algorithm} below.

Before presenting the fully discretized algorithms, we need to introduce some notation. We consider $M$ samples of the initial condition $X_0$, i.e.\ $(X_0(\omega_i)_{i=1,\dots, M}$), and of the Brownian motion $W_t$, i.e.\ $(W_t(\omega_i)_{i=1,\dots, M})$. 
We denote the fully discrete DLRA solution by $(\widehat{U}_n,\widehat{\mathbb{Y}}_n)$, where $\widehat{\mathbb{Y}}_n = \left(\widehat{Y}_n^1, \dots, \widehat{Y}_n^M\right)\in \mathbb{R}^{d \times M}$ is the collection of the $M$ particles in the reduced space, where the dynamics of each particle $\widehat{Y}^{i}$ is driven by an independent Brownian motion $W_t^{i}=W_t(\omega_i)$. In the ambient space, these algorithms propagate $M$ particles $\widehat{X}_n^{i}=\widehat{U}_n \widehat{Y}_n^{i}$, with independent initializations $\widehat{X}_0^{i} = X_0(\omega_i)$. We will refer to these algorithms as ``particle systems".

In our error analysis, we want to compare the $M$ particles $(\widehat{X}_n^{i})_{i=1,\dots,M}$ of the particle systems given by Algorithms \ref{alg: DLR EM SDE algorithm}, \ref{alg: Stoc Proj algorithm}, and \ref{alg: Eva Proj Splitt SDE algorithm} with the $M$ i.i.d.\ realizations $X_n^{i}=X_n(\omega_i)$, $i=1,\dots,M$ of the semi-discrete algorithms of Section \ref{sec: time integration} (without stochastic discretization).

More precisely, consider the sample space $\Omega^M := \Omega \times \dots \times \Omega$ and equip it with the product $\sigma$-algebra and the product measure $\mathbb{P}_M$, i.e.\ $\mathbb{P}_M := \sbigotimes\limits^{M}_{i=1} \mathbb{P}$.  
We are interested in estimating the mean squared error between $X_n(\omega_i)$ and $\widehat{X}_n^{i}$, averaged over all the particles, namely
\begin{equation}\label{eq:def-en}
	e_n : = \sqrt{\mathbb{E}_{\omega_1, \dots \omega_M}\left[\frac{1}{M}\sum_{i=1}^{M} |X_n(\omega_i)-\widehat{X}_n^{i}|^2\right]},
\end{equation}
where $\mathbb{E}_{\omega_1, \dots \omega_M}$ is the expectation with respect to the product measure $\mathbb{P}_M$ over the product space $\Omega^M$. Moreover, the symbol $\mathbb{E}_{\omega_i}[X]$ denotes the (conditional)  expectation of a random variable $X$ on $\Omega^M$ with respect to the sample $\omega_i$-only.
Notice that, for an integer $p\geq 1$, the quantity 
\begin{equation}\label{eq: Lp norm}
	\| \cdot \|_{L^p\left(\Omega^{M}, \mathbb{R}^{d \times M}\right)} :=\sqrt[p]{\mathbb{E}_{\omega_1, \dots \omega_M}\left[\frac{1}{M}\sum_{i=1}^{M} |\cdot|^p\right]}
\end{equation} 
is a norm in $L^p\left(\Omega^{M}, \mathbb{R}^{d \times M}\right)$.
For the sake of notation, we will write $\mathbb{E}$ in place of $\mathbb{E}_{\omega_1, \dots \omega_M}$, unless further precision is needed.

We will henceforth use the following notation: $Z^{i}=Z(\omega_i)$, $i=1,\dots,M$ will denote $M$ i.i.d.\ realizations of a random variable $Z \in L^{p}\left(\Omega; \mathbb{R}^k\right)$, and $\mathbb{Z} = [Z^i]_{i=1,\dots,M} \in \mathbb{R}^{k \times M}$ the matrix whose columns are the $M$ samples extracted from $Z$. Moreover, for any real matrices $a =[a^i]_{i=1,\dots,M}\in \mathbb{R}^{j \times M}$, $b =[b^i]_{i=1,\dots,M}\in \mathbb{R}^{\ell \times M}$ on $\Omega^M$, where $a^i \in \mathbb{R}^{j}$, $b^i \in \mathbb{R}^{\ell}$ are columns of $a$ and $b$, respectively, we denote $\widehat{\mathbb{E}}[ab^{\top}]= \frac{ \sum_{i=1}^{M} a^{i}(b^{i})^{\top} }{M}$ the empirical average of $(a^{i}(b^{i})^{\top})$ where the columns of $a$ and $b$ are interpreted as (not-necessarily independent) particles. 

Notice that all the semi-discrete algorithms presented in Section~\ref{sec: time integration} rely on the invertibility of the Gram matrix $C_{Y_n} = \mathbb{E}[Y_n Y_n^{\top}]$. 
It was shown in~\cite{kazashi2025dynamicalpartI} that under a uniform ellipticity noise assumption (and eventual restriction on $\Delta t_n$ for the DLR-EM scheme) the Gramian of the semi-discrete solution remains positive definite at every time step.
However, this result does not immediately carry over to the Monte Carlo discretization framework, where the expectation is approximated by an empirical average. 
As a result, 
\( \widehat{C}_{\widehat{Y}_n}= \widehat{\mathbb{E}}[\widehat{\mathbb{Y}}_n \widehat{\mathbb{Y}}_n^{\top}] \) may not be invertible.

To ensure the computability of the deterministic modes $\widehat{U}_n$, which implicitly depend on $(\omega_i)_{i=1,\dots, M}$, and the well-posedness of our algorithm, we first consider a regularized Gramian approach where the empirical Gramian $\widehat{C}_{\widehat{Y}_n}$ is replaced by
\begin{equation}\label{eq: reg}
	\widehat{C}_{\widehat{Y}_n}^{\alpha} = \widehat{\mathbb{E}}[\widehat{\mathbb{Y}}_n\widehat{\mathbb{Y}}_n^{\top}]+ \alpha I_{k \times k}=\frac{\sum_{i=1}^{M}\widehat{Y}_n^{i}(\widehat{Y}_n^{i})^{\top}}{M}+ \alpha I_{k \times k}.
\end{equation}
The parameter $\alpha$ may be chosen a priori as a function of $M$ and its optimal scaling with $M$ will be analyzed in the next sections. Obviously, for $\alpha=0$ no regularization is employed. The regularization may not be needed at every time step. This observation opens to other possible adaptive (a posteriori) regularization strategies that will be discussed in Section \ref{sec: error}.

We are now ready to introduce the fully discrete versions of the algorithms described in Section \ref{sec: time integration}, which are given in Algorithms \ref{alg: DLR EM SDE algorithm}, \ref{alg: Stoc Proj algorithm}, and \ref{alg: Eva Proj Splitt SDE algorithm}, respectively. From now on, for the sake of notation, we will write $\widehat{a}_n^{i} := a(t_n, \widehat{X}_n^{i})$, and $\widehat{b}_n^{i} := b(t_n, \widehat{X}_n^{i})$, with $\widehat{X}_n^{i}=\widehat{U}_n^{\top}\widehat{Y}_n^{i}$ the $i$-th particle of the fully discrete solution, and $\widehat{a}_n := (\widehat{a}_n^1, \dots, \widehat{a}_n^M)$, $\widehat{b}_n := (\widehat{b}_n^1, \dots, \widehat{b}_n^M)$.
\begin{algorithm}
	\caption{Monte Carlo DLR Euler--Maruyama}\label{alg: DLR EM SDE algorithm}
	
	\begin{flushleft}
		\textbf{Input}: initial data $U_0$, $\mathbb{Y}_0=U_0\mathbb{X}_0$, number of samples $M$, regularization parameter $\alpha>0$
		
		\textbf{Output:} approximation $\{ \widehat{\mathbb{X}}_n = \widehat{U}_n^{\top}\widehat{\mathbb{Y}}_n \}_{n=0,\ldots, N}$.
	\end{flushleft}
	
	\begin{algorithmic}[1]
		
		\ForAll {$n \in \{0, \ldots, N-1\}$}
		
		\State Generate $M$ Brownian increments $\Delta W_n^{i} \overset{\text{i.i.d.}}{\sim} \mathcal{N}(0,\Delta t_nI_{m \times m})$ with $i \in \{1,\dots,M\}.$
		
		\State Assemble $\widehat{C}_{\widehat{Y}_n}^{\alpha} = \widehat{\mathbb{E}}[\widehat{\mathbb{Y}}_n\widehat{\mathbb{Y}}_n^{\top}]+ \alpha I_{k \times k}$;
		
		\State Compute $ \tildhatY{i}{n+1} = \widehat{Y}_n^{i} + \widehat{U}_n \widehat{a}_n^{i} \Delta t_n +  \widehat{U}_n \widehat{b}_n^{i} \Delta W_n^{i}$, for all $i \in \{1,\dots,M\}.$
		
		\State Compute $ \widehat{C}_{\widehat{Y}_n}^{\alpha} \widetilde{\widehat{U}}_{n+1} =  \widehat{C}_{\widehat{Y}_n}^{\alpha} \widehat{U}_n +\widehat{\mathbb{E}}\left[\widehat{\mathbb{Y}}_n \widehat{a}_n^{\top}\right]\left(I_{d \times d} - P_{\widehat{U}_n} \right) \Delta t_n$
		\State Reorthonormalize the deterministic basis: find  $(\widehat{U}_{n+1}, \widehat{\mathbb{Y}}_{n+1})$ such that:
		\begin{equation*}
			\widehat{U}_{n+1}^{\top} \widehat{\mathbb{Y}}_{n+1} = \tildhatUT{n+1} \widetilde{\widehat{\mathbb{Y}}}_{n+1}, \quad \widehat{U}_{n+1}\widehat{U}_{n+1}^{\top} = I_{k\times k}.
		\end{equation*}
		via $(\widehat{U}_{n+1}^{\top}, \widehat{R}_{n+1}) = \texttt{QR}(\tildhatUT{n+1})$ and $\widehat{Y}_{n+1}^{i} = \widehat{R}_{n+1} \tildhatY{i}{n+1}$, for all $i \in \{1,\dots,M\}.$
		\EndFor
	\end{algorithmic}
\end{algorithm}

\begin{algorithm}[!h]
	\caption{Monte Carlo DLR Projector Splitting for SDEs }\label{alg: Stoc Proj algorithm}
	\begin{flushleft}
		\textbf{Input}: initial data $U_0$, $\mathbb{Y}_0=U_0\mathbb{X}_0$, number of samples $M$, regularization parameter $\alpha>0$
		
		\textbf{Output:} approximation $\{\widehat{\mathbb{X}}_n=\widehat{U}_n^{\top}\widehat{\mathbb{Y}}_n\}_{n=0,\ldots, N}$. 
	\end{flushleft}
	\begin{algorithmic}[1]
		
		\ForAll {$n \in  \{0, \ldots N-1\}$} 
		
		\State Generate $M$ Brownian increments $\Delta W_n^{i} \overset{\text{i.i.d.}}{\sim}  \mathcal{N}(0,\Delta t_n I_{m \times m})$ with $i \in \{1,\dots,M\}.$
		
		\State Compute $ \tildhatY{i}{n+1} = \widehat{Y}_n^{i} + \widehat{U}_n \widehat{a}_n^{i} \Delta t_n +  \widehat{U}_n \widehat{b}_n^{i} \Delta W_n^{i}$, for all $i \in \{1,\dots,M\}.$
		
		\State Assemble $\widehat{C}_{\widetilde{\widehat{Y}}_{n+1}}^{\alpha} = \widehat{\mathbb{E}}[\widetilde{\widehat{\mathbb{Y}}}_{n+1}(\widetilde{\widehat{\mathbb{Y}}}_{n+1})^{\top}]+ \alpha I_{k \times k}$
		\State Compute $\widetilde{\widehat{U}}_{n+1}$: 
		$$\widehat{C}_{\widetilde{\widehat{Y}}_{n+1}}^{\alpha} \widetilde{\widehat{U}}_{n+1} =  \widehat{C}_{\widetilde{\widehat{Y}}_{n+1}}^{\alpha} \widehat{U}_n + \widehat{\mathbb{E}}\left[\widetilde{\widehat{\mathbb{Y}}}_{n+1} \left(\widehat{a}_n^{\top}\Delta t_n\right)\right]\left(I_{d \times d} - P_{\widehat{U}_n} \right) $$
		\State Reorthonormalize the deterministic basis: find   $(\widehat{U}_{n+1}, \widehat{\mathbb{Y}}_{n+1})$ such that:
		\begin{equation*}
			\widehat{U}_{n+1}^{\top} \widehat{\mathbb{Y}}_{n+1} = \tildhatUT{n+1} \widetilde{\widehat{\mathbb{Y}}}_{n+1}, \quad \widehat{U}_{n+1}\widehat{U}_{n+1}^{\top} = I_{k\times k}.
		\end{equation*}
		via $(\widehat{U}_{n+1}^{\top}, \widehat{R}_{n+1}) = \texttt{QR}(\tildhatUT{n+1})$ and $\widehat{Y}_{n+1}^{i} = \widehat{R}_{n+1} \tildhatY{i}{n+1}$, for all $i \in \{1,\dots,M\}.$
		\EndFor
	\end{algorithmic}
\end{algorithm}

\begin{algorithm}[!h]
	\caption{Monte Carlo DLR Projector Splitting for EM}\label{alg: Eva Proj Splitt SDE algorithm}
	\begin{flushleft}
		\textbf{Input}: initial data $U_0$, $\mathbb{Y}_0=U_0\mathbb{X}_0$, number of samples $M$, regularization parameter $\alpha>0$
		
		\textbf{Output:} approximation $\{\widehat{X}_n=\widehat{U}_n^{\top}\widehat{Y}_n\}_{n=0,\ldots, N}$. 
	\end{flushleft}
	\begin{algorithmic}[1]
		
		\ForAll {$n \in  \{0, \ldots, N-1\}$} 
		
		\State Generate $M$ Brownian increments $\Delta W_n^{i} \overset{\text{i.i.d.}}{\sim} \mathcal{N}(0,\Delta t_n I_{m \times m})$ with $i \in \{1,\dots,M\}.$
		
		\State Compute $ \tildhatY{i}{n+1} = \widehat{Y}_n^{i} + \widehat{U}_n \widehat{a}_n^{i} \Delta t_n +  \widehat{U}_n \widehat{b}_n^{i} \Delta W_n^{i}$, for all $i \in \{1,\dots,M\}.$
		
		\State Assemble $\widehat{C}_{\widetilde{\widehat{Y}}_{n+1}}^{\alpha} = \widehat{\mathbb{E}}[\widetilde{\widehat{\mathbb{Y}}}_{n+1}(\widetilde{\widehat{\mathbb{Y}}}_{n+1})^{\top}]+ \alpha I_{k \times k}$
		\State Compute $\widetilde{\widehat{U}}_{n+1}$: 
		$$\widehat{C}_{\widetilde{\widehat{Y}}_{n+1}}^{\alpha} \widetilde{\widehat{U}}_{n+1} =  \widehat{C}_{\widetilde{\widehat{Y}}_{n+1}}^{\alpha} \widehat{U}_n + \widehat{\mathbb{E}}\left[\widetilde{\widehat{\mathbb{Y}}}_{n+1} \left(\widehat{a}_n^{\top}\Delta t_n+ (\Delta \mathbb{W}_n)^{\top}\widehat{b}_n^{\top}\right) \right]\left(I_{d \times d} - P_{\widehat{U}_n} \right),$$
		where $\Delta \mathbb{W}_n : = \bigl[\Delta W_n^1, \dots, \Delta W_n^M\bigr]$
		\State Reorthonormalize the deterministic basis: find   $(\widehat{U}_{n+1}, \widehat{Y}_{n+1})$ such that:
		\begin{equation*}
			\widehat{U}_{n+1}^{\top} \widehat{\mathbb{Y}}_{n+1} = \tildhatUT{n+1} \widetilde{\widehat{\mathbb{Y}}}_{n+1}, \quad \widehat{U}_{n+1}\widehat{U}_{n+1}^{\top} = I_{k\times k}.
		\end{equation*}
		via $(\widehat{U}_{n+1}^{\top}, \widehat{R}_{n+1}) = \texttt{QR}(\tildhatUT{n+1})$ and $\widehat{Y}_{n+1}^{i} = \widehat{R}_{n+1} \tildhatY{i}{n+1}$, for all $i \in \{1,\dots,M\}.$
		\EndFor
	\end{algorithmic}
\end{algorithm}

In Section \ref{sec: stoc proj} we will present a convergence analysis of Algorithm \ref{alg: Stoc Proj algorithm}. 
We defer the analysis of the convergence of Algorithms~\ref{alg: DLR EM SDE algorithm} and~\ref{alg: Eva Proj Splitt SDE algorithm} to Section~\ref{sec: DLR EM & KNV}, since, unlike Algorithm~\ref{alg: Stoc Proj algorithm}, it relies on additional assumptions which will be introduced there.

\subsection{General notation}\label{sec: notation}
For a vector $v \in \mathbb{R}^{m}$, its Euclidean norm is denoted by $|v|$. For a matrix $A\in \mathbb{R}^{m \times n}$, $| A |$ and $\| A \|_{\mathrm{F}}$ denote its spectral and Frobenius norm, respectively. We use the notation $A \succ B$ (respectively  $A \succeq B$) with $A,B$ square matrices, to denote that $A-B$ is positive definite (respectively positive semi-definite). 
For a general matrix $A \in \mathbb{R}^{m \times n}$, $\{\sigma^i(A)\}_{i\geq 1}^{\min \{n,m\}}$ denote its singular values listed in decreasing order, whereas, if $A$ is a square matrix, $\{ \lambda^i(A)\}_{i\geq 1}$ denotes its eigenvalues listed in decreasing order with respect to its real part. 
Moreover, we denote by $\mathbb{E}_n[X]$ the conditional expectation of the random variable $X$ with respect to $\mathcal{F}_{t_n}$, i.e. $\mathbb{E}_n[X]=\mathbb{E}[X | \mathcal{F}_{t_n}]$. 
Given a random vector $X \in L^p(\Omega, \mathbb{R}^d)$, we define the variance of $X$ as
\begin{equation}\label{eq: var}
\mathrm{Var}(X) =\mathbb{E}[|X-\mathbb{E}[X]|^2],
\end{equation}
which is actually the trace of the covariance matrix $C_X := \mathbb{E}[\left(X-\mathbb{E}[X]\right)\left(X-\mathbb{E}[X]\right)^{\top}]$.	
As we work on the sample space $\Omega^M$ with probability $\mathbb{P}_M$ introduced in the previous section, we denote by $\boldsymbol{\omega}:=(\omega_1,\dots, \omega_M)$ an arbitrary element of $\Omega^M$

In later sections, we will need to compare projections done with respect to the exact measure $\mathbb{P}$ and the empirical (product) measure $\mathbb{P}_M$ induced by the particle system. In the following paragraphs, we introduce specific notation and definitions for these projections.

\paragraph{Semi-discretized DLRA $\mathbf{X_n}$.} 
We consider first the semi-discrete DLRA of rank $k$, namely, $\{X_n(\omega)\}_n$ produced by any of the methods introduced in Section \ref{sec: time integration}, with $X_n (\omega)=  U_n^{\top}Y_n(\omega)$, where $U_n \in \mathbb{R}^{d \times k}$  has orthonormal rows and $Y_n \in L^{2}(\Omega, \mathbb{R}^{k})$.
By $Y_n^{i}=Y_n(\omega_i) \in  L^{2}(\Omega_i, \mathbb{R}^{k})$, $i=1,\dots, M$, we denote $M$ i.i.d. copies of $Y_n$, which we collect in a random matrix 
$$\mathbb{Y}_n(\boldsymbol{\omega})= [Y_n(\omega_{1}) \ | \ \dots \ | \ Y_n(\omega_{M})] = [Y^1_n \ | \ \dots \ | \ Y^M_n]  \in \mathbb{R}^{k \times M},$$
with  $\mathbb{Y}_n \in L^2(\Omega^M, \mathbb{R}^{k \times M})$. Likewise, $X_n^{i}= X_n(\omega_i)= U_n^{\top}Y_n^{i}$, $i=1,\dots,M$ denote the $M$ independent realizations (particles) of the semidiscrete solution, 
which we collect into the sample matrix 
\begin{equation}\label{eq: def X sample semidisc}
	\begin{aligned}
		\mathbb{X}_n(\boldsymbol{\omega}) :&= [X_n(\omega_{1}) \ | \ \dots \ | \ X_n(\omega_{M})] = [X^1_n \ | \ \dots \ | \ X^M_n]  = U_n^{\top}\mathbb{Y}_n(\boldsymbol{\omega}) \in \mathbb{R}^{d \times M}.            
	\end{aligned}          
\end{equation}
Notice that each particle $X_n^{i}$ depends only on the variable $\omega_i$ and not on the variables $\omega_j$ for $j \neq i$.
Since $Y_n^{i}$ are all i.i.d., the covariance matrix $C_{Y_n}=\mathbb{E}[Y_n^{i}(Y_n^{i})^{\top}] =\mathbb{E}[\frac{1}{M}\mathbb{Y}_n\mathbb{Y}_n^{\top}] \in\mathbb{R}^{k \times k}$ is the same for all $i$. Its singular values, which coincide with the first $k$ singular values of $\mathbb{E}[X_n^{i}(X_n^{i})^{\top}] \in\mathbb{R}^{d \times d}$ for all $i$, are denoted by $\sigma_n^1,\dots,\sigma_n^k$. 
They are deterministic quantities and if $Y_n$ is of rank $k$, then they are all positive. Likewise, the modes $U_n$ are deterministic and do not depend on $\boldsymbol{\omega}$.

We denote by $\widehat{C}_{Y_n}= \frac{1}{M} \sum_{i=1}^{M} Y_n^{i}(Y_n^{i})^{\top} = \frac{1}{M}\mathbb{Y}_n\mathbb{Y}_n^{\top}$ the empirical covariance of the sample $\mathbb{Y}_n$. Its singular values, which coincide with the $k$ singular values of $\frac{1}{\sqrt{M}} \mathbb{X}_n$, are denoted by $\sigma_{n,M}^1,\dots,\sigma_{n,M}^k$. Notice that these singular values are random variables, i.e.\ $\sigma_{n,M}^j=\sigma_{n,M}^j(\boldsymbol{\omega})$ depending on all variables $\omega_1,\dots, \omega_M$. We introduce analogous notations for the auxiliary stochastic basis $\widetilde{Y}_n$, namely $\widetilde{Y}_n^{i}=\widetilde{Y}_n(\omega_i)$ are i.i.d. copies, 
$\widetilde{\mathbb{Y}}_n = [\widetilde{Y}_n(\omega_{1}) \ | \ \dots \ | \ \widetilde{Y}_n(\omega_{M})],$
is the corresponding sample matrix and $\widetilde{\sigma}_{n,M}^{j}=\widetilde{\sigma}_{n,M}^{j}(\boldsymbol{\omega})$ denotes the $j$-th largest singular value of the empirical Gramian $\widehat{C}_{\widetilde{Y}_n}=  \frac{1}{M} \sum_{i=1}^{M} \widetilde{Y}_n^{i}(\widetilde{Y}_n^{i})^{\top}=\frac{1}{M} \widetilde{\mathbb{Y}}_n\widetilde{\mathbb{Y}}_n^{\top}$.

\paragraph{Fully-discretized DLRA $\mathbf{\widehat{X}_n}$.}
In the case of the fully discretize solutions produced by Algorithms \ref{alg: DLR EM SDE algorithm}, \ref{alg: Stoc Proj algorithm}, and \ref{alg: Eva Proj Splitt SDE algorithm}, the particles $\widehat{X}_n^{i}$ interact each other because of the empirical projection operators. It follows that each particle $\widehat{X}^i_n=\widehat{X}_n^{i}(\boldsymbol{\omega})$ is a function of all the variables $\boldsymbol{\omega}=(\omega_1,\dots, \omega_M)$ in contrast to the semi-discrete case. Likewise, the stochastic particle $\widehat{Y}_n^{i} = \widehat{Y}_n^{i}(\boldsymbol{\omega})$, $i=1, \dots, M$ as well as the ``deterministic modes" $\widehat{U}_n = \widehat{U}_n(\boldsymbol{\omega})$ are functions of $\boldsymbol{\omega}$ because of the use of the empirical averages instead of exact expectations.
	
For the particle system, we can also collect all the particles in sample matrices
\begin{equation*}
	\begin{aligned}
		\widehat{\mathbb{X}}_n &:=  [\widehat{X}^{1}_n(\omega_{1},\dots,\omega_{M}) \ | \ \dots \ | \ \widehat{X}^{M}_n(\omega_{1},\dots,\omega_{M})]  = [\widehat{X}^1_n \ | \ \dots \ | \ \widehat{X}^M_n]  \in \mathbb{R}^{d \times M}, \\
		\widehat{\mathbb{Y}}_n &:= [\widehat{Y}_n^{1}(\omega_{1},\dots,\omega_{M}) \ | \ \dots \ | \ \widehat{Y}_n^{M}(\omega_{1},\dots,\omega_{M})] = [\widehat{Y}_n^1 \ | \ \dots \ | \ \widehat{Y}_n^M \ ]  \in \mathbb{R}^{k \times M}.      
	\end{aligned}          
\end{equation*}
We denote the singular values of the empirical Gramian $\widehat{C}_{\widehat{Y}_n}(\boldsymbol{\omega}) = \frac{1}{M} \sum_{j=1}^M \widehat{Y}_n^{j}(\boldsymbol{\omega})(\widehat{Y}_n^{j}(\boldsymbol{\omega}))^{\top} = \frac{1}{M} \widehat{\mathbb{Y}}_n(\boldsymbol{\omega}) \widehat{\mathbb{Y}}_n(\boldsymbol{\omega})^{\top}$ by $\widehat{\sigma}_{n,M}^1,\dots,\widehat{\sigma}_{n,M}^k$, which coincide with the first $k$ singular values of $\frac{1}{\sqrt{M}} \widehat{\mathbb{X}}_n(\boldsymbol{\omega})$.
Similarly
\begin{equation}\label{eq: def Y_tilde hat sample semidisc}
	\begin{aligned}
		\widetilde{\widehat{\mathbb{Y}}}_n :&=  [\tildhatY{1}{n} \ | \ \dots \ | \ \tildhatY{M}{n} \ ]  \in \mathbb{R}^{k \times M},         
	\end{aligned}          
\end{equation}
denotes the sample matrix of the intermediate stochastic particles and $\tildhatsigma{j}{n}{M}$ the $j$-th largest singular value of the associated empirical Gramian $\frac{1}{M}\widetilde{\widehat{\mathbb{Y}}}_n\widetilde{\widehat{\mathbb{Y}}}_n^{\top}$.

\paragraph{Projectors.} 
The convergence analysis with respect to the stochastic discretization will involve the comparison between the semi- and the fully-discrete DLRA. Consequently, we will need to evaluate differences between projectors appearing in the update of the semi- and fully-discrete solutions. For the sake of clarity, we introduce here the various projector operators that will appear in Section \ref{sec: error}. For this discussion it is helpful to notice that an element $f \in L^p(\Omega^{M};\mathbb{R}^{d \times M})$ can be written as 
	\begin{equation*}
		\begin{aligned}
			f(\boldsymbol{\omega}) &:=  [f^{1}(\omega_{1},\dots,\omega_{M}) \ | \ \dots \ | \ f^{M}(\omega_{1},\dots,\omega_{M})]  = [f^1 \ | \ \dots \ | \ f^M]  \in \mathbb{R}^{d \times M}, \\
		\end{aligned}          
	\end{equation*}
	where $f^{i} \in L^p\left(\Omega^M,\mathbb{R}^d\right)$ is the $i$-th column of $f$. 
	 Furthermore, we define the set of random variables that are $i.i.d.$ and have bounded $p$-moments as
	\begin{equation}\label{eq: Lp norm iid}
			L^p_{iid}\left(\Omega^{M}, \mathbb{R}^{d \times M}\right) :=\left\{f : \Omega^M \to \mathbb{R}^{d \times M} ,  \text{ s.t. } \exists g \in L^p\left(\Omega^{M}, \mathbb{R}^{d}\right) :  f(\boldsymbol{\omega})= \left[g(\omega_1) , \dots , g(\omega_M)\right] \right\}
	\end{equation}
Notice that $L^p_{iid}\left(\Omega^{M}, \mathbb{R}^{d \times M}\right)$ is isomorphic to $L^{p}\left(\Omega, \mathbb{R}^d\right)$ and for $f \in L^{p}\left(\Omega, \mathbb{R}^d\right)$, with little abuse of notation, we will still denote by $f$ the corresponding element in $L^p_{iid}\left(\Omega^{M}, \mathbb{R}^{d \times M}\right)$ i.e.\ $f(\boldsymbol{\omega})= \left[f(\omega_1) , \dots , f(\omega_M)\right]$.

Given $U_n$ and $\widehat{U}_n$, matrices with orthogonal rows associated to $\mathbb{X}_n$ and $\widehat{\mathbb{X}}_n$, respectively, 
\begin{equation}\label{eq: P_Us}
	P_{U_n}:= U_n^{\top} U_n, \text{ and } P_{\widehat{U}_n}:= \widehat{U}_n^{\top} \widehat{U}_n,
\end{equation}
denote the orthogonal projectors onto the range of $\mathbb{X}_n$, $\widehat{\mathbb{X}}_n$, respectively. 
For the treatment of Section \ref{sec: bound proj}, it is important to observe that, in contrast to $P_{{U}_n}$, \textit{the projector $P_{\widehat{U}_n}$ depends on $\boldsymbol{\omega}$}, since $\widehat{U}_n$ depends on $\boldsymbol{\omega}$.

Recall that $X_{n}^i= U_n^{\top}Y_n^{i}$
for $i=1,\dots,M$. 
Let us assume that the Gramian $\mathbb{E}\left[ Y_nY_n^{\top}\right]$ of the semi-discrete solution has smallest singular value $\sigma_n^k>0$.  
Then, we define the following projectors $P_{Y_n^{i}}: L^2(\Omega^{M};\mathbb{R}^{d \times M}) \to L^2(\Omega,\mathbb{R}^{d})$ as
\begin{equation}\label{eq: Y projector}
	\begin{aligned}
		P_{Y_n^{i}}[f] =	& \mathbb{E}\left[ \frac{1}{M} f \mathbb{Y}_n^{\top} \right]\left(\mathbb{E}\left[ \frac{1}{M} \mathbb{Y}_n\mathbb{Y}_n^{\top}\right]\right)^{-1}Y_n^{i} \\
		=& \mathbb{E}\left[ \frac{1}{M} \sum_{i=1}^{M} f^{i} (Y_n^{i})^{\top} \right] C_{Y_n}^{-1} Y_n^{i},\quad f \in L^2(\Omega^{M};\mathbb{R}^{d \times M}),\\
	\end{aligned}
\end{equation}
and $P_{\mathbb{Y}_n}: L^2(\Omega^{M};\mathbb{R}^{d \times M}) \to  L^2(\Omega^{M};\mathbb{R}^{d \times M})$ as 
\begin{equation}\label{eq: P mathbb Y}
	P_{\mathbb{Y}_n}[f]:=[P_{Y_n^{i}}[f]]_{i=1,\dots,M}.
\end{equation}
Notice that $P_{\mathbb{Y}_n}$ is an orthogonal projector with respect to the empirical Frobenius inner product in $\mathbb{R}^{ d \times M}$. 
Moreover if $f \in L^2(\Omega, \mathbb{R}^d)$ and still denote $f$ the corresponding element in $L^{p}_{iid}(\Omega^M,\mathbb{R}^{d \times M})$, then \eqref{eq: Y projector} is equivalent to
\begin{equation}\label{eq: Y projector 2}
	\begin{aligned}
		P_{Y_n^{i}}[f] =&  \mathbb{E}\left[ \frac{1}{M} \sum_{i=1}^{M} f(\omega_{i}) (Y_n(\omega_i))^{\top} \right] C_{Y_n}^{-1} Y_n^{i}	\\
		=& \mathbb{E}\left[ f Y_n^{\top} \right] C_{Y_n}^{-1} Y_n^{i},\quad f \in L^2(\Omega;\mathbb{R}^{d}).\\
	\end{aligned}
\end{equation}
For a positive constant $\alpha$, we define the regularized empirical Gramian $\widehat{C}_{Y_n}^{\alpha} := \left(\frac{\sum_{j=1}^{M} Y_n^j (Y_n^j)^{\top}}{M} + \alpha I_{k \times k} \right)$ and the following regularized projectors
\begin{equation}\label{eq: Y discretized projectors}
	\begin{aligned}
		\widehat{P}^{\alpha}_{Y_n^{i}}[f] =	&\left(\frac{1}{M} f \mathbb{Y}_n^{\top} \right) \left(\frac{1}{M} \mathbb{Y}_n \mathbb{Y}_n^{\top} + \alpha I_{k \times k}\right)^{-1} Y_n^{i} =   \widehat{\mathbb{E}}\left[ f \mathbb{Y}_n^{\top} \right] (\widehat{C}_{Y_n}^{\alpha})^{-1} Y_n^{i},\\
	\end{aligned}
\end{equation}
and $\widehat{P}^{\alpha}_{\mathbb{Y}_n}: L^2(\Omega^{M};\mathbb{R}^{d \times M}) \to  L^2(\Omega^{M};\mathbb{R}^{d \times M})$ as 
\begin{equation*}
	\widehat{P}^{\alpha}_{\mathbb{Y}_n}[f]:=[\widehat{P}^{\alpha}_{Y_n^{i}}[f]]_{i=1,\dots,M}, \text{ for } f \in L^2(\Omega^{M};\mathbb{R}^{d \times M}).
\end{equation*} 
Finally, we define the ``tangent space" projectors 
\begin{itemize}
	\item $P_{U_n^{\top}Y_n^{i}} : L^2\left(\Omega^M, \mathbb{R}^{d \times M} \right) \to L^2\left(\Omega^M, \mathbb{R}^{d} \right)$,
\begin{equation*}
	\begin{aligned}
	P_{U_n^{\top}Y_n^{i}} [f] = & P_{U_n}^{\perp}P_{Y^{i}_{n}}[\,f \,]+  P_{U_n}[\,f^{i} \,] \\
	=&\left(I_{d \times d}-U^{\top}_nU_n\right)\mathbb{E}\left[\frac{1}{M}\,f\, \mathbb{Y}_{n}^{\top}\right]C_{Y_{n}}^{-1}Y^{i}_{n} + U^{\top}_nU_n[\,f^i\,],
	\end{aligned}
\end{equation*}
\item $P_{U_n^{\top}\mathbb{Y}_n} : L^2\left(\Omega^M, \mathbb{R}^{d \times M} \right) \to L^2\left(\Omega^M, \mathbb{R}^{d \times M} \right)$, 
\begin{equation*}
P_{U_n^{\top}\mathbb{Y}_n}[f]:=[P_{U_n^{\top}Y_n^{i}} [f]]_{i=1,\dots,M}, \text{ for } f \in L^2(\Omega^{M};\mathbb{R}^{d \times M}),
\end{equation*}
\item $\widehat{P}_{U_n^{\top}Y_n^{i}}^{\alpha} : L^2\left(\Omega^M, \mathbb{R}^{d \times M} \right) \to L^2\left(\Omega^M, \mathbb{R}^{d} \right)$,
\begin{equation*}
	\begin{aligned}
		\widehat{P}_{U_n^{\top}Y_n^{i}}^{\alpha} [f] = & P_{U_n}^{\perp}\widehat{P}_{Y^{i}_{n}}^{\alpha}[\,f \,]+  P_{U_n}[\,f^{i} \,] \\
		=&\left(I_{d \times d}-U^{\top}_nU_n\right)\widehat{\mathbb{E}}\left[\,f\, \mathbb{Y}_{n+1}^{\top}\right]\left(\widehat{C}_{Y_{n+1}}^{\alpha}\right)^{-1}Y^{i}_{n} + U^{\top}_nU_n[\,f^i\,],
	\end{aligned}
\end{equation*}
and 
\item $\widehat{P}_{U_n^{\top}\mathbb{Y}_n}^{\alpha} : L^2\left(\Omega^M, \mathbb{R}^{d \times M} \right) \to L^2\left(\Omega^M, \mathbb{R}^{d \times M} \right)$ defined as 
\begin{equation*}
	\widehat{P}_{U_n^{\top}\mathbb{Y}_n}^{\alpha}[f]:=[\widehat{P}_{U_n^{\top}Y_n^{i}}^{\alpha} [f]]_{i=1,\dots,M}, \text{ for } f \in L^2(\Omega^{M};\mathbb{R}^{d \times M}).
\end{equation*}
\end{itemize}}
Analogous definitions of $\widehat{P}^{\alpha}_{\widehat{Y}_{n}}$, $\widehat{P}^{\alpha}_{\widehat{U}^{\top}_n\widehat{Y}_{n}}$ can be given for the fully-discrete solution. In all previous projectors, the stochastic basis $\mathbb{Y}_n$ (respectively $\widehat{\mathbb{Y}}_n$) can be replaced by the intermediate one $\widetilde{\mathbb{Y}}_n$ (respectively $\widetilde{\widehat{\mathbb{Y}}}_n$).

\section{Monte Carlo Convergence of the DLR Projector Splitting for SDEs (Algorithm \ref{alg: Stoc Proj algorithm})}\label{sec: stoc proj}
In this section, we will analyze the convergence of Algorithm~\ref{alg: Stoc Proj algorithm}.
 Therefore, here $\widehat{\mathbb{X}}_n$ will always refer to the fully-discrete numerical solution of Algorithm \ref{alg: Stoc Proj algorithm}, whereas $X_n$ will denote its semi-discrete counterpart by relations \eqref{eq: Stoc Proj U & Y} and \eqref{eq: Stoc Proj X} and $\mathbb{X}_n$ its particle version.
We will derive useful Lipschitz bounds between projectors, as well as probability bounds on the smallest singular value of the Gram matrix.
 Depending on the regularization used, we obtain different convergence results with varying stochastic convergence rates.
 To establish these results, we need first to establish bounds on the moments, smallest singular value, and tangent space projectors of the fully-discrete and semi-discrete solutions. These are established in the following subsections.
 
 \subsection{Moments bounds}
\begin{Lemma}[$L^{2p}(\Omega)$-bound on the semi-discrete DLRA $X_n$]\label{lem: L2p norm semidiscretized solution}
	
	For any sequence of $\{\Delta t_n\}_n$, if $\mathbb{E}[|Y_0|^{2}]< + \infty$, then there exists a positive constant $K_2(T)$ such that 
	$$\mathbb{E}[|X_n|^2] \leq K_2(T):=\left(\mathbb{E}[|X_0|^2]+1\right)\exp\{(1 + C_{\mathrm{lgb}}(2+T))T\};$$
	moreover, if $\mathbb{E}[|Y_0|^{2p}] < + \infty$ for some integer $p>1$, and if
	\begin{equation}\label{eq: dt semidiscrete condition}
		\Delta t_n \leq \frac{\widetilde{\sigma}_{n+1}^k}{3C_{\mathrm{lgb}}\left(1 + K_2(T)\right)}, \quad \forall n \in \{1, \dots, N\},
	\end{equation}
	then there exists a positive constant $K_{2p}(T)$ such that for all $n\in\{1,\dots,N\}$ one has
	\begin{equation}\label{eq: Lp moment bound semidiscretized X}
		\mathbb{E}[|\widetilde{Y}_{n}|^{2p}] \leq	\mathbb{E}[|X_n|^{2p}] \leq K_{2p}(T). 
	\end{equation}
	\begin{proof}The proof is deferred to Appendix \ref{app: proof norm bounds}.
	\end{proof}
\end{Lemma}
Notice that the uniform bound in $L^2$ holds for any sequence of time steps $\{\Delta t_n\}_n$. On the other hand, we were not able to obtain $L^{2p}$ bounds with $p>1$ independent of the smallest singular value $\widetilde{\sigma}_n^{k}$ and the time step condition \eqref{eq: dt semidiscrete condition} is needed to compensate such dependence and recover a uniform bound. One may ask whether condition \eqref{eq: dt semidiscrete condition} is ever satisfied or not. In the next subsection, Proposition \ref{prop: lowerbound exp}, we provide a uniform lower bound on $\widetilde{\sigma}_n^{k}$ that holds for any $n$ and any $\{\Delta t_n\}_n$ under uniform ellipticity assumptions on the diffusion. Hence, under Assumption \ref{ass: diff} below, condition \eqref{eq: dt semidiscrete condition} can be replaced by the uniform condition.
\begin{equation}\label{eq: dt semidiscrete condition 2}
	\Delta t_n \leq \left(\frac{\sigma_B}{2\sqrt{3} C_{\mathrm{lgb}}(1+K_2(T))}\right)^{2}.
\end{equation}

\begin{Lemma}[$L^{2p}(\Omega)$-moment bound on the fully-discrete DLRA]\label{lem: Lp moment bound}
If $\mathbb{E}[|Y_0|^{2p}] < + \infty$ for some integer $p\ge 1$ and if
\begin{equation}\label{eq: M dt cond}
	M^{p-1} k^p \Delta t_n^{p}  \leq 1, \quad \forall n \in \{1, \dots, N\},
\end{equation}
then there exists a positive constant $\widehat{K}_{2p}(T)$ such that for all $n\in\{1,\dots,N\}$ and $i \in \{1,\dots,M\}$ one has
\begin{equation}\label{eq: Lp moment bound X_hat Y_hat}
	\mathbb{E}[|\tildhatY{i}{n}|^{2p}] \leq  \mathbb{E}[|\widehat{X}_n^{i}|^{2p}] \leq \widehat{K}_{2p}(T). 
\end{equation}
\begin{proof} The proof is deferred to Appendix \ref{app: proof norm bounds}.
\end{proof}
\end{Lemma}
\begin{Remark}[On the condition on $\Delta t_n$]
One can notice that condition \eqref{eq: M dt cond} does not depend on the smallest singular value $\widetilde{\sigma}_{n+1}^k$, unlike condition \eqref{eq: dt semidiscrete condition}. This is not surprising, as in the fully discretized algorithm a regularization of the Gramian is employed at every time step, allowing for a uniform pathwise norm bound on the regularized projector. 
\end{Remark}
\subsection{Bounds for the smallest singular value of the Gramian} \label{sec: singular values}
In this section, we present some results concerning positive lower bounds on the smallest singular value of the empirical Gramians presented in Section \ref{sec: notation}.

Estimates of the smallest singular value of the Gramian are of primal interest in the DLRA literature and, consequently, also for DLRA for SDEs. Indeed, in Part I \cite{kazashi2025dynamicalpartI}, we investigated the influence of the smallest singular value of the Gramian on the accuracy and stability of the semi-discrete numerical schemes recalled in Section \ref{sec: time integration}. In particular, we showed that the semi-discrete DLR PS SDE scheme is insensitive to the smallest singular value. However, this might not be the case when introducing the stochastic discretization: indeed, our error estimates in this work require comparing different projectors on the same manifold, which inevitably brings a dependence on the curvature of the manifold, hence, on the smallest singular value of the DLR solution, see e.g. \cite[Appendix A]{kieri2019projection}.

An additional difficulty in this analysis comes from the fact that the singular values of any empirical Gramian may depend on the number of samples $M$ employed in the approximation. In this regard, in Section \ref{sec: bound proj} we derive Lipschitz bounds on empirical projection-type operators which involve only the smallest singular values $\sigma_{n}^k$ or $\widetilde{\sigma}_{n}^k$ of the semi-discrete Gramian, quantities that are independent of $M$, and are uniformly lower bounded for all $n$ under the following assumption \cite[Proposition 5.2]{kazashi2025dynamical}. 
\begin{Assumption}[Uniform elliptic diffusion]\label{ass: diff}
	There exists a positive constant $\sigma_{B}$ such that $$b(t, x)b(t,x)^{\top}\succeq\sigma_{B} \cdot I_{d \times d} \succ 0,$$
	for all $t \in [0, \infty)$, and for all $ x \in \mathbb{R}^d$. 
\end{Assumption}

The $M$-independence of our estimates of the smallest singular value $\tildhatsigma{k}{n+1}{M}$ comes from the regularization applied at each time step $t_n$. On the one hand, having $M$-independent lower bounds on $\tildhatsigma{k}{n+1}{M}$ simplifies the convergence analysis. On the other hand, when employing the regularization of the Gramian at each $t_n$, in our theoretical analysis presented in Section \ref{sec: error} we are not able to recover the standard $\frac12$ Monte Carlo convergence rate, although our result might not be sharp. We do recover, however, Monte Carlo rates when activating the regularization adaptively, only when needed, as we will show in Section \ref{sec: adap reg}. We would like to point out that all the results presented here can be proven for DLR EM and DLR PS EM, as well, under additional time-step conditions if needed.

First, we provide a condition that ensures the invertibility of the expectation of the empirical Gramian of the semi-discrete intermediate stochastic basis $\widetilde{\mathbb{Y}}_n$ under Assumption \ref{ass: diff}.

\begin{Proposition}[Lower-boundedness in expectation of the empirical Gramian of $\widetilde{\mathbb{Y}}_{n+1}$]\label{prop: lowerbound exp}
	Suppose that Assumption \ref{ass: diff} is satisfied. Then for any sequence of time-steps $\{\Delta t_n\}_n$ with $\Delta t_n = t_{n+1} - t_n$ and for all $M$, we have:
	\begin{equation}\label{eq: cov 1 Consistent Projector}
		\mathbb{E}\left[\frac{\sum_{j=1}^{M}  \widetilde{Y}^{j}_{n+1}  ( \widetilde{Y}^{j}_{n+1})^{\top} }{M}\right] \succeq \sigma_{B} \Delta t_n \cdot I_{k \times k}.
	\end{equation}
	Moreover, if $C_{Y_0}= \mathbb{E}[Y_0Y_0^{\top}] \succeq \sigma_{Y_0} \cdot I_{k \times k} \succ 0$, where $\sigma_{Y_0}$ is the smallest eigenvalue of $C_{Y_0}$, then for a uniform step-size $\Delta t$ one has
	\begin{equation}\label{eq: cov 2 Consistent Projector}
		\begin{aligned}
		C_{Y_{n+1}}	\succeq \mathbb{E}\left[\frac{\sum_{j=1}^{M}  \widetilde{Y}^{j}_{n+1}  ( \widetilde{Y}^{j}_{n+1})^{\top} }{M}\right] \succeq  & \min \{\sigma_{Y_0}, \frac{\sigma_{B}^2}{4C_{\mathrm{lgb}}(1 + K_2(T))}\}.
		\end{aligned}
	\end{equation}
	\begin{proof}
		Let us recall that $(U_{n+1}^{\top}, R_{n+1}) = \texttt{QR}(\widetilde{U}_{n+1}^{\top})$ and $Y_{n+1} = R_{n+1} \widetilde{Y}_{n+1}$, where \texttt{QR} is the QR decomposition. Then, the proof follows verbatim the one of \cite[Proposition 5.2]{kazashi2025dynamicalpartI}, using the fact that $R_n$ is always of full rank $k$ and by Ostrowski's Theorem \cite[Theorem 4.5.9]{horn2012matrix}. 
	\end{proof}
\end{Proposition} 
\begin{Remark}[Relation between the singular values of intermediate steps]\label{rmk: interm cov}
	One can notice that thanks to the orthogonality of the rows of the deterministic modes $U_n$, the following equality holds
	\begin{equation}
		\begin{aligned}
			\sigma^i\left(\mathbb{E}\left[\frac{\sum_{j=1}^{M}  \widetilde{Y}^{j}_{n+1}  ( \widetilde{Y}^{j}_{n+1})^{\top} }{M}\right]\right)& = \sigma^i\left(U_nU_n^{\top}\mathbb{E}\left[\frac{\sum_{j=1}^{M} \widetilde{Y}^{j}_{n+1} (\widetilde{Y}^{j}_{n+1})^{\top}}{M}\right]U_nU_n^{\top}\right)\\
			&= \sigma^i\left(\mathbb{E}\left[\frac{\sum_{j=1}^{M} U_n^{\top}\widetilde{Y}^{j}_{n+1} (\widetilde{Y}^{j}_{n+1})^{\top}U_n}{M}\right]\right),
		\end{aligned}
	\end{equation}
	for all $i=1, \dots k$, i.e.\ the first $k$ largest singular values of the Gramian $\mathbb{E}\left[\frac{\sum_{j=1}^{M}
		\widetilde{Y}^{j}_{n+1} (\widetilde{Y}^{j}_{n+1})^{\top}}{M}\right]$ are equivalent to the first $k$ largest ones of $\mathbb{E}\left[\frac{\sum_{j=1}^{M} U_n^{\top}\widetilde{Y}^{j}_{n+1} (\widetilde{Y}^{j}_{n+1})^{\top}U_n}{M}\right]$. Therefore, thanks to Proposition \ref{prop: lowerbound exp}, we have
	\begin{equation}\label{eq: cov 1 intermediate}
		\mathbb{E}\left[\frac{\sum_{j=1}^{M} U_n^{\top}\widetilde{Y}^{j}_{n+1} (\widetilde{Y}^{j}_{n+1})^{\top}U_n}{M}\right] \succeq \sigma_{B} \Delta t_n \cdot I_{k \times k},
	\end{equation}
	and, for a uniform time-step $\Delta t$,
	\begin{equation}
		\mathbb{E}\left[\frac{\sum_{j=1}^{M} U_n^{\top}\widetilde{Y}^{j}_{n+1} (\widetilde{Y}^{j}_{n+1})^{\top}U_n}{M}\right] \succeq \min \{\sigma_{Y_0}, \frac{\sigma_{B}^2}{4C_{\mathrm{lgb}}(1 + K_2(T))}\}.
	\end{equation}
	
\end{Remark}

When Assumption \ref{ass: diff} holds, relation \eqref{eq: cov 2 Consistent Projector} tells us that there exists a positive constant $\eta$, independent on $n$ and $M$, such that $C_{Y_{n+1}} \succeq C_{\widetilde{Y}_{n+1}} \succeq \eta$. Under the same assumption we can obtain a similar lower bound in probabilty on the empirical Gramian of the intermediate stochastic basis of the fully-discrete DLRA $\widetilde{\widehat{Y}}_n$. For this goal, the following quantity will be found in several lower bounds of the empirical Gramian
\begin{equation}\label{eq: rho}
	\varrho:=\frac{\sigma_{B}^2}{8 C_{\mathrm{lgb}}(1+2K_2(T))},
\end{equation}
where $K_2(T)$ is defined as in Lemma \ref{lem: L2p norm semidiscretized solution}.

\begin{Proposition}[Lower-bound on the empirical Gramian of $\widetilde{\widehat{\mathbb{X}}}_{n+1}$]\label{prop: lower-bound Gramian}
	Suppose that $\mathbb{E}[|Y_0|^{4p}] < \infty$ for some $p\geq1$, that 
	$\sigma_{Y_0}= \sigma^k\left(\mathbb{E}[Y_0Y_0^{\top}] \right)>0$, that Assumption \ref{ass: diff} holds and assume that $\Delta t_n$ satisfies \eqref{eq: dt semidiscrete condition 2} and \eqref{eq: M dt cond}.  Then, there exists a sequence of nested events $\{E_n\}_n \subset \mathcal{F}$ with $E_n \subset E_{n+1}$ such that for any $n=1, \dots, N$ $\tildhatsigma{k}{n}{M} \geq \min\{ \frac{\varrho}{2};  \frac{\sigma_{Y_0}}{2} \}$ on $E_n$. Moreover, there exists a positive constant $C_5$ independent of $\widetilde{\sigma}^k_n$ and of $n$ such that
	\begin{equation}\label{eq: E_n compl}
	\mathbb{P}(E_n^C) \leq C_5  \left( (1+\frac{1}{\sigma_B^{3p}})\frac{1}{M^p} + \frac{1 }{\sigma_B^{4p}M^{2p-1}} +  \frac{\sigma_B^{p}}{M^p} +  \frac{1 }{M^{2p-1}} \right).
	\end{equation}
	\begin{proof}
		The proof can be found in Appendix \ref{app: Gramian bounds}.
	\end{proof}
\end{Proposition}
Notice that the dependence on $\sigma_B$ in \eqref{eq: E_n compl} derives from the dependence on $\varrho$ defined in \eqref{eq: rho} of the upper bound of $\mathbb{P}(E_n^C) $. 

For numerical stability purposes, it is beneficial to understand whether the covariance matrix of the stochastic basis can be ill-conditioned or not so that one can employ a regularization procedure only if the smallest singular value of the Gramian is below a certain positive quantity. One hopes that this regularization strategy can be beneficial for convergence purposes, obtaining a better convergence rate in the number of samples $M$ with respect to the strategy of regularizing at each $t_n$. Proposition \ref{prop: lower-bound Gramian} will be of help for this strategy (see Algorithm \ref{alg: mod Stoc Proj algorithm}).

Finally, in the proofs of the projectors bounds and the error analysis in the next sections, we need to assure the invertibility of the un-regularized empirical Gramian. Therefore, we define the following set, which we will extensively use in the proofs in the Appendix, as well as in Theorem \ref{thm: convergence of the mod Stoc Proj}, where the empirical Gramians are strictly greater than a prescribed positive constant in a matrix sense:
\begin{equation}\label{eq: def A_n}
	A_n := E_n \cap \{ \widetilde{\sigma}^{k}_{n,M} \geq \frac{\varrho}{2} \},
\end{equation}
where $E_n$ is the event appearing in Proposition \ref{prop: lower-bound Gramian} (its precise definition is given in the proof of the Proposition in Appendix \ref{app: Gramian bounds}, see \eqref{eq: E_n}).

\subsection{Projector bounds}\label{sec: bound proj}
We collect here different Lipschitz-type bounds between the various projectors that will be employed in the error estimates of Section \ref{sec: error}. 

In order to prove these bounds, we need a stronger control on the tails of the distribution of the semi-discrete DLRA $X_n$, namely that they are subgaussian. We therefore introduce the following assumption on the diffusion term $b$ and the initial condition $X_0$, which will allow us to prove in Proposition \ref{Proposition: exp tails Xn} the required subgaussian behavior.

\begin{Assumption}[Boundedness of the diffusion and subgaussian tails of $X_0$]\label{ass: bound diff and tails}
	There exists a positive constant $K_{b}$ such that $\|b(t,x)\|_{\mathrm{F}} \leq K_{b}$, for all $(t,x) \in \mathbb{R}_{+}\times \mathbb{R}^{d}$. Furthermore, there exist positive constants $C_0,\mu$ such that
	\begin{equation}\label{eq:y0-subgaus-tail}
		\begin{aligned}
			\mathbb{P}(|X_0|\ge r)\le C_0e^{-\mu r^2},\qquad \text{ for all } r> 0.
		\end{aligned}
	\end{equation}
\end{Assumption}

\begin{Proposition}[Gaussian tails of the semidiscrete DLRA]\label{Proposition: exp tails Xn}
	Suppose that Assumptions \ref{ass: linear-growth bound}, \ref{ass: initial value}, \ref{ass: diff}, and \ref{ass: bound diff and tails} hold. Then, there exist constants $C_{T}:=C(T)$, $c_{T}:=c(T)>0$ independent of $N$ such that for all $R>0$ we have
	\begin{equation}
		\mathbb{P}(\sup\limits_{0 \leq n \leq N} |X_n| > R) \leq C_{T}e^{-c_{T}R^2}
	\end{equation}
	\begin{proof}
		The proof can be found in Appendix \ref{app: proof proj}.
	\end{proof}
\end{Proposition}
Notice that Assumption \ref{ass: diff} is employed in the proof of Proposition \ref{Proposition: exp tails Xn} in order to obtain a uniform positive lower bound on the smallest singular value of the semidiscrete solution, essential to obtain the estimate on tail of $X_n$ due to the presence of the projector $P_{\widetilde{Y}_{n+1}}$. Any other assumption that provides such a lower bound is suitable to obtain this statement.

We first analyze
the difference between the projector onto the $U$ modes obtained without stochastic discretization and the one $\widehat{U}$ provided via the Monte Carlo method.
\begin{Lemma}[$L^2(\Omega^{M};\mathbb{R}^{d \times M})$-Lipschitz bound on the difference between $P_{U_n}$ and $P_{\widehat{U}_n}$]\label{lem: diff Proj U}
Consider the semi-discrete DLRA solution \( X_n(\boldsymbol{\omega}) \) and the fully-discrete DLRA solution \( \widehat{X}_n(\boldsymbol{\omega}) \), respectively, at the same time \( t_n \).
Consider the operators $P_{U_n}$, $P_{\widehat{U}_n}$ defined in \eqref{eq: P_Us}, i.e.\ the orthogonal projectors onto the range of $X_n(\boldsymbol{\omega})$ and onto the range of $\widehat{\mathbb{X}}_n(\boldsymbol{\omega})$ (as a matrix for a fixed realization $\boldsymbol{\omega}$), respectively. Consider $f \in L^4_{iid}\left(\Omega^{M}, \mathbb{R}^{d \times M}\right)$ such that for positive constants $\tilde{C},\tilde{c}$, one has $\mathbb{P}(|f^{i}| > R) \leq \tilde{C} e^{-\tilde{c}R^2} $ for all $R>0$. 

If $\mathbb{E}[|Y_0|^{4}] < \infty$ and relation \eqref{eq: dt semidiscrete condition} holds, then there exists a positive constant $G_1:=G_1(\sigma_n^k)=O(\frac{1}{\sigma_n^k})$ dependent on the inverse of the smallest singular value $\sigma_n^k$, such that for positive constants $c,C$, it holds
\begin{equation}\label{eq: bound between U proj}
	\begin{aligned}
	\sqrt{\mathbb{E}\left[\frac{1}{M}\sum_{i=1}^{M} \left|\left(P_{U_n}-P_{\widehat{U}_n}\right)f^{i}\right|^2\right]} \leq &
	\left(e_n + \frac{1}{\sqrt{M}}\right) G_1 R +  \sqrt[4]{C} e^{-\frac{cMR^2}{4}} \|f\|_{L^4(\Omega^M;\mathbb{R}^{d \times M})},
	\end{aligned}
\end{equation}
 for all  $R > \|f\|_{L^2(\Omega^M;\mathbb{R}^{d \times M})}$, where $e_n$ is defined in \eqref{eq:def-en}.
\begin{proof}
	The proof is deferred to Appendix \ref{app: proof proj}.
\end{proof}
\end{Lemma}
In the next result we analyze the difference between the $L^2(\Omega)$-projector-type operator obtained from the semi-discrete algorithm and the completely discretized one by the Monte Carlo method, appearing in the particle system, both applied to the same random vector. One can notice that the latter is not an orthogonal projector, as a regularization of the Gramian is used. We show, in particular, how the optimal regularization parameter can be chosen as a function of the number of samples in order to have the best possible convergence rate.
\begin{Lemma}[$L^2(\Omega^{M};\mathbb{R}^{d \times M})$-Lipschitz bound on the difference between $P_{\widetilde{\mathbb{Y}}_{n+1}}$ and $\widehat{P}^{\alpha}_{\widetilde{\mathbb{Y}}_{n+1}}$]\label{lem: error between stochastic proj}
	Consider the semi-discrete intermediate stochastic basis \( \widetilde{Y}_n \) at time \(t_n\).
	Assume that the semi-discrete Gramian $\mathbb{E}\left[\widetilde{Y}_{n+1}\widetilde{Y}_{n+1}^{\top}\right]$ has smallest singular value $\widetilde{\sigma}_{n+1}^k>0$. Recall the operators $P_{\widetilde{\mathbb{Y}}_{n+1}}: L^2(\Omega^{M};\mathbb{R}^{d \times M}) \to  L^2(\Omega^{M};\mathbb{R}^{d \times M})$ defined in \eqref{eq: P mathbb Y} and $\widehat{P}^{\alpha}_{\widetilde{\mathbb{Y}}_{n+1}}: L^2(\Omega^{M};\mathbb{R}^{d \times M}) \to  L^2(\Omega^{M};\mathbb{R}^{d \times M})$ defined in \eqref{eq: Y discretized projectors} with 
a positive regularization parameter $\alpha>0$.

If $\mathbb{E}[|Y_0|^{4}] < \infty$ and relation \eqref{eq: dt semidiscrete condition} holds, then for $\alpha=O(\frac{1}{\sqrt[3]{M}})$ there exists a positive constant $G_2:=G_2(\widetilde{\sigma}_{n+1}^k)=O(\frac{1}{\widetilde{\sigma}_{n+1}^k})$ dependent on the inverse of the smallest singular value $\widetilde{\sigma}_{n+1}^k$ such that the difference between the exact projector $P_{\widetilde{\mathbb{Y}}_{n+1}}$ and its empirical, regularized approximation $\widehat{P}_{\widetilde{\mathbb{Y}}_{n+1}}^{\alpha}$ satisfies
\begin{equation}\label{eq: Y projectors error}
	\begin{aligned}
		\sqrt{\mathbb{E}\left[\frac{1}{M}\sum_{i=1}^{M}
		\left|\left(P_{\widetilde{Y}^{i}_{n+1}}-\widehat{P}^{\alpha}_{\widetilde{Y}^{i}_{n+1}}\right)f\right|^2\right]}  \leq &\frac{1}{\sqrt[3]{M}}G_2 \|f\|_{L^4(\Omega^{M};\mathbb{R}^{d \times M})},
	\end{aligned}
\end{equation} 
for all $f \in L^4_{iid}\left(\Omega^{M}, \mathbb{R}^{d \times M}\right)$. 
\begin{proof}
	The proof is deferred to Appendix \ref{app: proof proj}.
\end{proof}
\end{Lemma}
Finally, we derive a $L^2(\Omega)$-type bound on the difference of two discretized stochastic projectors applied to two different random vectors, respectively. To succeed in this goal, we exploit Proposition \ref{prop: lower-bound Gramian}.
\begin{Lemma}[$L^2(\Omega^{M};\mathbb{R}^{d \times M})$-error between the discretized stochastic projectors $\widehat{P}^{\alpha}_{\widetilde{\widehat{\mathbb{Y}}}_{n+1}}$ and $\widehat{P}^{\alpha}_{\widetilde{\mathbb{Y}}_{n+1}}$]\label{lem: discr proj Y}
	
Suppose that $\mathbb{E}[|Y_0|^{8}] < \infty$, that $\sigma_{Y_0}= \mathbb{E}[Y_0Y_0^{\top}] >0$, and that Assumption \ref{ass: diff} holds. Assume that $\Delta t_n$ satisfies \eqref{eq: dt semidiscrete condition} and \eqref{eq: M dt cond} with $p=4$, and consider $\varrho>0$ as defined in \eqref{eq: rho}.
Furthermore, consider the operators $\widehat{P}^{\alpha}_{\widetilde{Y}_{n+1}^{i}},\widehat{P}^{\alpha}_{\tildhatY{i}{n+1}}: L^2(\Omega^{M};\mathbb{R}^{d \times M}) \to L^2(\Omega^M,\mathbb{R}^{d})$ defined as in \eqref{eq: Y discretized projectors} for a positive regularization parameter $\alpha$. 

Consider $f \in L^4_{iid}\left(\Omega^{M}, \mathbb{R}^{d \times M}\right)$ such that for positive constants $\tilde{C},\tilde{c}$, one has $\mathbb{P}(|f^{i}| > R) \leq \tilde{C} e^{-\tilde{c}MR^2} $ for all $R>0$. 
Then, choosing $\alpha=O(\frac{1}{\sqrt[3]{M}})$, there exists a positive constant $G_3:=G_3(\widetilde{\sigma}_{n+1}^k,\varrho)=O(\max\{\frac{1}{\widetilde{\sigma}_{n+1}^k}; \frac{1}{\varrho}\})$ dependent on the inverse of $\widetilde{\sigma}_n^k$ and the inverse of $\varrho$ such that for positive constants $c,C$ it holds
\begin{equation}\label{eq: discr proj Y}
	\begin{aligned}
	\sqrt{\mathbb{E}\left[\frac{1}{M}\sum_{i=1}^{M}
	\left|\left(\widehat{P}^{\alpha}_{\tildhatY{i}{n+1}}-\widehat{P}^{\alpha}_{\widetilde{Y}_{n+1}^{i}}\right)f\right|^2\right]} \leq &  \left(e_n + \frac{1}{\sqrt{M}}\right)   G_3 R +  \left(\frac{1}{\sqrt[3]{M}} +  \sqrt[4]{C} e^{-\frac{cMR^2}{4}}\right) \|f\|_{L^4(\Omega^M;\mathbb{R}^{d \times M})}, 
\end{aligned}
\end{equation}
	for all $R > \|f\|_{L^2(\Omega^M;\mathbb{R}^{d \times M})}$.
\begin{proof}
	The proof is deferred to Appendix \ref{app: proof proj}.
\end{proof}
\end{Lemma}

\begin{Remark}[On the choice of $p$ in relation \eqref{eq: M dt cond}]
		Relation \eqref{eq: M dt cond} links the time-step size to the number of samples through the parameter $p$, which characterizes the order of moment boundedness of the initial condition. At first sight, one would therefore choose $p$ as small as possible. However, in order to satisfy the assumptions of Lemma \ref{lem: discr proj Y}, one must take at least $p=4$, which leads to a condition $\Delta t \propto \frac{1}{M^{\frac34}}$. Indeed, we need $\mathbb{E}[|Y_0|^{8}] < \infty$ in order to guarantee the sought bond in $M$ in relation \eqref{eq: bound requiring 8th} by Proposition \ref{prop: lower-bound Gramian}.
\end{Remark}
\subsection{Error analysis}\label{sec: error}
In this section, we study the convergence analysis of Algorithm \ref{alg: Stoc Proj algorithm}. 
More precisely, we analyze the $L^2$ error over all the particles between the DLR approximation obtained without the stochastic discretization and the one with the Monte Carlo method. 
By regularizing the discretized Gramian at each time step, with parameter $\alpha_{\mathrm{opt}} \propto M^{-1/3}$, we show that the method converges with a rate close to the one of the Monte Carlo estimator. The resulting error bound depends only on the smallest singular values of the semi-discrete Gramians, which is uniformly lower-bounded under Assumption \ref{ass: diff}.

\begin{Theorem}[Monte Carlo Convergence of the DLR Projector Splitting for SDEs]\label{thm: convergence of the Stoc Proj}
 Assume that $\mathbb{E}[|Y_0|^{8}] < + \infty$, Assumptions \ref{ass: diff} and \ref{ass: bound diff and tails}, and relations \eqref{eq: dt semidiscrete condition}-\eqref{eq: M dt cond} with $p=4$ hold. Moreover, let $\gamma$ be a positive constant such that $\sigma_n^k\geq \gamma >0$  for all $n=1,\dots,N$,  $\widetilde{\sigma}_{n+1}^k\geq \gamma >0$ for all $n$ such that $n=0,\dots,N-1$, and $\varrho\geq \gamma$, where $\varrho$ is defined in \eqref{eq: rho} (the existence of such $\gamma$ is guaranteed by \cite[Proposition 5.2]{kazashi2025dynamicalpartI} and Proposition \ref{prop: lowerbound exp} under Assumptions \ref{ass: diff} and \ref{ass: bound diff and tails}, and conditions \eqref{eq: dt semidiscrete condition}-\eqref{eq: M dt cond}). Then
	\begin{equation}\label{eq: error for the Stoc Proj}
		e_{n+1} \leq M^{-\frac{1}{3}+ \frac{\log \log M}{2\log M}CT} \exp \left(C T\right), \quad  \text{for all } \ 0 \leq n \leq N-1,
	\end{equation}
	with positive constant $C:=C(\gamma,T)$ independent of $M$, $\Delta t$, but dependent on the inverse of $\gamma$.
	\begin{proof}
		The strategy of the proof is the following: we want to derive a recursion on the error $e_n$ so that we can conclude via Gronwall's lemma tracking explicitly the dependence of all constants on $\frac{1}{M}$. In order to obtain such bound, we compare the quantities concerning the updates of $X_n$ and $\widehat{X}_n$ and we bound them according to the results of Section \ref{sec: bound proj}.
		
		Let us recall that the error we want to estimate is
		\begin{equation*}
			e_{n+1} =  \sqrt{\mathbb{E}\left[\frac{1}{M}\sum_{i=1}^{M} |X_{n+1}^{i}-\widehat{X}_{n+1}^{i}|^2\right]}.
		\end{equation*}
		First, notice that by construction 
		\begin{equation}\label{eq: e_0}
			e_0 = \sqrt{\mathbb{E}\left[\frac{1}{M}\sum_{i=1}^{M} |X_{0}^{i}-\widehat{X}_{0}^{i}|^2\right]} = \sqrt{\mathbb{E}\left[\frac{1}{M}\sum_{i=1}^{M} |X_{0}^{i}-X_{0}^{i}|^2\right]} =0.
		\end{equation}
		Then, expanding the definition of $e_{n+1}^2$ and using the scheme \eqref{eq: Stoc Proj X}, we obtain
		\begin{equation}\label{eq: intermediate error}
			\begin{aligned}
				e^2_{n+1} = &  \mathbb{E}\left[\frac{1}{M}\sum_{i=1}^{M} |X_{n+1}^{i}-\widehat{X}_{n+1}^{i}|^2\right]\\
				= &  \mathbb{E}\left[\frac{1}{M}\sum_{i=1}^{M} |X_{n}^{i}+P_{{U}^{\top}_n\widetilde{Y}_{n+1}^{i}}[ a_n] \Delta t_n+ P_{{U}_n}[b_n^{i} \Delta W_n^{i}] -\widehat{X}_{n}^{i}-\widehat{P}^{\alpha}_{\widehat{U}^{\top}_n\tildhatY{i}{n+1}}[\widehat{a}_n] \Delta t_n- P_{\widehat{U}_n}[\widehat{b}_n^{i} \Delta W_n^{i}] |^2\right]\\
				=& e_{n}^2 + 2\mathbb{E}\left[\frac{1}{M}\sum_{i=1}^{M} \langle X_{n}^{i}-\widehat{X}_{n}^{i},P_{{U}^{\top}_n\widetilde{Y}_{n+1}^{i}}[ a_n] \Delta t_n -\widehat{P}^{\alpha}_{\widehat{U}^{\top}_n\tildhatY{i}{n+1}}[\widehat{a}_n] \Delta t_n \rangle \right]\\
				& + \mathbb{E}\left[\frac{1}{M}\sum_{i=1}^{M} |P_{{U}^{\top}_n\widetilde{Y}_{n+1}^{i}}[ a_n] \Delta t_n+ P_{{U}_n}[b_n^{i} \Delta W_n^{i}] -\widehat{P}^{\alpha}_{\widehat{U}^{\top}_n\tildhatY{i}{n+1}}[\widehat{a}_n] \Delta t_n- P_{\widehat{U}_n}[\widehat{b}_n^{i} \Delta W_n^{i}] |^2\right],
			\end{aligned}
		\end{equation}
		where in the second to last line we used the property of independence of the Brownian increments with respect to the discretized solutions. Now, we give reasonable bounds on the crossed term and on the last quadratic term in \eqref{eq: intermediate error} in order to obtain a relation on $e_n$ to conclude via Gronwall's lemma. These computations will involve projectors that have been introducted in Section \ref{sec: notation} and whose differences have been analyzed in Section \ref{sec: bound proj}.
		
		Concerning the cross-term in the right-hand side of \eqref{eq: intermediate error}, using Cauchy-Schwarz inequality and Lipschitz property of the drift, i.e.\
		\begin{equation}\label{eq: lip drift}
			|a_n^i-\widehat{a}_n^i|\leq C_{\mathrm{Lip}}|X_n^i-\widehat{X}_n^i|, \quad \forall i=1,\dots, M,
		\end{equation}
		we decompose it into a sum of differences between projections:
		\begin{equation*}
			\begin{aligned}
				&\mathbb{E}\left[\frac{1}{M}\sum_{i=1}^{M} \langle X_{n}^{i}-\widehat{X}_{n}^{i},P_{{U}^{\top}_n\widetilde{Y}_{n+1}^{i}}[ a_n] \Delta t_n -\widehat{P}^{\alpha}_{\widehat{U}^{\top}_n\tildhatY{i}{n+1}}[\widehat{a}_n] \Delta t_n \rangle \right] \\
				\leq & e_n\sqrt{\mathbb{E}\left[\frac{1}{M}\sum_{i=1}^{M} \left|P_{{U}^{\top}_n\widetilde{Y}_{n+1}^{i}}[ a_n] \Delta t_n -\widehat{P}^{\alpha}_{\widehat{U}^{\top}_n\tildhatY{i}{n+1}}[\widehat{a}_n] \Delta t_n \right|^2 \right]}\\
				\leq & e_n\sqrt{\mathbb{E}\left[\frac{1}{M}\sum_{i=1}^{M} \left|P_{{U}^{\top}_n\widetilde{Y}_{n+1}^{i}}[ a_n] \Delta t_n - \widehat{P}^{\alpha}_{\widehat{U}^{\top}_n\tildhatY{i}{n+1}}[ a_n] \Delta t_n +\widehat{P}^{\alpha}_{\widehat{U}^{\top}_n\tildhatY{i}{n+1}}[ a_n] \Delta t_n  -\widehat{P}^{\alpha}_{\widehat{U}^{\top}_n\tildhatY{i}{n+1}}[\widehat{a}_n] \Delta t_n \right|^2 \right]}\\
				\leq & e_n \sqrt{\mathbb{E}\left[\frac{1}{M}\sum_{i=1}^{M} \left|P_{{U}^{\top}_n\widetilde{Y}_{n+1}^{i}}[ a_n] \Delta t_n -\widehat{P}^{\alpha}_{{U}^{\top}_n\widetilde{Y}_{n+1}^{i}}[ a_n] \Delta t_n + \widehat{P}^{\alpha}_{{U}^{\top}_n\widetilde{Y}_{n+1}^{i}}[ a_n] \Delta t_n  -\widehat{P}^{\alpha}_{\widehat{U}^{\top}_n\tildhatY{i}{n+1}}[a_n] \Delta t_n \right|^2 \right]}\\
				&+C_{\mathrm{Lip}}e_n^2\Delta t_n\\
				\leq & e_n \sqrt{\mathbb{E}\left[\frac{1}{M}\sum_{i=1}^{M} \left|P_{{U}_n}^{\perp}\left(P_{\widetilde{Y}_{n+1}^{i}} -\widehat{P}^{\alpha}_{\widetilde{Y}_{n+1}^{i}}\right)[ a_n] \Delta t_n \right|^2 \right]} \\
				&+e_n \sqrt{\mathbb{E}\left[\frac{1}{M}\sum_{i=1}^{M}\left| \left(P_{{U}_n}-P_{\widehat{U}_n}\right)\left[ \left(I_{d \times d} -  \widehat{P}^{\alpha}_{\widetilde{Y}_{n+1}^{i}}\right) a_n\right] \Delta t_n  +P_{\widehat{U}_n}^{\perp}\left(\widehat{P}^{\alpha}_{\widetilde{Y}_{n+1}^{i}}- \widehat{P}^{\alpha}_{\tildhatY{i}{n+1}}\right)[a_n] \Delta t_n \right|^2 \right]} \\
				&+ C_{\mathrm{Lip}}e_n^2\Delta t_n\\
				\leq & C_{\mathrm{Lip}}e_n^2\Delta t_n+e_n \sqrt{\mathbb{E}\left[\frac{1}{M}\sum_{i=1}^{M} \left|P_{{U}_n}^{\perp}\left(P_{\widetilde{Y}_{n+1}^{i}} -\widehat{P}^{\alpha}_{\widetilde{Y}_{n+1}^{i}}\right)[ a_n] \Delta t_n \right|^2 \right]} \\
				&+e_n \sqrt{\mathbb{E}\left[\frac{1}{M}\sum_{i=1}^{M}\left| \left(P_{{U}_n}-P_{\widehat{U}_n}\right)\left[ \left(I_{d \times d} -  \widehat{P}^{\alpha}_{\widetilde{Y}_{n+1}^{i}}\right) a_n\right] \Delta t_n \right|^2 \right]} \\
				&+e_n \sqrt{\mathbb{E}\left[ \frac{1}{M}\sum_{i=1}^{M}\left| \left(\widehat{P}^{\alpha}_{\widetilde{Y}_{n+1}^{i}}- \widehat{P}^{\alpha}_{\tildhatY{i}{n+1}}\right)[a_n] \Delta t_n \right|^2 \right]}.
\end{aligned}
\end{equation*}
We now use Lemmata \ref{lem: diff Proj U}, \ref{lem: error between stochastic proj}, and \ref{lem: discr proj Y} to bound the differences of projectors and simplify the resulting expression by grouping together terms of the same order in $e_n$ and $M$. Then, using the unitary norm of orthogonal projectors, linear-growth bound, Lemma \ref{lem: L2p norm semidiscretized solution} and monotonicity of the moments of a random variable, for all $R > \max\{ \| a_n\|_{L^2\left(\Omega^M, \mathbb{R}^{d\times M}\right)}, \| b_n\|_{L^2\left(\Omega^M, \mathbb{R}^{d \times m \times M}\right)} \}$ one obtains
\begin{equation*}
	\begin{aligned}
		&\mathbb{E}\left[\frac{1}{M}\sum_{i=1}^{M} \langle X_{n}^{i}-\widehat{X}_{n}^{i},P_{{U}^{\top}_n\widetilde{Y}_{n+1}^{i}}[ a_n] \Delta t_n -\widehat{P}^{\alpha}_{\widehat{U}^{\top}_n\tildhatY{i}{n+1}}[\widehat{a}_n] \Delta t_n \rangle \right] \\
		\leq & C_{\mathrm{Lip}}e_n^2\Delta t_n+  e_n \frac{1}{\sqrt[3]{M}}G_2 \|a_n\|_{L^2(\Omega^{M};\mathbb{R}^{d \times M})} \Delta t_n + e_n \left(e_n + \frac{1}{\sqrt{M}}\right) G_1 R \Delta t_n +  \|a_n\|_{L^4(\Omega^{M};\mathbb{R}^{d \times M})}   \sqrt[4]{C} e^{-\frac{cMR^2}{4}} \Delta t_n \\
		& +e_n \frac{1}{\sqrt[3]{M}} \|a_n\|_{L^4(\Omega^{M};\mathbb{R}^{d \times M})}  \Delta t_n +e_{n}\left(e_n + \frac{1}{\sqrt{M}}\right)   G_3 R \Delta t_n +  \|a_n\|_{L^4(\Omega^{M};\mathbb{R}^{d \times M})}   \sqrt[4]{C} e^{-\frac{cMR^2}{4}}  \Delta t_n \\
			\leq & C_{\mathrm{Lip}}e_n^2\Delta t_n+  e_n \frac{1}{\sqrt[3]{M}}G_2 \sqrt{C_{\mathrm{lgb}}}(1+\sqrt{K_2(T)}) \Delta t_n + e_n \left(e_n + \frac{1}{\sqrt{M}}\right) G_1 R \Delta t_n   \\
		&+ e_n \frac{1}{\sqrt[3]{M}}\sqrt{C_{\mathrm{lgb}}}\sqrt[4]{2}(1+\sqrt[4]{K_4(T)}) \Delta t_n +e_{n}\left(e_n + \frac{1}{\sqrt{M}}\right)   G_3 R \Delta t_n + 2\sqrt[4]{2}\sqrt{C_{\mathrm{lgb}}}(1+\sqrt[4]{K_4(T)}) \sqrt[4]{C} e^{-\frac{cMR^2}{4}}  \Delta t_n \\	
				\leq & \left(C_{\mathrm{Lip}}+ (G_1 + G_3) R \right)e_n^2\Delta t_n+  e_n \frac{1}{\sqrt[3]{M}}((G_2+\sqrt[4]{2}) \sqrt{C_{\mathrm{lgb}}}(1+\sqrt[4]{K_4(T)}) \Delta t_n \\
			&+ e_n \left( \frac{1}{\sqrt{M}}\right) (G_1  + G_3) R \Delta t_n + 2^{\frac54}\sqrt{C_{\mathrm{lgb}}}(1+\sqrt[4]{K_4(T)}) \sqrt[4]{C} e^{-\frac{cMR^2}{4}}  \Delta t_n.
			\end{aligned}
		\end{equation*}
	Note that via Young's inequality, one has that 
	\begin{equation}\label{eq: young ineq quan}
		\frac{1}{\sqrt[3]{M}} e_n \leq \frac{1}{2}\frac{1}{\sqrt[3]{M^2}}+ \frac{1}{2}e_n^2,  \text{ and } e_n \frac{1}{\sqrt{M}} \leq \frac{1}{2}e_n^2 + \frac{1}{2} \frac{1}{M}.
	\end{equation}
	Then, using relations \eqref{eq: young ineq quan} on cross terms and the fact that 
	$\frac{1}{M} \leq \frac{1}{\sqrt[3]{M^2}}$ for $M \geq 1$, there exists a positive constant $C:=C(\gamma,T)$ independent of $M$, $\Delta t$, but dependent on the inverse of $\gamma$, such that one has
	\begin{equation*}
		\begin{aligned}
			&\mathbb{E}\left[\frac{1}{M}\sum_{i=1}^{M} \langle X_{n}^{i}-\widehat{X}_{n}^{i},P_{{U}^{\top}_n\widetilde{Y}_{n+1}^{i}}[ a_n] \Delta t_n -\widehat{P}_{\widehat{U}^{\top}_n\tildhatY{i}{n+1}}[\widehat{a}_n] \Delta t_n \rangle \right] \\
			\leq & \left(C_{\mathrm{Lip}}+ (G_1+G_3)R \right)e_n^2\Delta t_n+  \frac{1}{2}(\frac{1}{\sqrt[3]{M^2}}+ e_n^2)((G_2 +\sqrt[4]{2})\sqrt{C_{\mathrm{lgb}}}(1+\sqrt[4]{K_4(T)})) \Delta t_n  \\
			&+ \left(\frac{1}{2}  \frac{1}{M} +  \frac{1}{2} e_n^2 \right) (G_1 + G_3) R \Delta t_n + 2^{\frac{5}{4}}\sqrt{C_{\mathrm{lgb}}}(1+\sqrt[4]{K_4(T)}) \sqrt[4]{C} e^{-\frac{cMR^2}{4}}  \Delta t_n\\
		\leq & e_n^2  C(\gamma,T) \left(1+ R\right)\Delta t_n + C(\gamma,T) \left( \frac{1}{\sqrt[3]{M^2}} + \frac{R}{M} +  e^{-MR^2}  \right)\Delta t_n.
	\end{aligned}
\end{equation*}

We now turn to the last term in the right-hand side in \eqref{eq: intermediate error}, which contains the quadratic increments.
\begin{equation*}
	\begin{aligned}
&\mathbb{E}\left[\frac{1}{M}\sum_{i=1}^{M} |P_{{U}^{\top}_n\widetilde{Y}_{n+1}^{i}}[ a_n] \Delta t_n+ P_{{U}_n}[b_n^{i} \Delta W_n^{i}] -\widehat{P}^{\alpha}_{\widehat{U}^{\top}_n\tildhatY{i}{n+1}}[\widehat{a}_n] \Delta t_n- P_{\widehat{U}_n}[\widehat{b}_n^{i} \Delta W_n^{i}] |^2\right]\\
\leq  & \mathbb{E}\Bigg[\frac{1}{M}\sum_{i=1}^{M} |P_{{U}^{\top}_n\widetilde{Y}_{n+1}^{i}}[ a_n] \Delta t_n-\widehat{P}^{\alpha}_{\widehat{U}^{\top}_n\tildhatY{i}{n+1}}[ a_n] \Delta t_n+ P_{{U}_n}[b_n^{i} \Delta W_n^{i}] - P_{\widehat{U}_n}[b_n^{i} \Delta W_n^{i}] \\
& + \widehat{P}^{\alpha}_{\widehat{U}^{\top}_n\tildhatY{i}{n+1}}[ a_n]  \Delta t_n  -\widehat{P}^{\alpha}_{\widehat{U}^{\top}_n\tildhatY{i}{n+1}}[\widehat{a}_n] \Delta t_n+ P_{\widehat{U}_n}[b_n^{i} \Delta W_n^{i}]- P_{\widehat{U}_n}[\widehat{b}_n^{i} \Delta W_n^{i}] |^2\Bigg]\\
\leq & 4C_{\mathrm{Lip}} e_n^2 \left((\Delta t_n)^2 + \Delta t_n \right) \\
& + 2 \mathbb{E}\left[\frac{1}{M}\sum_{i=1}^{M} |P_{{U}^{\top}_n\widetilde{Y}_{n+1}^{i}}[ a_n] \Delta t_n  -\widehat{P}^{\alpha}_{\widehat{U}^{\top}_n\tildhatY{i}{n+1}}[a_n] \Delta t_n+ P_{{U}_n}[b_n^{i} \Delta W_n^{i}]- P_{\widehat{U}_n}[b_n^{i} \Delta W_n^{i}] |^2\right].\\
\end{aligned}
\end{equation*}
Therefore, via similar computations employed for the second term in \eqref{eq: intermediate error}, with an abuse of notation on constants there exists a positive constant $C:=C(\gamma,T)$ independent of $M$, $\Delta t$, but dependent on the inverse of $\gamma$, such that one has 
\begin{equation*}
	\begin{aligned}
	&\mathbb{E}\left[\frac{1}{M}\sum_{i=1}^{M} |P_{{U}^{\top}_n\widetilde{Y}_{n+1}^{i}}[ a_n] \Delta t_n+ P_{{U}_n}[b_n^{i} \Delta W_n^{i}] -\widehat{P}^{\alpha}_{\widehat{U}^{\top}_n\tildhatY{i}{n+1}}[\widehat{a}_n] \Delta t_n- P_{\widehat{U}_n}[\widehat{b}_n^{i} \Delta W_n^{i}] |^2\right]\\
 \leq & e_n^2  C(\gamma,T) \left(1+ R\right)\Delta t_n + C(\gamma,T)  \left( \frac{1}{\sqrt[3]{M^2}} + \frac{R}{M} +  e^{-MR^2}  \right)\Delta t_n.
\end{aligned}
\end{equation*}

Putting all these intermediate computations together, from \eqref{eq: intermediate error} one finally obtains
	\begin{equation*}
	\begin{aligned}
		e^2_{n+1} \leq &  e_{n}^2 + e_n^2 2 C(\gamma,T) \left(1+ R\right) \Delta t_n +2 C(\gamma,T) \left( \frac{1}{\sqrt[3]{M^2}} + \frac{R}{M} +  e^{-MR^2}  \right)\Delta t_n,
	\end{aligned}
\end{equation*}
and via a discrete-type Gronwall's lemma \cite[Theorem 1.19]{sanz2023inverse} one has that for all $n$ and for all $R > \max\{ \| a_n\|_{L^2\left(\Omega^M, \mathbb{R}^{d\times M}\right)},$ $\| b_n\|_{L^2\left(\Omega^M, \mathbb{R}^{d \times m \times M}\right)} \}$ we have for a positive constant $C$ that 
\begin{equation}\label{eq: quasi final err}
	\begin{aligned}
			e^2_{n+1} \leq & \left( \frac{1}{\sqrt[3]{M^2}}  + \frac{R}{M} +  e^{-MR^2}   \right) \exp\{ C \left(1+ R\right) T \}.
	\end{aligned}
\end{equation}
In order to obtain the best possible rate of convergence with respect to $M$, we balance all the terms depending on $R$, by choosing $R$ as function of $M$.  Taking $R$ as $R(M)= \log \log M$ for $M$ large enough one has that
\begin{equation*}
	\begin{aligned}
		e^2_{n+1} \leq & \left( \frac{1}{\sqrt[3]{M^2}}  (\log M)^{CT}  \right) \exp\{ C T \} \\
		\leq & \left( M^{-\frac{2}{3}+ \frac{\log \log M}{\log M}CT} \right) \exp\{ C T \},
	\end{aligned}
\end{equation*}
which implies the thesis.
	\end{proof}
\end{Theorem}

\subsection{The case of adaptive regularization}\label{sec: adap reg}
The convergence rate established in Theorem \ref{thm: convergence of the Stoc Proj} for the stochastic error of Algorithm \ref{alg: Stoc Proj algorithm} is unfortunately lower than the standard Monte Carlo rate $\frac{1}{2}$. This is essentially due to the regularization on the empirical Gramian, whose parameter scales as $\alpha \propto \frac{1}{\sqrt[3]{M}}$. On the other hand, Proposition \ref{prop: lower-bound Gramian} establishes a lower bound on the probability that the smaller singular value $\widehat{\sigma}^k_{n,M}$ remains larger than a fixed positive constant. Such lower bound is small and controlled by $M$. We do expect, therefore, that in most of the time steps the smallest singular value is sufficiently far from zero and the regularization of the Gramian is not needed. 

In this section, we explore a modification of Algorithm \ref{alg: Stoc Proj algorithm}, detailed in Algorithm \ref{alg: mod Stoc Proj algorithm} below, in which the regularization is activated only when $\widehat{\sigma}^k_{n,M}$ falls below a prescribed value $\eta$. We will show in Theorem \ref{thm: convergence of the mod Stoc Proj} that this strategy allows us to recover the standard Monte Carlo rate, up to logarithmic terms, when the threshold $\eta$ is carefully chosen.

\begin{algorithm}
	\caption{Modified Monte Carlo DLR Projector Splitting for SDEs}\label{alg: mod Stoc Proj algorithm}
	\begin{flushleft}
		\textbf{Input}: initial data $U_0$, $Y_0$, number of samples $M$, positive quantity $\eta>0$.
		
		\textbf{Output:} approximation $\{\widehat{\mathbb{X}}_n=\widehat{U}_n^{\top}\widehat{\mathbb{Y}}_n\}_{n=0,\ldots, N}$. 
	\end{flushleft}
	\begin{algorithmic}[1]
		
		\ForAll {$n \in \{0, \ldots N-1\}$} 
		
		\State Generate $M$ Brownian increments $\Delta W_n^{i} \overset{\text{i.i.d.}}{\sim} \mathcal{N}(0,\Delta t_n)$ with $i \in \{1,\dots,M\}.$
		
		\State Compute $ \tildhatY{i}{n+1} = \widehat{Y}_n^{i} + \widehat{U}_n a(t_n,\widehat{U}_{n}^{\top}\widehat{Y}_n^{i}) \Delta t_n +  \widehat{U}_n b(t_n,\widehat{U}_{n}^{\top}\widehat{Y}_n^{i}) \Delta W_n^{i}$, for all $i \in \{1,\dots,M\}.$
		
		\State If $\sigma^k(\widehat{C}_{\widetilde{\widehat{Y}}_{n+1}}) > \eta$, then set
		
		$\qquad \widehat{C}_{\widetilde{\widehat{Y}}_{n+1}}^{\ast} = \widehat{\mathbb{E}}[\widetilde{\widehat{\mathbb{Y}}}_{n+1}(\widetilde{\widehat{\mathbb{Y}}}_{n+1})^{\top}]$
		
		else
		
		$\qquad \widehat{C}_{\widetilde{\widehat{Y}}_{n+1}}^{\ast}  = \widehat{\mathbb{E}}[\widetilde{\widehat{\mathbb{Y}}}_{n+1}(\widetilde{\widehat{\mathbb{Y}}}_{n+1})^{\top}]+ \frac{1}{\sqrt[3]{M}} I_{k \times k}$
		
		\State Compute $\widetilde{\widehat{U}}_{n+1}$: 
		$$\widehat{C}_{\widetilde{\widehat{Y}}_{n+1}}^{\ast} \widetilde{\widehat{U}}_{n+1} =  \widehat{C}_{\widetilde{\widehat{Y}}_{n+1}}^{\ast}  \widehat{U}_n + \widehat{\mathbb{E}}\left[\widetilde{\widehat{\mathbb{Y}}}_{n+1} \left(a(t_n,\widehat{U}_{n}^{\top}\widehat{\mathbb{Y}}_{n})^{\top}\Delta t_n\right)\right]\left(I_{d \times d} - P^{\text{row} }_{\widehat{U}_n} \right) $$
		\State Reorthonormalize the deterministic basis: find   $(\widehat{U}_{n+1}, \widehat{\mathbb{Y}}_{n+1})$ such that:
		\begin{equation*}
			\widehat{U}_{n+1}^{\top} \widehat{\mathbb{Y}}_{n+1} = \tildhatUT{n+1} \widetilde{\widehat{\mathbb{Y}}}_{n+1}, \quad \widehat{U}_{n+1}\widehat{U}_{n+1}^{\top} = I_{k\times k}.
		\end{equation*}
		For example, $(\widehat{U}_{n+1}^{\top}, \widehat{R}_{n+1}) = \texttt{QR}(\tildhatUT{n+1})$ and $\widehat{Y}_{n+1}^{i} = \widehat{R}_{n+1} \tildhatY{i}{n+1}$, for all $i \in \{1,\dots,M\}.$
		\EndFor
	\end{algorithmic}
\end{algorithm}

To analyze the convergence of Algorithm \ref{alg: mod Stoc Proj algorithm}, we introduce the un-regularized orthogonal empirical projectors $\widehat{P}_{\widetilde{Y}_{n+1}^{i}}: L^2(\Omega^{M};\mathbb{R}^{d \times M}) \to L^2(\Omega^M,\mathbb{R}^{d})$ as
\begin{equation}\label{eq: Y discretized projectors - full rank}
	\begin{aligned}
		\widehat{P}_{\widetilde{Y}_{n+1}^{i}}[f] =	& (\widetilde{Y}_{n+1}^{i})^{\top}\left(\frac{\sum_{j=1}^{M} \widetilde{Y}_{n+1}^j (\widetilde{Y}_{n+1}^j)^{\top}}{M}\right)^{-1}\widehat{\mathbb{E}}\left[ \widetilde{\mathbb{Y}}_{n+1} f^{\top} \right] =  (\widetilde{Y}_{n+1}^{i})^{\top}\widehat{C}_{\widetilde{Y}_{n+1}}^{-1}\widehat{\mathbb{E}}\left[ \widetilde{\mathbb{Y}}_{n+1} f^{\top} \right],\quad f\in L^2(\Omega^{M};\mathbb{R}^{d \times M}),\\
	\end{aligned}
\end{equation}
and $\widehat{P}_{\widetilde{\mathbb{Y}}_{n+1}}: L^2(\Omega^{M};\mathbb{R}^{d \times M}) \to  L^2(\Omega^{M};\mathbb{R}^{d \times M})$ as $\widehat{P}_{\widetilde{\mathbb{Y}}_{n+1}}:=[\widehat{P}_{\widetilde{Y}_{n+1}^{i}}[f]]$, for $f \in L^2(\Omega^{M};\mathbb{R}^{d \times M})$, which have already been used in Appendix \ref{app: proof proj} (see \eqref{eq: bar proj Y}) to prove Lemma \ref{lem: discr proj Y}. Notice that $\widehat{P}_{\widetilde{Y}_{n+1}^{i}}[f] = \widehat{P}^{\alpha}_{\widetilde{Y}_{n+1}^{i}}[f]$ and $\widehat{P}_{\tildhatY{i}{n+1}}[f] = \widehat{P}^{\alpha}_{\tildhatY{i}{n+1}}[f]$ for $\alpha=0$.
Then, the following estimate 
concerning the difference between non-degenerate projectors is useful.

\begin{Lemma}[Error between the stochastic projectors - full-rank case]\label{lem: error between stochastic proj - full rank}

Let us recall the definition of the semi-discrete Gramians $C_{\widetilde{Y}_{n+1}}=\mathbb{E}\left[\frac{1}{M} \widetilde{\mathbb{Y}}_{n+1}\widetilde{\mathbb{Y}}_{n+1}^{\top}\right]$, $C_{\widetilde{\widehat{Y}}_{n+1}}=\mathbb{E}\left[ \frac{1}{M} \widetilde{\widehat{\mathbb{Y}}}_{n+1}\widetilde{\widehat{\mathbb{Y}}}_{n+1}^{\top}\right]$, and empirical Gramians $\widehat{C}_{\widetilde{Y}_{n+1}}=\frac{\sum_{j=1}^{M} \widetilde{Y}_{n+1}^j (\widetilde{Y}_{n+1}^j)^{\top}}{M}$, $\widehat{C}_{\widetilde{\widehat{Y}}_{n+1}}=\frac{\sum_{j=1}^{M} \tildhatY{j}{n+1} (\tildhatY{j}{n+1})^{\top}}{M}$. Take $\gamma>0$ and define 
\begin{equation}\label{eq: sigma_gamma}
	\Sigma_{n}^{\gamma} = \{ \boldsymbol{\omega} \in \Omega^M \ : \ C_{\widehat{Y}_n}, C_{\widehat{\widehat{Y}}_n}, \widehat{C}_{\widehat{Y}_{n}}, \widehat{C}_{\widehat{\widehat{Y}}_{n}} \succeq \gamma I_{k \times k}\}.
\end{equation}
	
	If $\mathbb{E}[|Y_0|^{4}] < \infty$ and relation \eqref{eq: M dt cond} holds with $p=2$, then there exist constants $G_4:=G_4(\gamma)=O(\frac{1}{\gamma})$ and $G_5:=G_5(\gamma)=O(\frac{1}{\gamma})$ dependent on the inverse of $\gamma$ such that
	\begin{equation}\label{eq: Y projectors error - full rank}
		\sqrt{\mathbb{E}\left[\frac{1}{M}\sum_{i=1}^{M}
		\left|\left(P_{\widetilde{Y}^{i}_{n+1}}-\widehat{P}_{\widetilde{Y}^{i}_{n+1}}\right)f \right|^2 \mathbbm{1}_{\Sigma_{n+1}^{\gamma}}\right]} \leq \frac{1}{\sqrt{M}} G_4 \|f\|_{L^4(\Omega^{M};\mathbb{R}^{d \times M})}
	\end{equation} 
	for all $f \in L^4_{iid}\left(\Omega^{M}, \mathbb{R}^{d \times M}\right)$. Moreover, for all $R>\|f\|_{L^2(\Omega^M;\mathbb{R}^{d \times M})}$ one has for positive constants $c,C$ that
		\begin{equation}
		\begin{aligned}\label{eq: Y projectors error - full rank 2}
			\sqrt{\mathbb{E}\left[\frac{1}{M}\sum_{i=1}^{M}
			\left|\left(\widehat{P}_{\tildhatY{i}{n+1}}-\widehat{P}_{\widetilde{Y}_{n+1}^{i}}\right)f \right|^2 \mathbbm{1}_{\Sigma_{n+1}^{\gamma}}\right]} \leq \left(e_n + \frac{1}{\sqrt{M}}\right) G_5 R +  \left( C e^{-cMR^2} \right) \|f\|_{L^2(\Omega^M;\mathbb{R}^{d \times M})},
	\end{aligned}
	\end{equation} 
	for all $f \in L^2_{iid}\left(\Omega^{M}, \mathbb{R}^{d \times M}\right)$ such that for positive constants $\tilde{C},\tilde{c}$, for all $R>0$ one has $\mathbb{P}(|f^{i}| > R) \leq \tilde{C} e^{-\tilde{c}R^2} $. 
	\begin{proof}
		The proof is deferred to Appendix \ref{app: proof adap}.
	\end{proof}
\end{Lemma}

Lemma \ref{lem: error between stochastic proj - full rank} should be compared with Lemma \ref{lem: discr proj Y}. It shows that if we restrict to the set $\Sigma_{n+1}^{\gamma}$ where all the covariance matrices are uniformly lower bounded, and do not apply any regularization to the empirical covariances, the estimate \eqref{eq: Y projectors error - full rank 2} does not contain the slow converging term $O(\frac{1}{\sqrt[3]{M}})$ which appears, instead, in \eqref{eq: discr proj Y}. Furthermore, notice that, unlike Lemma \ref{lem: discr proj Y}, here we only require $\mathbb{E}[|Y_0|^{4}] < \infty$ as we  evaluate our difference in a set where all the empirical Gramians have a positive lower bound. Thanks to this improved estimate, we are able to show in the next Theorem \ref{thm: convergence of the mod Stoc Proj} for Algorithm \ref{alg: mod Stoc Proj algorithm} a better convergence rate that what obtained in Theorem \ref{thm: convergence of the Stoc Proj} for Algorithm \ref{alg: Stoc Proj algorithm}. To prove this estimate, we still need the condition $\mathbb{E}[|Y_0|^{8}] < \infty$, as it is required to bound the probability that the Gramian is below a certain positive threshold.

\begin{Theorem}[Stochastic Convergence of the Modified DLR Projector Splitting for SDEs]\label{thm: convergence of the mod Stoc Proj}
	Suppose that $\mathbb{E}[|Y_0|^{8}] < \infty$, that $\sigma_{Y_0}=\sigma^k\left( \mathbb{E}[ Y_0 Y_0^{\top}] \right)>0$, that Assumptions \ref{ass: diff}-\ref{ass: bound diff and tails} hold and assume that $\Delta t_n$ satisfies \eqref{eq: dt semidiscrete condition} and \eqref{eq: M dt cond} with $p=4$.  
	Let $\gamma = \min \{\frac{\varrho}{2}, \frac{\sigma_{Y_0}}{4}\}$, where $\varrho$ is defined in \eqref{eq: rho}. 
	If the constant $\eta$ in Algorithm \ref{alg: mod Stoc Proj algorithm} satisfies $\eta < \gamma$, then for a positive constant $C:=C(\eta, \gamma)=O(\max\{\frac{1}{\eta};\frac{1}{\gamma}\})$ dependent on the inverses of $\eta$ and $\gamma$ we have 
	\begin{equation}\label{eq: error for the mod Stoc Proj}
		e_{n+1} \leq  \left( M^{-\frac{1}{2}+ \frac{\log \log M}{2 \log M} CT}\right) \exp \left(C  T\right), \quad \text{ for all } \ 0 \leq n \leq N-1.
	\end{equation}
	
	\begin{proof}
	 The proof is similar to the one of Theorem \ref{thm: convergence of the Stoc Proj} and we only detail the bounds for the terms that are different.
		
		To exploit the regularization strategy, we decompose the error according to whether the Gramian remains well-conditioned or not. In the former case, we retrieve a Monte-Carlo-type error of order $O(\frac{1}{\sqrt{M}})$, whereas in the latter we want to bound everything by the probability of the bad event times the regularization rate of $\frac{1}{\sqrt[3]{M}}$.
		For all $n$, consider
		$A_{n}$ defined as in \eqref{eq: def A_n}. Notice that Lemma \ref{lem: Lp moment bound}, Proposition \ref{prop: lower-bound Gramian}, as well as all of the bounds in Lemmata \ref{lem: prop bound sigma tilde}-\ref{lem: p-mom cond covariance} in Appendix \ref{app: Gramian bounds} do apply to the solution produced by Algorithm \ref{alg: mod Stoc Proj algorithm} as well. From Proposition \ref{prop: lower-bound Gramian}, we have that $\widetilde{\widehat{\sigma}}_j^{k} \geq \min \{ \frac{\varrho}{2}, \frac{\sigma_{Y_0}}{4}\} > \eta$ for all $j=1,\dots,n$, i.e.\ $A_n$ contains only paths for which no regularization has been applied (at least until $t_n$).
		Then, one has 
		\begin{equation*}
			\begin{aligned}
				e^2_{n+1} = &  \mathbb{E}\left[\frac{1}{M}\sum_{i=1}^{M} |X_{n+1}^{i}-\widehat{X}_{n+1}^{i}|^2\right]\\
				=& e_{n}^2 +
					2\mathbb{E}\left[\frac{1}{M}\sum_{i=1}^{M} \langle X_{n}^{i}-\widehat{X}_{n}^{i},P_{{U}^{\top}_n\widetilde{Y}_{n+1}^{i}}[ a_n] \Delta t_n -\widehat{P}_{\widehat{U}^{\top}_n\tildhatY{i}{n+1}}[\widehat{a}_n] \Delta t_n \rangle \mathbbm{1}_{A_{n+1}} \right]\
				\\
				&+ 2\mathbb{E}\left[\frac{1}{M}\sum_{i=1}^{M} \langle X_{n}^{i}-\widehat{X}_{n}^{i},P_{{U}^{\top}_n\widetilde{Y}_{n+1}^{i}}[ a_n] \Delta t_n -\widehat{P}^{\alpha}_{\widehat{U}^{\top}_n\tildhatY{i}{n+1}}[\widehat{a}_n] \Delta t_n \rangle \mathbbm{1}_{A_{n+1}^{C}} \right]  \\
				& + \mathbb{E}\left[\frac{1}{M}\sum_{i=1}^{M} |P_{{U}^{\top}_n\widetilde{Y}_{n+1}^{i}}[ a_n] \Delta t_n+ P_{U_n}[b_n^{i} \Delta W_n^{i}] -\widehat{P}_{\widehat{U}^{\top}_n\tildhatY{i}{n+1}}[\widehat{a}_n] \Delta t_n- P_{\widehat{U}_n}[\widehat{b}_n^{i} \Delta W_n^{i}] |^2 \mathbbm{1}_{A_{n+1}}\right] \\
				& +  \mathbb{E}\left[\frac{1}{M}\sum_{i=1}^{M} |P_{{U}^{\top}_n\widetilde{Y}_{n+1}^{i}}[ a_n] \Delta t_n+ P_{U_n}[b_n^{i} \Delta W_n^{i}] -\widehat{P}^{\alpha}_{\widehat{U}^{\top}_n\tildhatY{i}{n+1}}[\widehat{a}_n] \Delta t_n- P_{\widehat{U}_n}[\widehat{b}_n^{i} \Delta W_n^{i}] |^2\mathbbm{1}_{A_{n+1}^{C}}\right],
			\end{aligned}
		\end{equation*}
	We now estimate each contribution separately.
	
		We first consider the contribution on the set $A_{n+1}$, where the Gramian is uniformly non-degenerate. Using Lemma \ref{lem: error between stochastic proj - full rank} with $\gamma = \min \{ \frac{\varrho}{2}, \frac{\sigma_{Y_0}}{4}\}$ and observing that $A_{n+1} \subset \Sigma_{n+1}^{\gamma}$, where $\Sigma_{n+1}^{\gamma}$ is defined in \eqref{eq: sigma_gamma}, 
			\begin{equation*}
			\begin{aligned}
				\mathbb{E}\left[\frac{1}{M}\sum_{i=1}^{M}  \left| \left(P_{\widetilde{Y}_{n+1}^{i}}-\widehat{P}_{\widetilde{Y}_{n+1}^{i}}\right)[ a_n]  \right|^2 \mathbbm{1}_{A_{n+1}} \right]^{\frac{1}{2}} 
				\leq  &\mathbb{E}\left[\frac{1}{M}\sum_{i=1}^{M}  \left| \left(P_{\widetilde{Y}_{n+1}^{i}}-\widehat{P}_{\widetilde{Y}_{n+1}^{i}}\right)[ a_n]\right|^2 \mathbbm{1}_{\Sigma_{n+1}^{\gamma}} \right]^{\frac{1}{2}} \\
				   \leq & G_4 \frac{1}{\sqrt{M}} \sqrt[4]{K_4(T)} \|a_n \|_{L^4(\Omega^M;\mathbb{R}^{d \times M})}.\\
			\end{aligned}
		\end{equation*} 
		On the other hand, using the same proof of Proposition \ref{prop: lower-bound Gramian} one can prove that $\mathbb{P}(A_{n+1}^{C}) \leq C_p 1/M^{p}$, where $C_p:=C_p(\varrho)$ is a constant dependent on $\frac{1}{\varrho^{2p}}$. Then, choosing $\alpha = \frac{1}{\sqrt[3]{M}}$ and proceeding similarly to the proof of Lemma \ref{lem: error between stochastic proj}, we bound the difference between the projector and the regularized empirical one applied to the semidiscretized basis in order to retrieve a rate with respect to the number of samples of order $\frac{1}{\sqrt[3]{M}}$ times the probability of the bad event.
		Indeed, by decomposing the contribution of 
				\(P_{\widetilde{Y}_{n+1}^{i}}-\widehat P_{\widetilde{Y}_{n+1}^{i}}^\alpha
				\) on the event $A_{n+1}^C$
				into the deterministic regularization error, the empirical Gramian error, and the empirical moment error, we have
			\begin{equation*}
			\begin{aligned}
			&	\|\left(P_{\widetilde{Y}_{n+1}^{i}}-\widehat{P}^{\alpha}_{\widetilde{Y}_{n+1}^{i}}\right)[ a_n] \Delta t\mathbbm{1}_{A_{n+1}^{C}}\|_{L^2(\Omega^{M};\mathbb{R}^{d \times M})} \\
				\leq& \left(\mathbb{E}\left[\frac{1}{M} \sum_{i=1}^{M} \left|(\widetilde{Y}_{n+1}^{i})^{\top}\left(\left(\mathbb{E}\left[ \widetilde{Y}_{n+1}\widetilde{Y}_{n+1}^{\top}\right]\right)^{-1}-\left(\mathbb{E}\left[ \widetilde{Y}_{n+1}\widetilde{Y}_{n+1}^{\top}+\alpha I_{k \times k}\right]\right)^{-1} \right)\mathbb{E}\left[ \widetilde{Y}_{n+1} a_n^{\top} \Delta t_n\right]  \mathbbm{1}_{A_{n+1}^{C}} \right|^2\right]\right)^{\frac{1}{2}}\\
			&\resizebox{1\linewidth}{!}{$+ \left(\mathbb{E}\left[\frac{1}{M} \sum_{i=1}^{M} \left|(\widetilde{Y}_{n+1}^{i})^{\top}\left(\mathbb{E}\left[ \widetilde{Y}_{n+1}\widetilde{Y}_{n+1}^{\top}+\alpha I_{k \times k}\right]^{-1}-\left(\frac{\sum_{j=1}^{M} \widetilde{Y}_{n+1}^j (\widetilde{Y}_{n+1}^j)^{\top}}{M} + \alpha I_{k \times k}\right)^{-1}\right)\mathbb{E}\left[ \widetilde{Y}_{n+1} a_n^{\top} \Delta t_n \right]  \mathbbm{1}_{A_{n+1}^{C}} \right|^2\right]\right)^{\frac{1}{2}}   $} \\
			&+ \left(\mathbb{E}\left[\frac{1}{M} \sum_{i=1}^{M} \left|(\widetilde{Y}_{n+1}^{i})^{\top}\left(\frac{\sum_{j=1}^{M} \widetilde{Y}_{n+1}^j (\widetilde{Y}_{n+1}^j)^{\top}}{M} + \alpha I_{k \times k}\right)^{-1}\left(\mathbb{E}\left[ \widetilde{Y}_{n+1} a_n^{\top} \Delta t\right] -\widehat{\mathbb{E}}\left[ \widetilde{Y}_{n+1} a_n^{\top} \Delta t_n\right]  \right)\right|^2  \mathbbm{1}_{A_{n+1}^{C}} \right]\right)^{\frac{1}{2}}.\\
				\leq &  \frac{k}{\gamma} \mathbb{E}\left(\left[\left|\mathbb{E}\left[ \widetilde{Y}_{n+1}\widetilde{Y}_{n+1}^{\top} \right]-\mathbb{E}\left[ \widetilde{Y}_{n+1}\widetilde{Y}_{n+1}^{\top}+\frac{1}{\sqrt[3]{M}}  I_{k \times k}  \right]\right|^2 \mathbbm{1}_{A_{n+1}^{C}} \right]\right)^{\frac{1}{2}} \|a_n\|_{L^2(\Omega^{M};\mathbb{R}^{d \times M})} \cdot \Delta t_n\\
				&+  \frac{ k}{\sqrt{\gamma}}  \sqrt[6]{M}  \Big(\mathbb{E} \Big[ \Big|\frac{1}{M} \sum_{j=1}^{M} \frac{\widetilde{Y}_{n+1}^j (\widetilde{Y}_{n+1}^j)^{\top}}{M}- \mathbb{E}\Big[ \widetilde{Y}_{n+1}\widetilde{Y}_{n+1}^{\top} \Big]\Big|^2 \mathbbm{1}_{A_{n+1}^{C}}\Big]\Big)^{\frac{1}{2}} \|a_n\|_{L^2(\Omega^{M};\mathbb{R}^{d \times M})}\Delta t_n \\
				& + \sqrt{k}\sqrt[6]{M} \left(\mathbb{E}\left[ \left|\widehat{\mathbb{E}}\left[ \widetilde{Y}_{n+1} (a_n \Delta t_n)^{\top}  \right]-\mathbb{E}\left[ \widetilde{Y}_{n+1} (a_n \Delta t_n)^{\top}  \right] \right|^2  \mathbbm{1}_{A_{n+1}^{C}} \right]\right)^{\frac{1}{2}}.\\
			\end{aligned}
		\end{equation*}
		
		Now, notice that $\mathbb{R}^{k\times k}$ with the Frobenium norm is a Hilbert space, i.e.\ it is of Rademacher type $2$ \cite[Definition 1.1]{woyczynski1980marcinkiewicz}. Then, the Monte Carlo estimator of size $M$ applied to the random variable $\widetilde{Y}_{n+1}(\omega) (\widetilde{Y}_{n+1}(\omega))^{\top}: \Omega \to \mathbb{R}^{k \times k}$ converges in $L^4(\Omega, \mathbb{R}^{k \times k})$ to its average $\mathbb{E}\Big[ \widetilde{Y}_{n+1}\widetilde{Y}_{n+1}^{\top}\Big]$ with a rate of $\frac{1}{\sqrt{M}}$ \cite[Corollary 2.1]{woyczynski1980marcinkiewicz}. The same holds for the random variable $\widetilde{Y}_{n+1}(\omega) (a_{n}(\omega))^{\top}: \Omega \to \mathbb{R}^{k \times k}$. Via this fact and using Hölder inequality with exponents $\frac{1}{2}$, we obtain that there exists a constant $c_4$ independent on $\gamma$ and $d$ such that it holds
		\begin{equation*}
			\begin{aligned}
				&	\|\left(P_{\widetilde{Y}_{n+1}^{i}}-\widehat{P}^{\alpha}_{\widetilde{Y}_{n+1}^{i}}\right)[ a_n] \Delta t_n\mathbbm{1}_{A_{n+1}^{C}}\|_{L^2(\Omega^{M};\mathbb{R}^{d \times M})} \\
				\leq &  \frac{k}{\gamma}  \frac{1}{\sqrt[3]{M}} \mathbb{P}(A_{n+1}^{C})^{\frac{1}{2}}  \|a_n\|_{L^2(\Omega^{M};\mathbb{R}^{d \times M})}\Delta t_n \\
				&+  \frac{ k}{\sqrt{\gamma}} \sqrt[6]{M} \Big(\mathbb{E} \Big[ \Big|\frac{\sum_{j=1}^{M} \widetilde{Y}_{n+1}^j (\widetilde{Y}_{n+1}^j)^{\top}}{M}  - \mathbb{E}\Big[ \widetilde{Y}_{n+1}\widetilde{Y}_{n+1}^{\top} \Big]\Big|^4\Big]\Big)^{\frac{1}{4}} \mathbb{P}(A_{n+1}^{C})^{\frac{1}{4}}  \|a_n\|_{L^2(\Omega^{M};\mathbb{R}^{d \times M})}\Delta t_n \\
				& +\sqrt{k} \sqrt[6]{M} \left(\mathbb{E}\left[ \left|\widehat{\mathbb{E}}\left[ \widetilde{Y}_{n+1} (a_n \Delta t_n)^{\top} \right]-\mathbb{E}\left[ \widetilde{Y}_{n+1} (a_n \Delta t_n)^{\top}  \right] \right|^4   \right]\right)^{\frac{1}{4}} \mathbb{P}(A_{n+1}^{C})^{\frac{1}{4}} \\
				\leq &  \frac{k}{\gamma}  \frac{1}{\sqrt[3]{M}} \mathbb{P}(A_{n+1}^{C})^{\frac{1}{2}} \|a_n\|_{L^2(A_{n+1}^{C};\mathbb{R}^{d \times M})}\Delta t_n +  \frac{k}{\sqrt{\gamma}}  \frac{c_4}{\sqrt[3]{M}} \mathbb{P}(A_{n+1}^{C})^{\frac{1}{4}} \sqrt[4]{K_4(T)}   \|a_n\|_{L^2(\Omega^{M};\mathbb{R}^{d \times M})}\Delta t_n  \\
				& +\frac{k c_4}{\sqrt[3]{M}}  \mathbb{E}[\|\widetilde{Y}_{n+1} (a_n)^{\top}\|_{\mathrm{F}}^4]^{\frac{1}{4}} \mathbb{P}(A_{n+1}^{C})^{\frac{1}{4}} \Delta t_n \\
				\leq &\frac{k}{\sqrt[3]{M}} \left( \frac{1}{\gamma} \mathbb{P}(A_{n+1}^{C})^{\frac{1}{2}} + c_4 \sqrt[4]{K_4(T)} \frac{1}{\sqrt{\gamma}} \mathbb{P}(A_{n+1}^{C})^{\frac{1}{4}} + c_4 \sqrt[8]{K_8(T)} \mathbb{P}(A_{n+1}^{C})^{\frac{1}{4}} \right) \|a_n\|_{L^8(\Omega^{M};\mathbb{R}^{d \times M})}\Delta t_n \\
				 \leq& \frac{k}{\sqrt[3]{M}} \mathbb{P}(A_{n+1}^{C})^{\frac{1}{4}} \left( \frac{1}{\gamma} + c_4 \sqrt[4]{K_4(T)} \frac{1}{\sqrt{\gamma}} +  c_4 \sqrt[8]{K_8(T)}  \right) \|a_n\|_{L^8(\Omega^{M};\mathbb{R}^{d \times M})}\Delta t_n,
			\end{aligned}
		\end{equation*}
		 where in the second inequality we exploit the convergence of the Monte Carlo estimator, and in the last line we use the fact that $L^p$ norms are monotone with respect to the exponent $p\geq1$, 
		and $\mathbb{P}(A_{n+1}^{C}) \leq 1$.
		  
		One applies a similar treatment to similar differences arisen in the bound of the error. Then, there exists a positive constant $C_1:=C_1(T,\gamma)$ dependent on the inverse of $\gamma$, but
		independent on $M$ and $\Delta t$, such that
			\begin{equation*}
			\begin{aligned}
				e^2_{n+1} \leq &   e_n^2  \left(1 +2 C_1 (1+R) \Delta t_n  \right) +2 C_1 \left(\frac{R}{M} +e^{-MR^2}   + \frac{1}{\sqrt[3]{M^2}} \mathbb{P}(A_{n+1}^{C})^{\frac{1}{4}}\right)\Delta t_n .
			\end{aligned}
		\end{equation*}
		Using Propositions \ref{prop: lowerbound exp} and \ref{prop: lower-bound Gramian}, and Hölder's inequality, for positive constants $C_2:=C_2(T, \gamma)$ and $\tilde{C}_5:=\tilde{C}_5(T,p, \gamma)$ independent of $M$, $\Delta t_n$, but dependent on the inverse of $\gamma$, one gets
			\begin{equation*}
			\begin{aligned}
				e^2_{n+1} \leq &   e_n^2  \left(1 +2 C_2 (1+R) \Delta t_n \right) +2 C_2  \left(\frac{R}{M} +e^{-MR^2}  + \frac{1}{\sqrt[3]{M^2}} \left( \frac{\sqrt{\tilde{C}_5}}{\sqrt{M}} \right)\right)\Delta t_n. \\
			\end{aligned}
		\end{equation*}
		Via Gronwall's lemma for a positive constant $C:=C(\gamma,\eta)$ independent of $M$ and $\Delta t_n$ it holds that
	   	\begin{equation*}
	   	\begin{aligned}
	   e^2_{n+1} \leq & \left( \frac{1}{M^{\frac{7}{6}}}  + \frac{R}{M} +  e^{-MR^2}   \right) \exp\{ C \left(1+ R\right) T \}, \quad \text{ for all } n \ : \ 1 \leq n \leq N.
	   	\end{aligned}
	   \end{equation*}
	   Similarly to the proof of Theorem \ref{thm: convergence of the Stoc Proj}, choosing $R = \log \log M$ for $M$ large we obtain the thesis.
	   \end{proof}
\end{Theorem}
 
We finally estimate the probability that no regularization is employed over the full time-trajectory discretization.
\begin{Corollary}
	Assume that the same hypotheses of Theorem \ref{thm: convergence of the mod Stoc Proj} hold, but for an initial condition $\mathbb{E}[|Y_0|^{4p}] < \infty$, $p\geq 1$. Furthermore, assume a uniform time step $\Delta t$. Then there exist an event
	$D_{N} \in \mathcal{F}$ and a positive constant $C$ dependent on the inverse of $\gamma = \min \{\frac{\varrho}{2}, \frac{\sigma_{Y_0}}{4}\}$ such that 
	\begin{equation*}
		\mathbb{P}(D_N)
		\geq 1 -  \tilde{C}\frac{1}{\Delta t M^{p}}\\
	\end{equation*}
	and 
	\begin{equation}\label{eq: error for the mod Stoc Proj - trajectories}
		\mathbb{E}\left[\frac{1}{M}\sum_{i=1}^{M} |X_{n+1}^{i}-\widehat{X}_{n+1}^{i}|^2 \mathbbm{1}_{D_N}\right] \leq  \left( M^{-\frac{1}{2}+ \frac{\log \log M}{2 \log M} CT}\right)   \exp \left(C  T\right), \quad \text{ for all } \ 0 \leq n \leq N-1.
	\end{equation}
	\begin{proof}
		Let us define the set $D_N$ as follows
		\begin{equation*}
			D_N := E_N \cap  \Bigg\{ \bigcap_{0 \leq n \leq N-1}  \Big\{ \widetilde{\sigma}^{k}_{n+1,M} > \frac{\varrho}{2} \Big\}\Bigg\},
		\end{equation*}
		where $E_N$ is defined using relation \eqref{eq: E_n} for $n=N$. By definition of $D_N$, if $\boldsymbol{\omega} \in D_N$, then all the empirical Gramians involved in the analysis of  Algorithm \ref{alg: mod Stoc Proj algorithm} are lower bounded by a positive constant (see also Lemma \ref{lem: prop bound sigma tilde} and Proposition \ref{prop: lowerbound exp}). Therefore, there is no need to employ a regularization of the Gramians inherent to the stochastic basis for $\omega \in D_N$, allowing us to retrieve directly the same convergence rate of Theorem \ref{thm: convergence of the mod Stoc Proj} closed to Monte Carlo one if the expectation is restricted to the set $D_N$. Finally, Lemma \ref{lem: prop bound sigma tilde} and Proposition \ref{prop: lower-bound Gramian} give us a lower bound on the probability that $D_N$ happens. 
	\end{proof}
\end{Corollary}

\subsection{The case of degenerate diffusion}\label{sec: diff reg}
In Sections \ref{sec: bound proj} and \ref{sec: error}, as well in Lemma \ref{lem: discr proj Y}, Assumption \ref{ass: diff} is essential to retrieve convergence. However, in several cases this hypothesis can lead to systems that do not present any low dimensional structure, and, therefore, the employment of DLRA could lead to a well-posed, but inaccurate, approximation.

In order to deal with a possible degenerate diffusion $b$, we propose a modification of the diffusion $b(t, x)$ in the computation of Algorithm \ref{alg: Stoc Proj algorithm} by adding a positive parameter $\beta$, independent of $i$, for all $i=1,\dots,M$ and $n$, times the identity matrix, giving a new diffusion. This operation translates into the following stochastic basis update
\begin{equation}\label{eq: Y beta}
	\tildhatY{i}{n+1} = \widehat{Y}_n^{i} + \widehat{U}_n \widehat{a}_n^{i} \Delta t_n +  \widehat{U}_n \underbrace{\left( \widehat{b}_n^{i} + \beta I_{d \times m} \right)}_{\widehat{b}_n^{\beta,i}} \Delta W_n^{i}, \quad i=1, \dots, M,
\end{equation}
where $I_{d \times m}$ is a rectangular matrix such that $(I_{d \times m})_{ij}:= 1$, if $i=j$ and $0$ otherwise. Hereafter, we focus on the case $m=d$ and a positive diffusion matrix. 
\begin{Assumption}\label{ass: b sspd}
	The matrix $b(t,x) \in  \mathbb{R}^{d \times d}$ is positive semi-definite for all $t \in [0,\infty]$ and for all $x \in \mathbb{R}^d$.
\end{Assumption}

Notice that under this Assumption, \eqref{eq: Y beta} implies that
\begin{equation*}
	\left(b(t, x)+ \beta I_{d \times d}\right) \left(b(t, x)+ \beta I_{d \times d}\right)^{\top}\succeq \beta^2 I_{d \times d}, \quad \text{for all } t   \in [0,\infty), \text{ for all } x \in \mathbb{R}^d,
\end{equation*}
and the modified diffusion $b^{\beta}(t,x)=b(t,x) + \beta I_{d \times d}$ satisfies Assumption \ref{ass: diff}. We call this new version of Algorithm \ref{alg: Stoc Proj algorithm} as \emph{Diffusion-regularized DLR PS SDE}.

We want to show convergence of this new algorithm with respect to the semidiscretized DLRA $X_n$. In order for our study to make sense, we need $X_n$ to be well-defined. We then work with the following assumption for the semi-discrete DLRA solution.
\begin{Assumption}\label{ass: gamma cov}
	There exists a positive constant $\gamma$ independent of $n$ for $1 \leq n \leq N$, such that $\sigma^{k}_n \geq \gamma$ for all $n=1, \dots, N$. 
\end{Assumption}
Notice that Assumption \ref{ass: gamma cov} implies that $\widetilde{\sigma}^{k}_n \geq \gamma$, as well. Furthermore, notice that $\gamma$ is independent of the parameter of regularization $\beta$.
Under this new framework, we want to choose $\beta$ according to the number of samples $M$ so that we can recover the best convergence with respect to the Monte Carlo error. We observe that in Proposition \ref{prop: lower-bound Gramian}, the Assumption \ref{ass: diff} can be replaced with Assumption \ref{ass: gamma cov}, modulo replacing $\sigma_B$ with $\beta^2$. The parameter $\varrho$ defined in \eqref{eq: rho} now becomes
\begin{equation}\label{eq: varrho - beta}
	\varrho := \frac{\beta^2}{4 C_{\mathrm{lgb}}(1 + 2 K_2(T))}.
\end{equation}
and the probability of the complement of the set $E_n$ defined in \eqref{eq: E_n} is now bounded by
\begin{equation}
	\begin{aligned}
		\mathbb{P}(E_n^C) \leq & \tilde{C}_5 \left( \left(1 +\frac{1}{\beta^{3p}}\right) \frac{1}{M^p} + \frac{1}{\beta^{4p}} \frac{1}{M^{2p-1}} +\frac{\beta^p}{M^{p}} + \frac{1}{M^{2p-1}} \right) \\
		\leq &\mathring{C}_5 (\frac{1}{M^p} + \frac{1}{\beta^{3p}} \frac{1}{M^p} + \frac{1}{\beta^{4p}} \frac{1}{M^{2p-1}}) \\
	\end{aligned}
\end{equation} for positive constants $\tilde{C}_5$ and $\mathring{C}_5$ dependent on $\gamma$, but independent of $\beta$, if the initial datum $Y_0$ has at least $4p$-moments. In $E_n$, $\tildhatsigma{k}{n}{M}$ is always strictly positive, and hence $P_{\tildhatY{i}{n+1}}$ does not need a regularization of the Gramian. In this setting, the event $A_n$ becomes
\begin{equation}\label{eq: A_n reg diff}
	A_n := E_n \cap \{ \widetilde{\sigma}^k_{n,M} \geq \frac{\gamma}{2}\}.
\end{equation}

To prove convergence, we want to compare the non-regularized empirical projector defined in \eqref{eq: Y discretized projectors - full rank} in the favorable event $A_n$ with the empirical projector of the semi-discrete solution, and bound their differences in terms of the error between the semi-discrete DLRA solution and the regularized fully discrete one.

\begin{Lemma}[Bound on the diffusion-regularized empirical stochastic projectors]\label{lemma: regularized diffusion empir Proj}
	Suppose that $\mathbb{E}[|Y_0|^{4}] < \infty$, that $\sigma_{Y_0}= \mathbb{E}[|\widehat{Y}_0\widehat{Y}_0^{\top}|^2] >0$, that Assumptions \ref{ass: b sspd} and \ref{ass: gamma cov} hold, and assume that $\Delta t_n$ satisfies \eqref{eq: dt semidiscrete condition} and \eqref{eq: M dt cond} with $p=2$. 
	
	Furthermore, consider the operators $\widehat{P}_{\widetilde{Y}_{n+1}^{i}},\widehat{P}_{\tildhatY{i}{n+1}}: L^2(\Omega^{M};\mathbb{R}^{d \times M}) \to L^2(\Omega^M,\mathbb{R}^{d})$ defined as in \eqref{eq: Y discretized projectors - full rank}.
	
	Consider $f\in L^4_{iid}\left(\Omega^{M}, \mathbb{R}^{d \times M}\right)$ such that for positive constants $\tilde{C},\tilde{c}$, one has $\mathbb{P}(|f^{i}| > R) \leq \tilde{C} e^{-\tilde{c}R^2} $ for all $R>0$. Then, there exist positive constants  $G_6:=G_6(\gamma)=O(\frac{1}{\gamma})$ dependent on the inverse of $\gamma$, $c$, and $C$, such that 
	\begin{equation}\label{eq: discr proj Y - A_n}
		\begin{aligned}
			\sqrt{\mathbb{E}\left[\frac{1}{M}\sum_{i=1}^{M}
				\left|\left(\widehat{P}_{\tildhatY{i}{n+1}}-\widehat{P}_{\widetilde{Y}_{n+1}^{i}}\right)f\right|^2 \mathbbm{1}_{A_{n+1}}\right]} \leq &  \left(e_n + \frac{1}{\sqrt{M}} + \beta \sqrt{\Delta t} \right)   G_6 R + (C e^{-cMR^2} ) \|f\|_{L^4(\Omega^M;\mathbb{R}^{d \times M})},
		\end{aligned}
	\end{equation}
	for all $R > \|f\|_{L^2(\Omega^M;\mathbb{R}^{d \times M})}$.
	\begin{proof}
		The proof can be found in Appendix \ref{app: proof deg diff}.
	\end{proof}
\end{Lemma} 

\begin{Theorem}[Stochastic Convergence of the Diffusion-regularized modified DLR PS SDE]
	
	Suppose that $\mathbb{E}[|Y_0|^{8}] < \infty$, that $\sigma_{Y_0}= \mathbb{E}[\widehat{Y}_0\widehat{Y}_0^{\top}] >0$, that Assumptions \ref{ass: b sspd} and \ref{ass: gamma cov} hold, and assume that $\Delta t_n$ satisfies \eqref{eq: dt semidiscrete condition} and \eqref{eq: M dt cond} for $p=4$.  
	
	Choosing $\beta=M^{-\frac{p}{(4+3p)}}$, then there exists a positive constant $C:=C( \gamma)=O(\frac{1}{\gamma})$ dependent on the inverse of $\gamma$ such that it holds
	\begin{equation}\label{eq: error diff-mod Stoc Proj}
		e_{n+1} \leq  \left(M^{-\frac{p}{(4+3p)}+ \frac{\log \log M}{2 \log M}CT}  \right) \exp \left(C  T\right), \quad \text{ for all } \ 0 \leq n \leq N-1.
	\end{equation}
	\begin{proof}
		We compare the semi-discrete DLRA solution evaluated at the samples $(\omega_i)_{i=1,\dots,M}$, i.e.\
		\(X_{n+1}^{i}=X_{n+1}(\omega_i)\), with the diffusion-regularized fully
		discretized solution \(\widehat X_{n+1}^{i}\).	
		As in the proof of Theorem~\ref{thm: convergence of the mod Stoc Proj},
		we split the drift contribution in
		\(X_{n+1}^{i}-\widehat X_{n+1}^{i}\) over \(A_{n+1}\) and \(A_{n+1}^C\), where \(A_{n+1}\) is defined in \eqref{eq: A_n reg diff}.
		On \(A_{n+1}\), the empirical Gramian is bounded from below, so the empirical
		projector in the fully discrete drift term is stable; on \(A_{n+1}^C\), the
		loss of stability is compensated by the probability estimate for
		\(A_{n+1}^C\).  The diffusion contribution is treated as in the
		non-regularized case, except for the additional difference
		\(\widehat b_n^{\beta,i}-\widehat b_n^i\), which is controlled by the
		regularization parameter \(\beta\). Namely, we consider
		
		\begin{equation}\label{eq: intermediate error diff-reg}
			\begin{aligned}
				e^2_{n+1} = &  \mathbb{E}\left[\frac{1}{M}\sum_{i=1}^{M} |X_{n+1}^{i}-\widehat{X}_{n+1}^{i}|^2\right]\\
				= &  \mathbb{E}\left[\frac{1}{M}\sum_{i=1}^{M} |X_{n}^{i}+P_{{U}^{\top}_n\widetilde{Y}_{n+1}^{i}}[ a_n] \Delta t_n+ P_{{U}_n}[b_n^{i} \Delta W_n^{i}] -\widehat{X}_{n}^{i}-\widehat{P}_{\widehat{U}^{\top}_n\tildhatY{i}{n+1}}[\widehat{a}_n] \Delta t_n- P_{\widehat{U}_n}[\widehat{b}_n^{\beta,i} \Delta W_n^{i}] |^2\right]\\
				=& e_{n}^2 + 2\mathbb{E}\left[\frac{1}{M}\sum_{i=1}^{M} \langle X_{n}^{i}-\widehat{X}_{n}^{i},P_{{U}^{\top}_n\widetilde{Y}_{n+1}^{i}}[ a_n] \Delta t_n -\widehat{P}_{\widehat{U}^{\top}_n\tildhatY{i}{n+1}}[\widehat{a}_n] \Delta t_n \rangle \right]\\
				& + \mathbb{E}\left[\frac{1}{M}\sum_{i=1}^{M} |P_{{U}^{\top}_n\widetilde{Y}_{n+1}^{i}}[ a_n] \Delta t_n+ P_{{U}_n}[b_n^{i} \Delta W_n^{i}] -\widehat{P}_{\widehat{U}^{\top}_n\tildhatY{i}{n+1}}[\widehat{a}_n] \Delta t_n- P_{\widehat{U}_n}[\widehat{b}_n^{\beta,i} \Delta W_n^{i}] |^2\right]\\
				\leq & 	e_{n}^2 + 2e_n \Big( \mathbb{E}\left[\frac{1}{M}\sum_{i=1}^{M} \left| P_{{U}^{\top}_n\widetilde{Y}_{n+1}^{i}}[ a_n]  -\widehat{P}_{\widehat{U}^{\top}_n\tildhatY{i}{n+1}}[\widehat{a}_n] \right|^2 \mathbbm{1}_{A_{n+1}}\right] \Big)^{\frac{1}{2}} \Delta t_n\\
				&+2e_n  \Big(\mathbb{E}\left[\frac{1}{M}\sum_{i=1}^{M} \left| P_{{U}^{\top}_n\widetilde{Y}_{n+1}^{i}}[ a_n]  -\widehat{P}_{\widehat{U}^{\top}_n\tildhatY{i}{n+1}}[\widehat{a}_n]  \right|^2 \mathbbm{1}_{A_{n+1}^C}\right]\Big)^{\frac{1}{2}} \Delta t_n \\
				& + \mathbb{E}\left[\frac{1}{M}\sum_{i=1}^{M} |P_{{U}^{\top}_n\widetilde{Y}_{n+1}^{i}}[ a_n] \Delta t_n+ P_{{U}_n}[b_n^{i} \Delta W_n^{i}] -\widehat{P}_{\widehat{U}^{\top}_n\tildhatY{i}{n+1}}[\widehat{a}_n] \Delta t_n - P_{\widehat{U}_n}[\widehat{b}_n^{i} \Delta W_n^{i}] |^2\right] \\
				&+ \mathbb{E}\left[\frac{1}{M}\sum_{i=1}^{M} |P_{\widehat{U}_n}[\left(\widehat{b}_n^{\beta,i} - \widehat{b}_n^{i}\right)\Delta W_n^{i}] |^2\right],
			\end{aligned}
		\end{equation}
		where in the second to last line we used the property of independence of the Brownian increments with respect to the discretized solutions. Now, we give reasonable bounds on the cross term and on the last quadratic term in \eqref{eq: intermediate error diff-reg} in order to obtain a relation on $e_n$ to conclude via discrete Gronwall's lemma. 
		
		Part of the proof is similar to the one of Theorem \ref{thm: convergence of the mod Stoc Proj}, as by Assumption \ref{ass: gamma cov} and by regularization of the diffusion the involved projectors are all of full rank $k$.
		For the first cross term in \eqref{eq: intermediate error diff-reg} via using Young's inequality, Lemmata \ref{lem: diff Proj U}, \ref{lem: error between stochastic proj - full rank}, the linear-growth bound, Lemmata \ref{lem: L2p norm semidiscretized solution} and \ref{lem: Lp moment bound}, one gets
		\begin{equation*}
			\begin{aligned}
				& e_n \Big( \mathbb{E}\left[\frac{1}{M}\sum_{i=1}^{M} \left| P_{{U}^{\top}_n\widetilde{Y}_{n+1}^{i}}[ a_n]  -\widehat{P}_{\widehat{U}^{\top}_n\tildhatY{i}{n+1}}[\widehat{a}_n] \right|^2 \mathbbm{1}_{A_{n+1}}\right] \Big)^{\frac{1}{2}} \Delta t_n \\
				\leq & e_n^2 C(\gamma,T) (1+R)\Delta t_n + \frac{1}{M} C(\gamma,T) (1+R)\Delta t_n  + \beta^2 \Delta t_n^2 (C(\gamma,T) (1+R) + 1 ) + C(\gamma,T) (1+R) e^{-MR^2}  \Delta t_n
			\end{aligned}
		\end{equation*}
		for a positive constant $C:=C(\gamma,T)$ independent of $M$, $\Delta t$, but dependent on the inverse of $\gamma$. For the second cross term, via Holder and Young's inequality one obtains
		\begin{equation*}
			\begin{aligned}
				&e_n\Big(\mathbb{E}\left[\frac{1}{M}\sum_{i=1}^{M} \left| P_{{U}^{\top}_n\widetilde{Y}_{n+1}^{i}}[ a_n]  -\widehat{P}_{\widehat{U}^{\top}_n\tildhatY{i}{n+1}}[\widehat{a}_n]  \right|^2 \mathbbm{1}_{A_{n+1}^C}\right]\Big)^{\frac{1}{2}} \Delta t_n \\
				\leq& \frac{e_n^2}{2} \Delta t_n  + \mathbb{E}\left[\frac{1}{M}\sum_{i=1}^{M} \left| P_{{U}^{\top}_n\widetilde{Y}_{n+1}^{i}}[ a_n]  -\widehat{P}_{\widehat{U}^{\top}_n\tildhatY{i}{n+1}}[\widehat{a}_n]  \right|^2 \mathbbm{1}_{A_{n+1}^C}\right] \Delta t_n \\
				\leq &\frac{e_n^2}{2} \Delta t_n + \mathbb{P}(A_{n+1}^C)^{\frac{1}{2}} \mathbb{E}\left[\left(\frac{1}{M}\sum_{i=1}^{M} \left| P_{{U}^{\top}_n\widetilde{Y}_{n+1}^{i}}[ a_n]  -\widehat{P}_{\widehat{U}^{\top}_n\tildhatY{i}{n+1}}[\widehat{a}_n]  \right|^2 \right)^2\right] ^{\frac{1}{2}} \Delta t_n \\
				\leq &\frac{e_n^2}{2} \Delta t_n + \mathbb{P}(A_{n+1}^C)^{\frac{1}{2}} \mathbb{E}\left[\left(\frac{2}{M}\sum_{i=1}^{M} \left| P_{{U}^{\top}_n\widetilde{Y}_{n+1}^{i}}[ a_n]  \right|^2 + \frac{2}{M}\sum_{i=1}^{M} \left|\widehat{P}_{\widehat{U}^{\top}_n\tildhatY{i}{n+1}}[\widehat{a}_n]  \right|^2 \right)^2\right] ^{\frac{1}{2}} \Delta t_n \\
				\leq &\frac{e_n^2}{2} \Delta t_n + \tilde{C} \left(\frac{1}{M^{p}}+ \frac{1}{\beta^{3p} M^{p}}+ \frac{1}{\beta^{4p}} \frac{\Delta t^p}{M^{p}} \right)^{\frac{1}{2}}  \Delta t_n,
			\end{aligned}
		\end{equation*}
		where in the second to last line we employed Lemma \ref{lem: P f < f} and \ref{lem: Lp moment bound}.
		Similar computations can be derived for the quadratic terms of the non-regularized diffusion.
		
		Finally, for the regularized-diffusion term one obtains that
		\begin{equation*}
			\mathbb{E}\left[\frac{1}{M}\sum_{i=1}^{M} |P_{\widehat{U}_n}[(\widehat{b}_n^{\beta,i} -\widehat{b}_n^{i})\Delta W_n^{i}] |^2\right] \leq \beta^2 \Delta t_n.
		\end{equation*}
		
		By grouping all these computations, there exists a positive constant $\widetilde{G}:=\widetilde{G}(\gamma, T)$ dependent on the inverse of $\gamma$, such that
		\begin{equation}\label{eq: gronw beta}
			\begin{aligned}
				e^2_{n+1} \leq & e_{n}^2 + \widetilde{G}(1+R) e_n^2 \Delta t_n + \widetilde{G} \left(\frac{R}{M} + R \beta^2 \Delta t + e^{-MR^2}  \right)\Delta t_n + \widetilde{G} \left(\frac{1}{\beta^{\frac{3p}{2}} M^{\frac{p}{2}}}+ \frac{1}{\beta^{2p}} \frac{\Delta t^{\frac{p}{2}}}{M^{\frac{p}{2}}} +  \beta^2\right) \Delta t_n.
			\end{aligned}
		\end{equation}
		Assume that $\Delta t \leq \beta$. Then, to have the best possible convergence rate with respect to the sample size $M$, we want the following relation to hold
		\begin{equation*}
			\frac{1}{\beta^{\frac{3p}{2}} M^{\frac{p}{2}} } \sim  \beta^2,
		\end{equation*}
		which gives the optimal $\beta$ relation with the following 
		\begin{equation*}\label{eq: optimal beta}
			\beta_{\mathrm{opt}} \propto M^{-\frac{p}{(4+3p)}}.
		\end{equation*}
		Then, \eqref{eq: gronw beta} becomes
		\begin{equation*}
			\begin{aligned}
				e^2_{n+1} \leq & e_{n}^2 +\widetilde{G}(1+R) e_n^2 \Delta t_n + \widetilde{G} \left(\frac{R}{M} + R M^{-\frac{2p}{(4+3p)}}+ e^{-MR^2}  \right)\Delta t_n + \widetilde{G} M^{-\frac{2p}{(4+3p)}} \Delta t_n,
			\end{aligned}
		\end{equation*}
		and the Gronwall's lemma implies that
		\begin{equation*}
			\begin{aligned}
				e^2_{n+1} \leq & \left( M^{-\frac{2p}{(4+3p)}} + R M^{-\frac{p}{(4+3p)}} + \frac{R}{M} +  e^{-MR^2}   \right) \exp\{ C \left(1+ R\right) T \}, \quad \text{ for all } \ 1 \leq n \leq N.
			\end{aligned}
		\end{equation*}
		Similarly to the proof of Theorem \ref{thm: convergence of the Stoc Proj}, via choosing $R= \log \log M$ for $M$ large one obtains the thesis.
	\end{proof}
\end{Theorem}

\section{On Monte-Carlo convergence of Algorithms \ref{alg: DLR EM SDE algorithm} and \ref{alg: Eva Proj Splitt SDE algorithm}}\label{sec: DLR EM & KNV}
In Section \ref{sec: stoc proj}, we established convergence results for the stochastic discretization of the DLR PS SDE. The analysis relied on Lipschitz-type estimates between projectors, which in turn depend on the boundedness of the $p$-moments of both the semi-discrete and fully discretized solutions.

However, similar convergence results cannot be obtained as easily for the DLR Euler--Maruyama scheme and the DLR Projector Splitting for EM (Algorithms \ref{alg: DLR EM SDE algorithm} and \ref{alg: Eva Proj Splitt SDE algorithm}).
In this section we explain the main difficulties that prevent to obtain similar results.

\paragraph{DLR Euler--Maruyama.}
Concerning the DLR Euler-Maruyama, the problem is due to the potential lack of boundedness of $p$-moments of the fully-discrete solution, which is crucial in the statements of Section \ref{sec: stoc proj}. Indeed, the usual update in the ambient space reads as
\begin{equation}\label{eq: MC DLR EM X}
	\begin{aligned}
		\widehat{X}_{n+1}^{i} = &  \tildhatUT{n+1} \tildhatY{i}{n+1} \\
		=& 	\widehat{X}_n ^{i}	+\widehat{P}_{\widehat{U}_n^{\top}\widehat{Y}_{n}^{i}}[\widehat{a}_n^{i}]\Delta t_n+  \widehat{P}_{\widehat{U}_n}[\widehat{b}_{n}^{i} \Delta W_n^{i}] \\
		&+\left(I_{d \times d} - \widehat{P}_{\widehat{U}_n} \right) \widehat{\mathbb{E}}\left[\widehat{a}_n(\widehat{Y}_{n})^{\top}\right]\widehat{C}^{-1}_{\widehat{Y}_n} \left[\widehat{U}_n \widehat{a}_n^{i} (\Delta t_n)^2 + \widehat{U}_n \widehat{b}_{n}^{i} (\Delta W_n^{i}) \Delta t_n \right].\\
	\end{aligned}
\end{equation}
Relation \eqref{eq: MC DLR EM X} consists in the sum of a forward in-time increment projected in the point $\widehat{X}_n$ and a cross term that does not seem to present any structure useful to establish a priori bounds.  It is not straighforward how to derive stability estimates for \eqref{eq: MC DLR EM X}, specifically bounds on the moments of the numerical solution, which are essential in proving the convergence of the DLR PS SDE. 

In \cite{kazashi2025dynamicalpartI}, the boundedness of the second moment (in strong sense) for \eqref{eq: DLR EM X} was proved by induction, exploiting the fact that the expectation present in the cross term was taken with respect to the true measure of the semi-discrete DLRA. It is not immediate how to proceed in \eqref{eq: MC DLR EM X}, unless additional hypotheses are made. For instance, in the case of linear deterministic drift, i.e. $a(t,x)= A(t)x$ with $A \in \mathbb{R}^{d \times d}$, we know that the DLR EM is independent of the smallest singular value of the Gramian, as the evolution of the deterministic modes is decoupled from the stochastic basis. This property of the drift implies that we can avoid regularization of the Gramian, i.e.\ we can choose $\alpha=0$, and, hence, \eqref{eq: MC DLR EM X} becomes
\begin{equation*}
	\begin{aligned}
		\widehat{X}_{n+1}^{i} = & 	\widehat{X}_{n}^{i}  + \left(I_{d \times d} - \widehat{P}_{\widehat{U}_n} \right) \widehat{\mathbb{E}}\left[A(t_n)\widehat{U}_n^{\top}\widehat{Y}_n\widehat{Y}_n ^{\top}\right]\widehat{C}^{-1}_{\widehat{Y}_n}\widehat{Y}_{n}^{i} \Delta t_n +  \widehat{P}_{\widehat{U}_n} A(t_n) \widehat{X}_{n}^{i}\Delta t_n\\
		&+\left(I_{d \times d} - \widehat{P}_{\widehat{U}_n} \right) \widehat{\mathbb{E}}\left[A(t_n)\widehat{U}_n^{\top}\widehat{Y}_n\widehat{Y}_n ^{\top}\right]\widehat{C}^{-1}_{\widehat{Y}_n} \left[\widehat{U}_n A(t_n)\widehat{X}_{n}^{i} (\Delta t_n)^2 + \widehat{U}_n \widehat{b}_{n}^{i} (\Delta W_n^{i}) \Delta t_n \right] \\
		=& 	\widehat{X}_{n}^{i}  + A(t_n)\widehat{X}_{n}^{i}\Delta t_n +\left(I_{d \times d} - \widehat{P}_{\widehat{U}_n} \right) A(t_n)\widehat{P}_{\widehat{U}_n}  \left[ A(t_n)\widehat{X}_{n}^{i} (\Delta t_n)^2 + \widehat{b}_{n}^{i} (\Delta W_n^{i}) \Delta t_n \right].\\
	\end{aligned}
\end{equation*}
Therefore, we bound the second moment effortlessly using the linear-growth bound, orthogonality of $\widehat{U}_n$, and Young's inequality, as
\begin{equation*}
	\begin{aligned}
		\mathbb{E}[	|\widehat{X}_{n+1}|^2] = & \mathbb{E}[|\widehat{X}_n |^2] + 2\mathbb{E}[\widehat{X}_n^{\top}A(t_n)\widehat{X}_n]\Delta t_n + \mathbb{E}[|A(t_n)\widehat{X}_n|^2](\Delta t_n)^2  \\
		&+\mathbb{E}\left[	\widehat{X}_n^{\top}A(t_n)^{\top}\left(I_{d \times d} - \widehat{P}_{\widehat{U}_n} \right) A(t_n)\widehat{P}_{\widehat{U}_n}  \left[ A(t_n)\widehat{X}_n (\Delta t_n)^2 + \widehat{b}_{n}^{i} (\Delta W_n^{i}) \Delta t_n \right]\right] \Delta t_n \\
		&+\mathbb{E}\left[	\left|\left(I_{d \times d} - \widehat{P}_{\widehat{U}_n} \right) A(t_n)\widehat{P}_{\widehat{U}_n}  \left[ A(t_n)\widehat{X}_n (\Delta t_n)^2 + \widehat{b}_{n}^{i} (\Delta W_n^{i}) \Delta t_n \right]\right|^2\right]\\
		\leq & \mathbb{E}[|\widehat{X}_n |^2] + A(t_n)\mathbb{E}[|\widehat{X}_n|^2]\Delta t_n  + |A(t_n)|^2 \mathbb{E}[|\widehat{X}_{n}|^2] (\Delta t_n)^2 +	|A(t_n)|^3 \mathbb{E}[|\widehat{X}_n |^2] (\Delta t_n)^3  \\
		&+		|A(t_n)|^4 \mathbb{E}[|\widehat{X}_n |^2] (\Delta t_n)^4  + \frac{C_{\mathrm{lgb}}}{2}(1 +\mathbb{E}[|\widehat{X}_{n}|^2]) (\Delta t_n)^3,\\
	\end{aligned}
\end{equation*}
where we also exploit the independence of $\Delta W_n^{i}$, Young's inequality, and the linear-growth bound. Finally, the claim follows by induction with a condition on $\Delta t_n$ without assuming any more hypothesis on the diffusion than its linear-growth boundedness (see \eqref{ass: linear-growth bound}). Similarly, one can derive norm-bounds for the semi-discrete solution that are independent of the smallest singular value $\sigma_n^{k}$.

\paragraph{DLR Projector Splitting for EM.}
Different issues arise when dealing with the stochastic convergence of the DLR PS EM. Indeed, pursuing the analysis of the error $e_{n+1}$ similarly to the proof of Theorem \ref{thm: convergence of the Stoc Proj}, the following cross term would appear in the recursion of \eqref{eq: intermediate error}:
\begin{equation}\label{eq: bias term}
	\mathfrak{B}_n:=\mathbb{E}\left[\frac{1}{M}\sum_{i=1}^{M} (X_{n}^{i}-\widehat{X}_{n}^{i})^{\top}\left(P_{{U}^{\top}_n\widetilde{Y}_{n+1}^{i}}[ b_n \Delta W_n]-\widehat{P}_{\widehat{U}^{\top}_n\tildhatY{i}{n+1}}[\widehat{b}_{n} \Delta W_n] \right)\right].
\end{equation}
Notice that no Lipschitz-type bound derived in Section \ref{sec: stoc proj} can imply a useful bound that yields the convergence of the algorithm in the same fashion of Algorithm \ref{alg: Stoc Proj algorithm}. A possible sufficient condition to prove convergence would be to bound $\mathfrak{B}_n \lesssim e_n^2 \Delta t_n$, which would lead to conclusion by Gronwall's lemma. However, find a reasonable upper bound of \eqref{eq: bias term} is not straightforward. Future developments of this work could also be oriented to finding strategies to prove stochastic convergence of the particle system generated by Algorithm \ref{alg: Eva Proj Splitt SDE algorithm}. 

\section{Numerical Experiments}
\label{sec: numerical experiments}
In this part, we present numerical experiments that validate the results obtained in the previous sections. 

We plot the so-called \emph{$L^2$-sup-in-time error} $\mathcal{E}_{\sup}$, which, in the discrete setting, is approximated as follows:
\begin{equation*}
	\begin{aligned}
		\mathcal{E}_{\sup}^2 := \frac{1}{M} \sum_{i=1}^{M} \max\limits_{0\leq n \leq N} |X_{n}(\omega_i)-\widehat{X}_n^i|^2.
	\end{aligned}
\end{equation*}
We also plot the \emph{$L^2$ error at a given time $t_n$}, denoted by $\mathcal{E}_n$, which is approximated by
\begin{equation*}
	\begin{aligned}
		\mathcal{E}_n^2 := \frac{1}{M} \sum_{i=1}^{M}  |X_{n}(\omega_i)-\widehat{X}_n^{i}|^2.
	\end{aligned}
\end{equation*}
Similar errors can be defined by comparing the fully discrete solution $\widehat{\mathbb{X}}_n$ to the true solution $X^{\mathrm{true}}$ or the time continuous DLRA solution $X(t_n)$. To mitigate the possible random error due to the generation of the involved random variables, we will plot a Monte Carlo estimator of the quantities $\mathcal{E}_{\sup}^2$ and $\mathcal{E}_{n}^2$, involving $M_{\mathrm{runs}}$ independent simulations of the same problem.

In the numerical experiments discussed next, the true solution is unavailable, as well as the time continuous DLRA.
Errors between a fully-discretized algorithm and the semi-discretized DLRA will be illustrated. Since a closed-form solution of the latter is unavailable, we approximate it using a DLR PS SDE (Algorithm \ref{alg: Stoc Proj algorithm}) for a large number of particles, on the same time mesh as the numerical surrogate being tested. This choice has been made accordingly to Theorem \ref{thm: convergence of the Stoc Proj} and results shown in \cite{kazashi2025dynamicalpartI}, implying that the DLR PS SDE is converging, both in time and in the stochastic space via Monte Carlo method, to the continuous DLRA. We approximate $X^{\mathrm{true}}$ in the same way. For the sake of notation, in the caption of the figures of this section we will denote the semi-discrete DLRA as $DLR$.

\subsection{Projected SDE example}\label{sec: proj sde}
In this simulation, we show numerical performance concerning the convergence under regularization of the Gramian with $\alpha = \frac{1}{\sqrt[3]{M}}$ at each time step (which is the benchmark for Theorem \ref{thm: convergence of the Stoc Proj} in case of elliptic diffusion for example), as well as the case without regularization, i.e.\ with $\alpha=0$.

We simulate the following stochastic process 
\begin{equation}\label{ex: toy example stoc proj sde}
	\begin{aligned}
		\hspace{-0.3cm}
		\begin{pmatrix}
			\mathrm{d}X_1^{\mathrm{true}} \\
			\mathrm{d}X_2^{\mathrm{true}} \\
			\mathrm{d}X_3^{\mathrm{true}} \\
		\end{pmatrix}
		=& V(t)^{\top} V(t)
		\begin{pmatrix}
			\lambda \sin(X_1^{\mathrm{true}}(t)) & - \lambda \sin( X_2^{\mathrm{true}}(t))& - \theta X_3^{\mathrm{true}}(t) \\
			\lambda X_1^{\mathrm{true}}(t) & -\lambda X_2^{\mathrm{true}}(t)& - \theta X_3^{\mathrm{true}}(t) \\
			\lambda X_1^{\mathrm{true}}(t) & -\lambda X_2^{\mathrm{true}}(t)& - \theta X_3^{\mathrm{true}}(t) \\
		\end{pmatrix}
		\mathrm{d}t\\
		&+  V(t)^{\top} V(t)
		\begin{pmatrix}
			\sqrt{1 + 2 |X^{\mathrm{true}}(t)|} &0 \\
			0 &\sqrt{1 + 2 |X^{\mathrm{true}}(t)|} \\
			0 &\sqrt{1 + 2 |X^{\mathrm{true}}(t)|} \\
		\end{pmatrix}\cdot
		\begin{pmatrix}
			\mathrm{d}W_1(t) \\
			\mathrm{d}W_2(t) \\
		\end{pmatrix}
	\end{aligned}
\end{equation}
where $V(t) \in \mathbb{R}^{2 \times 3}$ is a time-dependent orthogonal matrix defined as
\begin{equation}\label{eq: V(t) subspace}
	V(t) = \begin{pmatrix}
		\frac{1}{\sqrt{2}}\cos(t) &\frac{1}{\sqrt{2}} \sin(t) &\frac{1}{\sqrt{2}}  \\
		\sin(t) &-\cos(t) &0\\
	\end{pmatrix}
\end{equation}
for $t \in [0,1]$, and $\lambda=0.6$, $\theta=0.5$. The initial datum is given by the random vector $X^{\mathrm{true}}(0) = V(0)^{\top} V(0)\tilde{X}$, where we define $\tilde{X}_i = \text{Un}_{i}(1,4)$, for $i=1,2$, where $\text{Un}_{1}(1,4)$ and $\text{Un}_{2}(1,4)$ are independent uniform random variables in the interval $[1,4]$. In equation \eqref{ex: toy example stoc proj sde} the drift and the diffusion are projected into the same subspace of dimension $2$, which is evolving over time. In this sight, we always decide to employ a DLR approximation of rank $k=2$. For this numerical experiment, the true DLR solution is computed with $M=1.5\cdot 10^5$ paths and time step $\Delta t = 1 \cdot 10^{-4}$. For all the Monte Carlo approximations we implemented, we consider the same step-size $\Delta t$. Therefore, we expect to see the Monte Carlo rate of convergence if the other sources of error, such as the time discretization, are negligeable. We run $M_{\mathrm{runs}}=25$ independent simulations and then we average out the errors $\mathcal{E}_{n}^2$ and $\mathcal{E}_{\sup}^2$.

\begin{figure}[!h]
	\centering
	\includegraphics[scale=0.29]{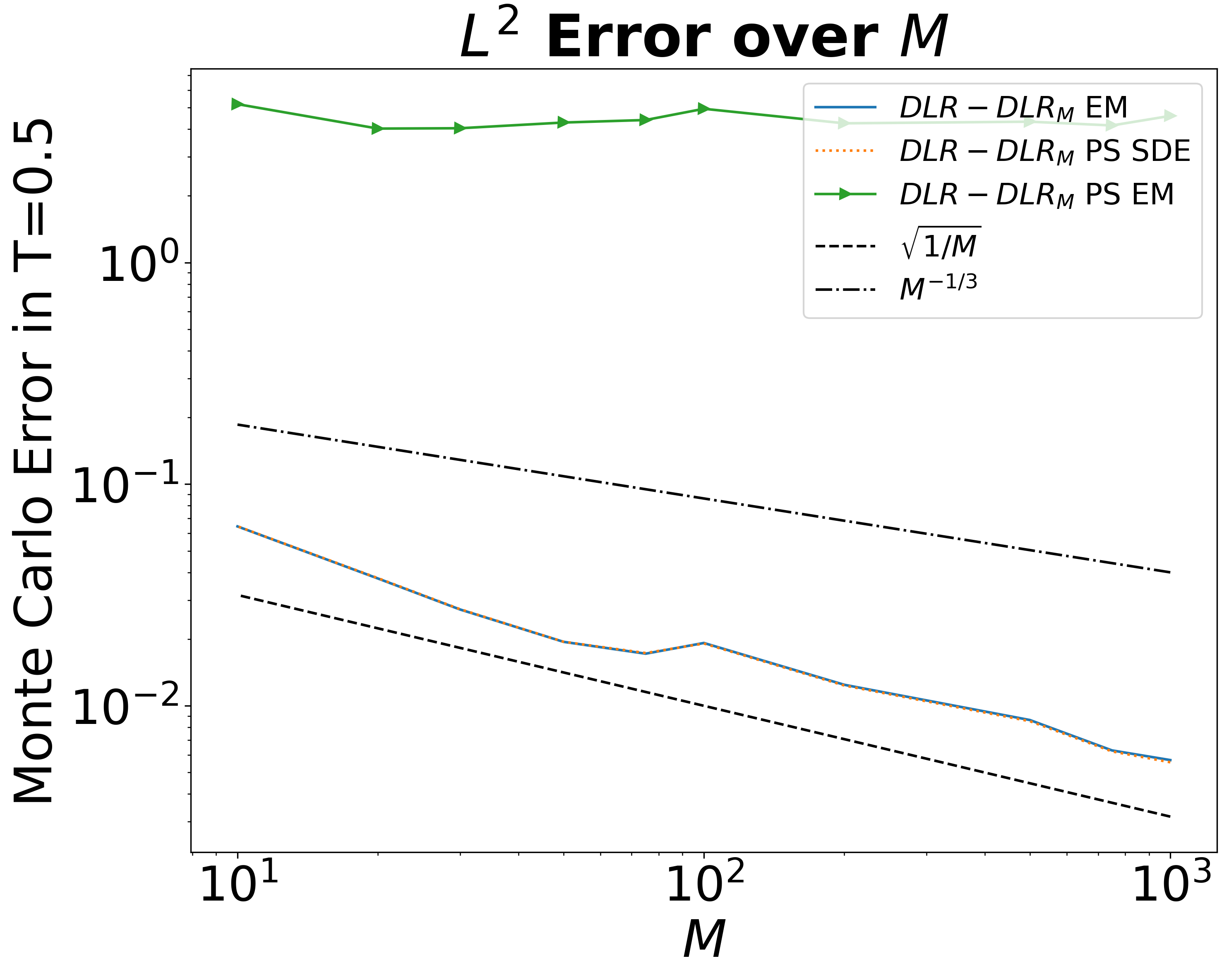}
	\includegraphics[scale=0.29]{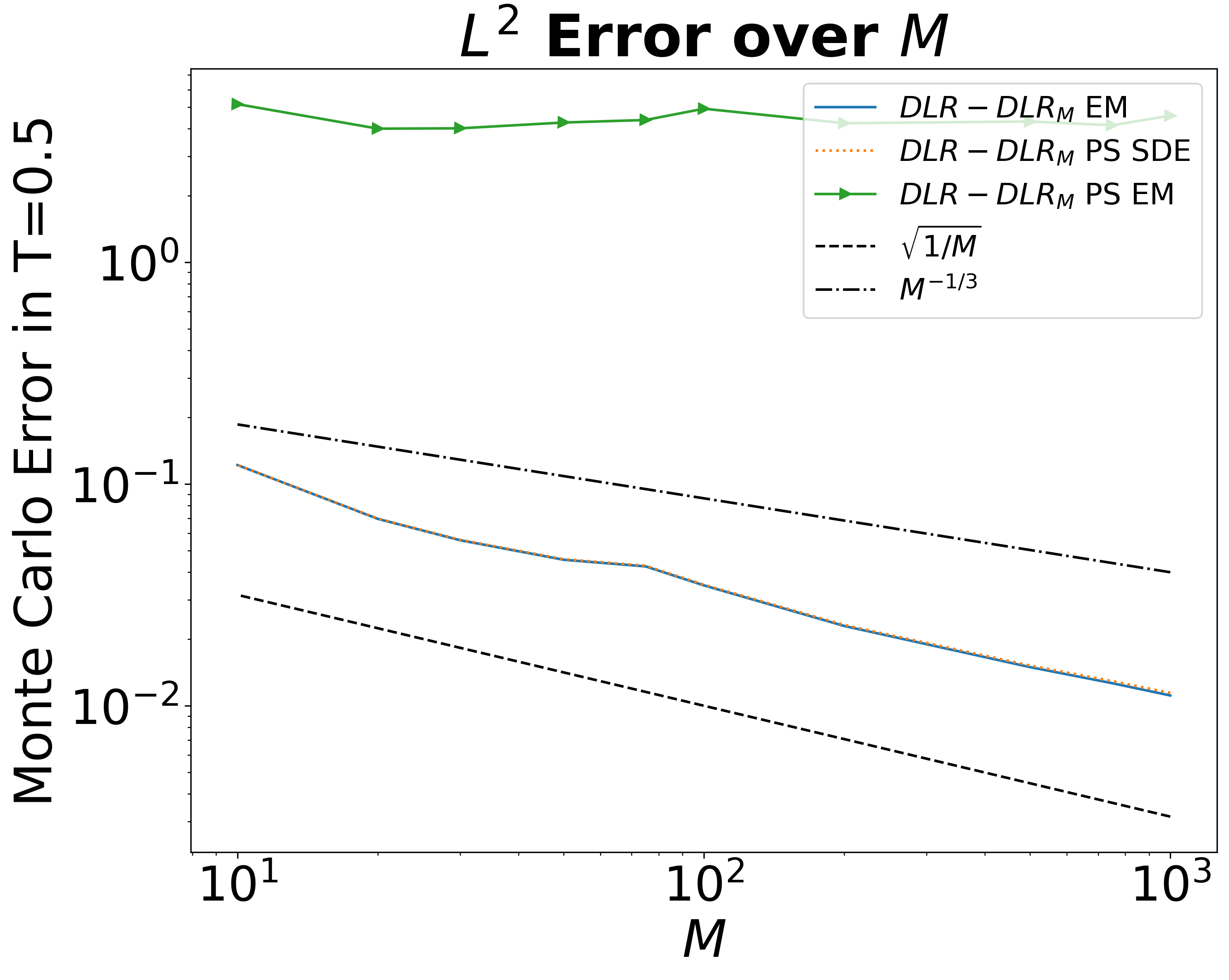}
	\caption{ $L^2$ errors for the DLR EM, DLR PS EM, and DLR PS SDE with respect to the $DLR$ solution at time $T=1$ (computed with $M=1.5\cdot 10^5$) without regularization $\alpha=0$ (left) and with regularization (right) of the Gramian with $\alpha = \frac{1}{\sqrt[3]{M}}$ for Problem \eqref{ex: toy example stoc proj sde}. $\Delta t=0.0001$, $k=2$.}
	\label{fig:toy example stoc proj sde - L2}
\end{figure}

\begin{figure}[!h]
	\centering
	\includegraphics[scale=0.29]{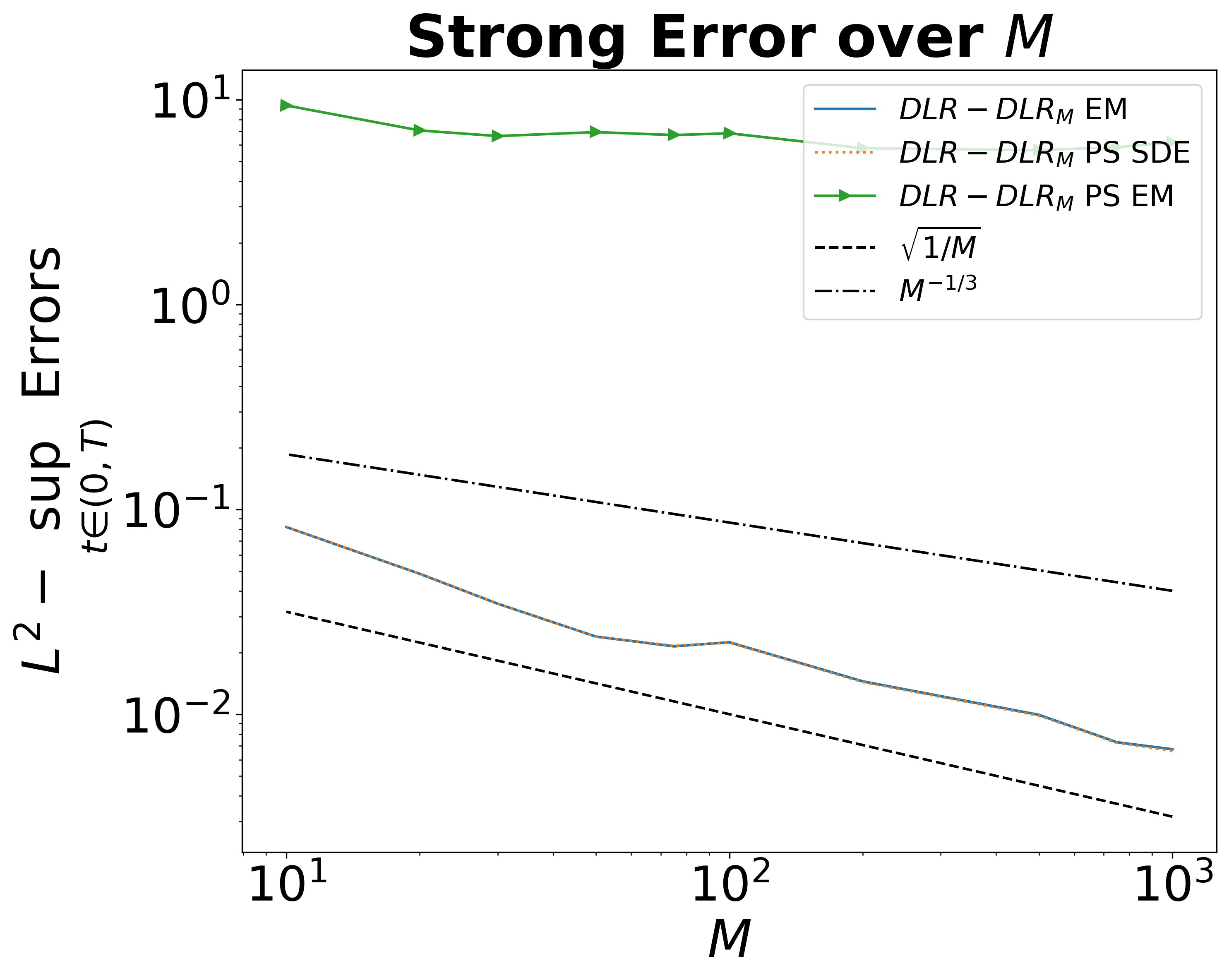}
	\includegraphics[scale=0.29]{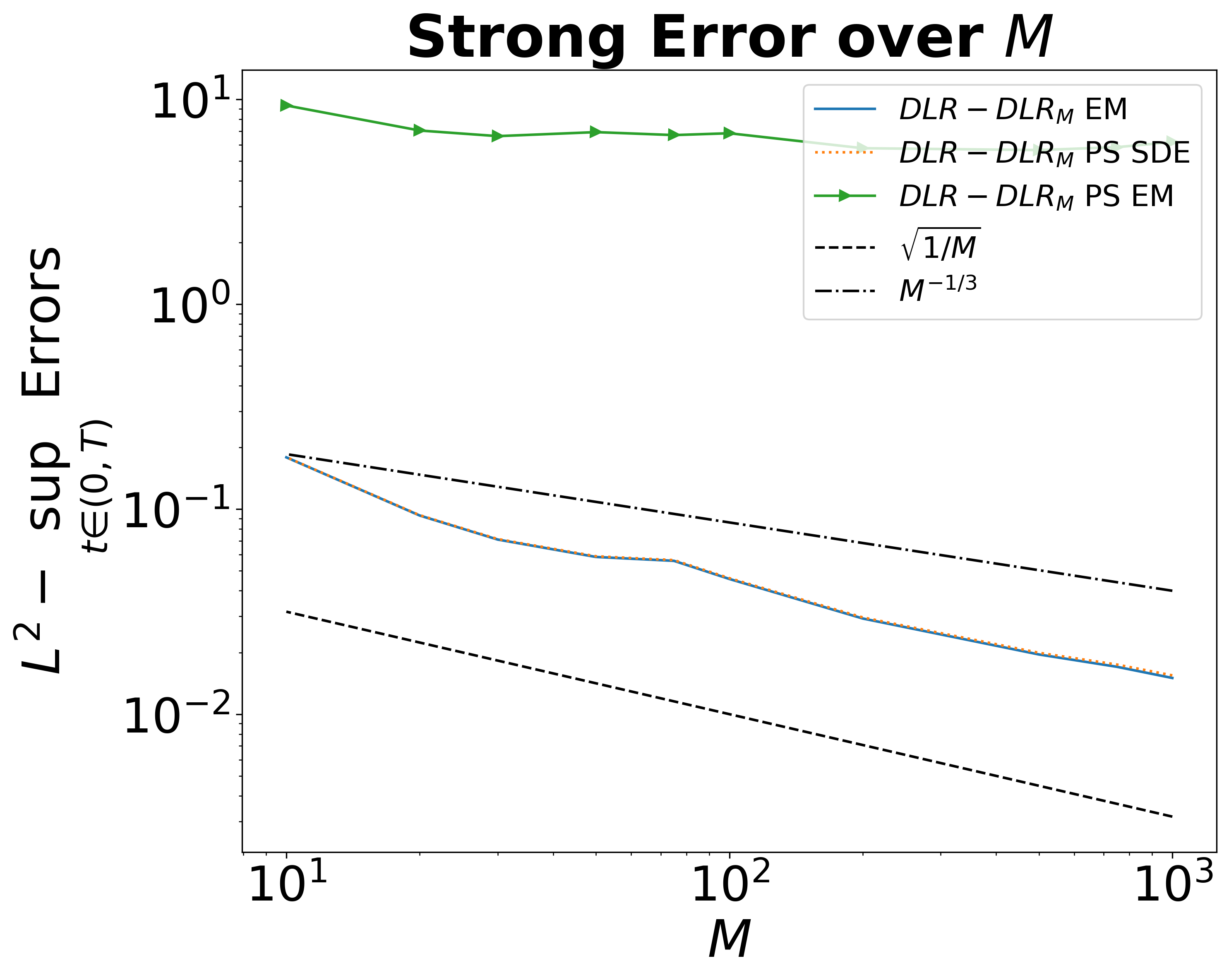}
	\caption{ $L^2$-sup-in-time errors for the DLR EM, DLR PS EM, and DLR PS SDE with respect to the semidiscrete $DLR$ solution at time $T=1$ (computed with $M=1.5\cdot 10^5$) without regularization, i.e.\ $\alpha=0$, (left) and with regularization (right) of the Gramian with $\alpha = \frac{1}{\sqrt[3]{M}}$ for Problem \eqref{ex: toy example stoc proj sde}. $\Delta t=0.0001$, $k=2$.}
	\label{fig:toy example stoc proj sde - L2 sup}
\end{figure}

To initiate all the DLRA algorithms, as $X^{\mathrm{true}}(0)=X(0)$ is already of rank equal to $2$, we consider a \texttt{QR} decomposition of the discretized initial datum. Precisely, the initial condition of Algorithms \ref{alg: DLR EM SDE algorithm}, \ref{alg: Stoc Proj algorithm}, and \ref{alg: Eva Proj Splitt SDE algorithm} is the rank-$k$ best approximation with respect to the Frobenius norm of $[(X^{\mathrm{true}})^{1}(0), \dots, (X^{\mathrm{true}})^{M}(0)] \in \mathbb{R}^{d \times M}$, i.e.\ the (rank-$k$) SVD of the matrix whose columns are the $M$ realizations of $X_0$. In this example the drift defined in \eqref{ex: toy example stoc proj sde} is non-linear and, hence, deterministic and stochastic basis of the DLRA are still coupled.
We retrieved the convergence over the stochastic discretization, as shown in Figure \ref{fig:toy example stoc proj sde - L2} for both the non-regularized with $\alpha=0$ and the regularized case. Notice that the DLR PS SDE saturates immediately, as it does not converges to the DLRA proposed in \cite{kazashi2025dynamical}, as pointed out in \cite{kazashi2025dynamicalpartI}. Concerning the $\mathcal{E}_{\sup}$ one can see in Figure \ref{fig:toy example stoc proj sde - L2 sup} that the non-regularized version approaches a Monte-Carlo rate of $O(\frac{1}{\sqrt{M}})$, whereas the Gramian-regularized version attains a rate of $O(\frac{1}{\sqrt[3]{M}})$ for a large number of particles, confirming the results of Theorem \ref{thm: convergence of the Stoc Proj} and \ref{thm: convergence of the mod Stoc Proj}.

\subsection{Stochastic Advection-Diffusion-Reaction PDE in 1D: Adaptive regularization}
Here, we present numerical performance of the adaptive regularization strategy (Algorithm \ref{alg: mod Stoc Proj algorithm}) via simulating an advection-diffusion-reaction problem under low-rank additive diffusion in 1D. The temporal domain is $[0,T]$, with $T=1$. The spatial domain is $[0,L]$ with length $L=1$, where we consider Neumann boundary conditions. More precisely, the problem reads:
\begin{equation}\label{ex: SADR}
	\begin{cases}
		\begin{aligned}
			\partial_t u(x,t,\omega)&= L u(x,t,\omega) + \sum_{i=0}^{m-1} \lambda^{i} \phi_i(x) \mathrm{d}W^i(\omega), \quad (x,t,\omega) \in [0,L] \times [0,T] \times \Omega, \\
			u_{0}(x,\omega)&=  \sum_{i = 0}^{m-1} \cos( \frac{i \pi x}{ L} )  \mathcal{N}_i(0,1), \quad (x,\omega) \in [0,L] \times \Omega\\
			\partial_x u(0,t, \omega)&= \partial_x u(1,t,\omega) = 0, \quad (t,\omega) \in [0,T] \times \Omega,
		\end{aligned}
	\end{cases}
\end{equation}
where in \eqref{ex: SADR} one has $L u(x,t,\omega) = a \Delta u(x,t,\omega) - v \nabla u(x,t,\omega) + r \sin (u(x,t,\omega))$,  $(\mathcal{N}_i(0,1))_{i=1,\dots,2m}$ are i.i.d.\ normal random variables independent of $(W^i)_{i=1,\dots,m}$, and $a=0.05$, $v= 0.01$, and $r= 0.01$. Concerning the diffusion term, we set $\phi_{i} = \cos(\frac{i \pi x}{L})$ for all $i=1,\cdots, m$, whereas $\lambda=\frac12$. Notice the scaling of the reaction term, for which we expect a stable and low-dimensional system.
The PDE is discretized with centered $2^{nd}$ order finite differences on a uniform grid with mesh size $\mathrm{d}x = 0.04$, hence leading to a $25$-dimensional problem. As usual, the time integration of the continuous DLR solution is made with the DLR PS SDE with mesh size $\Delta t=1 \cdot 10^{-4}$ and the number of paths used to simulate the reference DLRA is $M=10^5$. We keep the time-mesh fixed while simulating the various discrete DLR solutions. We choose $m=5$ and we employ a DLRA approximation of rank $k=5$, as the only term that can change the range of the solution, the advection one, has small contribution to the dynamics. Again, notice that as the drift is nonlinear, the evolution of $U_t$ and $Y_t$ is still coupled. We run $M_{\mathrm{runs}}=25$ independent simulations and then we average out the various errors.

\begin{figure}[!h]
	\centering
	\includegraphics[scale=0.29]{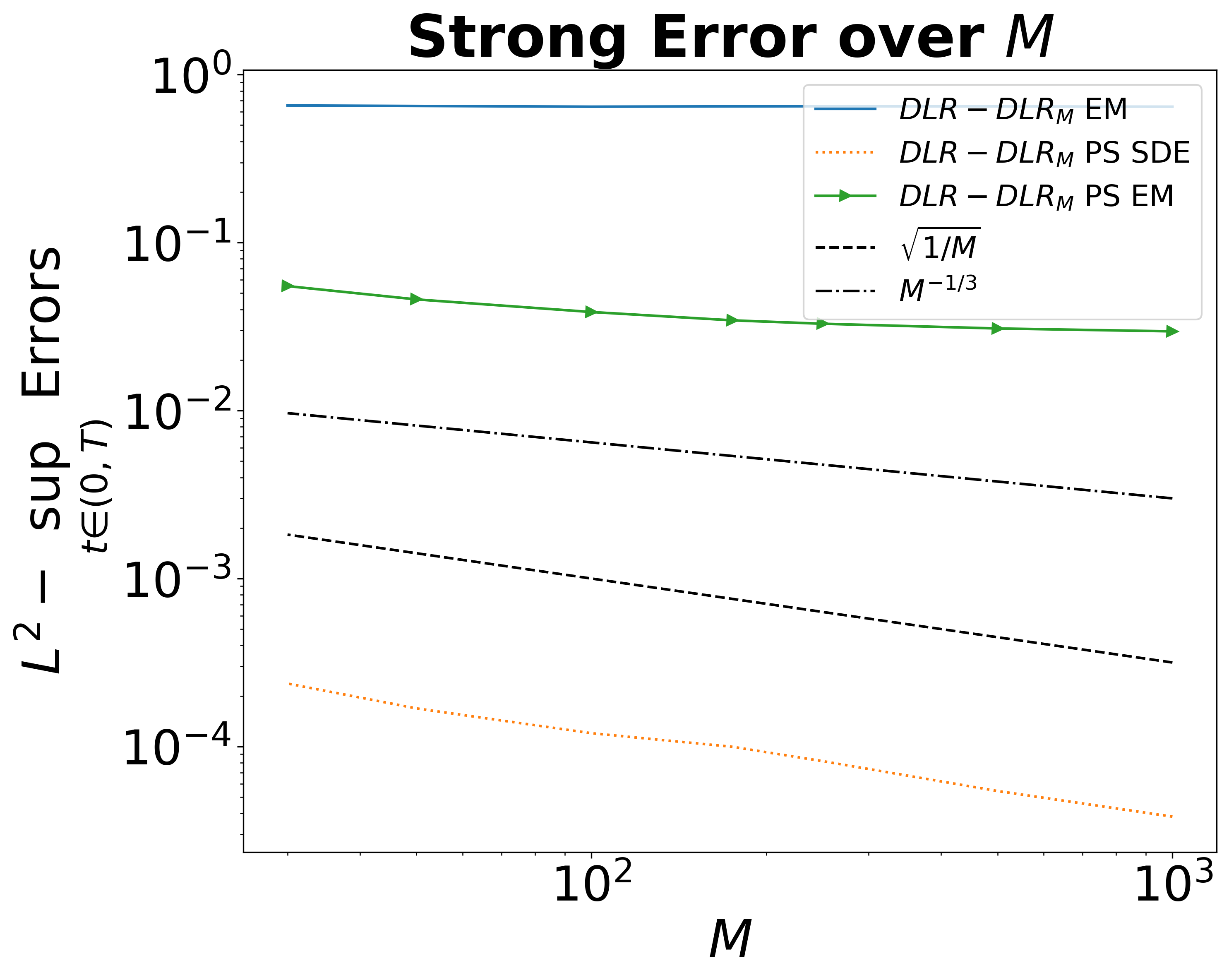}
	\includegraphics[scale=0.29]{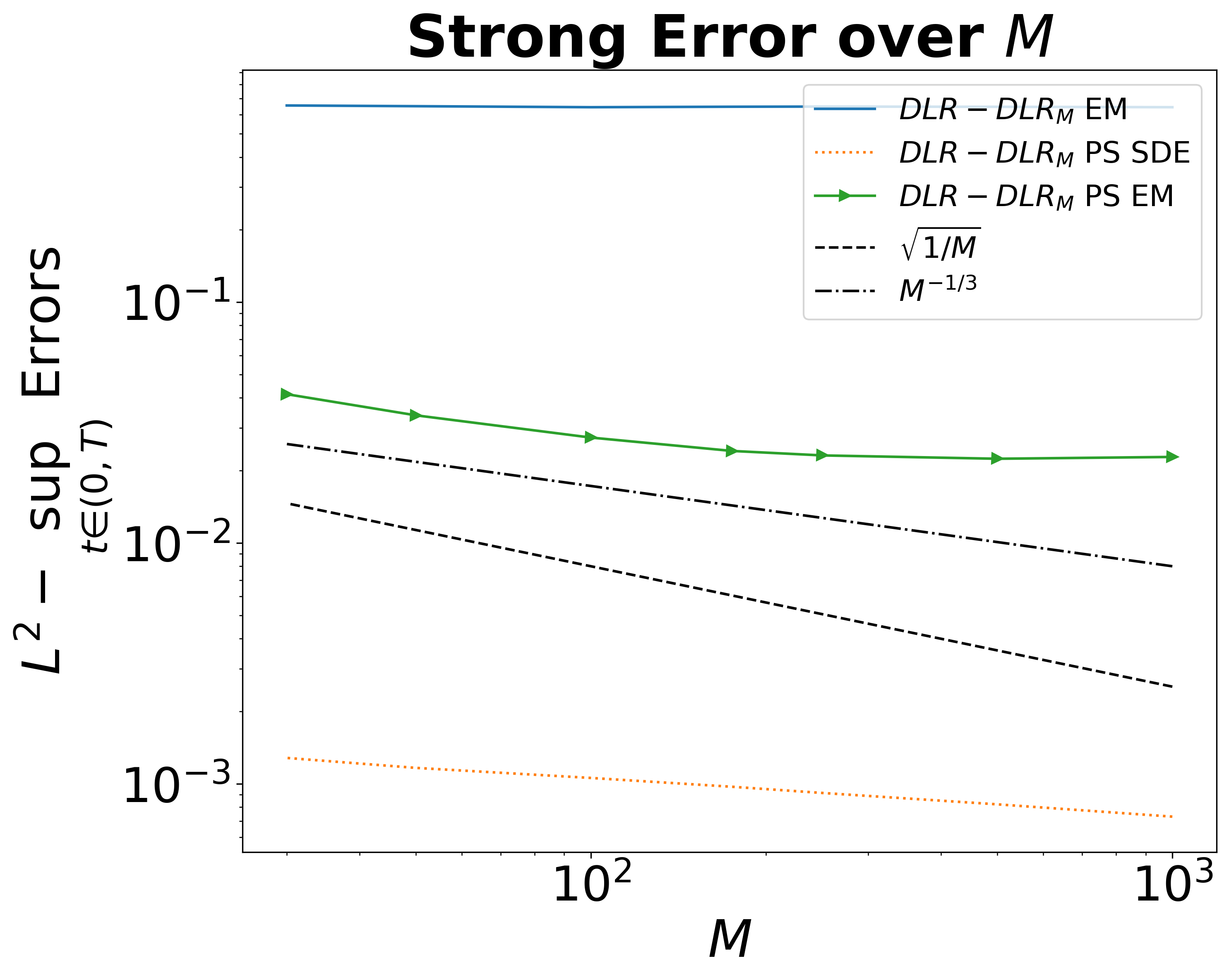}
	\caption{$L^2$-sup-in-time for $t \in [0,1]$ between the semidiscrete $DLR$ of the Problem \eqref{ex: SADR} and (left) un-regularized DLR algorithms, i.e.\ $\alpha=0$, and (right) their regularized version with $\alpha = \frac{1}{\sqrt[3]{M}}$, $k=5$, $\Delta t= 10^{-4}$, and reference DLR solution computed with $M=1\cdot 10^5$ paths.}
	\label{fig:SADR reg vs no reg errors}
\end{figure}

In Figure \ref{fig:SADR reg vs no reg errors} we plot the strong error of the particle-system algorithms versus the semi-discrete DLRA. Again, when we do not employ a regularization of the Gramian, DLR PS SDE achieves the standard Monte Carlo rate, whereas in the regularized case, the convergence rate reduces to $O(\frac{1}{\sqrt[3]{M}})$. Further, as noted in Part I \cite{kazashi2025dynamicalpartI}, DLR PS EM and DLR EM are less accurate than DLR PS SDE when compared against the true DLRA; for DLR EM, this is likely because time discretization is the dominant source of error.

\begin{figure}[!h]
	\centering
	\includegraphics[scale=0.11]{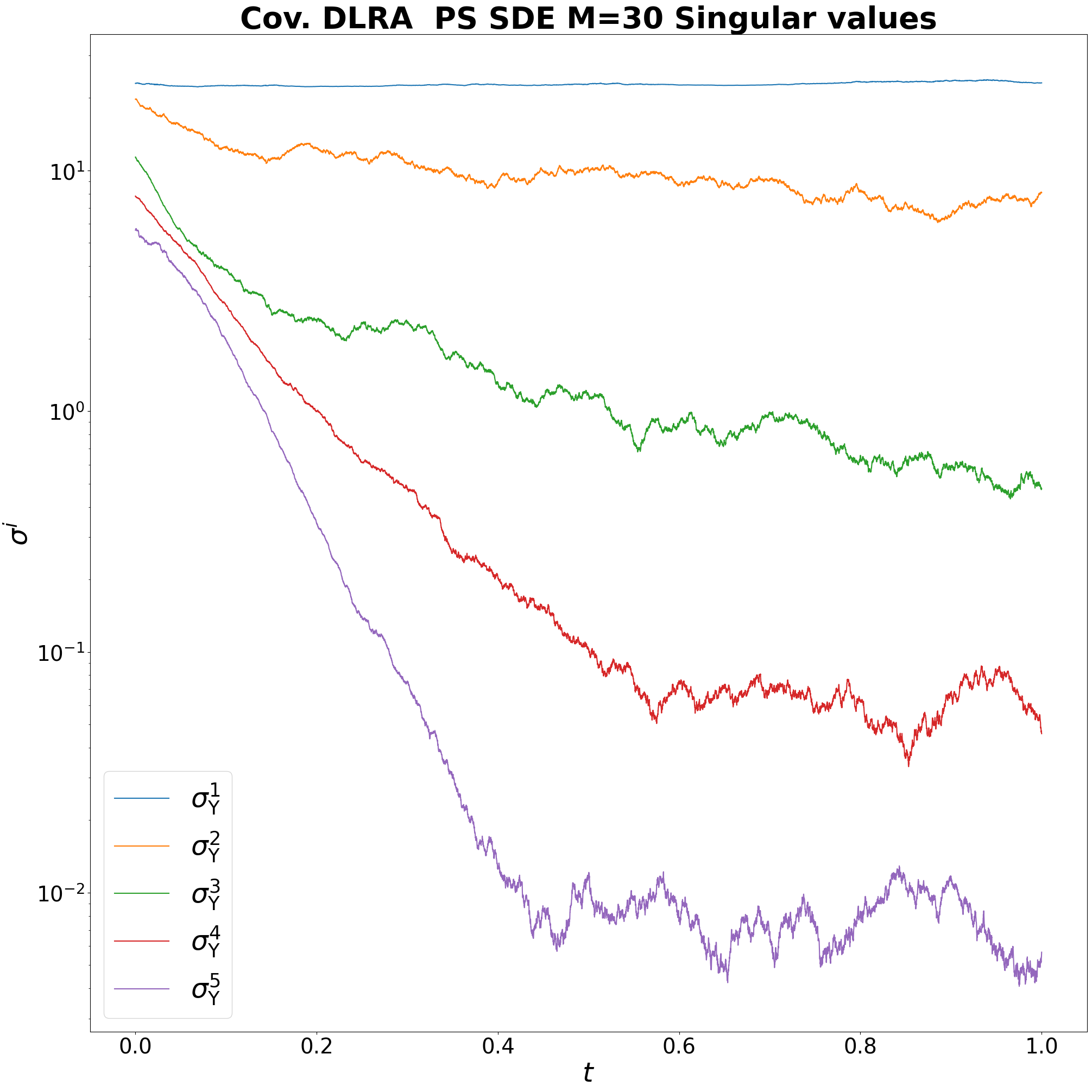}
	\includegraphics[scale=0.11]{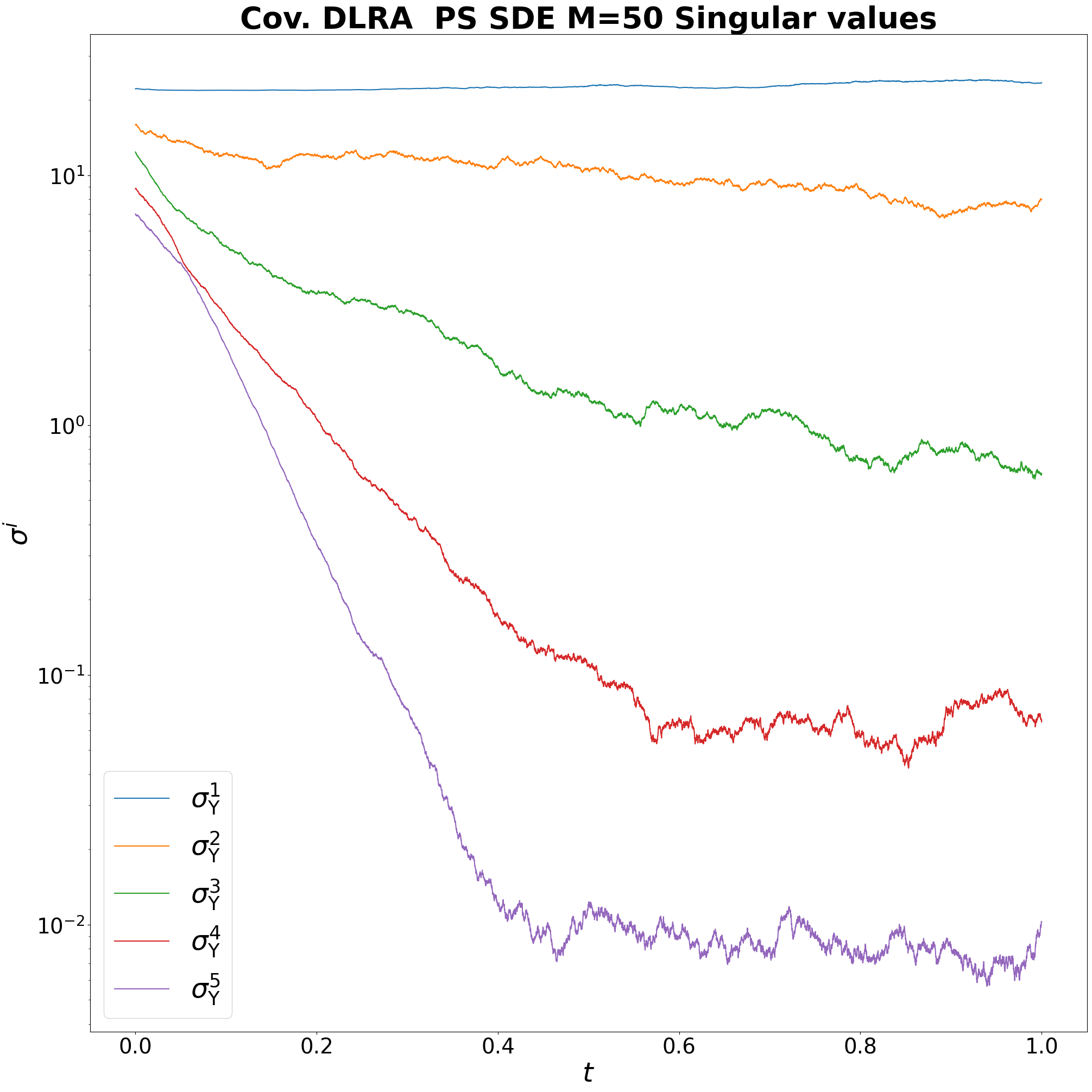}
	\caption{Singular values of $\mathbb{E}[Y_n Y_n^{\top}]$ of DLR PS SDE (Algorithm \ref{alg: Stoc Proj algorithm}) applied to Problem \eqref{ex: SADR} for $M=30$ (left) and $M=50$ (right). $k=5$, $\Delta t=  10^{-4}$. }
	\label{fig: sing val SADR example}
\end{figure}

In Figure \ref{fig: sing val SADR example} we show the evolution over time of the singular values of the Gramian for the stochastic basis of the DLR PS SDE for $M=30$ and $M=50$, respectively.

\begin{figure}[!h]
	\centering
	\includegraphics[scale=0.22]{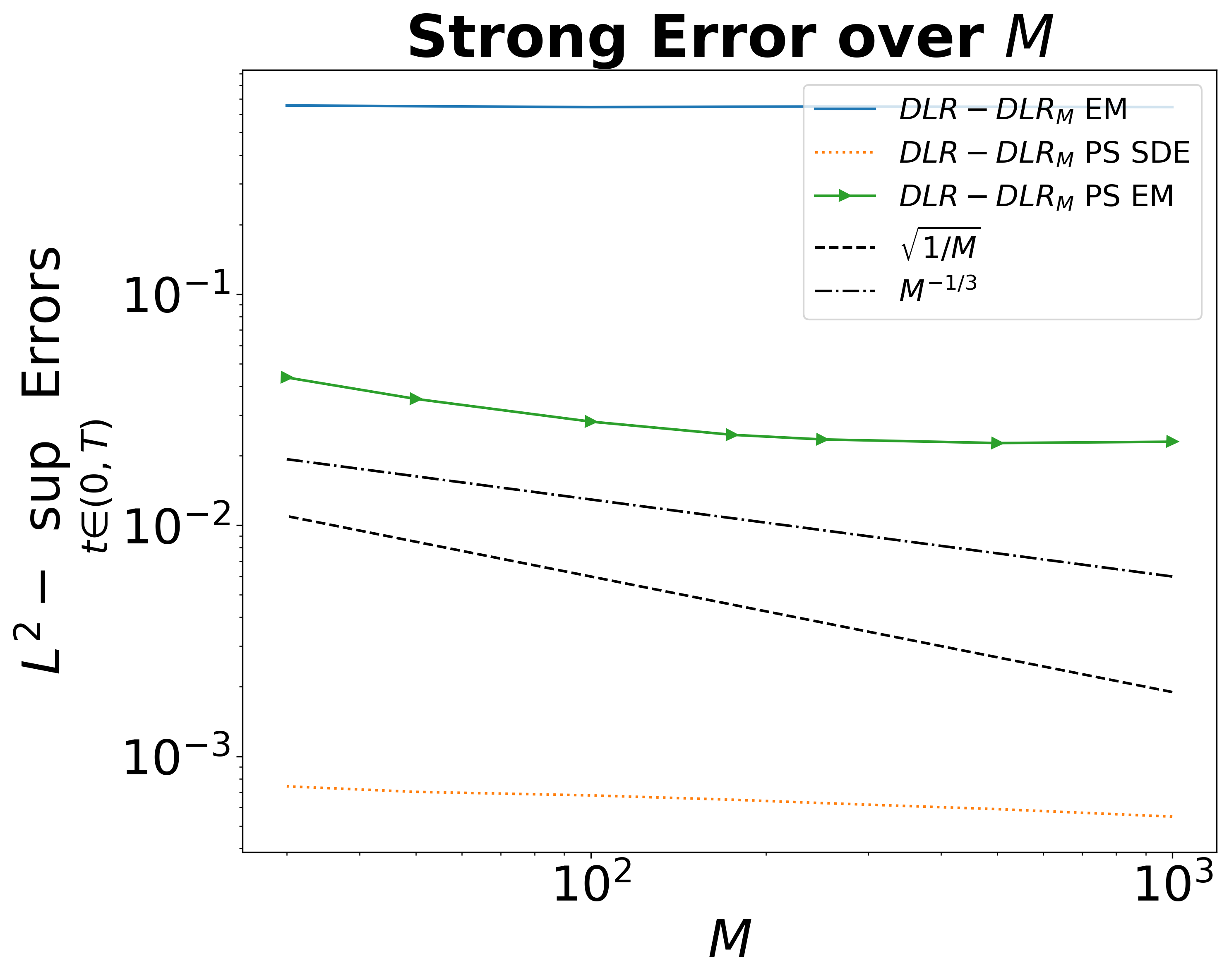}
	\includegraphics[scale=0.22]{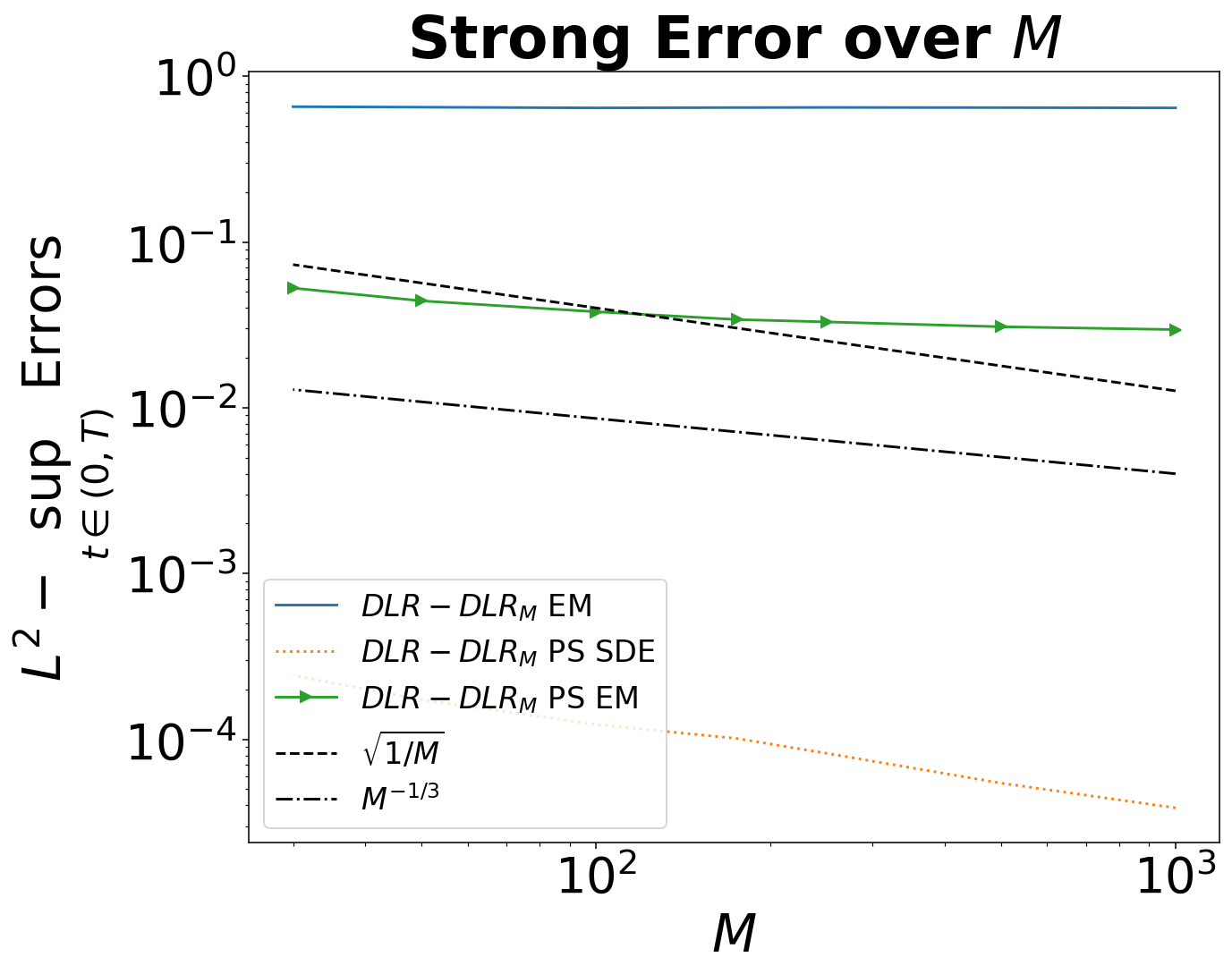}
	\includegraphics[scale=0.22]{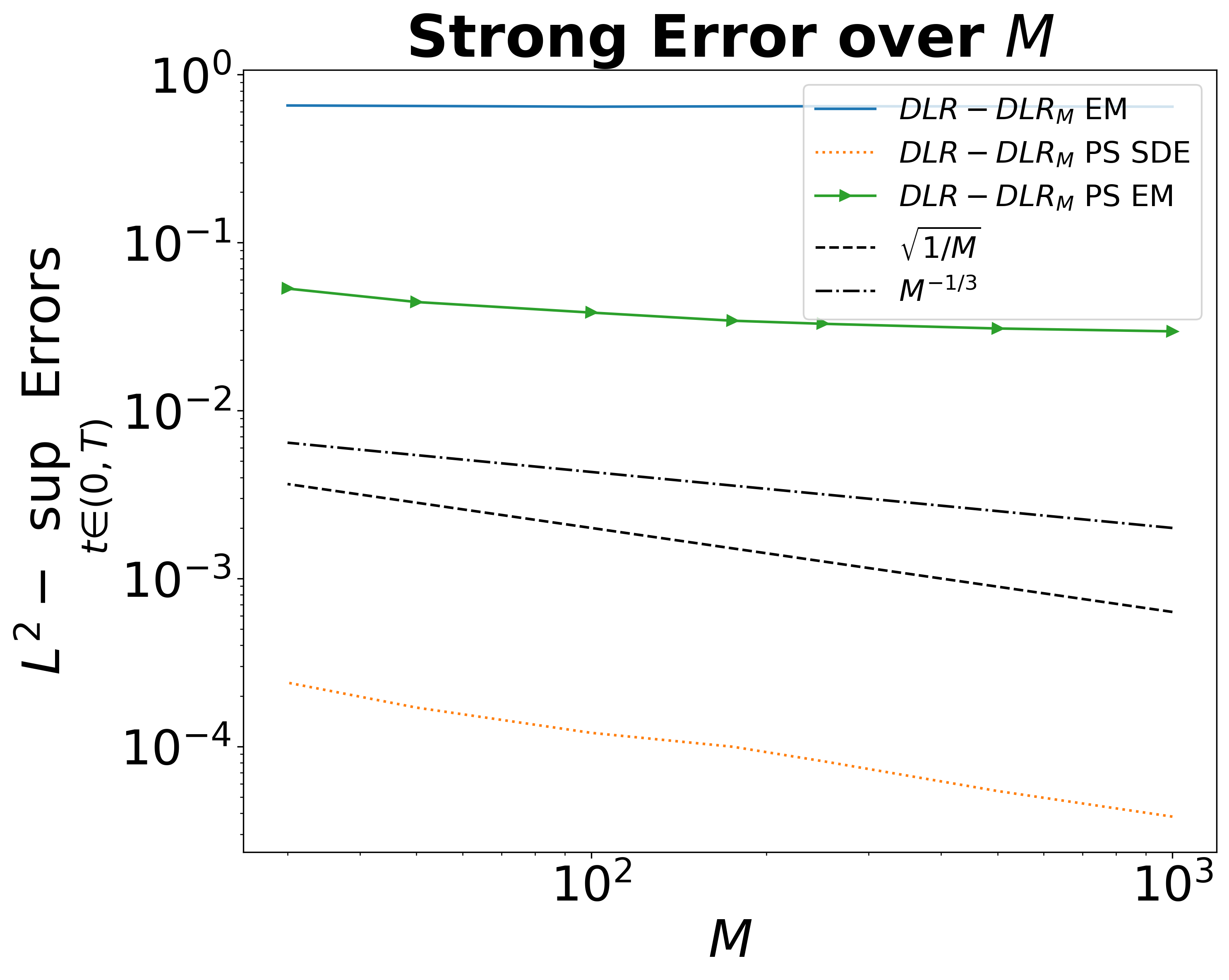}
	\caption{$L^2$-sup-in-time errors for $t \in [0,1]$ between the semidiscrete $DLR$ of the Problem \eqref{ex: SADR} and the DLR algorithms implemented using Algorithm \ref{alg: mod Stoc Proj algorithm}, for quantity $\eta=0.08$ (left), $\eta=0.0085$ (center),  $\eta=0.0085$ (right). In all the cases: $k=5$, $\Delta t= 10^{-4}$, and reference $DLR$ solution computed with $M=10^5$ paths.}
	\label{fig:SADR stoch alg example errors}
\end{figure}

In Figure \ref{fig:SADR stoch alg example errors}, we show the errors at final time $T=10$ between the modified stochastic version of the DLRA algorithms versus the continuous DLRA, varying the number of samples $M$ in the Monte Carlo discretization.
According to the trend of the singular values of $\mathbb{E}[Y_n Y_n^{\top}]$ presented in Figure \ref{fig: sing val SADR example}, we propose to implement Algorithm \ref{alg: mod Stoc Proj algorithm} with different thresholds $\eta$ for the regularization of the Gramian. 

Figure \ref{fig:SADR stoch alg example errors} illustrates that in this case the convergence with respect to the number of samples $M$ is between $O(\frac{1}{\sqrt{M}})$ and $O(\frac{1}{\sqrt[3]{M}})$ for the DLR PS SDE, a better result than the one where we employed the regularization procedure at each time step (see Figure \ref{fig:SADR reg vs no reg errors}). It can be observed that increasing the regularization parameter $\eta$ reduces the accuracy of the DLR PS EM method. In contrast, for smaller values of $\eta$, the convergence rate more closely approaches the theoretical rate of $O\left(\frac{1}{\sqrt{M}}\right)$ established in Theorem \ref{thm: convergence of the mod Stoc Proj}.

\subsection{Projected SDE example - diffusion regularization}
In this simulation, we show numerical performance concerning the convergence under diffusion-regularized algorithm presented in Section \ref{sec: diff reg}. We consider a slightly different SDE to the one employed in Section \ref{sec: proj sde}, described as follows
\begin{equation}\label{ex: toy example stoc proj sde - diff reg}
	\begin{aligned}
		\hspace{-0.3cm}
		\begin{pmatrix}
			\mathrm{d}X_1^{\mathrm{true}} \\
			\mathrm{d}X_2^{\mathrm{true}} \\
			\mathrm{d}X_3^{\mathrm{true}} \\
		\end{pmatrix}
		=& V(t)^{\top} V(t)
		\begin{pmatrix}
			\lambda \sin(X_1^{\mathrm{true}}(t)) & - \lambda \sin( X_2^{\mathrm{true}}(t))& - \theta X_3^{\mathrm{true}}(t) \\
			\lambda X_1^{\mathrm{true}}(t) & -\lambda X_2^{\mathrm{true}}(t)& - \theta X_3^{\mathrm{true}}(t) \\
			\lambda X_1^{\mathrm{true}}(t) & -\lambda X_2^{\mathrm{true}}(t)& - \theta X_3^{\mathrm{true}}(t) \\
		\end{pmatrix}
		\mathrm{d}t\\
		&+  (\sqrt{1 + 2 |X^{\mathrm{true}}(t)|})  V(t)^{\top} V(t)
		\begin{pmatrix}
			\mathrm{d}W_1(t) \\
			\mathrm{d}W_2(t) \\
			\mathrm{d}W_3(t) \\
		\end{pmatrix}
	\end{aligned}
\end{equation}
where $V(t) \in \mathbb{R}^{2 \times 3}$ is a time-dependent orthogonal matrix defined as in \eqref{eq: V(t) subspace}
for $t \in [0,0.5]$, and $\lambda=0.6$, $\theta=0.5$. Therefore, the diffusion matrix has actually rank $k=2$. The initial datum is given by the random vector $X^{\mathrm{true}}(0) = V(0)^{\top} V(0)\tilde{X}$, where we define component $\tilde{X}_i = \text{Un}_{i}(1,4)$, for $i=1,2$, where $\text{Un}_{1}(1,4)$ and $\text{Un}_{2}(1,4)$ are independent uniform random variables in the interval $[1,4]$. Also in this example, we employ a DLR approximation of rank $k=2$. The true DLR solution is computed with $M=10^5$ paths and time step $\Delta t = 1 \cdot 10^{-4}$. For all the Monte Carlo approximations we implemented, we consider the same step-size $\Delta t$. Therefore, we expect to see the Monte Carlo rate of convergence if the other sources of error, such as the time discretization, are negligeable. We run $M_{\mathrm{runs}}=25$ independent simulations and we average out all the error quantities.

\begin{figure}[!h]
	\centering
	\includegraphics[scale=0.29]{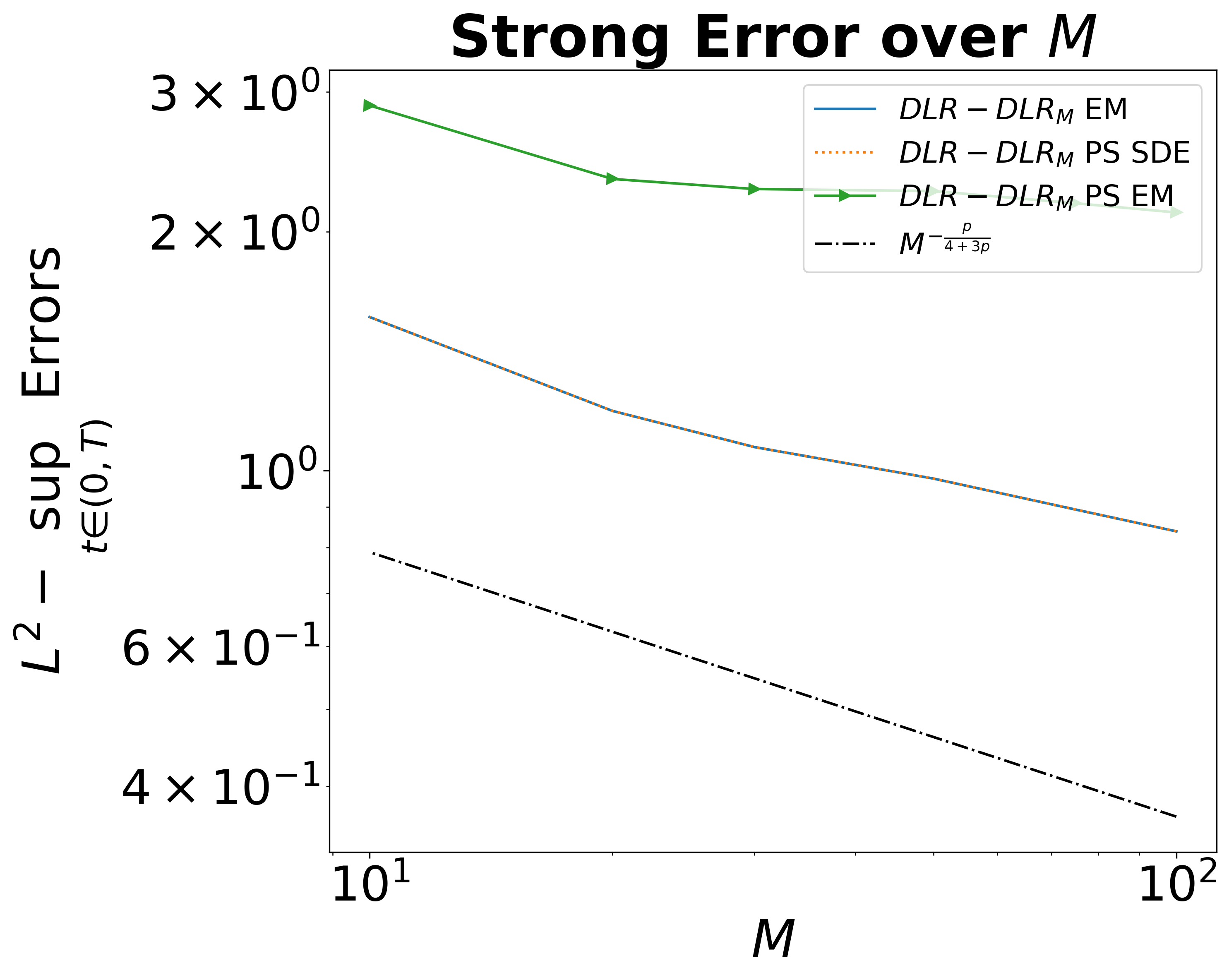}
	\includegraphics[scale=0.29]{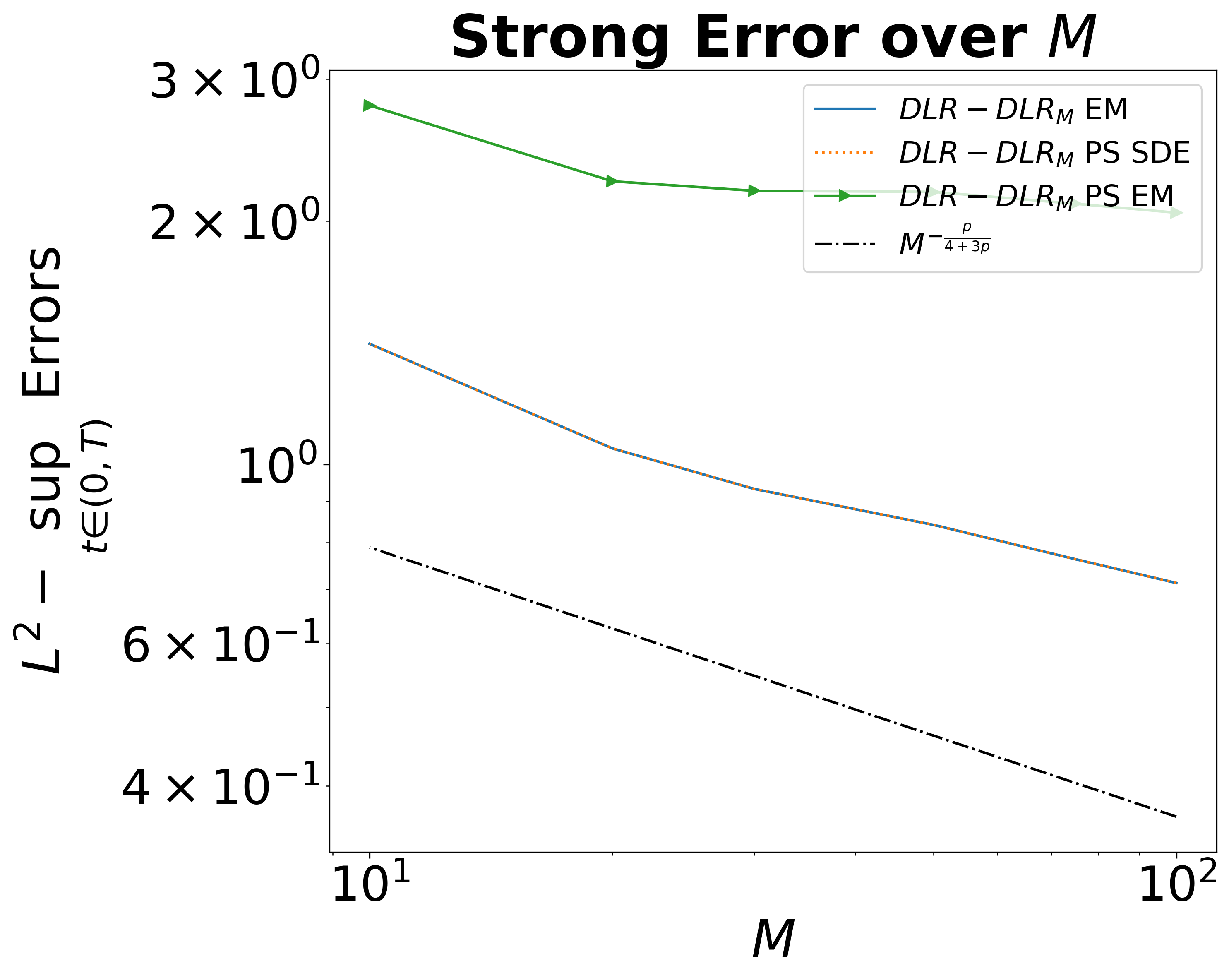}
	\caption{ $L^2$-sup-in-time errors for the DLR EM, DLR PS EM, and DLR PS SDE with respect to the semidiscrete $DLR$ solution at time $T=0.5$ (computed with $M=10^5$) with diffusion regularization parameter $\beta = M^{ -\frac{p}{4+3p}}$ with $p=8$, corresponding to $\beta=M^{-\frac{2}{7}}$ (left) and with $p=\infty$, corresponding to $\beta=M^{-\frac{1}{3}}$ (right) for Problem \eqref{ex: toy example stoc proj sde - diff reg}. In both cases: $\Delta t=0.001$, $k=2$.}
	\label{fig:toy example stoc proj sde - diff reg}
\end{figure}

Again, to initiate all the DLRA algorithms, as $X^{\mathrm{true}}(0)=X(0)$ is already of rank equal to $2$, we consider a \texttt{QR} decomposition of the discretized initial datum. Precisely, the initial condition of Algorithms \ref{alg: DLR EM SDE algorithm}, \ref{alg: Stoc Proj algorithm}, and \ref{alg: Eva Proj Splitt SDE algorithm} is the rank-$k$ best approximation with respect to the Frobenius norm of $[(X^{\mathrm{true}})^{1}(0), \dots, (X^{\mathrm{true}})^{M}(0)] \in \mathbb{R}^{d \times M}$, i.e.\ the (rank-$k$) SVD of the matrix whose columns are the $M$ realizations of $X_0$. 
Figure \ref{fig:toy example stoc proj sde - diff reg} shows convergence plots with respect to the number of samples $M$ for two different choices of the regularization parameter $\beta= M^{-\frac{p}{4+3p}}$, namely $p=8$, corresponding to $\beta=M^{-\frac{2}{7}}$, and of $p=\infty$, corresponding to $\beta=M^{-\frac{1}{3}}$.

\section{Conclusion}
\label{sec: conclusion}
In this second part, we provided a complete error analysis for the stochastic discretization of the three algorithms proposed in \cite{kazashi2025dynamicalpartI}. We employed a Monte Carlo method with $M$ samples to discretize the stochastic space. In the time discretized DO equations, resulting in a noisy interacting particle systems. The coupling that characterizes these algorithms presents several issues and, hence, its analysis is not standard. For instance, the components of the stochastic basis are not independent and, hence, probability bounds of the smallest singular value of its Gramian, which are essential for stability purposes, are challenging to prove. Moreover, the deterministic modes depend on all the paths, implying a not standard analysis in the computation of the error between the employed projectors. 

First, we achieved to prove convergence of the fully discretized DLR PS SDE to the semi-discretized version in the case of elliptic diffusion and subgaussian tails of the initial condition. In order to avoid possible ill-conditioning, we first proposed to regularize the Gramian. With the optimal choice of the regularization parameter with respect to the number of samples, we could prove a Monte Carlo convergence rate close to $O(\frac{1}{\sqrt[3]{M}})$ for large $M$. On the other hand, if we apply a regularization only when the smallest singular value is below a fixed threshold, we could prove a better rate close to $O(\frac{1}{\sqrt{M}})$ for large $M$ if high moments of the initial condition are available, which is the standard convergence rate for Monte Carlo. Lastly, we proposed a modification of the algorithm to deal with the case of degenerate diffusion and established convergence rate dependent on the number of available moments of the initial conditions.

The other two algorithms proposed in \cite{kazashi2025dynamicalpartI}, namely the DLR EM and the DLR PS EM, present several challenges and proving convergence of their stochastic discretization is still an open question. For the former, difficulties araise when proving the boundedness of the second moment of the numerical solution under linear-growth bound, unlike the time-discretization only setting. Concerning the latter, the presence of the diffusion term inside the expectation in the DO equations introduces a stochastic discrepance which is difficult to tackle. Future works should investigate these two methods more in-depth.

It would be desirable to obtain improved bounds on the discretized Gramian, both in probability and in expectation, as well as projector-type bounds with weaker assumptions than subgaussian tails of the semidiscrete DLRA. These estimates might allow us to improve the convergence rate, also in case of regularization. Furthermore, possible benefits could be found when approximations of the Gramian are taken from other similar numerical linear algebra algorithms, such as randomized SVD \cite{halko2011finding}. This might lead to a possible computational saving (but still remaining dependent on the smallest singular value in the error estimates though).

\section*{Acknowledgements}
FN and FZ were supported by the Swiss National Science Foundation under the
Project n. 200518 “Dynamical low rank methods for uncertainty quantification and data assimilation”.

\printbibliography

\appendix

\section{Proofs of $L^p$-norm bounds}\label{app: proof norm bounds}

\begin{proof}[\bfseries Proof of Lemma \ref{lem: L2p norm semidiscretized solution}] \label{proof: lem: L2 norm semidiscretized solution}
	The proof of the bound of the second moment can be found in \cite[Lemma 6.1]{kazashi2025dynamicalpartI}.
	
	Concerning the statement for $p>1$, the proof follows by bounding the increments of the update $X_{n+1}$ in a reasonable way to conclude by Gronwall's lemma. Consider the update of the DLR PS SDE given by \eqref{eq: Stoc Proj X}. Then, by means of\ orthogonality of the rows of $U_n$, one has $|P_{U_n}|=1$ and, by Young's inequality 
	\begin{equation}\label{eq: interm proj semidiscre X}
		\begin{aligned}
			|X_{n+1}|^{2p} = & \left(|X_n + P_{{U}^{\top}_n\widetilde{Y}_{n+1}}[ a_n \Delta t_n] + P_{{U}_n} [b_n] \Delta W_n| \right)^{2p} \\
			\leq & \left(|X_n| + |P_{{U}^{\top}_n\widetilde{Y}_{n+1}}[ a_n]| \Delta t_n + |P_{{U}_n} [b_n] \Delta W_n| \right)^{2p} \\
			\leq & \resizebox{1\linewidth}{!}{$ \left(|X_n|^2 + |P_{\widetilde{Y}_{n+1}}[ a_n]|^{2} \Delta t_n^2 + |P_{U_{n}}[ a_n]|^{2} \Delta t_n^2 + |b_n \Delta W_n|^2 + 2X_n^{\top} a_n \Delta t_n + 2 a_n^{\top}P_{U_{n}}b_n \Delta W_n \Delta t_n +2X_n^{\top}b_n \Delta W_n \right)^{p} $}  \\
				\leq & \resizebox{1\linewidth}{!}{$ \left(|X_n|^2 + |P_{\widetilde{Y}_{n+1}}[ a_n]|^{2} \Delta t_n^2 + |P_{U_{n}}[ a_n]|^{2} \Delta t_n^2 + 2|b_n \Delta W_n|^2 + |X_n|^2\Delta t_n +|a_n|^2 \Delta t_n + |a_n|^2 \Delta t_n^2   +2X_n^{\top}b_n \Delta W_n \right)^{p} $}  \\
			\leq & \resizebox{1\linewidth}{!}{$ \left(|X_n|^2 + |P_{\widetilde{Y}_{n+1}}[ a_n]|^{2} \Delta t_n^2 + C_{\mathrm{lgb}}(1 + |X_n|^2) \left(2 \Delta t_n^2 + 2| \Delta W_n|^2 + \Delta t_n \right) + |X_n|^2\Delta t_n +2X_n^{\top}b_n \Delta W_n \right)^{p} $} .\\
		\end{aligned}
	\end{equation}
	In \eqref{eq: interm proj semidiscre X}, the term requiring some care is $|P_{\widetilde{Y}_{n+1}}[ a_n]|^{2}$, as we cannot find a reasonable $L^p$-norm bound to it without involving the smallest singular value of the semi-discrete Gramian. Indeed, it holds that
	\begin{equation}\label{eq: proj Y delta t bound}
		\begin{aligned}
			|P_{\widetilde{Y}_{n+1}}[ a_n]|^{2} \Delta t_n^{2} =& \left|\mathbb{E}[a_n \widetilde{Y}_{n+1}^{\top}] C_{\widetilde{Y}_{n+1}}^{-1}\widetilde{Y}_{n+1}\right|^2 \Delta t_n^2\\
			\leq& \left|\mathbb{E}[a_n \widetilde{Y}_{n+1}^{\top}C_{\widetilde{Y}_{n+1}}^{-\frac{1}{2}}]\right|^2\left| C_{\widetilde{Y}_{n+1}}^{-\frac{1}{2}}\widetilde{Y}_{n+1}\right|^2 \Delta t_n^2 \\
			\leq & \mathbb{E}[|a_n|^2] \left| C_{\widetilde{Y}_{n+1}}^{-\frac{1}{2}}\right|^2 \left|Y_n + U_n a_n \Delta t_n + U_n b_n \Delta W_n \right|^2 \Delta t_n^2 \\
			\leq & C_{\mathrm{lgb}}\left(1 + K_2(T)\right) (\widetilde{\sigma}_{n+1}^k)^{-1} \left|Y_n + U_n a_n \Delta t_n + U_n b_n \Delta W_n \right|^2 \Delta t_n^2\\
			\leq & 3C_{\mathrm{lgb}}\left(1 + K_2(T)\right) (\widetilde{\sigma}_{n+1}^k)^{-1} \left(|Y_n|^2 + C_{\mathrm{lgb}}\left(1 + |X_n|^2\right) \left[(\Delta t_n)^2 + |\Delta W_n|^2\right] \right)\Delta t_n^2 \\
			= & 3C_{\mathrm{lgb}}\left(1 + K_2(T)\right) (\widetilde{\sigma}_{n+1}^k)^{-1} \left(|X_n|^2 + C_{\mathrm{lgb}}\left(1 + |X_n|^2\right)\left[(\Delta t_n)^2 + |\Delta W_n|^2\right]\right) \Delta t_n^2 \\
			= & 3C_{\mathrm{lgb}}\left(1 + K_2(T)\right) (\widetilde{\sigma}_{n+1}^k)^{-1} \Delta t_n \left(|X_n|^2 + C_{\mathrm{lgb}}\left(1 + |X_n|^2\right)\left[(\Delta t_n)^2 + |\Delta W_n|^2\right] \right) \Delta t_n\\
			\leq &\left( |X_n|^2 + C_{\mathrm{lgb}}\left(1 + |X_n|^2\right)\left[(\Delta t_n)^2 + |\Delta W_n|^2\right]\right) \Delta t_n
		\end{aligned}
	\end{equation}
	where we employed the linear-growth bound, the first statement of this lemma, Cauchy-Schwarz inequality, and the orthogonality of the rows of $U_n$, whereas in the last line we used condition \eqref{eq: dt semidiscrete condition}.  Notice that condition \eqref{eq: dt semidiscrete condition} allows us to obtain an upper bound of $|P_{\widetilde{Y}_{n+1}}[ a_n]|^{2}$ that does not depend on the smallest singular value $\widetilde{\sigma}_{n+1}^k$.
	
	Then, we apply the multinomial identity \cite[Section 1.2]{stanley2011enumerative}
	\begin{equation}\label{eq: multinomial relation}
		(x_1+x_2+\dots+x_d)^{p}
		=
		\sum\limits_{q_1+q_2+\dots+q_d=p}
		\frac{p!}{q_1!q_2!\dots q_d!}
		x_1^{q_1}x_2^{q_2}\dots x_d^{q_d},
	\end{equation}
	to the finite sum appearing in \eqref{eq: interm proj semidiscre X}. We then
	use \eqref{eq: proj Y delta t bound} to control the contributions containing
	\(
		\left(|P_{\widetilde{Y}_{n+1}}[a_n]|^2\Delta t_n^2\right)^{k} .
	\)
	for $k=1,\dots,M$. We
	then take conditional expectation with respect to \(\mathcal F_{t_n}\). Terms
	with an odd power of the increment \(\Delta W_n\) vanish whereas the remaining Brownian
	terms are estimated using
	\[
		\mathbb E_n[|\Delta W_n|^q]\leq C_q\Delta t_n^{q/2},
		\qquad q\geq1 .
	\]
	We then take expectation and use the tower property.
	Together with the linear-growth bound, this gives, for some \(C=C(p,T)>0\),
	\begin{equation*}
		\begin{aligned}
			\mathbb{E}[|X_{n+1}|^{2p}]
			\leq & \mathbb{E}[|X_{n}|^{2p}]	 + C ( 1 + \mathbb{E}[|X_{n}|^{2p}]) \Delta t_n,
		\end{aligned}
	\end{equation*}
	The desired bound follows from the discrete Gronwall lemma.
	
	Concerning the second part of the claim, one has that
	\begin{equation}\label{eq: X to Y 2p}
		|X_{n}|^{2p} = (X_{n}^{\top}X_{n})^{p} =  (\widetilde{Y}_{n}^{\top}\widetilde{U}_n\widetilde{U}_n^{\top}\widetilde{Y}_{n})^{p} \geq (\widetilde{Y}_{n}^{\top}\widetilde{Y}_{n})^{p} = |\widetilde{Y}_{n}|^{2p},
	\end{equation}
	where we use that fact that $\widetilde{U}_n\widetilde{U}_n^{\top}$ is symmetric positive definite for all $n$ (see \cite[Relation (27)]{kazashi2025dynamicalpartI}).
\end{proof}

\begin{proof}[\bfseries Proof of Lemma \ref{lem: Lp moment bound}] \label{proof: lem: Lp moment bound}

    To prove the right-hand-side inequality in \eqref{eq: Lp moment bound X_hat Y_hat}, we first derive suitable $p$-moment bounds for the terms appearing in the update of $\widehat{X}_{n+1}^i$. We then combine these estimates and use the time-step condition \eqref{eq: M dt cond} to apply Gronwall's lemma.
	
	Let us start from the analysis of the projectors. Consider $f=(f^1,\dots,f^M) \in L^2(\Omega^M,\mathbb{R}^{d \times M})$ with $f^i \in\mathbb{R}^d$ for $i=1,\dots,M$, and  
	define the matrix
	$$
	\widetilde{\widehat{\mathbb{Y}}}_{n+1}=\begin{bmatrix} \tildhatY{1}{n+1} & \tildhatY{2}{n+1} & \cdots & \tildhatY{M}{n+1} \end{bmatrix}\in\mathbb R^{k\times M}.
	$$
	We recall that the operators $\widehat{P}^{\alpha}_{\tildhatY{i}{n+1}}: (\mathbb{R}^{d \times M},\|\,\cdot\,\|_{\mathrm{F}}) \to (\mathbb{R}^{d},|\,\cdot\,|)$ and $\widehat{P}^{\alpha}_{ \widetilde{\widehat{\mathbb{Y}}}_{n+1}}: (\mathbb{R}^{d \times M},\|\,\cdot\,\|_{\mathrm{F}}) \to (\mathbb{R}^{d \times M},\|\,\cdot\,\|_{\mathrm{F}})$ satisfy
	$$\widehat{P}^{\alpha}_{\tildhatY{i}{n+1}}[f] =\widehat{\mathbb{E}}[f (\widetilde{\widehat{\mathbb{Y}}}_{n+1})^{\top} ] (\widehat{C}_{\widetilde{\widehat{Y}}_{n+1}}^{\alpha})^{-1} (\tildhatY{i}{n+1}),$$
	and 
	\begin{equation}\label{eq: P_alpha}
		\widehat{P}^{\alpha}_{\widetilde{\widehat{\mathbb{Y}}}_{n+1}}[f]:= \frac{1}{M} f\widetilde{\widehat{\mathbb{Y}}}_{n+1}^{\top} \Big(  \sum_{j=1}^{M} \frac{\tildhatY{j}{n+1}(\tildhatY{j}{n+1})^{\top}}{M} + \alpha I_{k\times k}  \Big)^{-1} \widetilde{\widehat{\mathbb{Y}}}_{n+1}  \in \mathbb{R}^{d \times M},
	\end{equation}
	respectively, and both operators depends on all the paths $\boldsymbol{\omega}=(\omega_1,\dots, \omega_M)$.
	The operation $	\widehat{P}^{\alpha}_{\widetilde{\widehat{\mathbb{Y}}}_{n+1}}[f]$ can be represented by (right) multiplying $f$ by the following stochastic ($\boldsymbol{\omega}$-dependent) matrix of size $M \times M$, which, with slight abuse of notation, we still denote by $\widehat{P}^{\alpha}_{\widetilde{\widehat{\mathbb{Y}}}_{n+1}}$:
	\begin{equation}\label{eq: matrix P_alpha}
		\begin{aligned}
	\left[\widehat{P}^{\alpha}_{\widetilde{\widehat{\mathbb{Y}}}_{n+1}}\right]_{ij} &=  \frac{(\tildhatY{i}{n+1})^{\top}}{M} \Big(  \sum_{\ell=1}^{M} \frac{\tildhatY{\ell}{n+1}(\tildhatY{\ell}{n+1})^{\top}}{M} + \alpha I_{k\times k}  \Big)^{-1} \tildhatY{j}{n+1}\\
	&= \frac{(\tildhatY{i}{n+1})^{\top}}{M} \Big(  \frac{\widetilde{\widehat{\mathbb{Y}}}_{n+1}(\widetilde{\widehat{\mathbb{Y}}}_{n+1})^{\top}}{M} + \alpha I_{k\times k}  \Big)^{-1} \tildhatY{j}{n+1}, \quad \text{ for } i,j = 1, \dots, M.
	\end{aligned}
	\end{equation}
	Using the convention $\widehat{a}_n=(\widehat{a}_n^1,\dots,\widehat{a}_n^M) \in \mathbb{R}^{d \times M}$, we notice that for all $i \in \{1,\dots,M\}$ it holds 
	\begin{equation*}
		\begin{aligned}
			\widehat{X}_{n+1}^{i} & = \widehat{X}_n^{i} + \widehat{P}^{\alpha}_{\widehat{U}^{\top}_n\tildhatY{i}{n+1}}[ \widehat{a}_n] \Delta t_n+ P_{\widehat{U}_n}[\widehat{b}_n^{i}] \Delta W_n^{i} \\
			& =  \widehat{X}_n^{i} + \widehat{P}_{\widehat{U}_n}^{\perp}\widehat{P}^{\alpha}_{\tildhatY{i}{n+1}}[ \widehat{a}_n] \Delta t_n +  P_{\widehat{U}_n}[ \widehat{a}_n^{i}] \Delta t_n + P_{\widehat{U}_n}[\widehat{b}_n^{i}] \Delta W_n^{i},
		\end{aligned}
	\end{equation*}
	where the operator $\widehat{P}^{\alpha}_{\widehat{U}^{\top}_n\tildhatY{i}{n+1}}: L^2(\Omega^M,\mathbb{R}^{d \times M}) \to L^2(\Omega^M,\mathbb{R}^{d})$ is defined as the regularized stochastic discretization of \eqref{eq: disc proj} (see \eqref{eq: P_alpha})
	\begin{equation*}
		\begin{aligned}
			\widehat{P}^{\alpha}_{\widehat{U}^{\top}_n\tildhatY{i}{n+1}}[\,f \,] := & \widehat{P}_{\widehat{U}_n}^{\perp}\widehat{P}^{\alpha}_{\tildhatY{i}{n+1}}[\,f \,] +  P_{\widehat{U}_n}[\,f^i \,] \\
			=&\left(I_{d \times d}-\widehat{U}^{\top}_n\widehat{U}_n\right)\widehat{\mathbb{E}}\left[\,f \, \widetilde{\widehat{\mathbb{Y}}}_{n+1}^{\top}\right]\left(\widehat{C}_{\widetilde{\widehat{Y}}_{n+1}}^{\alpha}\right)^{-1}\tildhatY{i}{n+1} + \widehat{U}^{\top}_n\widehat{U}_n[\,f^i \,],
		\end{aligned}
	\end{equation*}
	where $f^i(\boldsymbol{\omega})$ is the $i$-th column of $f(\boldsymbol{\omega})\in \mathbb{R}^{d \times M}$.
	
	By exchangeability of the particles, it holds that
	\begin{equation*}
		\mathbb{E}[|\widehat{X}_{n+1}^{i}|^{2p}] = \mathbb{E}[\frac{1}{M} \sum_{j=1}^{M}|\widehat{X}_{n+1}^{j}|^{2p}], \quad \text{ for } i \in \{1,\dots,M\}.
	\end{equation*}
	We will use the properties of projectors to get a bound of this quantity.
	Indeed, by triangle inequality we have
	\begin{equation*}
		\begin{aligned}
			| \widehat{X}_{n+1}^{i} |^2 
			=& |  \widehat{X}_n^{i} + P_{\widehat{U}_n}^{\perp}\widehat{P}^{\alpha}_{\tildhatY{i}{n+1}}[ \widehat{a}_n] \Delta t_n +  P_{\widehat{U}_n}[ \widehat{a}_n^{i}] \Delta t_n + P_{\widehat{U}_n}[\widehat{b}_n^{i}] \Delta W_n^{i} |^2 \\
			\leq & |  \widehat{X}_n^{i}|^2 + |P_{\widehat{U}_n}^{\perp}\widehat{P}^{\alpha}_{\tildhatY{i}{n+1}}[ \widehat{a}_n]|^2 \Delta t_n^2 +  |P_{\widehat{U}_n}[ \widehat{a}_n^{i}]|^2 \Delta t_n^2 + |P_{\widehat{U}_n}[\widehat{b}_n^{i}] \Delta W_n^{i} |^2 \\
			& + 2 (\widehat{X}_n^{i})^{\top}\widehat{a}_n^{i}\Delta t_n + 2 (\widehat{X}_n^{i})^{\top}\widehat{b}_n^{i}\Delta W_n^{i} + 2 (\widehat{a}_n^{i})^{\top}P_{\widehat{U}_n}\widehat{b}_n^{i}\Delta W_n^{i}\Delta t_n\\
			\leq & |  \widehat{X}_n^{i}|^2 + |\widehat{P}^{\alpha}_{\tildhatY{i}{n+1}}[ \widehat{a}_n]|^2 \Delta t_n^2 +  |\widehat{a}_n^{i}|^2 \Delta t_n^2 + |\widehat{b}_n^{i} \Delta W_n^{i} |^2 \\
			& + 2 (\widehat{X}_n^{i})^{\top}\widehat{a}_n^{i}\Delta t_n + 2 (\widehat{X}_n^{i})^{\top}\widehat{b}_n^{i}\Delta W_n^{i} + 2 (\widehat{a}_n^{i})^{\top}P_{\widehat{U}_n}\widehat{b}_n^{i}\Delta W_n^{i}\Delta t_n,
		\end{aligned}
	\end{equation*}
      where in the last line we employ the unitary norm of orthogonal projectors. Similarly to the proof of Lemma \ref{lem: L2p norm semidiscretized solution}, we use the multinomial relation \eqref{eq: multinomial relation} 
	\begin{equation}\label{eq: norm comb}
		\begin{aligned}
			\mathbb{E}[| \widehat{X}_{n+1}^{i} |^{2p}] 
			\leq &  \mathbb{E}[\Big(  |  \widehat{X}_n^{i}|^2 + |\widehat{P}^{\alpha}_{\tildhatY{i}{n+1}}[ \widehat{a}_n]|^2 \Delta t_n^2 +  |\widehat{a}_n^{i}|^2 \Delta t_n^2 + |\widehat{b}_n^{i} \Delta W_n^{i} |^2 \\
			& \quad + 2 (\widehat{X}_n^{i})^{\top}\widehat{a}_n^{i}\Delta t_n + 2 (\widehat{X}_n^{i})^{\top}\widehat{b}_n^{i}\Delta W_n^{i} + 2 (\widehat{a}_n^{i})^{\top}P_{\widehat{U}_n}\widehat{b}_n^{i}\Delta W_n^{i}\Delta t_n \Big)^{p} ]\\
			= &   \mathbb{E}\left[\left( \sum\limits_{q_1+q_2+\dots+q_7=p} \frac{p!}{q_1!q_2!\dots q_7!}  |  \widehat{X}_n^{i}  |^{2q_1} |  \widehat{P}^{\alpha}_{\tildhatY{i}{n+1}}[ \widehat{a}_n] |^{2q_2}  \Delta t_n^{2 q_2}   \dots ... \right) \right] \\
			= &  \mathbb{E}[|  \widehat{X}_n^{i}  |^{2p}]  + \mathbb{E}\left[\left( \sum\limits_{\substack{q_1+q_2+\dots+q_7=p  \\ q_1 \neq p, \ q_6+ q_7 \text{ even}} } \frac{p!}{q_1!q_2!\dots q_7!}  |  \widehat{X}_n^{i}  |^{2q_1} |  \widehat{P}^{\alpha}_{\tildhatY{i}{n+1}}[ \widehat{a}_n] |^{2q_2}  \Delta t_n^{2q_2}  ... \right) \right],\\
		\end{aligned}
	\end{equation}
	where the power terms whose sum $q_6 +q_7$ which is odd are null, whereas all the other power terms involving $q_6$ and/or $q_7$ are at least of order $O(\Delta t_n)$ due to the moments of Brownian increments. 

	In \eqref{eq: norm comb}, the terms requiring careful treatment are those involving
	\[
	|\widehat{P}^{\alpha}_{\tildhatY{i}{n+1}}[\widehat{a}_n]|^{2q_2}\Delta t_n^{2q_2}.
	\]
	Among them, the term that implies the most restrictive condition with respect to the number of sample $M$ is the one corresponding to $q_2=1$ arising in the cross term $ 	\frac{1}{M} \sum_{i=1}^{M} |  \widehat{X}_n^{i}  |^{2(p-1)} |  \widehat{P}^{\alpha}_{\tildhatY{i}{n+1}}[ \widehat{a}_n] |^{2}  \Delta t_n^{2}$. With this in mind, using Young inequality with exponents $\frac{p}{p-1}$ and $p$ one obtains that
	\begin{equation}
		\begin{aligned}
		&\frac{1}{M} \sum_{i=1}^{M}  |  \widehat{X}_n^{i}  |^{2(p-1)} |  \widehat{P}^{\alpha}_{\tildhatY{i}{n+1}}[ \widehat{a}_n] |^{2}  \Delta t_n^{2} \\
		 \leq & \frac{1}{M} \sum_{i=1}^{M} \frac{p-1}{p} \left(  |  \widehat{X}_n^{i}  |^{2(p-1)}   \Delta t^{\frac{p-1}{p}}  \right)^{\frac{p}{p-1}}  +\frac{1}{M} \sum_{i=1}^{M} \frac{1}{p} |  \widehat{P}^{\alpha}_{\tildhatY{i}{n+1}}[ \widehat{a}_n] |^{2p} (\Delta t ^{2-\frac{p-1}{p}})^{p} \\
		 = & \frac{p-1}{p}  \frac{1}{M} \sum_{i=1}^{M}   |  \widehat{X}_n^{i}  |^{2p}   \Delta t   +  \frac{1}{p}  \frac{1}{M} \sum_{i=1}^{M} |  \widehat{P}^{\alpha}_{\tildhatY{i}{n+1}}[ \widehat{a}_n] |^{2p} \Delta t ^{p+1} \\
				\end{aligned}
	\end{equation}
	We now focus on the term $\frac{1}{M}\sum_{i=1}^M
	|\widehat{P}^{\alpha}_{\tildhatY{i}{n+1}}[ \widehat{a}_n]|^{2p}$.
	The matrix $(\widehat{P}^{\alpha}_{\widetilde{\widehat{Y}}_{n+1}})$ in \eqref{eq: matrix P_alpha} is symmetric and positive semidefinite with eigenvalues
	$\lambda^r(\widehat{P}^{\alpha}_{\widetilde{\widehat{Y}}_{n+1}})= \frac{\lambda^r(\widetilde{\widehat{\mathbb{Y}}}_{n+1}(\widetilde{\widehat{\mathbb{Y}}}_{n+1})^{\top}/M)}{\lambda^r(\widetilde{\widehat{\mathbb{Y}}}_{n+1}(\widetilde{\widehat{\mathbb{Y}}}_{n+1})^{\top}/M) + \alpha} \in [0,1),$ for all $r=1,\dots,k$.
		
    Denote by $v_i$ the $i$-th column of the matrix $\widehat{P}^{\alpha}_{\widetilde{\widehat{Y}}_{n+1}}$ for $i=1,\dots,M$, i.e.\
    	$$v_i := \frac{1}{\sqrt{M}} (\widetilde{\widehat{\mathbb{Y}}}_{n+1})^\top  \Big( \frac{\widetilde{\widehat{\mathbb{Y}}}_{n+1}(\widetilde{\widehat{\mathbb{Y}}}_{n+1})^\top }{M} + \alpha I_{k\times k} \Big)^{-1} \frac{\tildhatY{i}{n+1}}{\sqrt{M}} \in \mathbb{R}^M.$$
    	We are interested in bounding the quantity $\frac{1}{M} \sum_{i=1}^{M} |  \widehat{P}_{\tildhatY{i}{n+1}}[ f ] |^{2p}$.
    	We can write
    	$$f v_i = \sum_{j=1}^{M} (v_i)_j f^j,$$
    	where $(v_i)_j$ is the $j$-th entry of $v_i\in \mathbb{R}^M$ and $f^j \in \mathbb{R}^d$ is the $j$-th column of $f\in \mathbb{R}^{d \times M}$.
    	Then,
    	for $p \ge 1$, via Hölder's inequality  one has that 
    	$$|f v_i| \le \Big( \sum_{j=1}^M |(v_i)_j|^{\frac{2p}{2p-1}} \Big)^{\frac{2p-1}{2p}} \Big( \sum_{j=1}^M |f^j|^{2p} \Big)^{1/(2p)}.$$
    	Raising both sides to the power $2p$ gives
    	$$|f v_i|^{2p} \le \Big( \sum_{j=1}^M |(v_i)_j|^{\frac{2p}{2p-1}} \Big)^{2p-1} \sum_{j=1}^M |f^j|^{2p},$$ 
    	and summing over $i=1,\dots,M$ we get
    	\begin{equation}\label{eq: p bound on vi}
    		\sum_{i=1}^M |f v_i|^{2p} \le \sum_{i=1}^M \Big( \sum_{j=1}^M |(v_i)_j|^{\frac{2p}{2p-1}} \Big)^{2p-1} \sum_{j=1}^M |f^j|^{2p}
    		= \Bigg( \sum_{i=1}^M \Big( \sum_{j=1}^M |(v_i)_j|^{\frac{2p}{2p-1}} \Big)^{2p-1} \Bigg) \sum_{j=1}^M |f^j|^{2p}.
    \end{equation}
	
	Now, via Hölder's inequality with exponents $\frac{2p-1}{p}$ and $\frac{2p-1}{p-1}$ and standard norm inequalities one has the following bound 
		\begin{equation}
			\begin{aligned}
				\Bigg( \sum_{i=1}^M \Big( \sum_{j=1}^M |(v_i)_j|^{\frac{2p}{2p-1}} \Big)^{2p-1} \Bigg) \leq & \Bigg( \sum_{i=1}^M \Big( M^{\frac{p-1}{2p-1}} |v_i|^{\frac{2p}{2p-1}} \Big)^{2p-1} \Bigg)\\
				\leq &  M^{p-1} \Bigg(\sum_{i=1}^M |v_i|^{2p} \Bigg) \\
				\leq &  M^{p-1} \Bigg(\sum_{i=1}^M |v_i|^2 \Bigg)^p = M^{p-1} \Bigg(\| \widehat{P}^{\alpha}_{\widetilde{\widehat{Y}}_{n+1}}]\|_{\mathrm{F}}^2 \Bigg)^p
			\end{aligned}
		\end{equation}
		and, hence, \eqref{eq: p bound on vi} can be further bounded as
		\begin{equation}\label{eq: p bound on vi 2}
			\begin{aligned}
				\sum_{i=1}^M |f v_i|^{2p} \le&  M^{p-1} \Bigg(\| \widehat{P}^{\alpha}_{\widetilde{\widehat{Y}}_{n+1}}\|_{\mathrm{F}}^2 \Bigg)^p \sum_{j=1}^M |f^j|^{2p}\\
				=& M^{p-1} \Bigg( \sum_{r = 1}^k \left(\frac{\lambda^r(\widetilde{\widehat{\mathbb{Y}}}_{n+1}(\widetilde{\widehat{\mathbb{Y}}}_{n+1})^{\top}/M)}{\lambda^r(\widetilde{\widehat{\mathbb{Y}}}_{n+1}(\widetilde{\widehat{\mathbb{Y}}}_{n+1})^{\top}/M) + \alpha}\right)^2 \Bigg)^p \sum_{j=1}^M |f^j|^{2p} \leq M^{p-1} k^{p} \sum_{j=1}^M |f^j|^{2p}
			\end{aligned}
		\end{equation}
		Relation \eqref{eq: p bound on vi 2} implies that
		\begin{equation}\label{eq: f M norm bound}
			\begin{aligned}
				\frac{1}{M}	\sum_{i=1}^M |\widehat{P}^{\alpha}_{\tildhatY{i}{n+1}}[ \widehat{a}_n] |^{2p} (\Delta t_n)^{p+1} \le& M^{p-1} k^{p} \frac{1}{M} \sum_{i=1}^M |\widehat{a}_n^i|^{2p}(\Delta t_n)^{p+1} \leq \frac{1}{M} \sum_{i=1}^M |\widehat{a}_n^i|^{2p}\Delta t_n,
			\end{aligned}
		\end{equation}
		where we have used condition \eqref{eq: M dt cond} in the last inequality. Notice that under this condition, $|\widehat{P}^{\alpha}_{\tildhatY{i}{n+1}}[\widehat{a}_n]|^{2(p+1)}\Delta t_n^{2p}$ has a bound that does not depend on $M$ and is proportional to $O(\Delta t_n)$. The other terms involving  $|\widehat{P}^{\alpha}_{\tildhatY{i}{n+1}}[\widehat{a}_n]|^{2q_2}\Delta t_n^{2q_2}$ are multiplied by at least one factor $\Delta t_n$ and can be treated similarly.
	
	We proceed similarly to the proof of Lemma \ref{lem: L2p norm semidiscretized solution} to conclude by Gronwall's lemma. 
	Summing \eqref{eq: norm comb} over $i=1,\dots,M$, dividing by $M$, using the interchangeability of particles, and using the linear-growth bound together with \eqref{eq: f M norm bound}, we obtain that there exists a positive constant $C=C(p,T)$ such that
	\begin{equation*}
		\begin{aligned}
			\mathbb{E}[|\widehat{X}_{n+1}^{i}|^{2p}]
			\leq & \mathbb{E}[|\widehat{X}_{n}^{i}|^{2p}]	 + C ( 1 + \mathbb{E}[|\widehat{X}_{n}^{i}|^{2p}]) \Delta t_n
			.
		\end{aligned}
	\end{equation*}
The conclusion follows by an application of the discrete Gronwall's lemma.

	On the other hand, to prove the left-hand-side inequality in \eqref{eq: Lp moment bound X_hat Y_hat}, notice that 
	\begin{equation}\label{eq: UnUn psd}
		\begin{aligned}
		\widetilde{\widehat{U}}_n \left(\widetilde{\widehat{U}}_n\right)^{\top} &\succeq I_{k \times k} + \left((\widehat{C}^{\alpha}_{\widetilde{\widehat{Y}}_{n}})^{-1}\widehat{\mathbb{E}}\left[\widetilde{\widehat{\mathbb{Y}}}_{n} \left(\widehat{a}_{n-1}^{\top}\Delta t_{n-1}\right) \right]P_{\widehat{U}_{n-1}}^{\perp} \right) \left((\widehat{C}^{\alpha}_{\widetilde{\widehat{Y}}_{n}})^{-1}\widehat{\mathbb{E}}\left[\widetilde{\widehat{\mathbb{Y}}}_{n}  \left(\widehat{a}_{n-1}^{\top}\Delta t_{n-1}\right) \right]P_{\widehat{U}_{n-1}}^{\perp}\right)^{\top} \\
		&\succeq I_{k \times k},
		\end{aligned}
	\end{equation}
	for all $n \in \{1,\dots,N\}$. Then, a similar relation to \eqref{eq: X to Y 2p} can be derived and the proof follows verbatim.  
\end{proof}

\section{Bound on the smallest singular value of the empirical covariances}\label{app: Gramian bounds}

In this appendix, we collect auxiliary results that provide lower bounds on the smallest singular value of the empirical Gramian of the semi-discrete $\mathbb{Y}_n$ and the fully-discrete $\widehat{\mathbb{Y}}_n$ stochastic basis. While the former case, involving i.i.d.-particles, reverts to a standard discussion on Monte-Carlo estimators, the latter case is more cumbersome. Indeed, the non-independence of the particles $(\widehat{Y}_n^{i})_{i=1,\dots, M}$ does not allow us to use Monte-Carlo estimates straight away. Dealing with the intermediate points $(\tildhatY{i}{n+1})_{i=1,\dots, M}$ via standard martingale inequalities conditioned on $\mathcal{F}_{t_n}$ is not straightforward either. Indeed, the presence of the regularized projector $(\widehat{P}^{\alpha}_{\tildhatY{i}{n+1}})_{i=1,\dots, M}$ which depends on the Brownian increments, does not make these inequalities directly applicable and represent an additional non-trivial difficulty in proving those estimates.

In order to derive our estimates, our strategy is to compare the empirical Gramians with the semidiscretized one, assuming that the latter is always strictly positive definite, which holds true under Assumption \ref{ass: diff} (see \cite{kazashi2025dynamicalpartI}).
	
We start providing an estimate on the probability that the empirical Gramian of the semidiscretized intermediate basis $\widetilde{\mathbb{Y}}_{n+1}$ is greater than a positive constant.
\begin{Lemma}[Probability lower bound on the semidiscretized singular value]\label{lem: prop bound sigma tilde}
	Assume that the initial datum $X_0$ satisfies $\mathbb{E}[|Y_0|^{4p}] < \infty$ for some $p\geq 1$. 
	Suppose that $\Delta t_n$ satisfies \eqref{eq: dt semidiscrete condition}.
	
	Recall that the $k$-th largest singular value of the intermediate point of the semidiscretized empirical Gramian $U_n^{\top}\frac{\sum_{j=1}^{M}\widetilde{Y}_{n+1}^{j}(\widetilde{Y}_{n+1}^{j})^{\top}}{M}U_n$ is $\widetilde{\sigma}^{k}_{n+1,M}$. Let $\vartheta>0$ be a uniform lower bound on $\widetilde{\sigma}^{k}_{n}$ for any $n$ and $k$. Then, there exists a constant $C_1:=C_1(p)$ independent of $\widetilde{\sigma}^{k}_{n+1}$ such that 
	\begin{equation}\label{eq: sigma_tilde_k_M}
		\mathbb{P}\left(\widetilde{\sigma}^{k}_{n+1,M} \leq \frac{\vartheta}{2}\right) \leq C_1 \frac{1}{\vartheta^{2p}}\frac{1}{M^p}.
    \end{equation}
	
	Furthermore, when employing a uniform time step $\Delta t$ satisfying \eqref{eq: dt semidiscrete condition}, we have
	\begin{equation}\label{eq: int sigma_tilde_k_M}
		\mathbb{P}\left(\bigcap_{0 \leq n \leq N-1}  \Big\{ \widetilde{\sigma}^{k}_{n+1,M} > \frac{\vartheta}{2} \Big\} \right) \geq 1- C_1 \frac{1}{\vartheta^{2p}} \frac{T}{\Delta t M^p}.
	\end{equation}
	\begin{proof}
		The proof follows by using the Weyl's inequality to obtain a useful expression for the sought probability dependent on the Monte-Carlo error.
				
		By assumption and via Ostrowski's Theorem \cite[Theorem 4.5.9]{horn2012matrix}, we have that $\mathbb{E}\left[U_n^{\top}\frac{\sum_{j=1}^{M}\widetilde{Y}_{n+1}^{j}(\widetilde{Y}_{n+1}^{j})^{\top}}{M}U_n\right] \succeq  \vartheta$. Via Weyl's inequality \cite[Theorem 4.3.1]{horn2012matrix} and equivalence of the largest singular value of a matrix with its spectral norm, one has that
		\begin{equation*}
			\begin{aligned}
				\widetilde{\sigma}^{k}_{n+1,M} \geq & \sigma^k\left(\mathbb{E}\left[U_n^{\top}\frac{\sum_{j=1}^{M}\widetilde{Y}_{n+1}^{j}(\widetilde{Y}_{n+1}^{j})^{\top}}{M}U_n\right]\right) \\
				&- \left|U_n^{\top}\frac{\sum_{j=1}^{M}\widetilde{Y}_{n+1}^{j}(\widetilde{Y}_{n+1}^{j})^{\top}}{M}U_n- \mathbb{E}\left[U_n^{\top}\frac{\sum_{j=1}^{M}\widetilde{Y}_{n+1}^{j}(\widetilde{Y}_{n+1}^{j})^{\top}}{M}U_n\right]\right| \\
				\geq & \vartheta - \left|\underbrace{\frac{\sum_{j=1}^{M}\widetilde{Y}_{n+1}^{j}(\widetilde{Y}_{n+1}^{j})^{\top}}{M}}_{=:S_{n+1}}- \mathbb{E}\left[\frac{\sum_{j=1}^{M}\widetilde{Y}_{n+1}^{j}(\widetilde{Y}_{n+1}^{j})^{\top}}{M}\right]\right|,
			\end{aligned}
		\end{equation*}
		where in the last line we employ the orthogonality of the rows of $U_n$. From this relation, one has that the event $\{\boldsymbol{\omega} \in \Omega^M : \frac{\vartheta}{2} \geq \widetilde{\sigma}^{k}_{n+1,M}(\boldsymbol{\omega})\}$ is included in $\{\boldsymbol{\omega} \in \Omega^M : \left|S_{n+1}(\boldsymbol{\omega})-\mathbb{E}\left[S_{n+1}(\boldsymbol{\omega})\right]\right| \geq \frac{\vartheta}{2}\}$. This last relation translates in the following inequality
		\begin{equation}\label{eq: low sigma}
			\begin{aligned}
				\mathbb{P}\left(\widetilde{\sigma}^{k}_{n+1,M}  \leq \frac{\vartheta}{2}\right) \leq &	\mathbb{P}\left(\left|S_{n+1}-\mathbb{E}\left[S_{n+1}\right]\right| \geq \frac{\vartheta}{2}\right). \\
			\end{aligned}
		\end{equation}
		Let us recall that $\{\widetilde{Y}_{n+1}^{j}\}_{j=1,\dots,M}$ are i.i.d.~particles and denote the $i$-th coordinate of the vector $\widetilde{Y}_{n+1}^{j}$ by $(\widetilde{Y}_{n+1}^{j})_i$. Then, we employ Markov inequality, $L^p$ convergence of the standard Monte-Carlo estimator \cite[Proposition 9.11]{ledoux1991probability}  as the samples $X^{1}_n,\dots,X^{M}_n$ of the semi-discrete DLRA are independent, and
		Lemma \ref{lem: L2p norm semidiscretized solution}, in order to obtain positive constants $C_{p}$ independent of $\widetilde{\sigma}^{k}_{n+1}$ and $\widetilde{\sigma}^{k}_{n+1,M} $, but dependent on $p$, such that one has
		\begin{equation}\label{eq: MC in k-dim}
			\begin{aligned}
				\mathbb{P}\left(\widetilde{\sigma}^{k}_{n+1,M}  \leq \frac{\vartheta}{2}\right)	\leq &	\frac{2^{2p}}{\vartheta^{2p}} \mathbb{E}\left[\left|S_{n+1}-\mathbb{E}\left[S_{n+1}\right]\right|^{2p}\right] \leq \frac{2^{2p}}{\vartheta^{2p}} \mathbb{E}\left[\left\|S_{n+1}-\mathbb{E}\left[S_{n+1}\right]\right\|^{2p}_{\mathrm{F}}\right] \\
				\leq & \frac{C_{p}}{\vartheta^{2p}} \frac{1}{M^p} \mathbb{E}\left[\left|\widetilde{Y}_{n+1}^{j}\right|^{4p}\right] \leq \frac{C_{p}}{\vartheta^{2p}} \frac{K_{4p}(T)}{M^p},  \\
			\end{aligned}
		\end{equation}
		which implies the sought relation. Furthermore, notice that this probability bound is independent of $n$.
		
		To prove statement \eqref{eq: int sigma_tilde_k_M}, using countable subadditivity of probability measures, it holds that
		\begin{equation*}
			\begin{aligned}
				\mathbb{P}\left(\bigcap_{0 \leq n \leq N-1}  \Big\{ \widetilde{\sigma}^{k}_{n+1,M} > \frac{\vartheta}{2} \Big\} \right)
				& = 1 - \mathbb{P}\left(  \bigcup_{0 \leq n \leq N-1}  \Big\{ \widetilde{\sigma}^{k}_{n+1,M} \leq \frac{\vartheta}{2} \Big\} \right) \\
				& \geq 1 - \sum_{n=0}^{N-1}\mathbb{P}\left( \Big\{ \widetilde{\sigma}^{k}_{n+1,M} \leq \frac{\vartheta}{2} \Big\} \right) \\
				& \geq 1 - N \frac{c_p}{\vartheta^{2p}}\frac{K_{4p}(T)}{M^p}
		\end{aligned}
	\end{equation*}
	and via using relation \eqref{eq: sigma_tilde_k_M} we get the thesis:
	\begin{equation*}
		\begin{aligned}
			\mathbb{P}\left(\bigcap_{0 \leq n \leq N-1}  \Big\{ \widetilde{\sigma}^{k}_{n+1,M} > \frac{\vartheta}{2} \Big\} \right) 
			& \geq 1 - \frac{T}{\Delta t} \frac{1}{\vartheta^{2p}} \frac{C_1}{M^p}.
		\end{aligned}
	\end{equation*}
	\end{proof}
\end{Lemma}

Now we proceed to analyze the empirical Gramian of the fully-discrete stochastic basis  $(\tildhatY{i}{n+1})_{i=1,\dots,M}$. In this discussion, we need to find a new discrete process whose empirical norm is greater than the fully-discrete solution but whose particles are independent conditioned on $\mathcal{F}_{t_n}$. Using this auxiliary process we can exploit convergence of the Monte-Carlo estimator. In addition, a Gronwall-type bound of this process is also needed. In order to prove this result, we need the following auxiliary lemma.

\begin{Lemma}\label{lem: P f < f}
	For any matrix $f = (f^1, \dots, f^M) \in \mathbb{R}^{d \times M}$ and any $\alpha >0$, one has that
	$$\frac{1}{M}\sum_{i = 1}^{M}|\widehat{P}^{\alpha}_{\widehat{U}_{n}^{\top}\tildhatY{i}{n+1}}[f]|^2 \leq \frac1M \sum_{i=1}^M |f^{i}|^2.$$
	\begin{proof}
		Thanks to the orthogonality of $U_n$ we have the following chain of relations
		\begin{equation*}
			\begin{aligned}
				\frac{1}{M}\sum_{i = 1}^{M}|\widehat{P}^{\alpha}_{\widehat{U}_{n}^{\top}\tildhatY{i}{n+1}}[f]|^2 & = \frac{1}{M}\sum_{i = 1}^{M}|P_{\widehat{U}_{n}}^{\perp}\widehat{P}^{\alpha}_{\tildhatY{i}{n+1}}[f]+ P_{\widehat{U}_{n}}f^i|^2 \\
				& = \frac{1}{M}\sum_{i = 1}^{M} \left(|P_{\widehat{U}_{n}}^{\perp}\widehat{P}^{\alpha}_{\tildhatY{i}{n+1}}[f] |^2+ |P_{\widehat{U}_{n}}f ^i|^2\right) \\
				& = \frac{1}{M}\sum_{i = 1}^{M} \left(| \frac{1}{M}(P_{\widehat{U}_{n}}^{\perp} f) \widetilde{\widehat{\mathbb{Y}}}_{n+1}^{\top} \left( \widehat{C}^{\alpha}_{\widetilde{\widehat{Y}}_{n+1}} \right)^{-1} \tildhatY{i}{n+1}|^2+ |P_{\widehat{U}_{n}}f ^i|^2\right). \\			
			\end{aligned}
		\end{equation*}
		Via the linearity and the cyclic property of the trace, we obtain that
		\begin{equation*}
			\begin{aligned}
				&\frac{1}{M}\sum_{i = 1}^{M}|\widehat{P}^{\alpha}_{\widehat{U}_{n}^{\top}\tildhatY{i}{n+1}}[f]|^2\\
				= &
				\frac{1}{M}\sum_{i = 1}^{M} \mathrm{Tr} \Bigg( \left(\frac{1}{M}(P_{\widehat{U}_{n}}^{\perp} f) \widetilde{\widehat{\mathbb{Y}}}_{n+1}^{\top} \left( \frac{1}{M} \widetilde{\widehat{\mathbb{Y}}}_{n+1} \widetilde{\widehat{\mathbb{Y}}}_{n+1}^{\top} + \alpha I_{k \times k} \right)^{-1} \tildhatY{i}{n+1}\right)^{\top} \\
				& \cdot \left(\left(\frac{1}{M}(P_{\widehat{U}_{n}}^{\perp} f) \widetilde{\widehat{\mathbb{Y}}}_{n+1}^{\top} \right)\left( \frac{1}{M} \widetilde{\widehat{\mathbb{Y}}}_{n+1} \widetilde{\widehat{\mathbb{Y}}}_{n+1}^{\top} + \alpha I_{k \times k} \right)^{-1} \tildhatY{i}{n+1}\right)\Bigg) \\
				&+  \frac{1}{M}\sum_{i = 1}^{M} \mathrm{Tr} \Bigg( (f^{i})^{\top}P_{\widehat{U}_{n}}P_{\widehat{U}_{n}}f ^i\Bigg)\\ 
				= &
				\frac{1}{M}\sum_{i = 1}^{M} \mathrm{Tr} \Bigg( \left(\left(\frac{1}{M}(P_{\widehat{U}_{n}}^{\perp} f) \widetilde{\widehat{\mathbb{Y}}}_{n+1}^{\top}\right) (\widehat{C}^{\alpha}_{\widetilde{\widehat{Y}}_{n+1}})^{-1}  \tildhatY{i}{n+1}\right) \left(\frac{1}{M}(P_{\widehat{U}_{n}}^{\perp} f) \widetilde{\widehat{\mathbb{Y}}}_{n+1}^{\top} (\widehat{C}^{\alpha}_{\widetilde{\widehat{Y}}_{n+1}})^{-1}  \tildhatY{i}{n+1}\right)^{\top} \\
				&+  \frac{1}{M}\sum_{i = 1}^{M} \mathrm{Tr} \Bigg( P_{\widehat{U}_{n}}f ^i (f^{i})^{\top}P_{\widehat{U}_{n}} \Bigg)\\ 
				= &
				\mathrm{Tr} \Bigg( \left( \frac{1}{M}(P_{\widehat{U}_{n}}^{\perp} f) \widetilde{\widehat{\mathbb{Y}}}_{n+1}^{\top} \right)(\widehat{C}^{\alpha}_{\widetilde{\widehat{Y}}_{n+1}})^{-1}  \left( 	\frac{1}{M}\sum_{i = 1}^{M} \tildhatY{i}{n+1}\tildhatY{i}{n+1}\right)^{\top} (\widehat{C}^{\alpha}_{\widetilde{\widehat{Y}}_{n+1}})^{-1}  \left(\frac{1}{M}(P_{\widehat{U}_{n}}^{\perp} f)  \widetilde{\widehat{\mathbb{Y}}}_{n+1}^{\top}\right)^{\top} \Bigg)\\
				&+  \frac{1}{M}\sum_{i = 1}^{M} \mathrm{Tr} \Bigg( P_{\widehat{U}_{n}}f ^i (f^{i})^{\top}P_{\widehat{U}_{n}} \Bigg).
			\end{aligned}
		\end{equation*}
		Then, due to the fact that $\alpha > 0$, one has that
		\begin{equation*}
			\begin{aligned}
				&\frac{1}{M}\sum_{i = 1}^{M}|\widehat{P}^{\alpha}_{\widehat{U}_{n}^{\top}\tildhatY{i}{n+1}}[f]|^2\\
				\leq &
				\mathrm{Tr} \Bigg( \left( \frac{1}{M}(P_{\widehat{U}_{n}}^{\perp} f) \widetilde{\widehat{\mathbb{Y}}}_{n+1}^{\top} \right) (\widehat{C}^{\alpha}_{\widetilde{\widehat{Y}}_{n+1}})^{-1}  \left(\frac{1}{M}(P_{\widehat{U}_{n}}^{\perp} f)  \widetilde{\widehat{\mathbb{Y}}}_{n+1}^{\top}\right)^{\top} \Bigg)+  \frac{1}{M}\sum_{i = 1}^{M} \mathrm{Tr} \Bigg( P_{\widehat{U}_{n}}f ^i (f^{i})^{\top}P_{\widehat{U}_{n}} \Bigg)\\ 
				\leq &
				\mathrm{Tr} \Bigg( \left( \frac{1}{M}(P_{\widehat{U}_{n}}^{\perp} f)\right) \widetilde{\widehat{\mathbb{Y}}}_{n+1}^{\top}  (\widehat{C}^{\alpha}_{\widetilde{\widehat{Y}}_{n+1}})^{-1} \widetilde{\widehat{\mathbb{Y}}}_{n+1}\left(\frac{1}{M}(P_{\widehat{U}_{n}}^{\perp} f)  \right)^{\top} \Bigg)+  \frac{1}{M}\sum_{i = 1}^{M} \mathrm{Tr} \Bigg( P_{\widehat{U}_{n}}f ^i (f^{i})^{\top}P_{\widehat{U}_{n}} \Bigg)\\ 
				\leq &
				\mathrm{Tr} \Bigg( \frac{1}{M} \left( (P_{\widehat{U}_{n}}^{\perp} f)\right) I_{M \times M} \left((P_{\widehat{U}_{n}}^{\perp} f)  \right)^{\top} \Bigg)+  \frac{1}{M}\sum_{i = 1}^{M} \mathrm{Tr} \Bigg( P_{\widehat{U}_{n}}f ^i (f^{i})^{\top}P_{\widehat{U}_{n}} \Bigg)\\ 
				= &
				\mathrm{Tr} \Bigg( \frac{1}{M}P_{\widehat{U}_{n}}^{\perp} \sum_{i = 1}^{M} f^i (f^i)^{\top}  P_{\widehat{U}_{n}}^{\perp}  \Bigg)+  \frac{1}{M}\sum_{i = 1}^{M} \mathrm{Tr} \Bigg( P_{\widehat{U}_{n}}f ^i (f^{i})^{\top}P_{\widehat{U}_{n}} \Bigg)\\ 
				\leq & \frac{1}{M} \sum_{i = 1}^{M} |P_{\widehat{U}_{n}}^{\perp}f ^i|^2 +  \frac{1}{M}\sum_{i = 1}^{M}|P_{\widehat{U}_{n}}f ^i|^2 = \frac{1}{M}\sum_{i = 1}^{M}|f ^i|^2,\\ 
			\end{aligned}
		\end{equation*}
		by orthogonality of $U_n$.
	\end{proof}
\end{Lemma}

Let us recall that $\mathbb{E}_n[X]$ is the conditional expectation of a random variable $X$ with respect to $\mathcal{F}_{t_n}$, i.e. $\mathbb{E}_n[X]=\mathbb{E}[X | \mathcal{F}_{t_n}]$.

\begin{Lemma}[Conditional bound for the empirical second moment of $\widetilde{\widehat{\mathbb{Y}}}_{n+1}$]\label{lem: cond bound Mn}
	Assume that $\Delta t_n$ satisfies \eqref{eq: M dt cond}.
	Let us define the stochastic process $(\widehat{Z}_{n}^i)_{n}$ as $\widehat{Z}_{0}^i = 	\widehat{X}_{0}^i$ and
	\begin{equation*}
		\widehat{Z}_{n+1}^i = 	\widehat{X}_{n}^i + a(t_n, 	\widehat{X}_{n}^i) \Delta t_n +  P_{\widehat{U}_n} b(t_n, 	\widehat{X}_{n}^i) \Delta W_n^{i}, \quad \text{ for } n \ : \ 0\leq n \leq N-1, \ i=\{1, \dots, M\}.
	\end{equation*}
	
	Furthermore, define $\widehat{M}_n:=\frac1M \sum_{i=1}^M |\widehat{X}_n^i|^2$ and $\widehat{N}_n:=\frac1M \sum_{i=1}^M |\widehat{Z}_n^i|^2$. Then, for all $1\leq n \leq N$ we have that $\widehat{M}_n \leq \widehat{N}_n$ and 
	\begin{equation}\label{eq: bound N_n}
		\begin{aligned}
			\mathbb{E}_{n}[\widehat{N}_{n+1}]
			\leq &  \overline{K}_2(\widehat{M}_0,t_{n+1})+\sum_{\ell=1}^{n} (1+ (1 +C_{\mathrm{lgb}} (2+T) )\Delta t_{\ell})^{n-\ell+1} \rho_{\ell,M}
		\end{aligned}
	\end{equation}
	where $\rho_{\ell,M}:= \widehat{N}_{\ell} - \mathbb{E}_{\ell-1}[\widehat{N}_{\ell}]$ and $\overline{K}_2(\widehat{M}_0,t):= \left( \widehat{M}_0 + 1\right) \exp\left((1 +C_{\mathrm{lgb}} (2+T) ) t\right)-1 + C_{\mathrm{lgb}}(2+T)\Delta t_0$.
	\begin{proof}
		We first prove the bound $\widehat{M}_n \leq \widehat{N}_n$, for all $n$.
		Exploiting the orthogonality of the rows of $\widehat{U}_n$, the fact that $ \widehat{X}_n^{i}= \widehat{U}_n^{\top}\widehat{Y}_n^{i}$, and $\frac{1}{M}\sum_{i = 1}^{M}|\widehat{P}^{\alpha}_{\widehat{U}_{n}^{\top}\tildhatY{i}{n+1}}[\widehat{a}_n]|^2 \leq \frac1M \sum_{i=1}^M |\widehat{a}_n^{i}|^2$ by Lemma \ref{lem: P f < f}, one has that
		\begin{equation*}
			\begin{aligned}
				\widehat{M}_{n+1}=&\frac1M \sum_{i=1}^M |\widehat{X}_{n+1}^i|^2 \\
				= & \widehat{M}_n + \frac2M \sum_{i=1}^M (\widehat{X}_{n}^i)^{\top}P_{\widehat{U}_{n}}[\widehat{a}_n^{i}] \Delta t_n + \frac2M \sum_{i=1}^M (\widehat{X}_{n}^i)^{\top}P_{\widehat{U}_{n}}[\widehat{b}_n^{i}] \Delta W_n^{i}  \\
				& + \frac1M \sum_{i=1}^M |\widehat{P}^{\alpha}_{\widehat{U}_{n}^{\top}\tildhatY{i}{n+1}}[\widehat{a}_n]|^2 (\Delta t_n)^2 + \frac2M \sum_{i=1}^M (\widehat{a}_n^{i})^{\top} P_{\widehat{U}_{n}}[\widehat{b}_n^{i}] \Delta W_n^{i} \Delta t_n +  \frac1M \sum_{i=1}^M |P_{\widehat{U}_{n}}[\widehat{b}_n^{i}] \Delta W_n^{i}|^2 \\
				\leq &  \widehat{M}_n+ \frac2M \sum_{i=1}^M (\widehat{X}_{n}^i)^{\top}\widehat{a}_n^{i} \Delta t_n + \frac2M \sum_{i=1}^M (\widehat{X}_{n}^i)^{\top}P_{\widehat{U}_{n}}[\widehat{b}_n^{i}] \Delta W_n^{i}  \\
				& + \frac1M \sum_{i=1}^M |\widehat{a}_n^{i}|^2 (\Delta t_n)^2 + \frac2M \sum_{i=1}^M (\widehat{a}_n^{i})^{\top} P_{\widehat{U}_{n}}[\widehat{b}_n^{i}] \Delta W_n^{i} \Delta t_n +  \frac1M \sum_{i=1}^M |P_{\widehat{U}_{n}}[\widehat{b}_n^{i}] \Delta W_n^{i}|^2 \\
				=&\frac1M \sum_{i=1}^M |\widehat{Z}_{n+1}^i|^2 = \widehat{N}_{n+1},
			\end{aligned}
		\end{equation*}
		which is the first part of the statement. 
		
		To prove the bound \eqref{eq: bound N_n}, we exploit the conditioning on $\mathcal{F}_{t_n}$, the orthogonality of the rows of $\widehat{U}_n$, and Young's inequality, to obtain that
		\begin{equation*}
			\begin{aligned}
				\mathbb{E}_n[\widehat{N}_{n+1}] = &  \widehat{M}_n+ \frac2M \sum_{i=1}^M (\widehat{X}_{n}^i)^{\top}\widehat{a}_n^{i} \Delta t_n + \frac1M \sum_{i=1}^M |\widehat{a}_n^{i}|^2 (\Delta t_n)^2  +  \frac1M \sum_{i=1}^M \|P_{\widehat{U}_{n}}[\widehat{b}_n^{i}]\|^2_{\mathrm{F}}  \Delta t \\
				\leq &  \widehat{M}_n + \widehat{M}_n  \Delta t_n  + \frac1M \sum_{i=1}^M |\widehat{a}_n^{i}|^2 \Delta t_n + \frac1M \sum_{i=1}^M |\widehat{a}_n^{i}|^2 (\Delta t_n)^2 +  \frac1M \sum_{i=1}^M \|\widehat{b}_n^{i}\|^2_{\mathrm{F}}  \Delta t_n\\
				\leq & \widehat{M}_n + \widehat{M}_n  \Delta t_n  + \frac1M \sum_{i=1}^M C_{\mathrm{lgb}}(1+|\widehat{X}_n^{i}|^2) \Delta t_n + \frac1M \sum_{i=1}^M C_{\mathrm{lgb}}(1+|\widehat{X}_n^{i}|^2) (\Delta t_n)^2 \\
				&+  \frac1M \sum_{i=1}^M C_{\mathrm{lgb}}(1+|\widehat{X}_n^{i}|^2)  \Delta t_n\\
				= & \widehat{M}_n + \widehat{M}_n  \Delta t_n  +  C_{\mathrm{lgb}}(1+\widehat{M}_n) (2\Delta t_n +  (\Delta t_n)^2 )\\
				\leq  &(1+ (1 +C_{\mathrm{lgb}} (2+T) )\Delta t_n ) \widehat{M}_n + C_{\mathrm{lgb}} (2+T) \Delta t_n.\\
				\leq  &(1+ (1 +C_{\mathrm{lgb}} (2+T) )\Delta t_n ) \widehat{N}_n + C_{\mathrm{lgb}} (2+T) \Delta t_n.\\
			\end{aligned}
		\end{equation*}
		where in the last line we used the first statement.
		We use now the decomposition $\widehat{N}_n = \mathbb{E}_{n-1}[\widehat{N}_{n}] + \widehat{N}_n - \mathbb{E}_{n-1}[\widehat{N}_{n}]$ and rewrite the previous inequality as
		\begin{equation*}
			\begin{aligned}
				\mathbb{E}_n[\widehat{N}_{n+1}] +1\leq & (1+ (1 +C_{\mathrm{lgb}} (2+T) )\Delta t_n )\mathbb{E}_{n-1}[\widehat{N}_{n}] + (1 + C_{\mathrm{lgb}} (2+T) \Delta t_n)\\
				& + (1+ (1 +C_{\mathrm{lgb}} (2+T) )\Delta t_n )\left(\widehat{N}_n - \mathbb{E}_{n-1}[\widehat{N}_{n}]\right) \\
				\leq & (1+ (1 +C_{\mathrm{lgb}} (2+T) )\Delta t_n ) (\mathbb{E}_{n-1}[\widehat{N}_{n}] +1)\\
				& + (1+ (1 +C_{\mathrm{lgb}} (2+T) )\Delta t_n )\left(\widehat{N}_n - \mathbb{E}_{n-1}[\widehat{N}_{n}]\right) \\
			\end{aligned}
		\end{equation*}
		
		Via Gronwall's lemma and exponential bound we conclude that
		\begin{equation*}
			\begin{aligned}
				\mathbb{E}_{n}[\widehat{N}_{n+1}] \leq & \left( \widehat{M}_0 +1\right) \exp\left( (1 +C_{\mathrm{lgb}} (2+T) ) t_{n+1} \right) -1 + C_{\mathrm{lgb}}(2+T)\Delta t_0 \\
				&+ \sum_{\ell=1}^{n} (1+ (1 +C_{\mathrm{lgb}} (2+T) )\Delta t_\ell)^{n-\ell+1} \rho_{\ell,M},
			\end{aligned}
		\end{equation*}
		where we stop the recursion at the step $\ell=1$.
	\end{proof}
\end{Lemma}

\begin{Remark}\label{rmk: link K_over and K}
	Notice that if $\widehat{M}_0 \leq 2 \mathbb{E}[|Y_0|^2]$, then $\overline{K}_2(\widehat{M}_0,t_n) \leq 2K_2$ for all $n$, where $K_2$ is defined in Lemma \ref{lem: L2p norm semidiscretized solution}. 
\end{Remark}

In order to prove a lower bound in probability on the empirical Gramian of the fully-discrete intermediate stochastic basis  $(\tildhatY{i}{n+1})_{i=1,\dots,M}$, we need Monte-Carlo type estimates of the conditional differences involving the auxiliary process.
\begin{Lemma}[Local Monte Carlo fluctuation of $\widehat{N}_n$]\label{lem: Local MC Mn}
	Suppose that we choose a uniform time-step $\Delta t$ satisfying \eqref{eq: M dt cond}. Further, suppose $\mathbb{E}[|Y_0|^{4p}] < \infty$ for some $p\geq 1$.
	For all $n$ such that $0\leq n \leq N-1$ we have that
	\begin{equation*}
		\begin{aligned}
			\mathbb{E}\left[\left|\widehat{N}_{n+1}-\mathbb{E}_{n}[\widehat{N}_{n+1}]\right|^{2p}\right]\leq C_2\frac{(\Delta t)^{p}}{M^{p}},
		\end{aligned}
	\end{equation*}
	for a positive constant $C_2 = C_2(p,T)$ independent of $M$ and $\Delta t$, where $\widehat{N}_n$ is defined in Lemma \ref{lem: cond bound Mn}. Furthermore, it holds that
	\begin{equation*}
		\mathbb{E}\left[  \sup_{1\leq h \leq n}\left| \sum_{\ell=1}^{h} \gamma^{-\ell}\rho_{\ell,M} \right|^{2p}\right] \leq C_{3} \frac{1}{M^{p}},
	\end{equation*}
	for any $ \gamma \geq 1$,
	for a positive constant $C_3 = C_3(p,T)$ independent of $M$ and $\Delta t$, where $\rho_{\ell,M}:= \widehat{N}_{\ell} - \mathbb{E}_{\ell-1}[\widehat{N}_{\ell}]$ as in Lemma \ref{lem: cond bound Mn}.
	\begin{proof}
		For brevity, define
		\begin{equation*}
			\widehat{Z}_{n+1}^i = m_i + \xi_i, \text{ with } m_i := \widehat{X}_n^i +  \widehat{a}_n^{i} \Delta t, \
			\xi_i := P_{\widehat{U}_n} \widehat{b}_n^{i} \Delta W_n^i.
		\end{equation*}
		Then,
		\begin{equation}\label{eq: Z_n hat}
			|\widehat{Z}_{n+1}^i|^2= |m_i + \xi_i|^2 = |m_i|^2 + 2 m_i^\top \xi_i + |\xi_i|^2.
		\end{equation} 
		Taking conditional expectation with respect to $\mathcal{F}_{t_n}$ and independence of Brownian increments, from \eqref{eq: Z_n hat} we obtain
		\begin{equation*}
			\begin{aligned}
				\mathbb{E}_n[|\widehat{Z}_{n+1}^i|^2]= |m_i|^2 + \mathbb{E}_n[|\xi_i|^2].
			\end{aligned}
		\end{equation*}
		which implies that $|\widehat{Z}^{i}_{n+1}|^2 - \mathbb{E}_n[|\widehat{Z}^{i}_{n+1}|^2]= 2 m_i^\top \xi_i + |\xi_i|^2 - \mathbb{E}_n[|\xi_i|^2].$
		
		Define $\zeta_i:=|\widehat Z_{n+1}^i|^2-\mathbb E_n[|\widehat Z_{n+1}^i|^2]$. Then $\mathbb E_n[\zeta_i]=0$. Moreover, the variables $(\zeta_i)_{i=1}^M$ are independent conditionally on $\mathcal F_{t_n}$.
		Notice that by triangle inequality and unitary norm of orthogonal projections we have
		\begin{equation*}
			\begin{aligned}
				|\zeta_i|\leq & 2|m_i||\xi_i|+\left||\xi_i|^2 - \mathbb{E}_n[|\xi_i|^2]\right| \\
				\leq &2 |\widehat{X}_n^i +  \widehat{a}_n^{i} \Delta t| |P_{\widehat{U}_n}  \widehat{b}_n^{i} \Delta W_n^i| + \left| |P_{\widehat{U}_n} \widehat{b}_n^{i} \Delta W_n^i| - \| P_{\widehat{U}_n} \widehat{b}_n^{i}\|^2_{\mathrm{F}} \Delta t\right|,
			\end{aligned}
		\end{equation*}
		and, thanks to the linear-growth bound, the tower property, Young's inequality, unitaty norm of orthogonal projectors, and Lemma \ref{lem: Lp moment bound}, for any $p \geq 1$ there exist positive constants $C_{p}$ and $C_1:=C_1(p,T)$ such that
		\begin{equation}\label{eq: bound eta_i at p}
			\begin{aligned}
				\mathbb{E}\left[|\zeta_i|^{2p}\right]\leq & \mathbb{E}\left[\left( 2|\widehat{X}_n^i +  \widehat{a}_n^{i} \Delta t|  |P_{\widehat{U}_n}  \widehat{b}_n^{i} \Delta W_n^i| + \left| | P_{\widehat{U}_n}\widehat{b}_n^{i} \Delta W_n^i|^2 - \| P_{\widehat{U}_n}\widehat{b}_n^{i}\|^2_{\mathrm{F}}  \Delta t\right|\right)^{2p}\right] \\
				\leq & 3^{2p-1} \Bigg(\mathbb{E}\left[\left(2 |\widehat{X}_n^i +  \widehat{a}_n^{i} \Delta t| |P_{\widehat{U}_n}  \widehat{b}_n^{i} \Delta W_n^i| \right)^{2p}\right]  + \mathbb{E}\left[\left(\| P_{\widehat{U}_n} \widehat{b}_n^{i} \Delta W_n^i\|^2_{\mathrm{F}}  \right)^{2p}\right] \\
				& +\mathbb{E}\left[\left(\| P_{\widehat{U}_n} \widehat{b}_n^{i}\|^2_{\mathrm{F}}  \Delta t\right)^{2p}\right] \Bigg)\\
				\leq & 3^{2p-1} \Bigg(\mathbb{E}\left[ \mathbb{E}_n\left[\left(2 |\widehat{X}_n^i +  \widehat{a}_n^{i} \Delta t|  |P_{\widehat{U}_n}  \widehat{b}_n^{i} \Delta W_n^i|\right)^{2p}\right] \right]  + \mathbb{E}\left[  \mathbb{E}_n\left[\left( \| P_{\widehat{U}_n} \widehat{b}_n^{i} \Delta W_n^i\|^2_{\mathrm{F}}  \right)^{2p}\right] \right]  \\
				&+ \mathbb{E}\left[\left(\|\widehat{b}_n^{i}\|_{\mathrm{F}}^2 \Delta t\right)^{2p}\right] \Bigg)\\
				\leq & 3^{2p-1} C_{p} \left(\mathbb{E}\left[\left( 2 |\widehat{X}_n^i +  \widehat{a}_n^{i} \Delta t| \|\widehat{b}_n^{i}\|_{\mathrm{F}} \right)^{2p}|\Delta t|^{p} \right]  + \mathbb{E}\left[\|\widehat{b}_n^{i}\|_{\mathrm{F}}^{4p} |\Delta t|^{2p} \right]  + \mathbb{E}\left[\left(\|\widehat{b}_n^{i}\|_{\mathrm{F}}^2 \Delta t\right)^{2p}\right] \right)\\
				\leq & 3^{2p-1} C_{p} \left( \mathbb{E}\left[ \left( |\widehat{X}_n^i +  \widehat{a}_n^{i} \Delta t|^2 +   \|\widehat{b}_n^{i}\|_{\mathrm{F}}^2 \right)^{2p} \right] (\Delta t)^{p} + \mathbb{E}\left[\|\widehat{b}_n^{i}\|_{\mathrm{F}}^{4p} \right] (\Delta t)^{2p}  + \mathbb{E}\left[\|\widehat{b}_n^{i}\|_{\mathrm{F}}^{4p} \right] \Delta t^{2p} \right)\\
				\leq & C_{1} (\Delta t)^p.
			\end{aligned}
		\end{equation}
		
		Using the tower property and a version of Rosenthal inequality for conditionally independent and centered random variables \cite[Theorem 1.1]{hu2011extension} similarly to relation \eqref{eq: MC in k-dim}, we get
		\begin{equation*}
			\begin{aligned}
				\mathbb{E}\left[\left|\frac{1}{M}\sum_{i=1}^M \zeta_i\right|^{2p}
				\right]=\mathbb{E}\left[\mathbb{E}_n\left[
				\left|\frac{1}{M}\sum_{i=1}^M \zeta_i\right|^{2p}\right]\right] \leq&\frac{C_{p}}{M^{2p}}\max\left\{\left(\sum_{i=1}^M 	\mathbb{E}\left[\mathbb{E}_n[|\zeta_i|^{2p}]\right], \mathbb{E}\left[\left(\sum_{i=1}^M\mathbb{E}_n[|\zeta_i|^{2}]\right)^p\right]\right)\right\} \\
				\leq & C_{p}\max\left\{\left(\frac{1}{M^{2p-1}}	\mathbb{E}\left[|\zeta_1|^{2p}\right], \frac{1}{M^{p}}\mathbb{E}\left[\left(\mathbb{E}_n[|\zeta_1|^{2}]\right)^p\right]\right)\right\}.
			\end{aligned}
		\end{equation*}
		for a positive constant $C_{p}$ dependent on $p$. Using relation \eqref{eq: bound eta_i at p} and Jensen's inequality, we obtain the first part of the statement:
		\begin{equation}\label{eq: eta_i}
			\begin{aligned}
				\mathbb{E}[|\widehat{N}_{n+1}-\mathbb{E}_n[\widehat{N}_{n+1}]|^{2p}]=\mathbb{E}\left[\left|\frac{1}{M}\sum_{i=1}^M \zeta_i\right|^{2p}
				\right]=\mathbb{E}\left[\mathbb{E}_n\left[
				\left|\frac{1}{M}\sum_{i=1}^M \zeta_i\right|^{2p}\right]\right] \leq C_{p,T}\frac{(\Delta t)^{p}}{M^{p}},
			\end{aligned}
		\end{equation}
		for a positive constant $C_{p,T}$ dependent on $p$ and $T$.
		
		To prove the second statement, notice that the sequence $\{ \rho_{\ell,M}:= \widehat{N}_{\ell} - \mathbb{E}_{\ell-1}[\widehat{N}_{\ell}] \}_{\ell}$ is a martingale difference sequence  with respect to the filtration $\mathcal{F}_{t_{\ell}}$, hence $\sum_{\ell=1}^{h} \gamma^{-\ell}\rho_{\ell,M}$ is a martingale. Then, by means of Doob's maximal inequality, Burkholder-Davis-Gundy inequality \cite[Section 2]{yaroslavtsev2020burkholder}, and relation \eqref{eq: eta_i}, there exist constants $C:=C(p,T)$, $\widetilde{C}:=\widetilde{C}(p,T)$ such that 
		\begin{equation*}
			\begin{aligned}
				\mathbb{E}\left[  \sup_{1\leq h \leq n} \left| \sum_{\ell=1}^{h}\gamma^{-\ell} \rho_{\ell,M} \right|^{2p}\right] \leq & \left(\frac{2p}{2p-1}\right)^{2p} \mathbb{E}\left[ \left| \sum_{\ell=1}^{n} \gamma^{-\ell}\rho_{\ell,M} \right|^{2p}\right] \\
				\leq &  C \mathbb{E}\left[ \left| \sum_{\ell=1}^{n}\gamma^{-2\ell}  \left| \rho_{\ell,M} \right|^{2}  \right|^{p}\right] \\
				\leq & C n^{p-1} \sum_{\ell=1}^{n}  \gamma^{-2p\ell} \mathbb{E}\left[ \left| \rho_{\ell,M} \right|^{2p} \right] \leq  C n^{p-1} \sum_{\ell=1}^{n} \mathbb{E}\left[ \left| \rho_{\ell,M} \right|^{2p}\right]\\
				\leq & C n^{p} C_{p,T}\frac{(\Delta t)^{p}}{M^{p}}  = \widetilde{C} \frac{1}{M^{p}},
			\end{aligned}
		\end{equation*}
		where in the second to last line we use the fact that $\gamma \geq1$, which concludes the proof.
	\end{proof}
\end{Lemma} 

In our estimates, we need to bound also the difference between the empirical Gramian of the fully-discrete stochastic basis and its conditional expectation with respect to the filtration at the previous time step.
\begin{Lemma}[$2p$-moment bound on centered conditional covariance]\label{lem: p-mom cond covariance}
	Assume that $\mathbb{E}[|Y_0|^{4p}] < \infty$ for some $p\geq 1$ and suppose that $\Delta t_n$ satisfies \eqref{eq: M dt cond} for every $n$.
	Consider a real number $q$ such that $0 \leq q \leq 1$ and define $\widetilde{\widehat{S}}_n:=	\frac1M\sum_{i=1}^{M}\tildhatY{i}{n}(\tildhatY{i}{n})^{\top}$. Then, there exists a constant $C_4:=C_4(p,T)$ independent of $M$ and $\{\Delta t_n\}_n$ such that 
	\begin{equation*}
		\mathbb{E}\left[\left|\sum_{\ell=0}^n  q^{n-\ell} \eta_{\ell+1,M}\right|^{2p}\right] \leq C_4 \frac{ (\Delta t)^p}{M^{p}}  \left(\sum_{\ell=0}^n q^{2(n-\ell)}\right)^{p},
	\end{equation*}
	where $\eta_{\ell+1,M}:= \left(\widetilde{\widehat{S}}_{\ell+1} - \mathbb{E}_\ell[\widetilde{\widehat{S}}_{\ell+1}]\right)$ and $\Delta t = \max\limits_{0\le \ell \le n} \Delta t_{\ell}$.
	\begin{proof}
		To prove the statement, we exploit an application of Burkholder-Davis-Gundy inequality.
		Indeed, the random variable $\widetilde{\widehat{S}}_{n+1} - \mathbb{E}_n[\widetilde{\widehat{S}}_{n+1}]$ is integrable and centered with respect to $\mathcal{F}_{t_n}$, since
		\begin{equation*}
			\mathbb{E}_n\left[\widetilde{\widehat{S}}_{n+1} - \mathbb{E}_n[\widetilde{\widehat{S}}_{n+1}]\right] = \mathbb{E}_n\left[\widetilde{\widehat{S}}_{n+1}\right]- \mathbb{E}_n[\widetilde{\widehat{S}}_{n+1}] = 0.
		\end{equation*}
		Hence, $\left(\widetilde{\widehat{S}}_{n+1} - \mathbb{E}_n[\widetilde{\widehat{S}}_{n+1}]\right)_{n}$ is martingale difference sequence  and the stochastic process 
		$\sum_{\ell=0}^n  q^{-\ell} \left(\widetilde{\widehat{S}}_{\ell+1} - \mathbb{E}_\ell[\widetilde{\widehat{S}}_{\ell+1}]\right)$ is a martingale. Then, by Burkholder-Davis-Gundy inequality \cite{hall2014martingale} there exists a positive constant $C_{p}$ dependent on $p$ such that one has
		\begin{equation*}
			\begin{aligned}
				\mathbb{E}	\left[\left|\sum_{\ell=0}^n  q^{n-\ell} \left(\widetilde{\widehat{S}}_{\ell+1} - \mathbb{E}_\ell[\widetilde{\widehat{S}}_{\ell+1}]\right)\right|^{2p}\right] = & q^{2pn}	\mathbb{E}	\left[\left|\sum_{\ell=0}^n  q^{-\ell} \left(\widetilde{\widehat{S}}_{\ell+1} - \mathbb{E}_\ell[\widetilde{\widehat{S}}_{\ell+1}]\right)\right|^{2p}\right]  \\
				\leq & q^{2pn}	\mathbb{E}	\left[\left\|\sum_{\ell=0}^n  q^{-\ell} \left(\widetilde{\widehat{S}}_{\ell+1} - \mathbb{E}_\ell[\widetilde{\widehat{S}}_{\ell+1}]\right)\right\|^{2p}_{\mathrm{F}}\right]  \\
				\leq 
				& q^{2pn} C_{p} \mathbb{E}	\left[ \left( \sum_{\ell=0}^n  q^{-2\ell} \left\|\left(\widetilde{\widehat{S}}_{\ell+1} - \mathbb{E}_\ell[\widetilde{\widehat{S}}_{\ell+1}]\right)\right\|^2_{\mathrm{F}} \right)^{p}\right] \\
				= 
				& C_{p} \mathbb{E}	\left[ \left( \sum_{\ell=0}^n  q^{2(n-\ell)} \left\|\left(\widetilde{\widehat{S}}_{\ell+1} - \mathbb{E}_\ell[\widetilde{\widehat{S}}_{\ell+1}]\right)\right\|^2_{\mathrm{F}} \right)^{p}\right] \\
				\leq  & C_{p}  \left( \sum_{\ell=0}^n  q^{2(n-\ell)} \right)^{p-1} \left( \sum_{\ell=0}^n  q^{2(n-\ell)}\mathbb{E}	\left[ \left\|\widetilde{\widehat{S}}_{\ell+1} - \mathbb{E}_\ell[\widetilde{\widehat{S}}_{\ell+1}] \right\|^{2p}_{\mathrm{F}}\right] \right) \\
				= & C_{p}  \left( \sum_{\ell=0}^n  q^{2(n-\ell)} \right)^{p-1} \left( \sum_{\ell=0}^n  q^{2(n-\ell)}\mathbb{E}	\left[ \mathbb{E}_{\ell} \left[ \left\|\widetilde{\widehat{S}}_{\ell+1} - \mathbb{E}_\ell[\widetilde{\widehat{S}}_{\ell+1}] \right\|^{2p}_{\mathrm{F}}\right] \right]\right) \\
			\end{aligned}
		\end{equation*}
		where we employed the weighted Holder inequality and the tower property.
		
		Now, we want to estimate $  \mathbb{E}_n\left[\left\|\widetilde{\widehat{S}}_{n+1} - \mathbb{E}_n[\widetilde{\widehat{S}}_{n+1}]\right\|^{2p}_{\mathrm{F}}\right]$. One has that
		\begin{equation*}
			\begin{aligned}
				\mathbb{E}_n\left[\left\|\widetilde{\widehat{S}}_{n+1} - \mathbb{E}_n[\widetilde{\widehat{S}}_{n+1}]\right\|^{2p}_{\mathrm{F}}\right]  
				\leq & C_{p} \frac{1}{M^{2p}} \sum_{i=1}^M \mathbb{E}_n\left[\left\|\tildhatY{i}{n+1}(\tildhatY{i}{n+1})^{\top}-\mathbb{E}_n\left[\tildhatY{i}{n+1}(\tildhatY{i}{n+1})^{\top}\right]\right\|^{2p}_{\mathrm{F}}\right] \\
				&+ C_{p} \frac{1}{M^{2p}} \left(\sum_{i=1}^M \mathbb{E}_n\left[\left\|\tildhatY{i}{n+1}(\tildhatY{i}{n+1})^{\top}-\mathbb{E}_n\left[\tildhatY{i}{n+1}(\tildhatY{i}{n+1})^{\top}\right]\right\|^{2}_{\mathrm{F}}\right] \right)^{p},
			\end{aligned}
		\end{equation*}
		where in the last line we employed the Rosenthal inequality \cite{oskekowski2012note} for conditionally independent and centered random variables similarly to the proof of \eqref{eq: MC in k-dim}, for a positive constant $C_{p}$ dependent on $p$. One gets
		\begin{equation*}
			\begin{aligned}
				\tildhatY{i}{n+1}(\tildhatY{i}{n+1})^{\top}-\mathbb{E}_n\left[\tildhatY{i}{n+1}(\tildhatY{i}{n+1})^{\top}\right]= & \widehat{Y}_n^i(\widehat{U}_n\widehat{b}_n^{i}\Delta W_n^{i})^{\top} + \widehat{U}_n\widehat{b}_n^{i}\Delta W_n^{i} (\widehat{Y}_n^i)^{\top} \\
				& + \widehat{U}_n\widehat{a}_n^{i} (\widehat{U}_n\widehat{b}_n^{i}\Delta W_n^{i})^{\top}\Delta t_n + \widehat{U}_n\widehat{b}_n^{i}\Delta W_n^{i} (\widehat{U}_n\widehat{a}_n^{i})^{\top}\Delta t_n \\
				&+  \widehat{U}_n\widehat{b}_n^{i}\Delta W_n^{i} (\widehat{U}_n\widehat{b}_n^{i}\Delta W_n^{i})^{\top}-  \widehat{U}_n\widehat{b}_n^{i}(\widehat{U}_n\widehat{b}_n^{i})^{\top} \Delta t_n
			\end{aligned}
		\end{equation*}
		and therefore, proceeding similarly to the proof of Lemma \ref{lem: Local MC Mn}, using Lemma \ref{lem: L2p norm semidiscretized solution} and the fact that $\Delta W_n^{i} \sim \mathcal{N}(0,\Delta t I_{m \times m})$, it follows that
		\begin{equation*}
			\begin{aligned}
				\mathbb{E}_n\left[\left\|\tildhatY{i}{n+1}(\tildhatY{i}{n+1})^{\top}-\mathbb{E}_n\left[\tildhatY{i}{n+1}(\tildhatY{i}{n+1})^{\top}\right]\right\|^{2p}_{\mathrm{F}}\right] \leq & \widetilde{C}_{p} \left(1 +  |\tildhatY{i}{n}|^{4p}\right) (\Delta t_n)^p,
			\end{aligned}
		\end{equation*}
		for a positive constant $\widetilde{C}_{p} = \widetilde{C}(p,T)$, which for $p=1$ is equal to
		\begin{equation*}
			\begin{aligned}
				\mathbb{E}_n\left[\left\|\tildhatY{i}{n+1}(\tildhatY{i}{n+1})^{\top}-\mathbb{E}_n\left[\tildhatY{i}{n+1}(\tildhatY{i}{n+1})^{\top}\right]\right\|^{2}_{\mathrm{F}}\right] \leq & \widetilde{C}_{1} \left(1 +  |\tildhatY{i}{n}|^{4}\right) \Delta t_n.
			\end{aligned}
		\end{equation*}
		
		Those estimates imply that 
		\begin{equation*}
			\begin{aligned}
				\mathbb{E}\left[ \mathbb{E}_{n}\left[\left\|\widetilde{\widehat{S}}_{n+1} - \mathbb{E}_n[\widetilde{\widehat{S}}_{n+1}]\right\|^{2p}_{\mathrm{F}}\right]\right]   
				\leq  &C \frac{1}{M^{p}}  \left(1 +  \mathbb{E}\left[ |\tildhatY{i}{n}|^{4p}\right]\right) (\Delta t_n)^p \leq \tilde{C} 	\frac{ (\Delta t_n)^p}{M^{p}},
			\end{aligned}
		\end{equation*}
		for positive constants $C,\widetilde{C}$, and therefore it holds that
		\begin{equation*}
			\begin{aligned}
				\mathbb{E}	\left[\left|\sum_{\ell=0}^n  q^{n-\ell} (\widetilde{\widehat{S}}_{\ell+1} - \mathbb{E}_\ell[\widetilde{\widehat{S}}_{\ell+1}])\right|^{2p}\right] \leq & C_{p}  \left( \sum_{\ell=0}^n  q^{2(n-\ell)} \right)^{p-1} \left( \sum_{\ell=0}^n  q^{2(n-\ell)}\mathbb{E}	\left[ \mathbb{E}_{\ell} \left[ \left\|\widetilde{\widehat{S}}_{\ell+1} - \mathbb{E}_\ell[\widetilde{\widehat{S}}_{\ell+1}] \right\|^{2p}_{\mathrm{F}}\right] \right]\right) \\
				\leq & C_{p}  \left( \sum_{\ell=0}^n  q^{2(n-\ell)} \right)^{p-1} \left( \sum_{\ell=0}^n  q^{2(n-\ell)} \tilde{C} 	\frac{ (\Delta t_\ell)^p}{M^{p}} \right) \\
				= & C_{p}\widetilde{C}\frac{ (\Delta t)^p}{M^{p}} \left( \left(\sum_{\ell=0}^n q^{2(n-\ell)}\right)^{p}\right),
			\end{aligned}
		\end{equation*}
		where $\Delta t = \max\limits_{0\le \ell \le n} \Delta t_{\ell}$.
	\end{proof}
\end{Lemma}

Finally, we can provide the proof of Proposition \ref{prop: lower-bound Gramian} regarding a lower-bound on the empirical Gramian of $\tildhatX{i}{n+1}=\widehat{U}_n^{\top}\tildhatY{i}{n+1}$ in probability.

\begin{proof}[\bfseries Proof of Proposition \ref{prop: lower-bound Gramian}]
	First, let us define $A:=\frac{\sigma_{B}}{2 C_{\mathrm{lgb}}(1+2K_2(T))}=\frac{\varrho}{4\sigma_{B}}$
	and the following event $E_n$
	\begin{equation}\label{eq: E_n}
		\begin{aligned}
			E_n:= &\left\{\frac{1}{M} \sum_{i=1}^M |Y_0^{i}|^2 < 2 \mathbb{E}[|Y_0|^2] \right\} \cap \left\{\frac{1}{M} \sum_{i=1}^M Y_0^{i}(Y_0^{i})^{\top} \succeq \frac{\sigma_{Y_0}}{2}  \right\} \\ &\cap  \left\{ \sup_{1\leq h \leq n}\left|\sum_{\ell=1}^h \gamma^{-\ell}\rho_{\ell,M}\right|  < \frac{1}{4}  \frac{2(1+2K_2(T))}{\exp\left((1 +C_{\mathrm{lgb}} (2+T))T\right) }\right\} \\
			& \cap \left\{  \left|\sum_{\ell=0}^n  (1- \frac{\Delta t}{A + \Delta t})^{n-\ell} \eta_{\ell+1,M} \right|  < \min\{ \frac{\varrho}{4},  \frac{\sigma_{Y_0}}{2} \} \right\},
		\end{aligned}
	\end{equation}
	where $\varrho$ is defined in \eqref{eq: rho} and $\gamma:=(1 +(1+C_{\mathrm{lgb}}(2+T))\Delta t)$. Notice that $E_n \subset E_{j}$ for all $j < n$, hence the sequence $\{E_n\}_n$ is nested.
	To prove the first statement, we retrieve a useful lower bound on the smallest singular value of the covariance of $\tildhatX{i}{n+1}$ at time $t_n$ in the event $E_n$. The remaining part of the proof will be dedicated to give an upper bound on the probability of the complement of $E_n$, aiming to find a Monte-Carlo type error.
	For the sake of notation, let us denote the empirical Gramians as follows
	\begin{equation*}
		\begin{aligned}
			\widehat{S}_n:=	\frac1M\sum_{i=1}^{M}\widehat{Y}_n^i(\widehat{Y}_n^i)^{\top}, \quad \widetilde{\widehat{S}}_n:=	\frac1M\sum_{i=1}^{M}\tildhatY{i}{n}(\tildhatY{i}{n})^{\top}.
		\end{aligned}
	\end{equation*}
	Notice that by Ostrowski's Theorem \cite[Theorem 4.5.9]{horn2012matrix} one obtains that 
	\begin{equation}\label{eq: equivalence S}
		\sigma^{k}\left(\widehat{S}_{n+1}\right)=	\sigma^{k}\left(\frac1M\sum_{i=1}^{M}\widehat{Y}_{n+1}^i(\widehat{Y}_{n+1}^i)^{\top}\right) = \sigma^{k}\left(\widehat{R}_{n+1}\frac1M\sum_{i=1}^{M}\tildhatY{i}{n+1}(\tildhatY{i}{n+1})^{\top} \widehat{R}_{n+1}^{\top}\right) \geq \sigma^{k}(\widetilde{\widehat{S}}_n),
	\end{equation}
	as $\widehat{R}_{n+1}\widehat{R}_{n+1}^{\top} \succeq I_{k \times k}$ by construction (see \eqref{eq: UnUn psd}).
	
	Let us now retrieve a recursion on $(\sigma^{k}(\widetilde{\widehat{S}}_n))_{n}$ via conditional expectation of $\widehat{S}_{n}$ and $\widetilde{\widehat{S}}_{n}$, and the Rayleigh quotient. We have
	\begin{equation*}
		\begin{aligned}
			\mathbb{E}_n\left[\widetilde{\widehat{S}}_{n+1}\right] = & \mathbb{E}_n\left[ \left(\frac1M\sum_{i=1}^{M} \tildhatY{i}{n+1}\left(\tildhatY{i}{n+1}\right)^{\top} \right)\right] \\
			= & \widehat{S}_n  + \frac1M\sum_{i=1}^{M}  \widehat{Y}_n^i(\widehat{a}_n^{i} )^{\top} \widehat{U}_n^{\top} \Delta t  + \frac1M\sum_{i=1}^{M} \widehat{U}_n\widehat{a}_n^{i} (\widehat{Y}_n^i)^{\top} \Delta t\\
			&+ \frac1M\sum_{i=1}^{M} \widehat{U}_n\widehat{a}_n^{i} (\widehat{a}_n^{i} )^{\top} \widehat{U}_n^{\top} (\Delta t)^2  + \frac1M\sum_{i=1}^{M} \widehat{U}_n\widehat{b}_n^{i}(\widehat{U}_n\widehat{b}_n^{i})^{\top} \Delta t.
		\end{aligned}
	\end{equation*}
	Using the following trivial identity
	\begin{equation*}
		\widetilde{\widehat{S}}_{n+1} = \mathbb{E}_n\left[\widetilde{\widehat{S}}_{n+1}\right]  +\widetilde{\widehat{S}}_{n+1} - \mathbb{E}_n\left[\widetilde{\widehat{S}}_{n+1}\right],
	\end{equation*}
	and taking $v \in \mathbb{R}^{k}$ with $|v|=1$, $v$ independent on $\boldsymbol{\omega}$, one has
	\begin{equation}\label{eq: inter S_n}
		\begin{aligned}
			v^{\top}\widetilde{\widehat{S}}_{n+1}v = &v^{\top}\mathbb{E}_n\left[\widetilde{\widehat{S}}_{n+1}\right]v  +  v^{\top}\left(\widetilde{\widehat{S}}_{n+1} - \mathbb{E}_n\left[\widetilde{\widehat{S}}_{n+1}\right]\right)v \\
			\geq & v^{\top}\widehat{S}_n v + \frac2M\sum_{i=1}^{M} v^{\top} \widehat{Y}_n^i(\widehat{a}_n^{i} )^{\top} \widehat{U}_n^{\top}v \Delta t + \frac1M\sum_{i=1}^{M} |v^{\top}\widehat{U}_n\widehat{a}_n^{i} |^2 (\Delta t)^2 \\
			& + \frac1M\sum_{i=1}^{M} v^{\top}\widehat{U}_n\widehat{b}_n^{i}(\widehat{U}_n\widehat{b}_n^{i})^{\top} v \Delta t + v^{\top} \left(\widetilde{\widehat{S}}_{n+1} - \mathbb{E}_n\left[\widetilde{\widehat{S}}_{n+1}\right]\right)v .\\
		\end{aligned}
	\end{equation}
	By Lemma \ref{lem: cond bound Mn}, we have that 
	$\frac1M\sum_{i=1}^{M} |\widehat{X}_n^{i} |^2 \leq \frac1M\sum_{i=1}^{M} |\widehat{Z}_n^{i} |^2$ and, therefore, on the event $E_n$, by the linear-growth bound, relation \eqref{eq: inter S_n} becomes 
	\begin{equation*}
		\begin{aligned}
			v^{\top}\widetilde{\widehat{S}}_{n+1}v 
			\geq & v^{\top}\widehat{S}_n v (1- \frac{\Delta t}{\varepsilon})  - \frac1M\sum_{i=1}^{M} C_{\mathrm{lgb}}(1+|\widehat{X}_n^{i} |^2) \left|(\Delta t)^2 - \varepsilon \Delta t\right| \\
			&+ \frac1M\sum_{i=1}^{M} v^{\top}\widehat{U}_n\widehat{b}_n^{i}(\widehat{U}_n\widehat{b}_n^{i})^{\top} v \Delta t + v^{\top}\left(\widetilde{\widehat{S}}_{n+1} - \mathbb{E}_n\left[\widetilde{\widehat{S}}_{n+1}\right]\right)v \\
			\geq & v^{\top}\widehat{S}_n v (1- \frac{\Delta t}{\varepsilon})  - \frac1M\sum_{i=1}^{M} C_{\mathrm{lgb}}(1+|\widehat{Z}_n^{i} |^2) \left|(\Delta t)^2 - \varepsilon \Delta t\right| \\
			&+ \frac1M\sum_{i=1}^{M} v^{\top}\widehat{U}_n\widehat{b}_n^{i}(\widehat{U}_n\widehat{b}_n^{i})^{\top} v \Delta t + v^{\top}\left(\widetilde{\widehat{S}}_{n+1} - \mathbb{E}_n\left[\widetilde{\widehat{S}}_{n+1}\right]\right)v.
		\end{aligned}
	\end{equation*}
	Moreover, by Remark \ref{rmk: link K_over and K} one has that $$\frac1M\sum_{i=1}^{M}\mathbb{E}_{n-1}\left[|\widehat{Z}_n^{i} |^2\right] \leq 2 K_2(T) + \left(\sum_{\ell =1}^{n-1} (1+ (1 +C_{\mathrm{lgb}} (2+T) )\Delta t )^{n-\ell}  \rho_{\ell,M}\right)$$
	if the event $\left\{\frac{1}{M} \sum_{i=1}^M |Y_0^i|^2 < 2 \mathbb{E}[|Y_0|^2] \right\}$ holds.
	For the sake of notation, let us recall that $\rho_{n,M}= \frac1M\sum_{i=1}^{M}|\widehat{Z}_n^{i} |^2- \mathbb{E}_{n-1}\left[|\widehat{Z}_n^{i} |^2\right]$, and $\eta_{n+1,M}=\widetilde{\widehat{S}}_{n+1} - \mathbb{E}_n\left[\widetilde{\widehat{S}}_{n+1}\right]$. Then, via choosing $\varepsilon = \frac{\sigma_{B}}{2(C_{\mathrm{lgb}}(1+2K_2(T))} + \Delta t$ and using Assumption \ref{ass: diff} one gets $\left|(\Delta t)^2 - \varepsilon \Delta t\right|  = A \Delta t$ and
	\begin{equation*}
		\begin{aligned}
			v^{\top}\widetilde{\widehat{S}}_{n+1}v \geq & v^{\top}\widehat{S}_n v (1- \frac{\Delta t}{\varepsilon})  - \left(\left( \frac1M\sum_{i=1}^{M} C_{\mathrm{lgb}}\left(1+\mathbb{E}_{n-1}\left[|\widehat{Z}_n^{i} |^2\right] \right)\right)+C_{\mathrm{lgb}} \rho_{n,M}\right)\left|(\Delta t)^2 - \varepsilon \Delta t\right| \\
			&+ \sigma_{B} \Delta t + v^{\top}\eta_{n+1,M}v\\
			\geq & v^{\top}\widehat{S}_n v (1- \frac{\Delta t}{\varepsilon})  - \frac1M\sum_{i=1}^{M} C_{\mathrm{lgb}}(1+2K_2(T)) \left|(\Delta t)^2 - \varepsilon \Delta t\right| \\
			&- C_{\mathrm{lgb}} \sum_{\ell=1}^{n}  (1+ (1 +C_{\mathrm{lgb}} (2+T) )\Delta t )^{n-\ell} \rho_{\ell,M}\left|(\Delta t)^2 - \varepsilon \Delta t\right| + \sigma_{B} \Delta t +v^{\top}\eta_{n+1,M}v\\
			\geq &  v^{\top}\widehat{S}_n v(1- \frac{\Delta t}{A + \Delta t}) - AC_{\mathrm{lgb}} \Delta t \left|\sum_{\ell=1}^{n} (1+ (1 +C_{\mathrm{lgb}} (2+T) )\Delta t )^{n-\ell}  \rho_{\ell,M}\right|  +  \frac{\sigma_{B}}{2}\Delta t +v^{\top}\eta_{n+1,M}v.\\
		\end{aligned}
	\end{equation*}
	Iterating over $n$ and recalling that $\gamma = 1 + (1+C_{\mathrm{lgb}}(2+T))\Delta t$ we obtain 
	\begin{equation*}
		\begin{aligned}
			v^{\top}\widetilde{\widehat{S}}_{n+1}v
			\geq & v^{\top}\widehat{S}_n v (1- \frac{\Delta t}{A + \Delta t})- AC_{\mathrm{lgb}}\exp\{(1 +C_{\mathrm{lgb}} (2+T) )n \Delta t\} \Delta t  \left|\sum_{\ell=1}^{n} \gamma^{-\ell}\rho_{\ell,M}\right|   \\
			&+  \frac{\sigma_{B}}{2}\Delta t +v^{\top}\eta_{n+1,M}v\\
			\geq & v^{\top} \widehat{S}_0 v (1- \frac{\Delta t}{A + \Delta t})^{n+1}+  \frac{\sigma_{B}}{2}\Delta t \sum_{\ell=0}^n  (1- \frac{\Delta t}{A + \Delta t})^{n-\ell} \\
			& -AC_{\mathrm{lgb}}  \Delta t \sum_{h=1}^{n}  \exp\{(1 +C_{\mathrm{lgb}} (2+T) )h \Delta t\} (1- \frac{\Delta t}{A + \Delta t})^{n-h}  \left|\sum_{\ell=1}^{h}  \gamma^{-\ell}\rho_{\ell,M}\right|  \\
			&- \left|\sum_{h=0}^n  (1- \frac{\Delta t}{A + \Delta t})^{n-h} v^{\top}\eta_{h+1,M}v \right| \\
			\geq & v^{\top} \widehat{S}_0 v (1- \frac{\Delta t}{A + \Delta t})^{n+1}+  \frac{\sigma_{B}}{2}\Delta t \sum_{\ell=0}^n  (1- \frac{\Delta t}{A + \Delta t})^{n-\ell} \\
			& -AC_{\mathrm{lgb}}  \Delta t \sum_{h=1}^{n} \exp\{(1 +C_{\mathrm{lgb}} (2+T) )h \Delta t\} (1- \frac{\Delta t}{A + \Delta t})^{n-h}  \left|\sum_{\ell=1}^{h}  \gamma^{-\ell} \rho_{\ell,M}\right|  \\
			&- \left|\sum_{h=0}^n  (1- \frac{\Delta t}{A + \Delta t})^{n-h} \eta_{h+1,M} \right| .
		\end{aligned}
	\end{equation*}
	Passing to the infimum with respect to $v$ in the term involving $\widehat{S}_n$ one gets
	\begin{equation*}
		\begin{aligned}
			\sigma^{k}\left(\widetilde{\widehat{S}}_{n+1}\right)
			\geq & \sigma_{\widehat{Y}_0}(1- \frac{\Delta t}{A + \Delta t})^{n+1}+  \frac{\sigma_{B}}{2}\Delta t \sum_{\ell=0}^n  (1- \frac{\Delta t}{A + \Delta t})^{n-\ell} - \left|\sum_{h=0}^n  (1- \frac{\Delta t}{A + \Delta t})^{n-h} \eta_{h+1,M}\right|\\
			& - AC_{\mathrm{lgb}}\Delta t \sum_{h=1}^{n}  \exp\{(1 +C_{\mathrm{lgb}} (2+T) )N\Delta t\} (1- \frac{\Delta t}{A + \Delta t})^{n-h}  \left| \sum_{\ell=1}^{h}  \gamma^{-\ell} \rho_{\ell,M}  \right| \\
			\geq & \sigma_{\widehat{Y}_0}(1- \frac{\Delta t}{A + \Delta t})^{n+1}+  \frac{\sigma_{B}}{2}\Delta t \sum_{\ell=0}^n  (1- \frac{\Delta t}{A + \Delta t})^{n-\ell} - \left|\sum_{h=0}^n  (1- \frac{\Delta t}{A + \Delta t})^{n-h} \eta_{h+1,M}\right|\\
			&- \frac{\sigma_{B}}{2 C_{\mathrm{lgb}}(1+2K_2(T))}C_{\mathrm{lgb}}\exp\{(1 +C_{\mathrm{lgb}} (2+T))T\} \Delta t \\
			& \quad \cdot  \sum_{h=0}^n (1- \frac{\Delta t}{A + \Delta t})^{n-h} \sup_{1\leq h \leq n}\left| \sum_{\ell=1}^h  \gamma^{-\ell} \rho_{\ell,M}  \right| . 
		\end{aligned}
	\end{equation*}
	Taking $\boldsymbol{\omega} \in E_n$ and using again the fact that $\left|1- \frac{\Delta t}{A + \Delta t}\right|<1$, where $A=\frac{\sigma_{B}}{2 C_{\mathrm{lgb}}(1+2K_2(T))}$, then we have
	\begin{equation*}
		\begin{aligned}
			\sigma^{k}\left(\widetilde{\widehat{S}}_{n+1}\right) \geq & \sigma_{Y_0}(1- \frac{\Delta t}{A + \Delta t})^{n+1}+  \frac{\sigma_{B}}{4}\Delta t \sum_{\ell=0}^n  (1- \frac{\Delta t}{A + \Delta t})^{n-\ell} - \frac{1}{2} \min\{\sigma_{Y_0}, \frac{\sigma_{B}^2}{8 C_{\mathrm{lgb}}(1+2K_2(T))} \}\\
			=& \sigma_{Y_0} \left(1 - \frac{1}{1+ \frac{A}{\Delta t}}\right)^{n+1} + \frac{\sigma_{B}}{4}\Delta t \sum_{s=0}^{n} \left(1 - \frac{1}{1+ \frac{A}{\Delta t}}\right)^{s} - \frac{1}{2} \min\{\sigma_{Y_0}, \frac{\sigma_{B}^2}{8 C_{\mathrm{lgb}}(1+2K_2(T))} \} \\
			= & \sigma_{Y_0} \left(1 - \frac{1}{1+ \frac{A}{\Delta t}}\right)^{n+1} + \frac{\sigma_{B}}{4}\Delta t \frac{1 - \left(1 - \frac{1}{1+ \frac{A}{\Delta t}}\right)^{n+1}}{1 -  \left(1 - \frac{1}{1+ \frac{A}{\Delta t}}\right)} - \frac{1}{2} \min\{\sigma_{Y_0}, \frac{\sigma_{B}^2}{8 C_{\mathrm{lgb}}(1+2K_2(T))} \} \\
			= & \sigma_{Y_0} \left(1 - \frac{1}{1+ \frac{A}{\Delta t}}\right)^{n+1} + \frac{\sigma_{B}}{4}(1+ \frac{A}{\Delta t})\Delta t \left[1 - \left(1 - \frac{1}{1+ \frac{A}{\Delta t}}\right)^{n+1}\right] \\
			&- \frac{1}{2} \min\{\sigma_{Y_0}, \frac{\sigma_{B}^2}{8 C_{\mathrm{lgb}}(1+2K_2(T))} \}\\
			= &\sigma_{Y_0}  \left(1 - \frac{1}{1+ \frac{A}{\Delta t}}\right)^{n+1} + \frac{\sigma_{B}}{4}A \left[1 - \left(1 - \frac{1}{1+ \frac{A}{\Delta t}}\right)^{n+1}\right]  + \frac{\sigma_{B}}{2}\Delta t \left[1 - \left(1 - \frac{1}{1+ \frac{A}{\Delta t}}\right)^{n+1}\right] \\
			& - \frac{1}{2} \min\{\sigma_{Y_0}, \frac{\sigma_{B}^2}{8 C_{\mathrm{lgb}}(1+2K_2(T))} \} \\
			\geq & \frac{1}{2} \min\{\sigma_{Y_0}, \frac{\sigma_{B}^2}{8 C_{\mathrm{lgb}}(1+2K_2(T))} \}+  \frac{\sigma_{B}}{2}\Delta t  \left[1-\left(1-\frac{1}{1+\frac{A}{\Delta t}}\right)^{n+1}\right]\\
			\geq & \min\{\frac{\sigma_{Y_0}}{2}, \frac{\sigma_{B}^2}{16C_{\mathrm{lgb}}(1+2K_2(T))}\},
		\end{aligned}
	\end{equation*}
	which implies the first statement.
	
	In order to prove the second statement, we show that we can bound the sought probability of the complement of $E_n$ by four components that can again be bounded by Monte-Carlo type errors. Indeed, we have
	\begin{equation*}
		\begin{aligned}
			\mathbb{P}\left(E_n^{C}\right) 
			\leq & \mathbb{P}\left( \frac{1}{M} \sum_{i=1}^M |Y_0^i|^2 > 2 \mathbb{E}[|Y_0|^2] \right) + \mathbb{P}\left( \frac{1}{M} \sum_{i=1}^M  Y_0^{i}(Y_0^{i})^{\top} \prec \frac{\sigma_{Y_0}}{2} \right)  \\
			& + \mathbb{P}\left( \sup_{1\leq h \leq n}\left|\sum_{\ell=0}^h  \gamma^{-\ell} \rho_{\ell,M} \right|  > \frac{1}{4}  \frac{2 (1+2K_2(T))}{\exp\left((1 +C_{\mathrm{lgb}} (2+T))T\right)}  \right) \\
			& + \mathbb{P}\left( \left|\sum_{\ell=0}^n  (1- \frac{\Delta t}{A + \Delta t})^{n-\ell} \eta_{\ell+1,M} \right| > \min\{ \frac{\varrho}{4}; \frac{\sigma_{Y_0}}{2} \} \right).
		\end{aligned}
	\end{equation*}
	For the first event, by construction we have by the Monte-Carlo estimator and Markov's inequality that
	\begin{equation*}
		\mathbb{P}\left( \frac{1}{M} \sum_{i=1}^M |Y_0^i|^2 > 2 \mathbb{E}[|Y_0|^2] \right) \leq 	\mathbb{P}\left( \left| \frac{1}{M} \sum_{i=1}^M |Y_0^i|^2 - \mathbb{E}[|Y_0|^2] \right| > \mathbb{E}[|Y_0|^2] \right) \leq C_p \frac{1}{M^{p}}
	\end{equation*}
	for a positive constant $C_p$ depending on $\mathbb{E}[|Y_0|^{4p}]$. For the second event, Lemma \ref{lem: prop bound sigma tilde} gives us the same upper bound in Monte-Carlo terms.
	
	We need to estimate the following probability
	\begin{equation*}
		\mathbb{P}\left( \sup_{1\leq h \leq n}\left|\sum_{\ell=1}^h \gamma^{-\ell} \rho_{\ell,M}\right|   > \frac{1}{4}  \frac{2 (1+2K_2(T))}{\exp\left((1 +C_{\mathrm{lgb}} (2T+T^2))\right)} \right) 
	\end{equation*}
	First notice that by Markov's inequality one has that 
	\begin{equation*}
		\begin{aligned}
			&\mathbb{P}\left(  \sup_{1\leq h \leq n}\left|\sum_{\ell=1}^h  \gamma^{-\ell} \rho_{\ell,M}\right|  > \frac{1}{4}  \frac{2 (1+2K_2(T))}{ \exp\left((1 +C_{\mathrm{lgb}} (2+T))T\right)} \right)\\ \leq&  \mathbb{E}\left[\sup_{1\leq h \leq n}\left|\sum_{\ell=1}^h  \gamma^{-\ell} \rho_{\ell,M} \right|^{2p}\right] \left(\frac{4 \exp\{(1 +C_{\mathrm{lgb}} (2+T))T\}}{2 (1+2K_2(T))}\right)^{2p}.\\
		\end{aligned}
	\end{equation*}
	Then, using Lemma \ref{lem: Local MC Mn}, we know that there exists a positive constant $C_3:=C_3(p,k,T)$ independent of the smallest singular value of the fully-discrete Gramian such that
	\begin{equation*}
		\mathbb{P}\left(  \sup_{1\leq h \leq n}\left|\sum_{\ell=1}^h  \gamma^{-\ell} \rho_{\ell,M} \right| > \frac{1}{4}  \frac{2 (1+2K_2(T))}{ \exp\left((1 +C_{\mathrm{lgb}} (2+T))T\right)} \right) \leq  C_3  \frac{1}{M^p}.\\
	\end{equation*}
	
	To conclude, we want to estimate $ \mathbb{P}\left(  \left| \sum_{\ell=0}^n  (1- \frac{\Delta t}{A + \Delta t})^{n-\ell} \eta_{n+1,M}\right| > \frac{\varrho}{4} \right)$. Again, via Markov's inequality one has that  
	\begin{equation*}
		\begin{aligned}
			&\mathbb{P}\left(   \left|\sum_{\ell=0}^n  (1- \frac{\Delta t}{A + \Delta t})^{n-\ell} \eta_{\ell+1,M} \right| > \min\{ \frac{\varrho}{4}, \frac{\sigma_{Y_0}}{2}\} \right)\\
			\leq & \mathbb{E}\left[\left|\sum_{\ell=0}^n  (1- \frac{\Delta t}{A + \Delta t})^{n-\ell} \eta_{\ell+1,M} \right|^{2p}\right] \left(\max\left\{\frac{4(4C_{\mathrm{lgb}}(1+2K_2(T)))}{\sigma_B^2}; \frac{2}{\sigma_{Y_0}}\right\}\right)^{2p}.\\
		\end{aligned}
	\end{equation*}
	Via Lemma \ref{lem: p-mom cond covariance}, there exists a positive constant $C_{p}$ independent of the smallest singular value of the fully-discrete Gramian such that
	\begin{equation*}
		\begin{aligned}
			\mathbb{E}\left[\left|\sum_{\ell=0}^n  (1- \frac{\Delta t}{A + \Delta t})^{n-\ell} \eta_{\ell+1,M} \right|^{2p}\right] \leq & C_{p}\widetilde{C}\frac{ (\Delta t)^p}{M^{p}} \left( \sum_{\ell=0}^n (1- \frac{\Delta t}{A + \Delta t})^{2(n-\ell)}\right)^{p} \\
			\leq & C_{p}\widetilde{C}\frac{ (\Delta t)^p}{M^{p}} \left(\sum_{\ell=0}^n (\frac{A}{A + \Delta t})^{2(n-\ell)}\right)^{p} \\
			\leq & C_{p}\widetilde{C}\frac{ (\Delta t)^p}{M^{p}} \left(\frac{1-(\frac{A}{A + \Delta t})^{2(n+1)}}{1-(\frac{A}{A + \Delta t})^{2}} \right)^{p} \\
			\leq & C_{p}\widetilde{C}\frac{ (\Delta t)^p}{M^{p}} \left(\frac{1}{1-(\frac{A}{A + \Delta t})^{2}} \right)^{p} \\
			\leq & C_{p}\widetilde{C}\frac{ (\Delta t)^p}{M^{p}} \left(\frac{1}{1-(\frac{A}{A + \Delta t})} \right)^{p} \\
			\leq & C_{p}\widetilde{C}\frac{ (\Delta t)^p}{M^{p}} \left( \left(\frac{A+ \Delta t}{\Delta t} \right)^{p}\right) =   \widetilde{C}_{p}\widetilde{C}\frac{ 1}{M^{p}} \left(A^p + \Delta t^p\right), \\
		\end{aligned}
	\end{equation*}
	for a positive constant $\widetilde{C}_{p}$. Then, recalling $A=\frac{\sigma_{B}}{2 C_{\mathrm{lgb}}(1+2K_2(T))}$, for a positive constant $C_4:=C_4(p,T)$  we finally obtain
	\begin{equation}\label{eq: last M bound}
		\begin{aligned}
			\mathbb{P}\left( \left|\sum_{\ell=0}^n  (1- \frac{\Delta t}{A + \Delta t})^{n-\ell} \eta_{\ell+1,M} \right| >  \min\{ \frac{\varrho}{4}, \frac{\sigma_{Y_0}}{2}\}  \right) \leq  C_4 \max\left\{\left( \frac{1}{\sigma_B^{3p}M^p} + \frac{\Delta t^p }{\sigma_B^{4p}M^p}\right); \frac{\sigma_B^{p} + \Delta t^p}{\sigma_{Y_0}^{2p}M^{p}} \right\}. \\
		\end{aligned}
	\end{equation}
	Finally, notice that \eqref{eq: M dt cond} implies $\Delta t^{p} \leq \frac{1}{M^{p-1}} \frac{1}{k^p}$. Using this relation in \eqref{eq: last M bound} yields the thesis.
\end{proof}

\section{Proofs of Projection-type bounds}\label{app: proof proj}

\begin{proof}[\bfseries Proof of Propositon \ref{Proposition: exp tails Xn}]\label{proof: Proposition: exp tails Xn}
	We want to prove that there exist two constants $C,c>0$, independent of $N$ such that for all $R>0$
	\begin{equation*}
		\begin{aligned}
			\mathbb{P}\left(\max_{0\le n\le N}|X_n|\ge R\right)\le Ce^{-cMR^2}.
		\end{aligned}
	\end{equation*}
	In order to achieve this result, we express $X_n$ as a sum of increments of drift and diffusion, then we prove that the diffusion-type terms have subgaussian tails, and we conclude using this tail relation to obtain a probability bound uniform in $N$.
	
	Iterating relation \eqref{eq: Stoc Proj X}, one can write the following relation on the semidiscrete DLRA:
	\begin{equation}\label{eq: sum X_n abs value}
		\begin{aligned}
			\max\limits_{1 \leq \ell \leq n+1} |X_\ell| &=	\max\limits_{0 \leq \ell \leq n} |X_{\ell+1}| \\
			& = \max\limits_{0 \leq \ell \leq n}|X_{\ell} +  P_{U^{\top}_\ell\widetilde{Y}_{\ell+1}}[a_\ell ] \Delta t_\ell + P_{U_\ell} [b_\ell] \Delta W_\ell |\\
			& = \max\limits_{0 \leq \ell \leq n}| X_{0} + \sum_{k=0}^{\ell} P_{U^{\top}_k\widetilde{Y}_{k+1}}[ a_k ] \Delta t_k + \sum_{k=0}^{\ell} P_{U_k} [b_k] \Delta W_k |\\
			& \leq |X_{0}| +  \max\limits_{0 \leq \ell \leq n} |\sum_{k=0}^{\ell} P_{U^{\top}_k\widetilde{Y}_{k+1}}[ a_k ] \Delta t_k | + \max\limits_{0 \leq \ell \leq n}  |\sum_{k=0}^{\ell} P_{U_k} [b_k] \Delta W_k | \\
			&=|X_{0}| +  \max\limits_{0 \leq \ell \leq n} \left( \left| \sum_{k=0}^{\ell} P_{U_k}[ a_k ] \Delta t_k+  P_{U_k}^{\perp}P_{\widetilde{Y}_{k+1}}[ a_k ] \Delta t_k \right|\right)  + \max\limits_{0 \leq \ell \leq n}  |\sum_{k=0}^{\ell} P_{U_k} [b_k] \Delta W_k | \\
			& \leq |X_{0}| +  \max\limits_{0 \leq \ell \leq n} \left( \sum_{k=0}^{\ell}   |a_k | \Delta t_k+ \left|\sum_{k=0}^{\ell}  P_{U_k}^{\perp}P_{\widetilde{Y}_{k+1}}[ a_k ] \Delta t_k \right|\right)  + \max\limits_{0 \leq \ell \leq n}  |\sum_{k=0}^{\ell} P_{U_k} [b_k] \Delta W_k |,
		\end{aligned}
	\end{equation}
	where we employed the subadditivity of the norm and the fact that the orthogonal projectors have operator norm equal to $1$. To deal with the term concerning the stochastic projector $P_{\widetilde{Y}_{k+1}}$ notices that
	\begin{equation}\label{eq: Proj abs val}
		\begin{aligned}
			P_{\widetilde{Y}_{k+1}}[ a_k ]  & = \mathbb{E}[a_k \widetilde{Y}_{k+1}] C_{\widetilde{Y}_{k+1}}^{-1} ( Y_k + U_k a_k \Delta t_k + U_k b_k \Delta W_k).
		\end{aligned}
	\end{equation}
	Assumption \ref{ass: diff} guarantees that the smallest singular value of the semidiscrete solution $\widetilde{\sigma}_{n+1}^k$ is lower bounded from below by a positive constant (see \cite[Proposition 5.2]{kazashi2025dynamicalpartI} and relation \eqref{eq: UnUn psd}). Notice that one can substitute Assumption \ref{ass: diff} with any other assumption that provides such a uniform positive lower bound. This observation allows us to control the term \eqref{eq: Proj abs val} uniformly.
	
	Plugging relation \eqref{eq: Proj abs val} into \eqref{eq: sum X_n abs value} and using the linear-growth bound, one obtains
	\begin{equation}\label{eq: sum X_n abs value 2}
		\begin{aligned}
			\max\limits_{1 \leq \ell \leq n+1} |X_\ell| 
			\leq& |X_{0}| +   \max\limits_{0 \leq \ell \leq n} \left( \sum_{k=0}^{\ell}   |a_k | \Delta t_k+ \left|\sum_{k=0}^{\ell}  P_{U_k}^{\perp} \left(\mathbb{E}[a_k \widetilde{Y}_{k+1}] C_{\widetilde{Y}_{k+1}}^{-1} ( Y_k + U_k a_k \Delta t_k + U_k b_k \Delta W_k)\right) \Delta t_k \right|\right)  \\
			&+ \max\limits_{0 \leq \ell \leq n}  |\sum_{k=0}^{\ell} P_{U_k} [b_k] \Delta W_k | \\
			\leq &  |X_{0}| +   \max\limits_{0 \leq \ell \leq n} \left( \sum_{k=0}^{\ell}   |a_k | \Delta t_k+ \sum_{k=0}^{\ell}  \left|\mathbb{E}[a_k \widetilde{Y}_{k+1}] C_{\widetilde{Y}_{k+1}}^{-1} ( Y_k + U_k a_k \Delta t_k) \right| \Delta t_k \right) \\
			&+  \max\limits_{0 \leq \ell \leq n}  \left| \sum_{k=0}^{\ell} \mathbb{E}[a_k \widetilde{Y}_{k+1}] C_{\widetilde{Y}_{k+1}}^{-1} U_k b_k \Delta W_k  \Delta t_k \right| + \max\limits_{0 \leq \ell \leq n}  |\sum_{k=0}^{\ell} P_{U_k} [b_k] \Delta W_k | \\
			\leq &  |X_{0}| + \max\limits_{0 \leq \ell \leq n} \Big( \sum_{k=0}^{\ell} C_{\mathrm{lgb}} (1+ |X_k |) \Delta t_k \\
			&+ \sum_{k=0}^{\ell}   (\widetilde{\sigma}_{n+1}^k)^{-\frac{1}{2}} \sqrt{C_{\mathrm{lgb}}(1+K_2(T))}  \left(|X_k| + C_{\mathrm{lgb}} (1+ |X_k |)  \Delta t_k \right) \Delta t_k  \Big)\\
			&+  \max\limits_{0 \leq \ell \leq n} \left(  \left| \sum_{k=0}^{\ell} \mathbb{E}[a_k \widetilde{Y}_{k+1}] C_{\widetilde{Y}_{k+1}}^{-1} U_k b_k \Delta W_k  \Delta t_k \right|\right) + \max\limits_{0 \leq \ell \leq n}  |\sum_{k=0}^{\ell} P_{U_k} [b_k] \Delta W_k | \\
		\end{aligned}
	\end{equation}
	where in the last inequality we employed the Cauchy-Schwarz inequality, and we used the fact that $U_n$ has orthogonal rows, and hence, $U_n^{\top}$ has orthogonal columns and preserves the Euclidean norm. 
	
	Now, for the sake of readability, let us define
	\begin{equation}\label{eq: martingales sum}
		\begin{aligned}
			M^1_n:=\sum_{k=0}^{n} \mathbb{E}[a_k \widetilde{Y}_{k+1}] C_{\widetilde{Y}_{k+1}}^{-1} U_k b_k \Delta W_k  \Delta t_k,\qquad 
			M^2_n:=\sum_{k=0}^{n}P_{U_k}b_k\Delta W_k.
		\end{aligned}
	\end{equation}
	Then, \(M^1_n\) and \(M^2_n\) are martingales with respect to the filtration \((\mathcal F_{t_n})_{n=0}^N\) as \(\mathbb{E}[a_n \widetilde{Y}_{n+1}] C_{\widetilde{Y}_{n+1}}^{-1}\),  \(U_nb_n\), and \(P_{U_n}b_n\) are \(\mathcal F_{t_n}\)-measurable, and \(\Delta W_n\) is independent of \(\mathcal F_{t_n}\), with mean zero. Thus
	\begin{equation*}
		\begin{aligned}
			\mathbb{E}[M^1_{n+1}-M^1_n\mid \mathcal F_{t_n}] =0, \qquad 	\mathbb{E}[M^2_{n+1}-M^2_n\mid \mathcal F_{t_n}] =0.
		\end{aligned}
	\end{equation*}
	
	Set
	\begin{equation*}
		\begin{aligned}
			R_n:=\max_{0\le j\le n}|X_j|,
			\qquad
			M_n^{1,*}:=\max_{0\le \ell \le n}|M^1_\ell|, \qquad M_n^{2,*}:=\max_{0\le \ell \le n}|M^2_\ell|.
		\end{aligned}
	\end{equation*}
	From \eqref{eq: sum X_n abs value 2} we then have 
	\begin{equation*}
		\begin{aligned}
			R_{n+1} \le& |X_0|+ C_{\mathrm{lgb}} \left(1 +\max_{0 \leq \ell \leq n+1} (\widetilde{\sigma}_{\ell}^k)^{-\frac{1}{2}} \sqrt{C_{\mathrm{lgb}}(1+K_2(T))}    \right) T \\
			&+  \left( C_{\mathrm{lgb}}  + (1 + C_{\mathrm{lgb}}  T)\max_{0 \leq \ell \leq n+1} (\widetilde{\sigma}_{\ell}^k)^{-\frac{1}{2}} \sqrt{C_{\mathrm{lgb}}(1+K_2(T))} \right) \sum_{k=0}^{n-1}R_k \Delta t_k +M_n^{1,*} +M_n^{2,*} \\
		\end{aligned}
	\end{equation*}
	Let
	\begin{equation*}
		\begin{aligned}
			A:=|X_0|+\underbrace{C_{\mathrm{lgb}} \left(1 +\sup_{0 \leq \ell \leq N} (\widetilde{\sigma}_{\ell}^k)^{-\frac{1}{2}} \sqrt{C_{\mathrm{lgb}}(1+K_2(T))}   \right)}_{=:M} T+M_N^{1,*} +M_N^{2,*}.
		\end{aligned}
	\end{equation*}
	Since \(M_n^{1,*}\le M_N^{1,*}\) and  \(M_n^{2,*}\le M_N^{2,*}\), for $K=\left(C_{\mathrm{lgb}}  + (1 + C_{\mathrm{lgb}}  T)\sup_{0 \leq \ell \leq N} (\widetilde{\sigma}_{\ell}^k)^{-\frac{1}{2}} \sqrt{C_{\mathrm{lgb}}(1+K_2(T))}  \right)$ we have
	\begin{equation}\label{eq:pre-gronwall}
		\begin{aligned}
			R_n\le A+K\sum_{k=0}^{n-1}R_k\Delta t_k,
		\end{aligned}
	\end{equation}
	and via the discrete Gronwall lemma we obtains
	\begin{equation}
		\begin{aligned}
			\max_{0\le n\le N}|X_n|
			\le e^{KT}\left(|X_0|+MT+\max_{0\le n\le N}|M_n^{1}|+ \max_{0\le n\le N}|M_n^{2}|\right).
			\label{eq:gronwall-bound}
		\end{aligned}
	\end{equation}
	
	We now estimate the martingale terms, similarly to the proof of \cite[Proposition 8.7]{baldi2017stochastic}. We prove this estimate for $M_n^2$, a similar result can be obtained for $M_n^1$ verbatim. Fix $\theta\in \mathbb{R}^{d}$ such that $|\theta|=1$, and define the scalar martingale $M_n^{2,\theta}:=\langle \theta,M_n\rangle.$
	Then one has that
	\begin{equation*}
		\begin{aligned}
			M_{n+1}^{2,\theta}-M_n^{2,\theta}=\langle \theta,P_{U_n}b_n\Delta W_n\rangle ,
		\end{aligned}
	\end{equation*}
	which conditionally on \(\mathcal F_{t_n}\), is a centered Gaussian with variance
	\begin{equation*}
		\begin{aligned}
			v_n^\theta
			&:=\mathbb E\left[\left(M_{n+1}^{2,\theta}-M_n^{2,\theta}\right)^2\mid \mathcal F_{t_n}\right] =\Delta t_n\,|\theta^\top P_{U_n}b_n|^2 \leq \Delta t_n\,\|b_n\|_{\mathrm{F}}^2  \le K_{b}^2 \Delta t_n,
		\end{aligned}
	\end{equation*}
	thanks to Assumption \ref{ass: bound diff and tails}. Hence its predictable quadratic variation satisfies
	\begin{equation} 	\label{eq:qv-bound}
		\begin{aligned}
			\langle M^{2,\theta}_N\rangle:=\sum_{n=0}^{N-1}v_n^\theta \le K_{b}^2T.
		\end{aligned}
	\end{equation}
	
	For \(\mu>0\), define
	\begin{equation*}
		\begin{aligned}
			Z_n^{\mu}:=\exp\left(\mu M_n^{2,\theta}-\frac{\mu^2}{2}\langle M^{2,\theta}_n\rangle\right).
		\end{aligned}
	\end{equation*}
	We prove that \(Z_n^{\mu}\) is a martingale. Since \(M_{n+1}^{2,\theta}-M_n^{2,\theta}\), conditionally on \(\mathcal F_{t_n}\), is Gaussian with mean \(0\) and variance \(v_n^\theta\), we have
	\begin{equation*}
		\begin{aligned}
			\mathbb E\left[\exp\left(\mu (M_{n+1}^{2,\theta}-M_n^{2,\theta})\right)
			\mid \mathcal F_{t_n}\right]=\exp\left(\frac{\mu^2}{2}v_n^\theta\right).
		\end{aligned}
	\end{equation*}
	Therefore
	\begin{equation*}
		\begin{aligned}
			\mathbb E[Z_{n+1}^{\mu}\mid \mathcal F_{t_n}]
			&=
			\mathbb E\left[
			\exp\left(
			\mu M_{n+1}^{2,\theta}-\frac{\mu^2}{2}\langle M^{2,\theta}_{n+1}\rangle
			\right)
			\mid \mathcal F_{t_n}
			\right] \\
			&=
			\exp\left(
			\mu M_n^{2,\theta}-\frac{\mu^2}{2}\langle M^{2,\theta}_n\rangle-\frac{\mu^2}{2}v_n^\theta
			\right)
			\mathbb E\left[
			\exp\left(\mu (M_{n+1}^{2,\theta}-M_n^{2,\theta})\right)
			\mid \mathcal F_{t_n}
			\right] \\
			&=
			\exp\left(
			\mu M_n^{2,\theta}-\frac{\mu^2}{2}\langle M^{2,\theta}_n\rangle
			\right) =Z_n^{\mu}.
		\end{aligned}
	\end{equation*}
	Thus \(Z_n^\mu\) is a non-negative martingale and, as \(\mathbb E[Z_0^\mu]=1\), one has \(\mathbb E[Z_n^\mu]=1\), for all $n$.
	Let \(r>0\). On the event
	\begin{equation*}
		\begin{aligned}
			\left\{\max_{0\le n\le N}M_n^{2,\theta}\ge r\right\},
		\end{aligned}
	\end{equation*}
	there exists \(n\le N\) such that \(M_n^{2,\theta}\ge r\). Since \(\langle M^{2,\theta}_n\rangle\le K_{b}^2T\), by \eqref{eq:qv-bound}, one has
	\begin{equation*}
		\begin{aligned}
			Z_n^{\mu}
			\ge
			\exp\left(\mu r-\frac{\mu^2}{2}K_b^2T\right).
		\end{aligned}
	\end{equation*}
	Hence
	\begin{equation*}
		\begin{aligned}
			\mathbb{P}\left(\max_{0\le n\le N}M_n^{2,\theta}\ge r\right)
			&\le
			\mathbb{P}\left(
			\max_{0\le n\le N}Z_n^{\mu}
			\ge
			\exp\left(\mu r-\frac{\mu^2}{2}K_b^2T\right)
			\right).
		\end{aligned}
	\end{equation*}
	By Doob's maximal inequality for non-negative submartingales \cite[Theorem 4.4.2]{durrett2019probability},
	\begin{equation*}
		\begin{aligned}
			\mathbb{P}\left(\max_{0\le n\le N}M_n^{2,\theta}\ge r\right) 
			&\le
			\exp\left(-\mu r+\frac{\mu^2}{2}K_b^2T\right)\mathbb E[Z_N^\mu]=
			\exp\left(-\mu r+\frac{\mu^2}{2}K_b^2T\right).
		\end{aligned}
	\end{equation*}
	As the first left-hand side is independent of $\mu$, by minimizing over $\mu$ the right-hand side choosing $\mu=\frac{r}{K_{b}^2T}$, we finally have
	\begin{equation}
		\begin{aligned}
			\mathbb{P}\left(\max_{0\le n\le N}M_n^{2,\theta}\ge r\right)
			\le
			\exp\left(-\frac{r^2}{2K_b^2T}\right).
			\label{eq:scalar-tail}
		\end{aligned}
	\end{equation}
	
	We now pass from scalar to vector estimates. Notice that if $\max_{0\le n\le N}|M_n^2| \geq r$, then necessarily there exists an $i$ in $1,\dots,d$, such that $\max_{0\le n\le N}|M^2_{n,i}|^{2} \geq \frac{r}{\sqrt{d}}$, where $M^2_{n,i}$ is the $i$-th component of $M^2_n$. Considering $\theta$ as a canonical vector basis in $\mathbb{R}^d$, we get 
	\begin{equation*}
		\begin{aligned}
			\mathbb{P}\left(\max_{0\le n\le N}|M^2_n| \ge r\right) \leq  \sum_{i=1}^{d} \mathbb{P}\left(\max_{0\le n\le N}|M^2_{n,i}|^{2} \geq \frac{r}{\sqrt{d}} \right) \leq 2d \exp\left(-\frac{r^2}{2K_b^2 d T}\right),
		\end{aligned}
	\end{equation*}
	where the right-hand side is independent of \(N\). Similar estimates hold for $M^{1}_n$.
	
	We now combine these estimates on the tails of the martingales $M_n^{1}$ and $M_n^{2}$ to conclude. From \eqref{eq:gronwall-bound},
	\begin{equation*}
		\begin{aligned}
			\left\{\max_{0\le n\le N}|X_n|\ge R\right\}
			&\subset	\left\{
			|X_0|+MT+\max_{0\le n\le N}|M_n^{1}| + \max_{0\le n\le N}|M_n^{2}|
			\ge e^{-KT}R
			\right\} \\
			&=
			\left\{
			|X_0|+\max_{0\le n\le N}|M_n^{1}| + \max_{0\le n\le N}|M_n^{2}|
			\ge e^{-KT}R-MT
			\right\}.
		\end{aligned}
	\end{equation*}
	Assume first that $R\ge 2MTe^{KT}=R_0$, then
	\begin{equation*}
		\begin{aligned}
			e^{-KT}R-MT\ge \frac12e^{-KT}R.
		\end{aligned}
	\end{equation*}
	Thus
	\begin{equation*}
		\begin{aligned}
			\left\{
			\max_{0\le n\le N}|X_n|\ge R
			\right\}
			&\subset
			\left\{
			|X_0|+\max_{0\le n\le N}|M_n^{1}| + \max_{0\le n\le N}|M_n^{2}|
			\ge \frac12e^{-KT}R
			\right\} \\
			&\subset
			\left\{
			|X_0|\ge \frac14e^{-KT}R
			\right\}
			\cup
			\left\{
			\max_{0\le n\le N}|M_n^{1}| \ge \frac18e^{-KT}R
			\right\} 	\cup
			\left\{
			\max_{0\le n\le N}|M_n^{2}|\ge \frac18e^{-KT}R
			\right\}.
		\end{aligned}
	\end{equation*}
	Therefore, by Assumption \ref{ass: bound diff and tails} we obtain
	\begin{equation}
		\begin{aligned}
			\mathbb{P}\left(\max_{0\le n\le N}|X_n|\ge R\right)
			&\le
			C_0\exp\left(-\frac{\mu}{16}e^{-2KT}R^2\right)
			+
			4d \exp\left(-\frac{e^{-2KT}R^2}{2K_b^2 64 d T}\right).
			\label{eq:large-R-tail}
		\end{aligned}
	\end{equation}
	Enlarging constants to also cover the case for \(0\le R\le R_0\), we conclude that there exist \(C_{T},c_{T}>0\), depending on $K_{b},T,d$, but not on $N$ such that
	\begin{equation*}
		\begin{aligned}
			\mathbb{P}\left(\max_{0\le n\le N}|X_n|\ge R\right)\le C_{T}e^{-c_{T}R^2},\qquad \text{for all } R\ge 0,
		\end{aligned}
	\end{equation*}
	i.e.\ $\max\limits_{0\le n\le N}|X_n|$ has subgaussian tails uniformly in $N$.
\end{proof}

Using Proposition \ref{Proposition: exp tails Xn}, we can now prove a Bernstein's type inequality, that will be of help to establish probability estimates on the tails of the particle system distribution.
\begin{Lemma}[Bernstein inequality]\label{cor: Bernstein ineq}
	Consider $f\in L^4_{iid}\left(\Omega^{M}, \mathbb{R}^{d \times M}\right)$ such that for positive constants $\tilde{C},\tilde{c}$, one has $\mathbb{P}(|f^{i}| > R) \leq \tilde{C} e^{-\tilde{c}R^2}$ for all $R>0$. Then there exist positive constants $c,C$ such that one has
	\begin{equation*}
		\mathbb{P}\left( \boldsymbol{\omega}:=(\omega_1, \dots, \omega_M)\ : \ \frac{1}{M} \sum_{i = 1}^M |f^i(\omega_i)|^2 > R^2\right) \leq C \exp\{ - c M R^2 \},
	\end{equation*}
	for all $R > \|f\|_{L^2(\Omega^M;\mathbb{R}^{d \times M})}$.
	\begin{proof}
		The proof follows from the Bernstein inequality \cite[Theorem 2.9.1]{vershynin2018high}. Indeed, as $f^i$ has subgaussian tails, then $|f^i|^2$ has exponential tails and, hence, one has that 
		\begin{equation*}
			\mathbb{P}\left( \boldsymbol{\omega} \ : \ \left|\frac{1}{M} \sum_{i = 1}^M |f^i|^2- \mathbb{E}[\frac{1}{M} \sum_{i = 1}^M |f^i|^2] \right| > R_1^2 \right) \leq C_1 \exp\{ - c_1 M \min \{R_1^4, R_1^2\} \},
		\end{equation*}
		for positive constant $c_1, C_1$. Then 
		\begin{equation*}
			\begin{aligned}
				&\mathbb{P}\left( \boldsymbol{\omega} \ : \frac{1}{M} \sum_{i = 1}^M |f^i|^2  > R_1^2+ \mathbb{E}[\frac{1}{M} \sum_{i = 1}^M |f^i|^2]\right)\\
				& \leq \mathbb{P}\left( \boldsymbol{\omega} \ : \ \left|\frac{1}{M} \sum_{i = 1}^M |f^i|^2- \mathbb{E}[\frac{1}{M} \sum_{i = 1}^M |f^i|^2] \right| > R_1^2\right) \\
				&\leq  C_2 \exp\{ - c_2 M \min \left\{(R_1^2 + \mathbb{E}[\frac{1}{M} \sum_{i = 1}^M |f^i|^2])^2 , R_1^2 + \mathbb{E}[\frac{1}{M} \sum_{i = 1}^M |f^i|^2]\} \right\},
			\end{aligned}
		\end{equation*}
		via enlarging the constants $c_1$ and $C_1$ to $c_2$ and $C_2$, respectively. This last inequality implies the thesis.
	\end{proof}
\end{Lemma}

\begin{proof}[\bfseries Proof of Lemma \ref{lem: diff Proj U}]\label{proof: lem: diff Proj U}
	To prove the claim, we will bound the error in \eqref{eq: bound between U proj} with $e_n= \sqrt{\mathbb{E}\left[\frac{1}{M}\sum_{i=1}^{M} |X_n(\omega_i)-\widehat{X}_n^{i}|^2\right]}$, i.e.\ the error between the semi-discrete DLRA solution and the fully discrete one, plus, via norm inequalities, the error between the  semi-discrete DLRA and a system of $M$-particles extracted from it, from which we expect a Monte-Carlo-type error. The proof will be based on a localization procedure in order to exploit the subgaussian tails of $f$.
	
	We prove the estimate by first reducing the projector error to an operator-norm bound. The first part of this proof follows closely the one developed in \cite[Lemmata A.2 and A.3]{kazashi2025dynamical}. 
	
	Let us consider the singular-value decomposition of 
	$\mathbb{X}_n(\boldsymbol{\omega})=U_n^{\top}\tilde{\Sigma}_n \tilde{V}_n^{\top}$, with $U_n$ deterministic, whereas $\tilde{\Sigma}_n$ and $\tilde{V}_n^{\top}$ depends on $\boldsymbol{\omega}$. We have that $V_n =  \sqrt{M} \tilde{V}_n \in \mathbb{R}^{M \times k}$ is a matrix whose columns are orthonormal in $\left(\mathbb{R}^{M}, \widehat{\mathbb{E}}[\ \cdot , \cdot \ ]\right)$, where $\widehat{\mathbb{E}}[xy] = \frac{1}{M} \langle x,y\rangle_{\mathbb{R}^{M}}$, with $x,y\in \mathbb{R}^M$. Moreover, let us define $\breve{\mathbb{X}}_n := U^{\top}_n \Sigma_n V^{\top}_n$, where $\Sigma_n= \mathrm{diag}(\sqrt{\sigma_n^1}, \dots, \sqrt{\sigma_n^k})$ is a diagonal matrix whose elements on the diagonal are the square root of the singular values of the Gramian for the semi-discrete process $X_n(\omega)$, i.e.\ the singular values of $\mathbb{E}[X_n(\omega)X_n(\omega)^{\top}]$. Notice that $\mathbb{X}_n$ and $\breve{\mathbb{X}}_n$ share the same ($\boldsymbol{\omega}-$independent) $U_n$, the same $V_n$ dependent on $\boldsymbol{\omega}$ up to a rescaling factor, but unlike for $\mathbb{X}_n$, the singular values of $\breve{\mathbb{X}}_n$ do not depend on $\boldsymbol{\omega}$.
	Then, following similar steps to the ones of the proof of \cite[Lemma A.2]{bachmayr2021existence} we have pathwise that
	\begin{equation*}
		|\left(P_{U_n}-P_{\widehat{U}_n}\right)f^{i}| \leq \|(P_{\widehat{U}_n})^{\perp}P_{U_n}\|_{\mathbb{R}^{d} \to \mathbb{R}^{d}}|f^{i}|.
	\end{equation*}
	In this regard, we want to bound $ \|(P_{\widehat{U}_n})^{\perp}P_{U_n}\|_{\mathbb{R}^{d} \to \mathbb{R}^{d}}$. It holds that
	\begin{equation}\label{eq: PuPu perp}
		\begin{aligned}
			\|(P_{\widehat{U}_n})^{\perp}P_{U_n}\|_{\mathbb{R}^{d} \to \mathbb{R}^{d}} 
			& = \sup\limits_{\substack{g \in \mathbb{R}^{d} \\ |g| = 1}} \left|\left(I_{d \times d}-P_{\widehat{U}_n}\right) \sum_{j=1}^{k } U_n^j \langle  U_n^j, g\rangle \right|\\
			& = \sup\limits_{\substack{g \in \mathbb{R}^{d} \\ |g| = 1}} \left|\left(I_{d \times d}-P_{\widehat{U}_n}\right) \sum_{j=1}^{k }\frac{\widehat{\mathbb{E}}[\breve{\mathbb{X}}_n^{j}V_n^j]}{\sqrt{\sigma_n^j}} \langle  U_n^j, g\rangle \right| \\
			& = \sup\limits_{\substack{g \in \mathbb{R}^{d} \\ |g| = 1}} \left|\left(I_{d \times d}-P_{\widehat{U}_n}\right) \sum_{j=1}^{k } \frac{1}{M}\frac{\breve{\mathbb{X}}_n V_n^j}{\sqrt{\sigma_n^j}} \langle  U_n^j, g\rangle \right| ,
		\end{aligned}
	\end{equation}
	where $\breve{\mathbb{X}}_n^{j} := \sqrt{\sigma_n^j}  U_n^j (V_n^j)^{\top}$ for all $j=1,\dots,k$ is a $d \times M$ matrix of rank $1$, \textit{$U_n^j$ denotes the $j$-th column of $U_n^{\top}$} and \textit{$V_n^j$ the $j$-th column of $V_n$}, and in the last line we use the fact that $(V_n^j)_j$ are orthogonal with respect to the empirical scalar product. 
	
	Similarly, we can consider the following singular value decomposition for the particle system $\widehat{\mathbb{X}}_n=\widehat{U}_n^{\top}\widehat{\Sigma}_n \widehat{V}_n^{\top}$, with $\widehat{U}_n$, $\widehat{\Sigma}_n$ and $\widehat{V}_n^{\top}$ all dependent on $\boldsymbol{\omega}$.
	Notice that from the orthogonality of the columns of $\widehat{U}^{\top}_n$ we have
	\begin{equation}\label{eq: orth P_U_hat}
		\left(I_{d \times d}-P_{\widehat{U}_n}\right)\widehat{\mathbb{X}}_n = \widehat{\mathbb{X}}_n-P_{\widehat{U}_n}\widehat{\mathbb{X}}_n = \widehat{\mathbb{X}}_n-\widehat{U}^{\top}_n\widehat{U}_n\widehat{U}^{\top}_n\widehat{\Sigma}_n \widehat{V}^{\top}_n = \widehat{\mathbb{X}}_n-\widehat{\mathbb{X}}_n= 0.
	\end{equation}
	After having defined the rank-$1$ matrices $\widehat{\mathbb{X}}^j_n:= \sqrt{\widehat{\sigma}_n^j}  \widehat{U}_n^j \sqrt{M}(\widehat{V}_n^j)^{\top} \in \mathbb{R}^{d \times M}$ for all $j=1,\dots,k$, where $\widehat{U}_n^j$ denotes the $j$-th column of $\widehat{U}_n^{\top}\in \mathbb{R}^{d \times k}$ and $\widehat{V}_n^j$ the $j$-th column of $\widehat{V}_n \in \mathbb{R}^{M \times k}$, \eqref{eq: orth P_U_hat} implies
	\begin{equation*}
		\begin{aligned}
			0 = \frac{1}{M}\left(I_{d \times d}-P_{\widehat{U}_n}\right)\widehat{\mathbb{X}}_nV_n^j& =\left(I_{d \times d}-P_{\widehat{U}_n}\right) \frac{1}{M}\widehat{\mathbb{X}}_nV_n^j, \quad \forall j=1,\dots,k.
		\end{aligned}
	\end{equation*}
	Using the following relation, we want to bound the quantity in \eqref{eq: PuPu perp} in terms of the difference between the auxiliary matrix $\breve{\mathbb X}_n$ and the particle matrix $\widehat{\mathbb X}_n$. Using orthogonality of $U^{\top}_n$, subadditivity of the norm, Cauchy-Schwarz inequality applied with respect to the empirical scalar product, and Bessel's inequality, one has that
	\begin{equation*}
		\begin{aligned}
			\|(P_{\widehat{U}_n})^{\perp}P_{U_n}\|_{\mathbb{R}^{d} \to \mathbb{R}^{d}}^2 & =\sup\limits_{\substack{g \in \mathbb{R}^{d} \\ |g| = 1}} \left|\left(I_{d \times d}-P_{\widehat{U}_n}\right) \sum_{j=1}^{k }\frac{(\breve{\mathbb{X}}_n-\widehat{\mathbb{X}}_n)V_n^j}{M \sqrt{\sigma_n^j}} \langle  U_n^j, g\rangle \right|^2 \\
			&  \leq \sup\limits_{\substack{g \in \mathbb{R}^{d} \\ |g| = 1}}  \left| \sum_{j=1}^{k} \frac{1}{M} \frac{(\breve{\mathbb{X}}_n-\widehat{\mathbb{X}}_n)V_n^j}{\sqrt{\sigma_n^j}} \langle  U_n^j, g\rangle \right|^2 \\
			&  \leq k \sup\limits_{\substack{g \in \mathbb{R}^{d} \\ |g| = 1}}  \sum_{j=1}^{k} \left| \frac{1}{M} \frac{(\breve{\mathbb{X}}_n-\widehat{\mathbb{X}}_n)V_n^j}{\sqrt{\sigma_n^j}} \langle  U_n^j, g\rangle \right|^2 \\
			& = k\sup\limits_{\substack{g \in \mathbb{R}^{d} \\ |g| = 1}}  \sum_{j=1}^{k} \frac{1}{M^2} \left\|  \frac{(\breve{\mathbb{X}}_n-\widehat{\mathbb{X}}_n)}{\sqrt{\sigma_n^j}} \right\|^2_{\mathrm{F}} | V_n^j|^2 | \langle  U_n^j, g\rangle |^2\\
			& = k\sup\limits_{\substack{g \in \mathbb{R}^{d} \\ |g| = 1}}  \frac{1}{M^2}   \frac{\left\|\breve{\mathbb{X}}_n-\widehat{\mathbb{X}}_n\right\|_{\mathrm{F}}^2}{\sigma_n^k} M \sum_{j=1}^{k} | \langle  U_n^j, g\rangle |^2 \\
			& \leq k\sup\limits_{\substack{g \in \mathbb{R}^{d} \\ |g| = 1}}  \frac{1}{M}  \frac{\left\|\breve{\mathbb{X}}_n-\widehat{\mathbb{X}}_n\right\|_{\mathrm{F}}^2 }{\sigma_n^k} |g|^2\\
			& \leq \frac{k}{\sigma_n^k}  \frac{1}{M}  \frac{\left\|\breve{\mathbb{X}}_n-\widehat{\mathbb{X}}_n\right\|_{\mathrm{F}}^2 }{\sigma_n^k}\\
		\end{aligned}
	\end{equation*}
	Therefore, we have
	\begin{equation}\label{eq: P_u_perpPu}
		\begin{aligned}
			\|(P_{\widehat{U}_n})^{\perp}P_{U_n}\|_{\mathbb{R}^{d} \to \mathbb{R}^{d}}^2 & \leq \frac{k}{\sigma_n^k} \frac{ \|\breve{\mathbb{X}}_n-\widehat{\mathbb{X}}_n\|_{\mathrm{F}}^2}{M}, \\
			& \leq \frac{2k}{\sigma_n^k} \left( \frac{ \|\breve{\mathbb{X}}_n-\mathbb{X}_n\|_{\mathrm{F}}^2}{M}+  \frac{ \|\mathbb{X}_n-\widehat{\mathbb{X}}_n\|_{\mathrm{F}}^2}{M}\right),
		\end{aligned}
	\end{equation}
	where we recall that $\mathbb{X}_n$ is defined as $\mathbb{X}_n(\omega) =  U^{\top}_n\tilde{\Sigma}_n (\omega)\tilde{V}^{\top}_n (\omega) = \sum_{j=1}^{k }
	\sqrt{\sigma_{n,M}^j} (\omega_j) U_n^j \sqrt{M} \tilde{V}_n^j (\omega_j)$ (see Section \ref{sec: notation}). On the other hand, we have the trivial bound $\|(P_{\widehat{U}_n})^{\perp}P_{U_n}\|_{\mathbb{R}^{d} \to \mathbb{R}^{d}}\leq 1$ as $\widehat{U}_n$ and $U_n$ are orthogonal in $\mathbb{R}^d$.
	
	We now return to bound the quantity in the statement \eqref{eq: bound between U proj} of the lemma using the previous operator-norm estimate. Define the event $\Omega_R = \{ \boldsymbol{\omega} \ : \ \frac{1}{M}\sum_{i=1}^{M} |f^{i}(\omega_i) |^2 \leq R^2\}$ for $R> \|f\|_{L^2(\Omega^M;\mathbb{R}^{d \times M})}$.
	Thus, through the interchangeability of particles, Lemma \ref{cor: Bernstein ineq}, Cauchy-Schwarz and triangle inequality, and relation \eqref{eq: P_u_perpPu} the $L^2$ error can be bounded in the following way
	\begin{equation}\label{eq: inter proj u bnd}
		\begin{aligned}
			&\sqrt{\mathbb{E}\left[\frac{1}{M}\sum_{i=1}^{M} \left|\left(P_{U_n}-P_{\widehat{U}_n}\right)f^{i}\right|^2\right]}\\
			\leq & \sqrt{\mathbb{E}\left[\frac{1}{M}\sum_{i=1}^{M} \|P_{U_n}-P_{\widehat{U}_n}\|^2_{\mathbb{R}^{d} \to \mathbb{R}^{d}} |f^{i}|^2\right]}\\
			\leq & \sqrt{\mathbb{E}\left[\frac{1}{M}\sum_{i=1}^{M} \|(P_{\widehat{U}_n})^{\perp}P_{U_n}\|^2_{\mathbb{R}^{d} \to \mathbb{R}^{d}}|f^{i}|^2 \mathbbm{1}_{\Omega_R}\right]} + \sqrt{\mathbb{E}\left[\frac{1}{M}\sum_{i=1}^{M} \|(P_{\widehat{U}_n})^{\perp}P_{U_n}\|^2_{\mathbb{R}^{d} \to \mathbb{R}^{d}}|f^{i}|^2 \mathbbm{1}_{\Omega_R^{C}}\right]}\\
			\leq & \sqrt{\mathbb{E}\left[\|(P_{\widehat{U}_n})^{\perp}P_{U_n}\|_{\mathbb{R}^{d} \to \mathbb{R}^{d}}^2\right]} R + \sqrt{\mathbb{E}\left[\frac{1}{M}\sum_{i=1}^{M} \|(P_{\widehat{U}_n})^{\perp}P_{U_n}\|^2_{\mathbb{R}^{d} \to \mathbb{R}^{d}}|f^{i}|^2 \mathbbm{1}_{\Omega_R^{C}}\right]}\\
			\leq & \sqrt{\mathbb{E}\left[\|(P_{\widehat{U}_n})^{\perp}P_{U_n}\|_{\mathbb{R}^{d} \to \mathbb{R}^{d}}^2\right]} R + \sqrt{\mathbb{E}\left[\frac{1}{M}\sum_{i=1}^{M}|f^{i}|^2 \mathbbm{1}_{\Omega_R^{C}}\right]}\\
			\leq & \sqrt{\mathbb{E}\left[\|(P_{\widehat{U}_n})^{\perp}P_{U_n}\|_{\mathbb{R}^{d} \to \mathbb{R}^{d}}^2\right]} R + \sqrt[4]{\mathbb{E}\left[\frac{1}{M^2}(\sum_{i=1}^{M}|f^{i}|^2 )^2\right]} \sqrt[4]{\mathbb{P}(\mathbbm{1}_{\Omega_R^{C}})}\\
			\leq &\sqrt{\mathbb{E}\left[\frac{2k}{\sigma_n^k} \left( \frac{ \|\mathbb{X}_n-\breve{\mathbb{X}}_n\|_{\mathrm{F}}^2}{M}+  \frac{ \|\mathbb{X}_n-\widehat{\mathbb{X}}_n\|_{\mathrm{F}}^2}{M}\right)\right]} R +  \sqrt[4]{\mathbb{E}\left[\frac{1}{M}\sum_{i=1}^{M}(|f^{i}|^2 )^2\right]} \sqrt[4]{C} e^{-\frac{cMR^2}{4}}\\
			= &\sqrt{ \frac{2k}{\sigma_n^k} }  \left(\sqrt{ \mathbb{E}\left[\frac{ \|\mathbb{X}_n-\widehat{\mathbb{X}}_n\|_{\mathrm{F}}^2}{M}\right]} + \sqrt{\mathbb{E}\left[ \frac{ \|\mathbb{X}_n-\breve{\mathbb{X}}_n\|_{\mathrm{F}}^2}{M} \right] } \right) R + \|f\|_{L^4(\Omega^M;\mathbb{R}^{d \times M})} \sqrt[4]{C} e^{-\frac{cMR^2}{4}} \\
			= &\sqrt{ \frac{2k}{\sigma_n^k}} \left(e_n+  \sqrt{\mathbb{E}\left[\frac{ \|\mathbb{X}_n-\breve{\mathbb{X}}_n\|_{\mathrm{F}}^2}{M} \right]}\right) R + \|f\|_{L^4(\Omega^M;\mathbb{R}^{d \times M})} \sqrt[4]{C} e^{-\frac{cMR^2}{4}},\\
		\end{aligned}
	\end{equation}
	where $\| \cdot \|_{L^4(\Omega^M;\mathbb{R}^{d \times M})}$ is defined in \eqref{eq: Lp norm}.
	It remains to estimate $
	\mathbb{E}\left[\frac{\|\mathbb X_n-\breve{\mathbb X}_n\|_{\mathrm F}^2}{M}\right]$,
	which quantifies the Monte Carlo error between the sampled semi-discrete matrix and its reconstruction from the exact Gramian.
	We proceed as follows using orthogonality of $U_n$ and $V_n$
	\begin{equation}
		\begin{aligned}
			\mathbb{E}\left[\frac{ \|\mathbb{X}_n-\breve{\mathbb{X}}_n\|_{\mathrm{F}}^2}{M} \right] & = \mathbb{E}\left[\frac{ \|U^{\top}_n\tilde{\Sigma}_n\tilde{V}^{\top}_n-U^{\top}_n\Sigma_n V^{\top}_n\|_{\mathrm{F}}^2}{M} \right]= \mathbb{E}\left[\frac{ \|\tilde{\Sigma}_n\tilde{V}^{\top}_n-\Sigma_n V^{\top}_n\|_{\mathrm{F}}^2}{M} \right]\\
			& = \mathbb{E}\left[\frac{ \|\tilde{\Sigma}_n\tilde{V}^{\top}_n -\Sigma_n \sqrt{M}  \tilde{V}^{\top}_n \|_{\mathrm{F}}^2}{M} \right] = \mathbb{E}\left[\frac{ \|\tilde{\Sigma}_n-\Sigma_n\sqrt{M}\|_{\mathrm{F}}^2}{M} \right]\\
			& = \mathbb{E}\left[\sum_{j=1}^{k}(\sqrt{M \sigma^j_{n,M}}- \sqrt{\sigma_n^j}\sqrt{M})^2\frac{ 1}{M} \right]= \mathbb{E}\left[\sum_{j=1}^{k}(\sqrt{\sigma^j_{n,M}}-\sqrt{\sigma^j_{n}})^2 \right]\\
			& =  \mathbb{E}\left[\sum_{j=1}^{k} \frac{(\sigma^j_{n,M} - \sigma^j_{n})^2}{(\sqrt{\sigma^j_{n,M}}+\sqrt{\sigma^j_{n}})^2}\right]  \\
			& =  \mathbb{E}\left[\sum_{j=1}^{k} \frac{(\sigma^j_{n,M} - \sigma^j_{n})^2}{\sigma^k_{n,M}+\sigma^k_{n} + 2 \sqrt{\sigma^k_{n,M}}\sqrt{\sigma^k_{n}}}\right] \\
			& \leq \frac{1}{\sigma^k_{n}} \mathbb{E}\left[\sum_{j=1}^{k} (\sigma^j_{n,M} - \sigma^j_{n})^2\right]  \\
			& \leq  \frac{1}{\sigma^k_{n}}  \mathbb{E}\left[ \|\frac{ 1}{M}\mathbb{X}_n\mathbb{X}^{\top}_n - \mathbb{E}[X_nX^{\top}_n]\|_{\mathrm{F}}^2 \right],\\
		\end{aligned}
	\end{equation}
	where in the last line we use Hoffman-Wielandt theorem \cite[Corollary 7.3.5]{horn2012matrix}.
	
	Note that $(\frac{ 1}{M}\mathbb{X}_n\mathbb{X}^{\top}_n)_{ij} = \frac{ 1}{M} \sum_{\ell = 1}^M (X_n(\omega_\ell))_{i}(X_n(\omega_\ell))_{j}$, i.e.\ the $(i,j)$-entry
	of $\frac{ 1}{M}\mathbb{X}_n\mathbb{X}_n^{\top}$ is the Monte Carlo estimator of $\mathbb{E}[(X_n(\omega))_{i}(X_n(\omega))_{j}]$, where $(X_n(\omega))_{i}$ denotes the $i$-th component of random vector $X_n(\omega)$.
	Then, one has
	\begin{equation}\label{eq: MC in xx_T}
		\begin{aligned}
			\mathbb{E}\left[\frac{ \| \mathbb{X}_n-\breve{\mathbb{X}}_n\|_{\mathrm{F}}^2}{M} \right] 
			& \leq \frac{1}{\sigma^k_{n}} \mathbb{E}\left[\|\frac{ 1}{M}\mathbb{X}_n\mathbb{X}^{\top}_n - \mathbb{E}[X_nX^{\top}_n]\|_{\mathrm{F}}^2 \right] = \frac{1}{\sigma^k_{n}} \frac{1}{M}  \mathbb{E}[\|X_nX^{\top}_n\|_{\mathrm{F}}^2]= \frac{1}{\sigma^k_{n}} \frac{1}{M} \mathbb{E} \left[|X_n|^{4} \right],\\
		\end{aligned}
	\end{equation}
	for some positive constant $C_4$, where we employed 
	the convergence of the Monte Carlo estimator in Hilbert spaces as the samples $X^{1}_n,\dots,X^{M}_n$ of the semi-discrete DLRA are independent, and the fact that thanks to the cycling property of the trace one has
	\begin{equation*}
		\begin{aligned}
			\mathbb{E}\big[\|X_n X_n^\top\|_{\mathrm{F}}^2\big]  
			&= \mathbb{E}\Big[ \mathrm{Tr}\big((X_n X_n^\top)^\top (X_n X_n^\top)\big) \Big] = \mathbb{E}\Big[ \mathrm{Tr}\big(X_n X_n^\top X_n X_n^\top\big) \Big] \\
			&= \mathbb{E}\Big[ \mathrm{Tr}\big(X_n (X_n^\top X_n) X_n^\top\big) \Big] 
			= \mathbb{E}\Big[ (X_n^\top X_n)\,\mathrm{Tr}(X_n X_n^\top)  \Big] 
			= \mathbb{E}\Big[ (X_n^\top X_n)^2 \Big] 
			= \mathbb{E}\big[|X_n|^4\big].
		\end{aligned}
	\end{equation*} 
	Finally, via using relation \eqref{eq: MC in xx_T} in \eqref{eq: inter proj u bnd} we finally obtain
	\begin{equation*}
		\begin{aligned}
			\sqrt{\mathbb{E}\left[\frac{1}{M}\sum_{i=1}^{M} \left|\left(P_{U_n}-P_{\widehat{U}_n}\right)f^{i}\right|^2\right]} \leq& \left(\sqrt{\frac{2k}{\sigma_n^k}} e_n + \frac{1}{\sqrt{M}} \frac{\sqrt{2k K_4(T)}}{\sigma_n^k} \right) R \\
			& + \|f\|_{L^4(\Omega^M;\mathbb{R}^{d \times M})} \sqrt[4]{C} e^{-\frac{cMR^2}{4}}. 
		\end{aligned}
	\end{equation*}
\end{proof}

\begin{proof}[\bfseries Proof of Lemma \ref{lem: error between stochastic proj}] \label{proof: lem: error between stochastic proj}
	The proof follows by adding and subtracting useful terms in the left-hand side of \eqref{eq: Y projectors error}, and then bounding them. For $f \in L^4(\Omega^{M};\mathbb{R}^{d \times M})$ with $f^i$ i.i.d.\ one aims to bound the following quantity with bounds linked to the Monte Carlo error and to the regularization parameter $\alpha$:
	\begin{equation*}
		\begin{aligned}
			& \left(\mathbb{E}\left[\frac{1}{M} \sum_{i=1}^{M} \left|\left(P_{\widetilde{Y}_{n+1}^{i}}- \widehat{P}^{\alpha}_{\widetilde{Y}_{n+1}^{i}} \right)f\right|^2\right]\right)^{\frac{1}{2}}  \\
			=&\resizebox{1\linewidth}{!}{$ \left(\mathbb{E}\left[\frac{1}{M} \sum_{i=1}^{M} \left|(\widetilde{Y}_{n+1}^{i})^{\top}\left(\mathbb{E}\left[ \frac{1}{M} \widetilde{\mathbb{Y}}_{n+1}\widetilde{\mathbb{Y}}_{n+1}^{\top}\right]\right)^{-1}\mathbb{E}\left[\frac{1}{M} \widetilde{\mathbb{Y}}_{n+1} f^{\top} \right] - (\widetilde{Y}_{n+1}^{i})^{\top}\left(\frac{\sum_{j=1}^{M} \widetilde{Y}_{n+1}^j (\widetilde{Y}_{n+1}^j)^{\top}}{M} + \alpha I_{k \times k}\right)^{-1}\widehat{\mathbb{E}}\left[ \widetilde{\mathbb{Y}}_{n+1} f^{\top} \right] \right|^2\right]\right)^{\frac{1}{2}}  $} \\
			\leq& \underbrace{\left(\mathbb{E}\left[\frac{1}{M} \sum_{i=1}^{M} \left|(\widetilde{Y}_{n+1}^{i})^{\top}\left(\left(\mathbb{E}\left[ \frac{1}{M} \widetilde{\mathbb{Y}}_{n+1}\widetilde{\mathbb{Y}}_{n+1}^{\top}\right]\right)^{-1}-\left(\mathbb{E}\left[ \frac{1}{M} \widetilde{\mathbb{Y}}_{n+1}\widetilde{\mathbb{Y}}_{n+1}^{\top}+\alpha I_{k \times k}\right]\right)^{-1} \right)\mathbb{E}\left[ \frac{1}{M} \widetilde{\mathbb{Y}}_{n+1} f^{\top} \right] \right|^2\right]\right)^{\frac{1}{2}}}_{=:L_1}\\
			&\resizebox{1\linewidth}{!}{$+ \underbrace{\left(\mathbb{E}\left[\frac{1}{M} \sum_{i=1}^{M} \left|(\widetilde{Y}_{n+1}^{i})^{\top}\left(\left(\mathbb{E}\left[ \frac{1}{M} \widetilde{\mathbb{Y}}_{n+1}\widetilde{\mathbb{Y}}_{n+1}^{\top}+\alpha I_{k \times k}\right]\right)^{-1}-\left(\frac{\sum_{j=1}^{M} \widetilde{Y}_{n+1}^j (\widetilde{Y}_{n+1}^j)^{\top}}{M} + \alpha I_{k \times k}\right)^{-1}\right)\mathbb{E}\left[ \frac{1}{M} \widetilde{\mathbb{Y}}_{n+1} f^{\top} \right] \right|^2\right]\right)^{\frac{1}{2}}}_{=:L_2} $} \\
			&+ \underbrace{\left(\mathbb{E}\left[\frac{1}{M} \sum_{i=1}^{M} \left|(\widetilde{Y}_{n+1}^{i})^{\top}\left(\frac{\sum_{j=1}^{M} \widetilde{Y}_{n+1}^j (\widetilde{Y}_{n+1}^j)^{\top}}{M} + \alpha I_{k \times k}\right)^{-1}\left(\mathbb{E}\left[ \frac{1}{M} \widetilde{\mathbb{Y}}_{n+1} f^{\top} \right] -\widehat{\mathbb{E}}\left[ \widetilde{\mathbb{Y}}_{n+1} f^{\top}\right]  \right)\right|^2\right]\right)^{\frac{1}{2}}}_{=:L_3},\\
		\end{aligned}
	\end{equation*} 
	where we exploited the invariance of the Euclidean norm under transposition and the triangle inequality.
	Let us recall the following useful matrix relation 
	for invertible matrices $A,B \in \mathbb{R}^{k \times k}$:
	\begin{equation}\label{eq: A-B}
		A^{-1}-B^{-1}= -B^{-1}\left(A-B\right)A^{-1}.
	\end{equation}
	This identity will be used to compare the inverse matrices appearing in the terms $L_1$ and $L_2$.
	Then, using also Cauchy-Schwarz inequality, the property of the product measure, the interchangeability of particles and the expectation $\mathbb{E}_i$ with respect to the $i$-th particle, relation \eqref{eq: A-B}, and $|\mathbb{E}\left[ \frac{1}{M} \widetilde{\mathbb{Y}}_{n+1}\widetilde{\mathbb{Y}}_{n+1}^{\top} \right]^{-1}|\leq (\widetilde{\sigma}_{n+1}^k)^{-1}$, we bound the term $L_1$ as
	\begin{equation*} 
		\begin{aligned}
			L_1  = & \left(\mathbb{E}\left[\frac{1}{M} \sum_{i=1}^{M} \left|(\widetilde{Y}_{n+1}^{i})^{\top}\left(\left(\mathbb{E}\left[ \frac{1}{M} \widetilde{\mathbb{Y}}_{n+1}\widetilde{\mathbb{Y}}_{n+1}^{\top} \right]\right)^{-1}-\left(\mathbb{E}\left[ \frac{1}{M} \widetilde{\mathbb{Y}}_{n+1}\widetilde{\mathbb{Y}}_{n+1}^{\top}+\alpha I_{k \times k}\right]\right)^{-1} \right)\mathbb{E}\left[ \frac{1}{M} \widetilde{\mathbb{Y}}_{n+1} f^{\top} \right] \right|^2\right]\right)^{\frac{1}{2}}\\
			= &\left( \frac{1}{M} \sum_{i=1}^{M} \mathbb{E}_{\omega_{1}}\left[ \left|(\widetilde{Y}_{n+1}^{1})^{\top}\left(\left(\mathbb{E}\left[ \frac{1}{M} \widetilde{\mathbb{Y}}_{n+1}\widetilde{\mathbb{Y}}_{n+1}^{\top}\right]\right)^{-1}-\left(\mathbb{E}\left[ \frac{1}{M} \widetilde{\mathbb{Y}}_{n+1}\widetilde{\mathbb{Y}}_{n+1}^{\top}+\alpha I_{k \times k}\right]\right)^{-1} \right)\mathbb{E}\left[ \frac{1}{M} \widetilde{\mathbb{Y}}_{n+1} f^{\top} \right] \right|^2\right]\right)^{\frac{1}{2}} \\
			= & \left(\mathbb{E}_{\omega_1}\left[ \left|(\widetilde{Y}_{n+1}^{1})^{\top}\left(\left(\mathbb{E}\left[ \frac{1}{M} \widetilde{\mathbb{Y}}_{n+1}\widetilde{\mathbb{Y}}_{n+1}^{\top}\right]\right)^{-1}-\left(\mathbb{E}\left[ \frac{1}{M} \widetilde{\mathbb{Y}}_{n+1}\widetilde{\mathbb{Y}}_{n+1}^{\top}+\alpha I_{k \times k}\right]\right)^{-1} \right)\mathbb{E}\left[ \frac{1}{M} \widetilde{\mathbb{Y}}_{n+1} f^{\top} \right] \right|^2\right]\right)^{\frac{1}{2}}\\
			= &\Big(\mathbb{E}_{\omega_1}\Big[ \Big|(\widetilde{Y}_{n+1}^{1})^{\top}\left(\mathbb{E}\left[ \frac{1}{M} \widetilde{\mathbb{Y}}_{n+1}\widetilde{\mathbb{Y}}_{n+1}^{\top}+\alpha I_{k \times k}\right]\right)^{-1}\left(\mathbb{E}\left[ \frac{1}{M} \widetilde{\mathbb{Y}}_{n+1}\widetilde{\mathbb{Y}}_{n+1}^{\top}\right]-\mathbb{E}\left[ \frac{1}{M} \widetilde{\mathbb{Y}}_{n+1}\widetilde{\mathbb{Y}}_{n+1}^{\top}+\alpha I_{k \times k}\right] \right) \\
			&\cdot \left(\mathbb{E}\left[ \frac{1}{M} \widetilde{\mathbb{Y}}_{n+1}\widetilde{\mathbb{Y}}_{n+1}^{\top}\right]\right)^{-1}\mathbb{E}\left[ \frac{1}{M} \widetilde{\mathbb{Y}}_{n+1} f^{\top}  \right] \Big|^2\Big]\Big)^{\frac{1}{2}}\\
			= & \left( \mathbb{E}_{\omega_1}\left[ \left|(\widetilde{Y}_{n+1}^{1})^{\top}\left(\mathbb{E}\left[ \frac{1}{M} \widetilde{\mathbb{Y}}_{n+1}\widetilde{\mathbb{Y}}_{n+1}^{\top}+\alpha I_{k \times k}\right]\right)^{-\frac{1}{2}} \right|^2\right]\right)^{\frac{1}{2}} \left| \left(\mathbb{E}\left[ \frac{1}{M} \widetilde{\mathbb{Y}}_{n+1}\widetilde{\mathbb{Y}}_{n+1}^{\top}+\alpha I_{k \times k}\right]\right)^{-\frac{1}{2}} \right| \\ 
			& \quad \cdot \left|\alpha I_{k \times k} \right| \left| \left(\mathbb{E}\left[ \frac{1}{M} \widetilde{\mathbb{Y}}_{n+1}\widetilde{\mathbb{Y}}_{n+1}^{\top}\right]\right)^{-\frac{1}{2}} \right| \left( \left|\mathbb{E}\left[ \left(\mathbb{E}\left[ \frac{1}{M} \widetilde{\mathbb{Y}}_{n+1}\widetilde{\mathbb{Y}}_{n+1}^{\top}\right]\right)^{-\frac{1}{2}} \frac{1}{M} \widetilde{\mathbb{Y}}_{n+1} f^{\top} \right] \right|^2\right)^{\frac{1}{2}}  \\
			\leq & \frac{\sqrt{k}}{\sqrt{\widetilde{\sigma}_{n+1}^k + \alpha}} \alpha \|f\|_{L^2(\Omega^{M};\mathbb{R}^{d \times M})} \frac{1}{\sqrt{\widetilde{\sigma}_{n+1}^k}} \\
			\leq & \frac{\sqrt{k}}{\widetilde{\sigma}_{n+1}^k}  \alpha  \|f\|_{L^2(\Omega^{M};\mathbb{R}^{d \times M})}, \\
		\end{aligned}
	\end{equation*} 
	where we employed repeatedly the following relation
	\begin{equation}\label{eq: Gramian cs}
		\begin{aligned}
			&\mathbb{E}\left[ \left|\frac{1}{\sqrt{M}} (\widetilde{\mathbb{Y}}_{n+1})^{\top}\left(\mathbb{E}\left[  \frac{1}{M} \widetilde{\mathbb{Y}}_{n+1}\widetilde{\mathbb{Y}}_{n+1}^{\top}+\alpha I_{k \times k}\right]\right)^{-\frac{1}{2}} \right|^2\right] \\
			=& 	\mathbb{E}\left[ \frac{1}{\sqrt{M}} (\widetilde{\mathbb{Y}}_{n+1})^{\top}\mathbb{E}\left[  \frac{1}{M} \widetilde{\mathbb{Y}}_{n+1}\widetilde{\mathbb{Y}}_{n+1}^{\top}+\alpha I_{k \times k}\right]^{-1}\widetilde{\mathbb{Y}}_{n+1} \frac{1}{\sqrt{M}}\right]\\
			=& \mathbb{E}\left[ \mathrm{Tr}\left(\frac{1}{\sqrt{M}}(\widetilde{\mathbb{Y}}_{n+1})^{\top}\mathbb{E}\left[ \frac{1}{M} \widetilde{\mathbb{Y}}_{n+1}\widetilde{\mathbb{Y}}_{n+1}^{\top}+\alpha I_{k \times k}\right]^{-1}\widetilde{\mathbb{Y}}_{n+1} \frac{1}{\sqrt{M}} \right)\right]\\
			=& \mathbb{E}\left[ \mathrm{Tr}\left(\frac{1}{M} \widetilde{\mathbb{Y}}_{n+1}\widetilde{\mathbb{Y}}_{n+1}^{\top}\mathbb{E}\left[ \frac{1}{M} \widetilde{\mathbb{Y}}_{n+1}\widetilde{\mathbb{Y}}_{n+1}^{\top}+\alpha I_{k \times k}\right]^{-1} \right)\right]\\
			=&  \mathrm{Tr}\left(\mathbb{E}\left[\frac{1}{M} \widetilde{\mathbb{Y}}_{n+1}\widetilde{\mathbb{Y}}_{n+1}^{\top}\right]\mathbb{E}\left[ \frac{1}{M} \widetilde{\mathbb{Y}}_{n+1}\widetilde{\mathbb{Y}}_{n+1}^{\top}+\alpha I_{k \times k}\right]^{-1} \right) \leq k.
		\end{aligned}
	\end{equation}
	Similarly, for $L_2$ we have
	\begin{equation*}
		\begin{aligned}
			L_2  =& \resizebox{1\linewidth}{!}{$ \left(\mathbb{E}\left[\frac{1}{M} \sum_{i=1}^{M} \left|(\widetilde{Y}_{n+1}^{i})^{\top}\left(\left(\mathbb{E}\left[  \frac{1}{M} \widetilde{\mathbb{Y}}_{n+1}\widetilde{\mathbb{Y}}_{n+1}^{\top}+\alpha I_{k \times k}\right]\right)^{-1}-\left(\frac{\sum_{j=1}^{M} \widetilde{Y}_{n+1}^j (\widetilde{Y}_{n+1}^j)^{\top}}{M} + \alpha I_{k \times k}\right)^{-1}\right)\mathbb{E}\left[ \frac{1}{M} \widetilde{\mathbb{Y}}_{n+1} f^{\top}  \right] \right|^2\right]\right)^{\frac{1}{2}} $} \\
			= & \resizebox{1\linewidth}{!}{$ \Big(\mathbb{E}\Big[\frac{1}{M} \sum_{i=1}^{M} \Big|(\widetilde{Y}_{n+1}^{i})^{\top}\Big(\frac{\sum_{j=1}^{M} \widetilde{Y}_{n+1}^j (\widetilde{Y}_{n+1}^j)^{\top}}{M} + \alpha I_{k \times k}\Big)^{-1}\Big(\frac{\sum_{j=1}^{M} \widetilde{Y}_{n+1}^j (\widetilde{Y}_{n+1}^j)^{\top}}{M}+ \alpha I_{k \times k} - \mathbb{E}\Big[  \frac{1}{M} \widetilde{\mathbb{Y}}_{n+1}\widetilde{\mathbb{Y}}_{n+1}^{\top}+\alpha I_{k \times k}\Big]\Big) $} \\
			&\cdot \mathbb{E}\Big[  \frac{1}{M} \widetilde{\mathbb{Y}}_{n+1}\widetilde{\mathbb{Y}}_{n+1}^{\top}+\alpha I_{k \times k}\Big]^{-1}\mathbb{E}\Big[ \frac{1}{M} \widetilde{\mathbb{Y}}_{n+1} f^{\top}  \Big] \Big|^2\Big]\Big)^{\frac{1}{2}} \\
			\leq & \Big(\mathbb{E} \Big[\frac{1}{M} \sum_{i=1}^{M}\Big|(\widetilde{Y}_{n+1}^{i})^{\top}\Big(\frac{\sum_{j=1}^{M} \widetilde{Y}_{n+1}^j (\widetilde{Y}_{n+1}^j)^{\top}}{M} + \alpha I_{k \times k}\Big)^{-1}\Big|^2\Big|\frac{\sum_{j=1}^{M} \widetilde{Y}_{n+1}^j (\widetilde{Y}_{n+1}^j)^{\top}}{M} - \mathbb{E}\Big[  \frac{1}{M} \widetilde{\mathbb{Y}}_{n+1}\widetilde{\mathbb{Y}}_{n+1}^{\top} \Big]\Big|^2\\
			&\cdot \Big|\mathbb{E}\Big[  \frac{1}{M} \widetilde{\mathbb{Y}}_{n+1}\widetilde{\mathbb{Y}}_{n+1}^{\top}+\alpha I_{k \times k}\Big]^{-1}\mathbb{E}\Big[ \frac{1}{M} \widetilde{\mathbb{Y}}_{n+1} f^{\top}  \Big] \Big|^2\Big]\Big)^{\frac{1}{2}} \\
			\leq & \frac{1}{\sqrt{\alpha}} \Big(\mathbb{E} \Big[ \Big( \frac{1}{M} \sum_{i=1}^{M}\Big|(\widetilde{Y}_{n+1}^{i})^{\top}\Big(\frac{\sum_{j=1}^{M} \widetilde{Y}_{n+1}^j (\widetilde{Y}_{n+1}^j)^{\top}}{M} + \alpha I_{k \times k}\Big)^{-\frac{1}{2}}\Big|^2\Big)\Big|\frac{\sum_{j=1}^{M} \widetilde{Y}_{n+1}^j (\widetilde{Y}_{n+1}^j)^{\top}}{M} - \mathbb{E}\Big[  \frac{1}{M} \widetilde{\mathbb{Y}}_{n+1}\widetilde{\mathbb{Y}}_{n+1}^{\top} \Big]\Big|^2\Big]\Big)^{\frac{1}{2}}\\
			&\cdot \Big( \Big|\mathbb{E}\Big[ \frac{1}{M} \widetilde{\mathbb{Y}}_{n+1}\widetilde{\mathbb{Y}}_{n+1}^{\top}+\alpha I_{k \times k}\Big]^{-\frac{1}{2}}\mathbb{E}\Big[ \mathbb{E}\Big[  \frac{1}{M} \widetilde{\mathbb{Y}}_{n+1}\widetilde{\mathbb{Y}}_{n+1}^{\top}+\alpha I_{k \times k}\Big]^{-\frac{1}{2}} \frac{1}{M} \widetilde{\mathbb{Y}}_{n+1} f^{\top} \Big] \Big|^2\Big)^{\frac{1}{2}}  \\
			\leq & 
			\frac{k}{\sqrt{\alpha}} \frac{1}{\sqrt{\widetilde{\sigma}_{n+1}^k}}  \|f\|_{L^2(\Omega^{M};\mathbb{R}^{d \times M})}  \Big(\mathbb{E}\Big[ \Big\|\frac{\sum_{j=1}^{M} \widetilde{Y}_{n+1}^j (\widetilde{Y}_{n+1}^j)^{\top}}{M} - \mathbb{E}\Big[  \frac{1}{M} \widetilde{\mathbb{Y}}_{n+1}\widetilde{\mathbb{Y}}_{n+1}^{\top}\Big]\Big\|^2_{\mathrm{F}} \Big]\Big)^{\frac{1}{2}}.\\
		\end{aligned}
	\end{equation*} 
	By construction, $(\widetilde{Y}_{n+1}^i )_i$ are independent with respect to $\mathbb{P}$, hence by convergence of the Monte Carlo estimator we find that
	\begin{equation*}
		\begin{aligned}
			L_2 
			\leq & \frac{1}{\sqrt{\alpha}} \frac{1}{\sqrt{\widetilde{\sigma}_{n+1}^k}}  \|f\|_{L^2(\Omega^{M};\mathbb{R}^{d \times M})} \frac{1}{\sqrt{M}} \mathbb{E}[|\widetilde{Y}_{n+1}\widetilde{Y}_{n+1}^{\top}|^2] \\
			\leq & \frac{k}{\sqrt{\alpha}} \frac{1}{\sqrt{\widetilde{\sigma}_{n+1}^k}}  \|f\|_{L^2(\Omega^{M};\mathbb{R}^{d \times M})} \frac{1}{\sqrt{M}} \sqrt{K_4(T)},
		\end{aligned}
	\end{equation*} 
	where $K_4(T)$ is defined in \eqref{eq: Lp moment bound semidiscretized X}.
	
	Let us finally analyze $L_3$, which corresponds to the error in approximating $\mathbb{E}[\widetilde{Y}_{n+1} f^\top]$ by its empirical counterpart.
	Using Cauchy-Schwarz inequality we get 
	\begin{equation*}
		\begin{aligned}
			L_3  =& \left(\mathbb{E}\left[\frac{1}{M} \sum_{i=1}^{M} \left|(\widetilde{Y}_{n+1}^{i})^{\top}\left(\frac{\sum_{j=1}^{M} \widetilde{Y}_{n+1}^j (\widetilde{Y}_{n+1}^j)^{\top}}{M} + \alpha I_{k \times k}\right)^{-1}\left(\widehat{\mathbb{E}}\left[ \widetilde{\mathbb{Y}}_{n+1} f^{\top} \right]-\mathbb{E}\left[ \frac{1}{M} \widetilde{\mathbb{Y}}_{n+1}  f^{\top}  \right]\right) \right|^2\right]\right)^{\frac{1}{2}} \\
			\leq & \frac{\sqrt{k}}{\sqrt{\alpha}} \left(\mathbb{E}\left[ \left|\widehat{\mathbb{E}}\left[ \widetilde{\mathbb{Y}}_{n+1} f^{\top} \right]-\mathbb{E}\left[ \frac{1}{M} \widetilde{\mathbb{Y}}_{n+1} f^{\top}  \right] \right|^2\right]\right)^{\frac{1}{2}} \\
			= & \frac{\sqrt{k}}{\sqrt{\alpha}} \left(\mathbb{E}\left[ \left|\widehat{\mathbb{E}}\left[ \widetilde{\mathbb{Y}}_{n+1} f^{\top} \right]-\mathbb{E}\left[ \frac{1}{M} \widetilde{\mathbb{Y}}_{n+1} f^{\top}  \right] \right|^2\right]\right)^{\frac{1}{2}} \leq \frac{\sqrt{k}}{\sqrt{\alpha}} \left(\mathbb{E}\left[ \left\| \sum_{j=1}^{M}\frac{ \widetilde{Y}_{n+1}^j (f^j)^{\top}}{M} -\mathbb{E}_{\omega_j}\left[  \widetilde{Y}_{n+1}^j (f^j)^{\top}  \right] \right\|^2_{\mathrm{F}}\right]\right)^{\frac{1}{2}}\\
			\leq &  \frac{\sqrt{k}}{\sqrt{\alpha}} \frac{1}{\sqrt{M}} \sqrt{\mathbb{E}[\| \widetilde{Y}_{n+1} (f^1)^{\top}\|^2_{\mathrm{F}}]} \\
			\leq & \frac{\sqrt{k}}{\sqrt{\alpha}} \frac{1}{\sqrt{M}} \mathbb{E}[| \widetilde{Y}_{n+1}|^4]^{\frac14} \mathbb{E}[|f^1|^4]^{\frac14} \leq  \frac{\sqrt{k}}{\sqrt{\alpha}} \frac{1}{\sqrt{M}}\sqrt{K_4(T)} \|f\|_{L^4(\Omega^{M};\mathbb{R}^{d \times M})}, \\
		\end{aligned}
	\end{equation*} 
	where we employ the interchangeability of particles, $\mathbb{E}_{\omega_{j}}$ indicates the expectation with respect to $\omega_{j}$ and in the third and the second to last line we used the interchangeability of i.i.d particles, and the convergence of the Monte Carlo estimator, and $K_4(T)$ is defined as in \eqref{eq: Lp moment bound semidiscretized X}.
	
	Collecting the bounds for $L_1$, $L_2$, and $L_3$, we obtain
	\begin{equation*}
		\begin{aligned}
			\left\|\left(P_{\widetilde{\mathbb{Y}}_{n+1}}-\widehat{P}^{\alpha}_{\widetilde{\mathbb{Y}}_{n+1}}\right)f\right\|_{L^2(\Omega^{M};\mathbb{R}^{d \times M})} \leq & L_1 + L_2 + L_3 \\
			\leq &  \frac{k}{\widetilde{\sigma}_{n+1}^k}  \alpha  \|f\|_{L^2(\Omega^{M};\mathbb{R}^{d \times M})} + \sqrt{K_4(T)} \frac{k}{\sqrt{\alpha}} \frac{1}{\sqrt{\widetilde{\sigma}_{n+1}^k}} \|f\|_{L^2(\Omega^{M};\mathbb{R}^{d \times M})} \frac{1}{\sqrt{M}} \\
			&+ \frac{\sqrt{k}}{\sqrt{\alpha}} \frac{1}{\sqrt{M}} \sqrt{K_2(T)} \|f\|_{L^2(\Omega^{M};\mathbb{R}^{d \times M})}.
		\end{aligned}
	\end{equation*}
	One finds the optimal value $\alpha_{\mathrm{opt}}$ via matching the order of convergence in the number of samples among the terms $L_1$, $L_2$, $L_3$:
	$$\alpha_{\mathrm{opt}} = \frac{1}{\sqrt{\alpha_{\mathrm{opt}}}} \frac{1}{\sqrt{M}},$$
	where we get $\alpha_{\mathrm{opt}} = \frac{1}{\sqrt[3]{M}}$ and via using the monotonicity of the $L^2(\Omega)$-norm the thesis follows.
\end{proof}

\begin{proof}[\bfseries Proof of Lemma \ref{lem: discr proj Y}] \label{proof: lem: discr proj Y}
	The proof is divided in two cases. The first concerns the event $A_n$ defined in \eqref{eq: def A_n} where both empirical Gramians are strictly lower bounded in a matrix sense by a positive constant. Then, we proceed by adding and subtracting useful terms in the left-hand side of \eqref{eq: discr proj Y}, bounding them, and putting all these computations together in order to obtain the sought bound. The second case concerns the complementary event $A_n^{C}$ where we bound the difference of the two projectors by a Monte-Carlo type error.
	Indeed, we have the following decomposition
	\begin{equation*}
		\begin{aligned}
			\sqrt{\mathbb{E}\left[\frac{1}{M}\sum_{i=1}^{M}
				\left|\left(\widehat{P}^{\alpha}_{\tildhatY{i}{n+1}}-\widehat{P}^{\alpha}_{\widetilde{Y}_{n+1}^{i}}\right)f\right|^2\right]} \leq &		\sqrt{\mathbb{E}\left[\frac{1}{M}\sum_{i=1}^{M}
				\left|\left(\widehat{P}^{\alpha}_{\tildhatY{i}{n+1}}-\widehat{P}^{\alpha}_{\widetilde{Y}_{n+1}^{i}}\right)f \right|^2 \mathbbm{1}_{A_{n+1}} \right]} \\
			&+ 	\sqrt{\mathbb{E}\left[\frac{1}{M}\sum_{i=1}^{M}
				\left|\left(\widehat{P}^{\alpha}_{\tildhatY{i}{n+1}}-\widehat{P}^{\alpha}_{\widetilde{Y}_{n+1}^{i}}\right)f \right|^2 \mathbbm{1}_{A_{n+1}^C} \right]} .
		\end{aligned}
	\end{equation*}
	
	Let us first analyze the first term on the right-hand side of the previous inequality.
	For given $(\widetilde{Y}_{n+1}^{i})_{\{i=1,...,M\}}$, with $\widetilde{Y}_{n+1}^{i} \in \mathbb{R}^k$,  let us define the following operator $\widehat{P}_{\widetilde{Y}_{n+1}^{i}} : L^2(\Omega^M;\mathbb{R}^{d \times M}) \to L^2(\Omega^M;\mathbb{R}^{d})$ as
	\begin{equation}\label{eq: bar proj Y}
		\widehat{P}_{\widetilde{Y}_{n+1}^{i}}[f] =(\widetilde{Y}_{n+1}^{i})^{\top}\left(\frac{\sum_{j=1}^{M} \widetilde{Y}_{n+1}^j (\widetilde{Y}_{n+1}^j)^{\top}}{M}\right)^{\dag}\widehat{\mathbb{E}}\left[ \widetilde{\mathbb{Y}}_{n+1} f^{\top} \right],
		\qquad i=1,\dots,M,
	\end{equation}
	where $\dag$ denotes the pseudoinverse of minimal norm \cite[Section 5.5.4]{golub2013matrix}, and consider also the operator $\widehat{P}_{\tildhatY{i}{n+1}}$ defined as in \eqref{eq: bar proj Y} for $\tildhatY{i}{n+1}$. Notice that if $\widetilde{\sigma}_{n+1,M}^k$ is positive, then the pseudoinverse in \eqref{eq: bar proj Y} is actually an inverse, and the same holds for $\widehat{P}_{\tildhatY{i}{n+1}}$ if $\tildhatsigma{k}{n+1}{M}>0$. Furthermore, notice that $\widehat{P}_{\widetilde{Y}_{n+1}^{i}}[f]$ depends on all the paths of $\boldsymbol{\omega}$.
	Consider $f=(f^i)_{i=1,\dots,M} \in L^4(\Omega^M;\mathbb{R}^{d \times M})$. Then, one bounds the difference of the projectors via triangle inequality in the following way:
	\begin{equation*}
		\begin{aligned}
			\sqrt{\mathbb{E}\left[\frac{1}{M}\sum_{i=1}^{M}
				\left|\left(\widehat{P}^{\alpha}_{\tildhatY{i}{n+1}}-\widehat{P}^{\alpha}_{\widetilde{Y}_{n+1}^{i}}\right)f\right|^2\mathbbm{1}_{A_{n+1}}\right]} \leq &	\sqrt{\mathbb{E}\left[\frac{1}{M}\sum_{i=1}^{M}
				\left|\left(\widehat{P}^{\alpha}_{\tildhatY{i}{n+1}}-\widehat{P}_{\tildhatY{i}{n+1}}\right)f\right|^2\mathbbm{1}_{A_{n+1}}\right]} \\
			&+ 	\sqrt{\mathbb{E}\left[\frac{1}{M}\sum_{i=1}^{M}
				\left|\left(\widehat{P}_{\tildhatY{i}{n+1}}-\widehat{P}_{\widetilde{Y}_{n+1}^{i}}\right)f \right|^2\mathbbm{1}_{A_{n+1}}\right]}\\
			& + \sqrt{\mathbb{E}\left[\frac{1}{M}\sum_{i=1}^{M}
				\left|\left(\widehat{P}_{\widetilde{Y}_{n+1}^{i}}-\widehat{P}^{\alpha}_{\widetilde{Y}_{n+1}^{i}}\right)f \right|^2\mathbbm{1}_{A_{n+1}}\right]}.\\
		\end{aligned}
	\end{equation*}
	Consider an SVD of $\widetilde{\widehat{\mathbb{Y}}}_{n+1} = [\tildhatY{1}{n+1} \ | \ \dots \ | \ \tildhatY{M}{n+1}\ ]=[\tildhatY{i}{n+1}]_i \in \mathbb{R}^{k \times M}$:
	$$\widetilde{\widehat{\mathbb{Y}}}_{n+1} = \widehat{W}^{\top}_{n+1} \widehat{\Sigma}_{n+1} \widehat{V}_{n+1} =  \widehat{W}^{\top}_{n+1} \tilde{\widehat{\Sigma}}_{n+1} \tilde{\widehat{V}}_{n+1}$$
	with $\tilde{\widehat{\Sigma}}_{n+1}:=\widehat{\Sigma}_{n+1} \frac{1}{\sqrt{M}}$ and  $\tilde{\widehat{V}}_{n+1}:= \sqrt{M}\widehat{V}_{n+1},$
	where $\widehat{W}^{\top}_{n+1} \in \mathbb{R}^{k \times k}$ has orthonormal rows, $\widehat{\Sigma}_{n+1} \in \mathbb{R}^{k \times k}$, with $$\widehat{\Sigma}_{n+1} = \mathrm{diag}(\sqrt{M \tildhatsigma{1}{n+1}{M}}, \sqrt{M \tildhatsigma{2}{n+1}{M}}, \dots, \sqrt{M \tildhatsigma{k}{n+1}{M}})$$ with $\tildhatsigma{i}{n+1}{M}$ singular values of $\frac{1}{M}\widetilde{\widehat{\mathbb{Y}}}_{n+1}(\widetilde{\widehat{\mathbb{Y}}}_{n+1})^{\top}$ taken in a decreasing order, and $\widehat{V}_{n+1}  \in \mathbb{R}^{k \times M}$ has orthonormal rows (notice the difference with the proof of Lemma \ref{lem: diff Proj U}). Note that $\tilde{\widehat{V}}_n$ has orthonormal rows with respect to the scalar product derived from the empirical average, namely $\widehat{\mathbb{E}}[xy] = \frac{1}{M} \langle x,y\rangle_{\mathbb{R}^{M}}$, with $x,y\in \mathbb{R}^M$.
	Analysing the first term in the right-hand side and denoting $\widetilde{\widehat{V}}^{i}_{n+1}$ the $i$-th column of $\widetilde{\widehat{V}}_{n+1}$, one obtains 
	\begin{equation*}
		\begin{aligned}
			& \sqrt{\mathbb{E}\left[\frac{1}{M}\sum_{i=1}^{M}
				\left|\left(\widehat{P}^{\alpha}_{\tildhatY{i}{n+1}}-\widehat{P}_{\tildhatY{i}{n+1}}\right)f \right|^2\mathbbm{1}_{A_{n+1}}\right]} \\
			& = \resizebox{1\linewidth}{!}{$ \sqrt{\mathbb{E}\left[\frac{1}{M}\sum_{i=1}^{M}
					\left|(\tildhatY{i}{n+1})^{\top} \left(\left(\frac{\sum_{j=1}^{M} \tildhatY{j}{n+1} (\tildhatY{j}{n+1})^{\top}}{M} + \alpha I_{k \times k}\right)^{-1} - \left(\frac{\sum_{j=1}^{M} \tildhatY{j}{n+1} (\tildhatY{j}{n+1})^{\top}}{M}\right)^{\dag} \right)\widehat{\mathbb{E}}\left[ \widetilde{\widehat{\mathbb{Y}}}_{n+1} f^{\top} \right] \right|^2 \mathbbm{1}_{A_{n+1}}\right]} $} \\
			& = \resizebox{1\linewidth}{!}{ $\sqrt{\mathbb{E}\left[\frac{1}{M}\sum\limits_{i=1}^{M}
					\left|(\tilde{\widehat{V}}_{n+1}^i)^{\top} \tilde{\widehat{\Sigma}}_{n+1}  \widehat{W}_{n+1} \left(\left(\widehat{W}_{n+1}^{\top} \tilde{\widehat{\Sigma}}_{n+1}^{2}\widehat{W}_{n+1}+ \alpha I_{k \times k}\right)^{-1} - \left(\widehat{W}_{n+1}^{\top} \tilde{\widehat{\Sigma}}_{n+1}^{2}\widehat{W}_{n+1}\right)^{\dag} \right) \widehat{W}_{n+1}^{\top} \tilde{\widehat{\Sigma}}_{n+1} \widehat{\mathbb{E}}\left[ \tilde{\widehat{V}}_{n+1} f^{\top}  \right] \right|^2 \mathbbm{1}_{A_{n+1}} \right]} $}  \\
			&= \resizebox{1\linewidth}{!}{$ \sqrt{\mathbb{E}\left[\frac{1}{M}\sum_{i=1}^{M}
					\left|(\tilde{\widehat{V}}_{n+1}^i)^{\top} \tilde{\widehat{\Sigma}}_{n+1}  \widehat{W}_{n+1} \widehat{W}_{n+1}^{\top}  \left(\left(\tilde{\widehat{\Sigma}}_{n+1}^{2}+ \alpha I_{k \times k}\right)^{-1} -  \left(\tilde{\widehat{\Sigma}}_{n+1}^{2}\right)^{\dag} \right)\widehat{W}_{n+1} \widehat{W}_{n+1}^{\top} \tilde{\widehat{\Sigma}}_{n+1} \widehat{\mathbb{E}}\left[ \tilde{\widehat{V}}_{n+1} f^{\top} \right] \right|^2 \mathbbm{1}_{A_{n+1}} \right]} $} \\
			&= \sqrt{\mathbb{E}\left[\frac{1}{M}\sum_{i=1}^{M}
				\left|(\tilde{\widehat{V}}_{n+1}^i)^{\top} \tilde{\widehat{\Sigma}}_{n+1}   \left(\left(\tilde{\widehat{\Sigma}}_{n+1}^{2}+ \alpha I_{k \times k}\right)^{-1} -  \left(\tilde{\widehat{\Sigma}}_{n+1}^{2}\right)^{\dag} \right) \tilde{\widehat{\Sigma}}_{n+1} \widehat{\mathbb{E}}\left[ \tilde{\widehat{V}}_{n+1} f^{\top} \right] \right|^2 \mathbbm{1}_{A_{n+1}} \right]}, \\
		\end{aligned}
	\end{equation*}
	where we exploited the invariance of the Euclidean norm under transposition.
	
	Notice that the event $E_{n+1}$ defined in \eqref{eq: E_n} satisfies $E_{n+1} \supset A_{n+1}$ and, hence,
	we have that $\mathrm{Rank}(\widetilde{\widehat{\mathbb{Y}}}_{n+1})=k$, and the diagonal matrix $$\tilde{\widehat{\Sigma}}_{n+1}  \left(\left(\tilde{\widehat{\Sigma}}_{n+1}^{2}+ \alpha I_{k \times k}\right)^{-1} -  \left(\tilde{\widehat{\Sigma}}_{n+1}^{2}\right)^{\dag} \right) \tilde{\widehat{\Sigma}}_{n+1}$$
	is given by:
	\begin{equation*}
		\left(\tilde{\widehat{\Sigma}}_{n+1}   \left(\left(\tilde{\widehat{\Sigma}}_{n+1}^{2}+ \alpha I_{k \times k}\right)^{-1} -  \left(\tilde{\widehat{\Sigma}}_{n+1}^{2}\right)^{\dag} \right) \tilde{\widehat{\Sigma}}_{n+1}\right)_{ii} = 
		\begin{cases}
			\frac{-\alpha}{(\tildhatsigma{i}{n+1}{M} +\alpha)} \text{ if } i=1,\dots,k.
		\end{cases} 
	\end{equation*}
	Therefore, by the orthogonality of $\tilde{\widehat{V}}_n$ with respect to the empirical scalar product and Cauchy-Schwarz inequality, one has 
	\begin{equation*}
		\begin{aligned}
			\sqrt{\mathbb{E}\left[\frac{1}{M}\sum_{i=1}^{M}
				\left|\left(\widehat{P}^{\alpha}_{\tildhatY{i}{n+1}}-\widehat{P}_{\tildhatY{i}{n+1}}\right)f\right|^2\mathbbm{1}_{A_{n+1}}\right]} &\leq  \sqrt{\mathbb{E}\left[  \frac{1}{M}\sum_{i=1}^{M}
				\left|(\tilde{\widehat{V}}_n^i)^{\top} \mathrm{diag}( \frac{-\alpha}{(\tildhatsigma{i}{n+1}{M} +\alpha)}) \widehat{\mathbb{E}}\left[ \tilde{\widehat{V}}_n f^{\top} \right] \right|^2 \mathbbm{1}_{A_{n+1}} \right]} \\
			&\leq  \sqrt{\mathbb{E}\left[  \frac{1}{M}\sum_{i=1}^{M}
				\left|(\tilde{\widehat{V}}_n^i)\right|^2 \left|\mathrm{diag}( \frac{-\alpha}{(\tildhatsigma{i}{n+1}{M} +\alpha)}) \widehat{\mathbb{E}}\left[ \tilde{\widehat{V}}_n f^{\top} \right] \right|^2 \mathbbm{1}_{A_{n+1}} \right]} \\
			&= \sqrt{\mathbb{E}\left[ 
				\left| \mathrm{diag}( \frac{\alpha}{(\tildhatsigma{i}{n+1}{M} +\alpha)}) \widehat{\mathbb{E}}\left[ \tilde{\widehat{V}}_n f^{\top} \right] \right|^2 \mathbbm{1}_{A_{n+1}} \right]} \\
			&\leq  \sqrt{\mathbb{E}\left[
				\left| \mathrm{diag}( \frac{\alpha}{(\tildhatsigma{i}{n+1}{M} +\alpha)})\right|^2 \left|\widehat{\mathbb{E}}\left[ \tilde{\widehat{V}}_n f^{\top} \right] \right|^2 \mathbbm{1}_{A_{n+1}} \right]} \\
			&\leq  \sqrt{\mathbb{E}\left[
				\left|\frac{\alpha}{\tildhatsigma{k}{n+1}{M}  +\alpha}\right|^2   \left( \frac{1}{M}\sum_{i=1}^{M} |\tilde{\widehat{V}}^i_n (f^i)^{\top}|  \right)^2 \mathbbm{1}_{A_{n+1}} \right]} \\
			&= \sqrt{\mathbb{E}\left[
				\left|\frac{\alpha}{\tildhatsigma{k}{n+1}{M}  +\alpha}\right|^2  \left( \sum_{i=1}^{M}  \frac{|\tilde{\widehat{V}}^i_n| }{\sqrt{M}} \frac{|f^i| }{\sqrt{M}}  \right)^2 \mathbbm{1}_{A_{n+1}} \right]} \\
			&\leq  \sqrt{\mathbb{E}\left[ \left|\frac{\alpha}{\tildhatsigma{k}{n+1}{M}  +\alpha}\right|^2 
				\frac{1}{M}\sum_{i=1}^{M}	\left|(\tilde{\widehat{V}}_n^i)\right|^2  \frac{1}{M}\sum_{i=1}^{M}
				\left| f^i \right|^2 \mathbbm{1}_{A_{n+1}} \right]} \\
			& \leq  \alpha \left(\min \left\{\frac{\varrho}{2}; \frac{\sigma_{Y_0}}{2} \right\}\right)^{-1} \|f\|_{L^2(\Omega^M;\mathbb{R}^{d \times M})}, \\
		\end{aligned}
	\end{equation*}
	where $\tildhatsigma{k}{n}{M}\geq \min \left\{\frac{\varrho}{2}; \frac{\sigma_{Y_0}}{2} \right\}$ in the set $E_n$ thanks to Proposition \ref{prop: lower-bound Gramian}.
	An analogous result holds for the term
	$$\sqrt{\mathbb{E}\left[\frac{1}{M}\sum_{i=1}^{M}
		\left|\left(\widehat{P}_{\widetilde{Y}_{n+1}^{i}}-\widehat{P}^{\alpha}_{\widetilde{Y}_{n+1}^{i}}\right)f\right|^2 \mathbbm{1}_{A_{n+1}}\right]}.$$ 
	
	Now, we want to evaluate the term $\sqrt{\mathbb{E}\left[\frac{1}{M}\sum_{i=1}^{M}
		\left|\left(\widehat{P}_{\tildhatY{i}{n+1}}-\widehat{P}_{\widetilde{Y}_{n+1}^{i}}\right)f\right|^2 \mathbbm{1}_{A_{n+1}} \right]}$. The matrices 
	$$\widetilde{\mathbb{X}}_{n+1} = [U_n^{\top} \widetilde{Y}^1_{n+1},\dots, U_n^{\top} \widetilde{Y}^M_{n+1}]\in \mathbb{R}^{d \times M} \text{ and } \widetilde{\widehat{\mathbb{X}}}_{n+1} = [ \widehat{U}_n^{\top}\tildhatY{1}{n+1},\dots, \widehat{U}_n^{\top}\tildhatY{M}{n+1}] \in \mathbb{R}^{d \times M}$$ 
	admit the following SVD decomposition (see Section \ref{sec: notation}):
	$$\widetilde{\mathbb{X}}_{n+1} = U^{\top}_n  \widetilde{\mathbb{Y}}_{n+1} = U^{\top}_n  \tilde{\Sigma}_{n+1} \tilde{V}_{n+1}, \quad \widetilde{\widehat{\mathbb{X}}}_{n+1} = \widehat{U}^{\top}_n  \widetilde{\widehat{\mathbb{Y}}}_{n+1}  = \widehat{U}^{\top}_n  \widehat{\Sigma}_{n+1} \widehat{V}_{n+1}.$$
	
	Then, notice that $\widehat{P}_{\widetilde{Y}_{n+1}^{i}}$ and $\widehat{P}_{\tildhatY{i}{n+1}}$ are the projectors onto the corange of $\widetilde{\mathbb{X}}_{n+1}$ and $\widetilde{\widehat{\mathbb{X}}}_{n+1}$, and, hence, by orthogonality on $\widehat{V}_{n+1}$ and $\widetilde{V}_{n+1}$ in $\mathbb{R}^{M}$ one has that
	\begin{equation*}
		\begin{aligned}
			&\frac{1}{M}\sum_{i=1}^{M}
			\left|\left(\widehat{P}_{\tildhatY{i}{n+1}}-\widehat{P}_{\widetilde{Y}_{n+1}^{i}}\right)f\right|^2\\
			= &	\frac{1}{M}\sum_{i=1}^{M}
			\left|(\tildhatY{i}{n+1})^{\top} \left(\frac{1}{M}  \widetilde{\widehat{\mathbb{Y}}}_{n+1}\widetilde{\widehat{\mathbb{Y}}}_{n+1}^{\top}\right)^{\dagger} \left(\frac{1}{M}\sum_{j=1}^{M} \tildhatY{j}{n+1} (f^{j})^{\top} \right)-(\widetilde{Y}^{j}_{n+1})^{\top} \left(\frac{1}{M}  \widetilde{\mathbb{Y}}_{n+1}\widetilde{\mathbb{Y}}_{n+1}^{\top}\right)^{\dagger} \left(\frac{1}{M}\sum_{j=1}^{M} \widetilde{Y}^{j}_{n+1} (f^{j})^{\top}\right)\right|^2  \\
			= &	\frac{1}{M}\sum_{i=1}^{M}
			\left| (\widehat{V}^{i}_{n+1})^{\top}  \widehat{\Sigma}_{n+1} \left(\frac{1}{M} \widehat{\Sigma}_{n+1}^2 \right)^{\dagger} \left(\frac{1}{M}\sum_{j=1}^{M} \widehat{\Sigma}_{n+1} \widehat{V}^{j}_{n+1} (f^{j})^{\top}\right)-(\widetilde{V}^{i}_{n+1})^{\top}  \tilde{\Sigma}_{n+1}  \left(\frac{1}{M}   \tilde{\Sigma}_{n+1}^2 \right)^{\dagger} \left(\frac{1}{M}\sum_{j=1}^{M}  \tilde{\Sigma}_{n+1} \widetilde{V}^{j}_{n+1} (f^{j})^{\top}\right)\right|^2  \\
			= &	\frac{1}{M}\sum_{i=1}^{M}
			\left| (\widehat{V}^{i}_{n+1})^{\top}  \sum_{j=1}^{M}  \widehat{V}^{j}_{n+1} (f^{j})^{\top}-(\widetilde{V}^{i}_{n+1})^{\top} \sum_{j=1}^{M}   \widetilde{V}^{j}_{n+1} (f^{j})^{\top} \right|^2 = \frac{1}{M}\sum_{i=1}^{M}
			\left| (\widehat{V}^{i}_{n+1})^{\top} \widehat{V}_{n+1} f^{\top}-(\widetilde{V}^{i}_{n+1})^{\top}   \widetilde{V}_{n+1} f^{\top}\right|^2 \\
			= & \frac{1}{M}\sum_{i=1}^{M}
			\left| e_i^{\top} \left(\widehat{V}_{n+1}^{\top} \widehat{V}_{n+1}-\widetilde{V}_{n+1}^{\top}  \widetilde{V}_{n+1}\right) f^{\top}\right|^2 \\
			= & \frac{1}{M}
			\left\|\left(\widehat{V}_{n+1}^{\top} \widehat{V}_{n+1}-\widetilde{V}_{n+1}^{\top}  \widetilde{V}_{n+1}\right) f^{\top}\right\|^2_{\mathrm{F}} \\
			\leq &
			\left| \widehat{V}_{n+1}^{\top} \widehat{V}_{n+1}-\widetilde{V}_{n+1}^{\top}  \widetilde{V}_{n+1} \right|^2 \frac{1}{M}\left\| f\right\|^2_{\mathrm{F}}
			= \left| \widehat{V}_{n+1}^{\top} \widehat{V}_{n+1}-\widetilde{V}_{n+1}^{\top}  \widetilde{V}_{n+1} \right|^2 \frac{1}{M} \sum_{i=1}^{M}| f^{i}|^2,
		\end{aligned}
	\end{equation*}
	where $e_i$ denotes the $i$-th canonical vector in $\mathbb{R}^M$.
	With the notation $(\tilde{V}_{n+1}^{\top}\tilde{V}_{n+1})^{\perp} = I_{M \times M} - (\tilde{V}_{n+1}^{\top}\tilde{V}_{n+1})^{\perp}$ it holds that
	\begin{equation*}
		\begin{aligned}
			& |(\widehat{V}_{n+1})^{\top}\widehat{V}_{n+1}-(\tilde{V}_{n+1})^{\top}\tilde{V}_{n+1}|^{2}\\
			=& | \left(\widehat{V}_{n+1}^{\top}\widehat{V}_{n+1}+(\widehat{V}_{n+1}^{\top}\widehat{V}_{n+1})^{\perp}\right)\left((\widehat{V}_{n+1})^{\top}\widehat{V}_{n+1}-(\tilde{V}_{n+1})^{\top}\tilde{V}_{n+1}\right)\left(\tilde{V}_{n+1}^{\top}\tilde{V}_{n+1}+ (\tilde{V}_{n+1}^{\top}\tilde{V}_{n+1})^{\perp}\right)|^{2}\\
			=& | \widehat{V}_{n+1}^{\top}\widehat{V}_{n+1} ((\tilde{V}_{n+1})^{\top}\tilde{V}_{n+1})^{\perp}- ((\widehat{V}_{n+1})^{\top}\widehat{V}_{n+1})^{\perp} \tilde{V}_{n+1}^{\top}\tilde{V}_{n+1} |^{2}\\
			\leq & \| \widehat{V}_{n+1}^{\top}\widehat{V}_{n+1} ((\tilde{V}_{n+1})^{\top}\tilde{V}_{n+1})^{\perp}- ((\widehat{V}_{n+1})^{\top}\widehat{V}_{n+1})^{\perp} \tilde{V}_{n+1}^{\top}\tilde{V}_{n+1} \|^{2}_{\mathrm{F}}\\
			\leq& 2 \| (\widehat{V}_{n+1}^{\top}\widehat{V}_{n+1})^{\perp} \tilde{V}_{n+1}^{\top}\tilde{V}_{n+1} \|^{2}_{\mathrm{F}},
		\end{aligned}
	\end{equation*}
	where in the first relation we interpret the projectors as matrices and in the second to last line we exploit the orthogonality of $\tilde{V}_{n+1}$ and $\widehat{V}_{n+1}$, using the same argument as developed in \cite[Theorem 2.6.1]{golub2013matrix}. With this intermediate estimate, we can prove that the error between the two projectors is bounded by the error between the two DLRA surrogates.
	Indeed, we proceed similarly to steps developed in the proof of Lemma \ref{lem: diff Proj U}. Let us define the matrix $\breve{\mathbb{X}}_{n+1} := U^{\top}_n \Sigma_{n+1} \sqrt{M} \tilde{V}_{n+1}$, where $\Sigma_n= \mathrm{diag}(\sqrt{\widetilde{\sigma}^1_{n+1}}, \dots, \sqrt{\widetilde{\sigma}^k_{n+1}})$ is the diagonal matrix whose elements on the diagonal are the square root of the singular values of the Gramian for the semi-discrete process $U_n^{\top}\widetilde{Y}_{n+1}(\omega)$, i.e.\ the singular values of $\mathbb{E}[U_n^{\top}\widetilde{Y}_{n+1}(\omega)\widetilde{Y}_{n+1}(\omega)^{\top}U_n]$. 
	Then, exploiting the orthogonality of $U_n$, $\tilde{V}_{n+1}$ and $\widehat{V}_{n+1}$, one has
	\begin{equation}\label{eq: vnvn 1}
		\begin{aligned}
			\|  (\widehat{V}_{n+1}^{\top}\widehat{V}_{n+1})^{\perp} \tilde{V}_{n+1}^{\top}\tilde{V}_{n+1} \|^{2}_{\mathrm{F}} & =  \|  (\widehat{V}_{n+1}^{\top}\widehat{V}_{n+1})^{\perp} \breve{\mathbb{X}}_{n+1}^{\top} U_n (\Sigma_{n+1} \sqrt{M})^{-1} \tilde{V}_{n+1} \|^{2}_{\mathrm{F}} \\
			&=  \|  (\widehat{V}_{n+1}^{\top}\widehat{V}_{n+1})^{\perp}  (\widetilde{\widehat{\mathbb{X}}}_{n+1}  -\breve{\mathbb{X}}_{n+1})^{\top} U_n  (\Sigma_{n+1} \sqrt{M})^{-1}  \tilde{V}_{n+1} \|^{2}_{\mathrm{F}} \\
			&\leq  |  (\widehat{V}_{n+1}^{\top}\widehat{V}_{n+1})^{\perp} |^2 \|(\widetilde{\widehat{\mathbb{X}}}_{n+1}  -\breve{\mathbb{X}}_{n+1})^{\top} U_n  (\Sigma_{n+1} \sqrt{M})^{-1}  \|^{2}_{\mathrm{F}}  | \tilde{V}_{n+1} |^2 \\
			&\leq  \|(\widetilde{\widehat{\mathbb{X}}}_{n+1}  -\breve{\mathbb{X}}_{n+1})^{\top}\|^{2}_{\mathrm{F}}  | (\Sigma_{n+1} \sqrt{M})^{-1}  |^2 \\
			& \leq  \frac{ \|  \widetilde{\widehat{\mathbb{X}}}_{n+1}   - \breve{\mathbb{X}}_{n+1} \|^{2}_{\mathrm{F}}}{(\widetilde{\sigma}^k_{n+1})M},
		\end{aligned}
	\end{equation}
	and it also holds that
	\begin{equation}\label{eq: vnvn 2}
		\begin{aligned}
			\|  (\widehat{V}_{n+1}^{\top}\widehat{V}_{n+1})^{\perp} \tilde{V}_{n+1}^{\top}\tilde{V}_{n+1} \|^{2}_{\mathrm{F}} & \leq	| (\widehat{V}_{n+1}^{\top}\widehat{V}_{n+1})^{\perp}| \|  \tilde{V}_{n+1}^{\top}\tilde{V}_{n+1} \|^{2}_{\mathrm{F}} \leq k.
		\end{aligned}
	\end{equation}
	
	Now, define the event  $\Omega_R = \{ \boldsymbol{\omega} \ : \ \frac{1}{M}\sum_{i=1}^{M} |f^{i}(\omega_i) |^2 \leq R^2\}$ for $R> \|f\|_{L^2(\Omega^M;\mathbb{R}^{d \times M})}$.
	Then, via Lemma \ref{cor: Bernstein ineq}, Cauchy-Schwarz inequality, property of operator norms for matrices, orthogonality of the projections, and relations \eqref{eq: vnvn 1} and \eqref{eq: vnvn 2} we retrieve that
	\begin{equation*}
		\begin{aligned}
			&	\sqrt{\mathbb{E}
				\left[\frac{1}{M}\sum_{i=1}^{M}
				\left|\left(\widehat{P}_{\widetilde{Y}_{n+1}^{i}}-\widehat{P}_{\tildhatY{i}{n+1}}\right)f\right|^2 \mathbbm{1}_{A_{n+1}} \right]} \\
			\leq & \sqrt{\mathbb{E}
				\left[\frac{1}{M}\sum_{i=1}^{M}
				\left|\left(\widehat{P}_{\widetilde{Y}_{n+1}^{i}}-\widehat{P}_{\tildhatY{i}{n+1}}\right)f\right|^2 \mathbbm{1}_{A_{n+1}} \mathbbm{1}_{\Omega_R}\right]} \\
			& + \sqrt{\mathbb{E}
				\left[\frac{1}{M}\sum_{i=1}^{M}
				\left|\left(\widehat{P}_{\widetilde{Y}_{n+1}^{i}}-\widehat{P}_{\tildhatY{i}{n+1}}\right)f\right|^2 \mathbbm{1}_{A_{n+1}}  \mathbbm{1}_{\Omega_R^C} \right]} \\
			\leq & \sqrt{\mathbb{E}
				\left[\frac{1}{M}\sum_{i=1}^{M}
				\left|\left(\widehat{P}_{\widetilde{Y}_{n+1}^{i}}-\widehat{P}_{\tildhatY{i}{n+1}}\right)f\right|^2 \mathbbm{1}_{A_{n+1}} \mathbbm{1}_{\Omega_R}\right]} \\
			& + \sqrt{\mathbb{E}
				\left[
				2\|  (\widehat{V}_{n+1}^{\top}\widehat{V}_{n+1})^{\perp} \tilde{V}_{n+1}^{\top}\tilde{V}_{n+1} \|^{2}_{\mathrm{F}} \frac{1}{M}\left\|f\right\|^2_{\mathrm{F}} \mathbbm{1}_{A_{n+1}} \mathbbm{1}_{\Omega_R^C} \right]} \\
			\leq&\sqrt{2}  \sqrt{\mathbb{E}\left[\|  (\widehat{V}_{n+1}^{\top}\widehat{V}_{n+1})^{\perp} \tilde{V}_{n+1}^{\top}\tilde{V}_{n+1} \|^2_{\mathrm{F}} \mathbbm{1}_{A_{n+1}}
				\right]} R + \sqrt{2k}  \sqrt{\mathbb{E}
				\left[\frac{1}{M}\sum_{i=1}^{M}
				\left|f^{i}\right|^2  \mathbbm{1}_{\Omega_R^C} \right]}   \\
			\leq&  \sqrt{2}  \sqrt{\mathbb{E}\left[\frac{1}{M}
				\frac{ \|  \widetilde{\widehat{\mathbb{X}}}_{n+1}  - \breve{\mathbb{X}}_{n+1} \|^{2}_{\mathrm{F}} }{ \widetilde{\sigma}^k_{n+1}}
				\right]} R + \sqrt{2k}  \|f\|_{L^4(\Omega^M;\mathbb{R}^{d \times M})} \sqrt[4]{C} e^{-\frac{cMR^2}{4}} .\\
		\end{aligned}
	\end{equation*}
	Then, we can proceed to bound the error in the following way
	\begin{equation*}
		\begin{aligned}
			&\sqrt{\mathbb{E}
				\left[\frac{1}{M}\sum_{i=1}^{M}
				\left|\left(\widehat{P}_{\widetilde{Y}_{n+1}^{i}}-\widehat{P}_{\tildhatY{i}{n+1}}\right)f\right|^2 \mathbbm{1}_{A_{n+1}} \right]} \\
			\leq&  \sqrt{\frac{2}{\widetilde{\sigma}^k_{n+1}}}  \left( \sqrt{\mathbb{E}\left[\frac{1}{M}
				\|  \widetilde{\widehat{\mathbb{X}}}_{n+1}  - \widetilde{\mathbb{X}}_{n+1} \|^{2}_{\mathrm{F}} \right]} +\sqrt{\mathbb{E}\left[\frac{1}{M}
				\|  \breve{\mathbb{X}}_{n+1}  -  \widetilde{\mathbb{X}}_{n+1} \|^{2}_{\mathrm{F}} 
				\right]} \right) R + \sqrt{2k}  \|f\|_{L^4(\Omega^M;\mathbb{R}^{d \times M})} \sqrt[4]{C} e^{-\frac{cMR^2}{4}}  \\
		\end{aligned}
	\end{equation*}
	Then, we can bound the difference between the sample matrices of the particle system and the semidiscretized DLRA in the following way
	\begin{equation*}
		\begin{aligned}
			&\sqrt{\mathbb{E}\left[\frac{1}{M}
				\|  \widetilde{\widehat{\mathbb{X}}}_{n+1}  - \widetilde{\mathbb{X}}_{n+1} \|^{2}_{\mathrm{F}} \right]} \\
			\leq &  \sqrt{\mathbb{E}\left[\frac{1}{M} 
				\|  \left(\widehat{\mathbb{X}}_{n}  + P_{\widehat{U}_{n}} \widehat{a}_n \Delta t + P_{\widehat{U}_{n}} \widehat{b}_n \Delta \mathbb{W}_n \right)- \left(\mathbb{X}_{n}  + P_{U_{n}} a_n \Delta t + P_{U_{n}} b_n \Delta  \mathbb{W}_n \right) \|^{2}_{\mathrm{F}} 
				\right]}  \\
			\leq& \resizebox{1\linewidth}{!}{$ \sqrt{\mathbb{E}\left[\frac{1}{M} 
					\|  \left(\widehat{\mathbb{X}}_{n} -\mathbb{X}_{n}\right) + P_{\widehat{U}_{n}} \left(\widehat{a}_n - a_n\right) \Delta t + P_{\widehat{U}_{n}} \left(\widehat{b}_n - b_n\right) \Delta  \mathbb{W}_n + \left(P_{\widehat{U}_{n}} - P_{U_{n}} \right)a_n \Delta t + \left(P_{\widehat{U}_{n}} - P_{U_{n}} \right) b_n \Delta  \mathbb{W}_n  \|^{2}_{\mathrm{F}} \right]} $}\\
		\end{aligned}
	\end{equation*}
	Now, by Lipschitzianity of the drift and the diffusion, by the orthogonality of $\widehat{U}_n$, as well as, thanks to the independence of the Brownian increments, one has that
	\begin{equation*}
		\begin{aligned}
			&\sqrt{\mathbb{E}\left[\frac{1}{M} 
				\| P_{\widehat{U}_{n}} \left(\widehat{a}_n - a_n\right) \Delta t + P_{\widehat{U}_{n}} \left(\widehat{b}_n - b_n\right) \Delta  \mathbb{W}_n \|^{2}_{\mathrm{F}} 
				\right]} \\
			\leq & \sqrt{\mathbb{E}\left[\frac{1}{M} 
				\| P_{\widehat{U}_{n}} \left(\widehat{a}_n - a_n\right) \Delta t  \|^{2}_{\mathrm{F}} 	\right]}  +  \sqrt{\mathbb{E}\left[\frac{1}{M} 
				\|  P_{\widehat{U}_{n}} \left(\widehat{b}_n - b_n\right) \Delta  \mathbb{W}_n \|^{2}_{\mathrm{F}} 
				\right]} \\
			\leq & \sqrt{\mathbb{E}\left[\frac{1}{M} 
				\| \left(\widehat{a}_n - a_n\right) \Delta t  \|^{2}_{\mathrm{F}} 	\right]}  +  \sqrt{\mathbb{E}\left[\frac{1}{M} 
				\|  \left(\widehat{b}_n - b_n\right) \Delta  \mathbb{W}_n \|^{2}_{\mathrm{F}} 
				\right]} \\
			= & \sqrt{\mathbb{E}\left[\frac{1}{M} \sum_{i=1}^{M}
				| \left(\widehat{a}^i_n - a^i_n\right) \Delta t  |^{2}	\right]}  +  \sqrt{\mathbb{E}\left[\frac{1}{M}  \sum_{i=1}^{M}
				\|  \left(\widehat{b}^i_n - b^i_n\right) \Delta  W^i_n \|^{2}_{\mathrm{F}} 
				\right]} \\
			\leq & C_{\mathrm{Lip}} \sqrt{\mathbb{E}\left[ \frac{1}{M} \|  \widehat{\mathbb{X}}_{n}  -  \widetilde{\mathbb{X}}_{n} \|^{2}_{\mathrm{F}} \right]}  (\Delta t + \sqrt{\Delta t}) =  \sqrt{2}C_{\mathrm{Lip}} e_n  (\Delta t + \sqrt{\Delta t}).
		\end{aligned}
	\end{equation*}
	Using this last result, one obtains
	\begin{equation*}
		\begin{aligned}
			&\sqrt{\mathbb{E}
				\left[\frac{1}{M}\sum_{i=1}^{M}
				\left|\left(\widehat{P}_{\widetilde{Y}_{n+1}^{i}}-\widehat{P}_{\tildhatY{i}{n+1}}\right)f\right|^2 \mathbbm{1}_{A_{n+1}} \right]} \\
			\leq & \resizebox{1\linewidth}{!}{$ \Bigg(\frac{2}{\sqrt{\widetilde{\sigma}^k_{n+1}}} \left( e_n + \tilde{G}_1 \left(e_n +  \frac{1}{\sqrt{M}} (\Delta t + \sqrt{\Delta t}) \right) + \left( \sqrt[4]{C} e^{-\frac{cMR^2}{4}} C_{\mathrm{lgb}} (1 + \sqrt[4]{K_4})\right) (\Delta t + \sqrt{\Delta t}) + C_{\mathrm{lip}} e_n (\Delta t + m \sqrt{\Delta t}) \right) $} \\
			&\quad + \frac{1}{\sqrt{M}}\frac{2\sqrt{k K_4(T)}}{\widetilde{\sigma}^k_{n+1}}\Bigg)R + \sqrt{2k}  \|f\|_{L^4(\Omega^M;\mathbb{R}^{d \times M})} \sqrt[4]{C} e^{-\frac{cMR^2}{4}},
		\end{aligned}
	\end{equation*}
	where in the last line we employ Lemma \ref{lem: diff Proj U} for the terms involving the differences between $P_{\widehat{U}_{n}}$ and $P_{U_{n}}$, where $\tilde{G}_1$ is a constant dependent on the smallest singular value $\widetilde{\sigma}_n^k$, the linear growth bound, and the convergence of the Monte Carlo estimator by proceeding as in \eqref{eq: MC in xx_T} in the proof of Lemma \ref{lem: diff Proj U} for the term $\sqrt{\mathbb{E}\left[\frac{1}{M} \|  \breve{\mathbb{X}}_{n+1}  -  \widetilde{\mathbb{X}}_{n+1} \|^{2}_{\mathrm{F}}\right]}$. This relation  simplifies to
	\begin{equation*}
		\begin{aligned}
			&\sqrt{\mathbb{E}
				\left[\frac{1}{M}\sum_{i=1}^{M}
				\left|\left(\widehat{P}_{\widetilde{Y}_{n+1}^{i}}-\widehat{P}_{\tildhatY{i}{n+1}}\right)f\right|^2 \mathbbm{1}_{A_{n+1}} \right]}  \\
			\leq &\tilde{G}_2 \left( \frac{1}{\sqrt{M}} +e_n\right) R +  \sqrt{2k}\sqrt[4]{\tilde{C}} e^{-\frac{cMR^2}{4}} \|f\|_{L^4(\Omega^M;\mathbb{R}^{d \times M})},
		\end{aligned}
	\end{equation*}
	for a positive constant $\tilde{G}_2$ dependent on the smallest singular value $\widetilde{\sigma}^k_{n+1}$ and a positive constant $\tilde{C}$.
	
	To conclude, we need to estimate this projector bound on the event $A_{n+1}^C$. Thanks to the fact that the empirical norm of regularized orthogonal projectors is smaller than $1$, Lemma \ref{lem: prop bound sigma tilde} and to Proposition \ref{prop: lower-bound Gramian} one obtains
	\begin{equation}\label{eq: bound requiring 8th}
		\begin{aligned}
			\sqrt{\mathbb{E}\left[\frac{1}{M}\sum_{i=1}^{M}
				\left|\left(\widehat{P}^{\alpha}_{\tildhatY{i}{n+1}}-\widehat{P}^{\alpha}_{\widetilde{Y}_{n+1}^{i}}\right)f \right|^2 \mathbbm{1}_{A_{n+1}^C}\right]} \leq &	\sqrt{ \mathbb{E}\left[\frac{1}{M}\sum_{i=1}^{M}
				\left|\widehat{P}^{\alpha}_{\tildhatY{i}{n+1}}[f]\right|^2\mathbbm{1}_{A_{n+1}^C}+\left|\widehat{P}^{\alpha}_{\widetilde{Y}_{n+1}^{i}}[f]\right|^2\mathbbm{1}_{A_{n+1}^C}\right]} \\
			\leq & \sqrt{\mathbb{E}\left[\frac{1}{M}\sum_{i=1}^{M}
				\left|f\right|^2\mathbbm{1}_{A_{n+1}^C}\right]} \leq  \mathbb{P}(\mathbbm{1}_{A_{n+1}^C})^{\frac{1}{4}}\|f\|_{L^4(\Omega^M;\mathbb{R}^{d \times M})} \\
			\leq & \sqrt[4]{C_5 + C_1 \frac{1}{\varrho^{8}}} \frac{1}{\sqrt{M}}\|f\|_{L^4(\Omega^M;\mathbb{R}^{d \times M})}.
		\end{aligned}
	\end{equation}
	
	Finally, putting all the computations together, we have
	\begin{equation*}
		\begin{aligned}
			&	\sqrt{\mathbb{E}\left[\frac{1}{M}\sum_{i=1}^{M}
				\left|\left(\widehat{P}^{\alpha}_{\tildhatY{i}{n+1}}-\widehat{P}^{\alpha}_{\widetilde{Y}_{n+1}^{i}}\right)f\right|^2\right]} \\
			\leq &	\sqrt{\mathbb{E}\left[\frac{1}{M}\sum_{i=1}^{M}
				\left|\left(\widehat{P}^{\alpha}_{\tildhatY{i}{n+1}}-\widehat{P}_{\tildhatY{i}{n+1}}\right)f\right|^2 \mathbbm{1}_{A_{n+1}} \right]} + 	\sqrt{\mathbb{E}\left[\frac{1}{M}\sum_{i=1}^{M}
				\left|\left(\widehat{P}_{\tildhatY{i}{n+1}}-\widehat{P}_{\widetilde{Y}_{n+1}^{i}}\right)f\right|^2\mathbbm{1}_{A_{n+1}}\right]}\\
			& + \sqrt{\mathbb{E}\left[\frac{1}{M}\sum_{i=1}^{M}
				\left|\left(\widehat{P}_{\widetilde{Y}_{n+1}^{i}}-\widehat{P}^{\alpha}_{\widetilde{Y}_{n+1}^{i}}\right)f\right|^2\mathbbm{1}_{A_{n+1}}\right]}+ 	\sqrt{\mathbb{E}\left[\frac{1}{M}\sum_{i=1}^{M}\left|\left(\widehat{P}^{\alpha}_{\tildhatY{i}{n+1}}-\widehat{P}^{\alpha}_{\widetilde{Y}_{n+1}^{i}}\right)f\right|^2\mathbbm{1}_{A_{n+1}^C}\right]} \\
			\leq  &   \left( 2 \alpha   \left(\min \left\{\frac{\varrho}{2}; \frac{\sigma_{Y_0}}{2} \right\}\right)^{-1}   + \sqrt{2k}  \sqrt[4]{\tilde{C}} e^{-\frac{\tilde{c}MR^2}{4}}  + \sqrt[4]{C_5 + C_1 \frac{1}{\varrho^{8}}} \frac{1}{\sqrt{M}} \right) \|f\|_{L^4(\Omega^M;\mathbb{R}^{d \times M})} + \tilde{G}_2 \left( e_n +  \frac{1}{\sqrt{M}}  \right) R,
		\end{aligned}
	\end{equation*}
	where we employ the monotonicity of the $L^p$ moments with respect to $p\geq1$, and via taking the optimal value $\alpha = \frac{1}{\sqrt[3]{M}}$ imposed by Lemma \ref{lem: discr proj Y} we get the thesis.
\end{proof} 

\section{The case of adaptive regularization}\label{app: proof adap}

\begin{proof}[\bfseries Proof of Lemma \ref{lem: error between stochastic proj - full rank}]\label{proof lem: error between stochastic proj - full rank}
	For $f \in L^4(\Omega^{M};\mathbb{R}^{d \times M})$, the proof follows similarly to Lemma \ref{lem: error between stochastic proj}. Indeed we will add and subtract suitable terms and we conclude via finding reasonable bounds for the appeared differences related to projector-type quantities. We split the error in two components, exploiting the invariance of the Euclidean norm under transposition and the triangular inequality:
	\begin{equation*}
		\begin{aligned}
			&\|\left(P_{\widetilde{Y}_{n+1}}-\widehat{P}_{\widetilde{Y}_{n+1}}\right)f \mathbbm{1}_{\Sigma_{n+1}^{\gamma}} \|_{L^2(\Omega^{M};\mathbb{R}^{d \times M})} =  \left(\mathbb{E}\left[\frac{1}{M} \sum_{i=1}^{M} \left|\left(P_{\widetilde{Y}_{n+1}^{i}}- \widehat{P}_{\widetilde{Y}_{n+1}^{i}} \right)f\right|^2 \mathbbm{1}_{\Sigma_{n+1}^{\gamma}} \right]\right)^{\frac{1}{2}}\\
			=& \resizebox{1\linewidth}{!}{$ \left(\mathbb{E}\left[\frac{1}{M} \sum_{i=1}^{M} \left|(\widetilde{Y}_{n+1}^{i})^{\top}\left(\mathbb{E}\left[ \frac{1}{M} \widetilde{\mathbb{Y}}_{n+1}\widetilde{\mathbb{Y}}_{n+1}^{\top}\right]\right)^{-1}\mathbb{E}\left[ \frac{1}{M} \widetilde{Y}_{n+1} f^{\top} \right] - (\widetilde{Y}_{n+1}^{i})^{\top}\left(\frac{\sum_{j=1}^{M} \widetilde{Y}_{n+1}^j (\widetilde{Y}_{n+1}^j)^{\top}}{M}\right)^{-1}\widehat{\mathbb{E}}\left[ \widetilde{Y}_{n+1} f^{\top} \right] \right|^2   \mathbbm{1}_{\Sigma_{n+1}^{\gamma}} \right]\right)^{\frac{1}{2}} $}\\
			\leq& \underbrace{\left(\mathbb{E}\left[\frac{1}{M} \sum_{i=1}^{M} \left|(\widetilde{Y}_{n+1}^{i})^{\top}\left(\left(\mathbb{E}\left[  \frac{1}{M} \widetilde{\mathbb{Y}}_{n+1}\widetilde{\mathbb{Y}}_{n+1}^{\top} \right]\right)^{-1}-\left(\frac{\sum_{j=1}^{M} \widetilde{Y}_{n+1}^j (\widetilde{Y}_{n+1}^j)^{\top}}{M} \right)^{-1}\right)\mathbb{E}\left[ \frac{1}{M} \widetilde{\mathbb{Y}}_{n+1} f^{\top}  \right] \right|^2   \mathbbm{1}_{\Sigma_{n+1}^{\gamma}} \right]\right)^{\frac{1}{2}}}_{=:L_1}\\
			&+ \underbrace{\left(\mathbb{E}\left[\frac{1}{M} \sum_{i=1}^{M} \left|(\widetilde{Y}_{n+1}^{i})^{\top}\left(\frac{\sum_{j=1}^{M} \widetilde{Y}_{n+1}^j (\widetilde{Y}_{n+1}^j)^{\top}}{M} \right)^{-1}\left(\mathbb{E}\left[ \frac{1}{M} \widetilde{\mathbb{Y}}_{n+1} f^{\top}  \right] -\widehat{\mathbb{E}}\left[ \frac{1}{M} \widetilde{\mathbb{Y}}_{n+1} f^{\top} \right]  \right)\right|^2   \mathbbm{1}_{\Sigma_{n+1}^{\gamma}} \right]\right)^{\frac{1}{2}}}_{=:L_2}.\\
		\end{aligned}
	\end{equation*} 
	By assumption we have that $\frac{\sum_{j=1}^{M} \widetilde{Y}_{n+1}^j (\widetilde{Y}_{n+1}^j)^{\top}}{M}\succeq \gamma I_{k \times k}$ and, hence, recalling $K_4(T)$ in relation \eqref{eq: Lp moment bound semidiscretized X}, we bound $L_1$ as follows
	\begin{equation*}
		\begin{aligned}
			L_1  =& \left(\mathbb{E}\left[\frac{1}{M} \sum_{i=1}^{M} \left|(\widetilde{Y}_{n+1}^{i})^{\top}\left(\left(\mathbb{E}\left[ \frac{1}{M} \widetilde{\mathbb{Y}}_{n+1}\widetilde{\mathbb{Y}}_{n+1}^{\top} \right]\right)^{-1}-\left(\frac{\sum_{j=1}^{M} \widetilde{Y}_{n+1}^j (\widetilde{Y}_{n+1}^j)^{\top}}{M} \right)^{-1}\right)\mathbb{E}\left[ \frac{1}{M} \widetilde{\mathbb{Y}}_{n+1} f^{\top}  \right] \right|^2 \mathbbm{1}_{\Sigma_{n+1}^{\gamma}} \right]\right)^{\frac{1}{2}} \\
			= & \resizebox{1\linewidth}{!}{$ \mathbb{E}\Big[\frac{1}{M} \sum_{i=1}^{M} \Big|(\widetilde{Y}_{n+1}^{i})^{\top}\Big(\frac{\sum_{j=1}^{M} \widetilde{Y}_{n+1}^j (\widetilde{Y}_{n+1}^j)^{\top}}{M} \Big)^{-1}\Big(\frac{\sum_{j=1}^{M} \widetilde{Y}_{n+1}^j (\widetilde{Y}_{n+1}^j)^{\top}}{M} - \mathbb{E}\Big[ \frac{1}{M} \widetilde{\mathbb{Y}}_{n+1}\widetilde{\mathbb{Y}}_{n+1}^{\top}\Big]\Big) \mathbb{E}\Big[ \frac{1}{M} \widetilde{\mathbb{Y}}_{n+1}\widetilde{\mathbb{Y}}_{n+1}^{\top} \mathbbm{1}_{\Sigma_{n+1}^{\gamma}} \Big]^{-1}\mathbb{E}\Big[ \frac{1}{M} \widetilde{\mathbb{Y}}_{n+1} f^{\top}  \Big] \Big|^2 \mathbbm{1}_{\Sigma_{n+1}^{\gamma}} \Big]^{\frac{1}{2}} $} \\
			= & \frac{1}{\gamma}  \Big(\mathbb{E}\Big[ \left(\frac{1}{M} \sum_{i=1}^{M} \Big|(\widetilde{Y}_{n+1}^{i})^{\top}\Big(\frac{\sum_{j=1}^{M} \widetilde{Y}_{n+1}^j (\widetilde{Y}_{n+1}^j)^{\top}}{M} \Big)^{-\frac{1}{2}} \big|^2\right) \Big|\Big(\frac{\sum_{j=1}^{M} \widetilde{Y}_{n+1}^j (\widetilde{Y}_{n+1}^j)^{\top}}{M} - \mathbb{E}\Big[ \frac{1}{M} \widetilde{\mathbb{Y}}_{n+1}\widetilde{\mathbb{Y}}_{n+1}^{\top}\Big]\Big)\Big|^2 \\
			& \cdot  \Big|\mathbb{E}\Big[  \mathbb{E}\Big[ \frac{1}{M} \widetilde{\mathbb{Y}}_{n+1}\widetilde{\mathbb{Y}}_{n+1}^{\top} \mathbbm{1}_{\Sigma_{n+1}^{\gamma}} \Big]^{-\frac{1}{2}} \frac{1}{M}\widetilde{\mathbb{Y}}_{n+1}f^{\top} \Big] \Big|^2 \mathbbm{1}_{\Sigma_{n+1}^{\gamma}} \Big]\Big)^{\frac{1}{2}} \\
			\leq & \frac{k}{\gamma} \|f\|_{L^2(\Omega^{M};\mathbb{R}^{d \times M})}  \Big(\mathbb{E} \left[\Big|\frac{\sum_{j=1}^{M} \widetilde{Y}_{n+1}^j (\widetilde{Y}_{n+1}^j)^{\top}}{M} - \mathbb{E}\Big[ \frac{1}{M} \widetilde{\mathbb{Y}}_{n+1}\widetilde{\mathbb{Y}}_{n+1}^{\top} \Big]\Big|^2\right]\Big)^{\frac{1}{2}}\\
			\leq & \frac{k}{\gamma} \|f\|_{L^2(\Omega^{M};\mathbb{R}^{d \times M})}  \sqrt{K_4(T)}\frac{1}{\sqrt{M}},\\
		\end{aligned}
	\end{equation*} 
	via independence of samples $\widetilde{Y}_{n+1}^{i}$ and Monte Carlo convergence.
	Similarly to the proof of Lemma \ref{lem: error between stochastic proj}, recalling the definition of $\mathrm{Var}$ in \eqref{eq: var} we get for $L_2$ 
	\begin{equation*}
		\begin{aligned}
			L_2  =& \left(\mathbb{E}\left[\frac{1}{M} \sum_{i=1}^{M} \left|(\widetilde{Y}_{n+1}^{i})^{\top}\left(\frac{\sum_{j=1}^{M} \widetilde{Y}_{n+1}^j (\widetilde{Y}_{n+1}^j)^{\top}}{M} \right)^{-1}\left(\widehat{\mathbb{E}}\left[\widetilde{\mathbb{Y}}_{n+1} f^{\top} \right]-\mathbb{E}\left[ \frac{1}{M} \widetilde{\mathbb{Y}}_{n+1} f^{\top}  \right]\right) \right|^2 \mathbbm{1}_{\Sigma_{n+1}^{\gamma}} \right]\right)^{\frac{1}{2}} \\
			= & \frac{\sqrt{k}}{\sqrt{\gamma}} \left(\mathbb{E}\left[ \left|\left(\widehat{\mathbb{E}}\left[ \widetilde{\mathbb{Y}}_{n+1} f^{\top} \right]-\mathbb{E}\left[ \frac{1}{M} \widetilde{\mathbb{Y}}_{n+1} f^{\top}  \right]\right) \right|^2\right]\right)^{\frac{1}{2}} \\
			\leq & \frac{\sqrt{k}}{\sqrt{\gamma}} \frac{1}{\sqrt{M}} \sqrt{\mathbb{E}[\|\widetilde{Y}_{n+1} (f^1)^{\top}\|_{\mathrm{F}}^2]} \\
			\leq & \frac{\sqrt{k}}{\sqrt{\gamma}} \frac{1}{\sqrt{M}} \sqrt{K_4(T)} \|f\|_{L^4(\Omega^{M};\mathbb{R}^{d \times M})}. \\
		\end{aligned}
	\end{equation*} 
	Putting all the intermediate steps together we get the statement \eqref{eq: Y projectors error - full rank}. Concerning \eqref{eq: Y projectors error - full rank 2}, computations are analogous to the first part and to the proof of Lemma \ref{lem: discr proj Y}.
\end{proof}

\section{The case of degenerate diffusion}\label{app: proof deg diff}
\begin{proof}[\bfseries Proof of Lemma \ref{lemma: regularized diffusion empir Proj}]
	The proof follows closely the proof of Lemma \ref{lem: discr proj Y}.
	Notice that in the event $A_{n+1}$ defined in \eqref{eq: A_n reg diff}, both $\widetilde{Y}_{n+1}$ and $\widetilde{\widehat{Y}}_{n+1}$ have rank $k$. Indeed, on \(A_{n+1}\), the empirical
	Gramians are uniformly non-degenerate, so the perturbation argument used in
	Lemma~\ref{lem: discr proj Y} applies.
	
	Part of the proof is similar to the one of Lemma \ref{lem: error between stochastic proj - full rank}, as by Assumption \ref{ass: gamma cov} and by regularization of the diffusion the involved projectors are all of full rank $k$.
	As proposed in the proof of Lemma \ref{lem: discr proj Y}, let us define the event  $\Omega_R = \{ \boldsymbol{\omega} \ : \ \frac{1}{M}\sum_{i=1}^{M} |f^{i}(\omega_i)|^2 \leq R^2 \}$ for $R> \|f\|_{L^2(\Omega^M;\mathbb{R}^{d \times M})}$. Then for a positive constant $\tilde{G}_2$ dependent on the smallest singular value $\widetilde{\sigma}^k_{n+1}$ we have
	\begin{equation*}
		\begin{aligned}
			\sqrt{\mathbb{E}\left[\frac{1}{M}\sum_{i=1}^{M}
				\left|\left(\widehat{P}_{\widetilde{Y}_{n+1}^{i}}- \widehat{P}_{\tildhatY{i}{n+1}}\right)f\right|^2 \mathbbm{1}_{A_{n+1}}\right]}  \leq & \sqrt{\mathbb{E}[\frac{1}{M} \sum_{i=1}^{M} \left| \left(\widehat{P}_{\widetilde{Y}_{n+1}^{i}}-\widehat{P}_{\tildhatY{i}{n+1}}\right)f\right|^2 \mathbbm{1}_{A_{n+1}} \mathbbm{1}_{\Omega_R} ]} \\
			&+ \sqrt{\mathbb{E}[\frac{1}{M} \sum_{i=1}^{M} \left| \left(\widehat{P}_{\widetilde{Y}_{n+1}^{i}}-\widehat{P}_{\tildhatY{i}{n+1}}\right)f\right|^2 \mathbbm{1}_{A_{n+1}} \mathbbm{1}_{\Omega_R^C} ]} \\
			\leq&  \sqrt{\frac{2}{\widetilde{\sigma}^k_{n+1}}}  \sqrt{\mathbb{E}\left[\frac{1}{M} \left(
				\|  \widetilde{\widehat{\mathbb{X}}}_{n+1}  - \widetilde{\mathbb{X}}_{n+1} \|^{2}_{\mathrm{F}} + \|  \breve{\mathbb{X}}_{n+1}  -  \widetilde{\mathbb{X}}_{n+1} \|^{2}_{\mathrm{F}} \right)
				\right]}R \\
			&+\sqrt{2k}  \sqrt[4]{C} e^{-\frac{cMR^2}{4}}  \|f\|_{L^4(\Omega^M;\mathbb{R}^{d \times M})} \\
			\leq  &   \tilde{G}_2 \left(e_n +  \frac{1}{\sqrt{M}} + \sqrt{\mathbb{E}\left[\frac{1}{M} 
				\| P_{\widehat{U}_n} \beta \Delta W_n \|^{2}_{\mathrm{F}} \right]}\right) R \\
			&+  \sqrt{2k}  \sqrt[4]{C} e^{-\frac{cMR^2}{4}} \|f\|_{L^4(\Omega^M;\mathbb{R}^{d \times M})}.
		\end{aligned}
	\end{equation*}
	Finally, the proof follows to the one of Lemma \ref{lem: discr proj Y} verbatim.
\end{proof} 

\end{document}